\documentclass[10pt,a4paper]{article}
\usepackage{amsmath,amsfonts,amssymb,amsthm}
\usepackage{paralist}

\newenvironment{mabstract}
{\begin{quote}\small {\bfseries Abstract. }}{\end{quote}\par}

\newenvironment{mkeywords}
{\begin{quote}\small {\bfseries Keywords. }}{\end{quote}\par}

\newenvironment{msubjclass}
{\begin{quote}\small {\bfseries AMS Subject Classification. }}{\end{quote}\par}

\newtheorem{theorem}{Theorem}[section]

\newcounter{appendixcounter}

\newtheorem{theoremappendix}{Theorem}[appendixcounter]

\newtheorem{lemma}{Lemma}[section]

\newtheorem{assumption}{Assumption}[section]

\newtheorem{assumptionappendix}{Assumption}[appendixcounter]

\newtheorem{corollary}{Corollary}[section]

\newtheorem{remark}{Remark}[section]

\newtheorem*{vremark}{Remark}

\newtheorem*{rproof}{\bfseries Proof}
\newenvironment{sproof}
{\begin{rproof}\rm}{\qed \end{rproof}}

\newtheorem{bproof}{}

\newtheorem{bcase}{}[section]

\newtheorem{step}{Step}

\makeatletter
\renewcommand\@seccntformat[1]{{\csname thesection\endcsname}.\hspace{0.5em}}
\makeatother

\begin{document}

\title{Convergence and Convergence Rate of Stochastic Gradient Search
in the Case of Multiple and Non-Isolated Extrema} 

\author{Vladislav B. Tadi\'{c}
\thanks{Department of Mathematics, University of Bristol, United Kingdom.  
(email: v.b.tadic@bristol.ac.uk). } }
\date{}

\maketitle

\begin{mabstract} 
The asymptotic behavior of stochastic gradient algorithms is studied. 
Relying on results from differential geometry 
(Lojasiewicz gradient inequality), 
the single limit-point convergence of the algorithm iterates is demonstrated 
and relatively tight bounds on the convergence rate 
are derived. 
In sharp contrast to the existing asymptotic results, 
the new results presented here allow 
the objective function to have multiple and non-isolated minima. 
The new results also offer new insights into the asymptotic properties of 
several classes of recursive algorithms which are 
routinely used in engineering, statistics, machine learning and operations research. \end{mabstract}

\begin{mkeywords}
Stochastic gradient search, 
single limit-point convergence, 
convergence rate, 
Lojasiewicz gradient inequality, 
supervised learning, reinforcement learning, recursive principal component analysis, 
recursive maximum likelihood estimation, 
recursive prediction error identification. 
\end{mkeywords}

\begin{msubjclass}
Primary 62L20; 
Secondary 90C15, 93E12, 93E35.
\end{msubjclass}

\section{Introduction} 

Stochastic gradient search (also known as stochastic gradient algorithm) 
is a stochastic optimization method of the stochastic approximation type. 
It computes minima (or maxima) of an objective function 
whose values are available only through noise-corrupted observations. 
Stochastic gradient search has found a wide range of applications in diverse areas 
such as statistical inference, signal processing, automatic control, communications, 
machine learning and operations research. 
Due to its practical importance, 
the asymptotic analysis of stochastic gradient search 
is the subject of a number of papers and books
(see \cite{benveniste}, \cite{borkar}, \cite{kushner&yin}, \cite{ljung2}, 
\cite{polyak}, \cite{spall} 
and references cited therein). 
The focus of the analysis has been on the single-point convergence and convergence rate
since these properties most precisely characterize 
the asymptotic behavior of stochastic gradient search. 
Unfortunately, the existing results on the asymptotic properties of stochastic gradient search 
hold under restrictive conditions. 
Unless each minimum of the objective function (minimized by stochastic gradient search) 
is isolated, 
the existing results do not guarantee that stochastic gradient search converges to a single point. 
In addition to this, the existing results do not provide any information on the convergence rate
unless the objective function is strongly convex. 
In the case of many stochastic gradient algorithms met in practice, 
these conditions are not only hard to verify, but likely not to hold at all. 

In this paper, the convergence and convergence rate of stochastic gradient search 
are studied under conditions which allow the objective function 
to be non-convex and 
to have multiple and non-isolated minima. 
Relying on results from differential geometry (Lojasiewicz gradient inequality), 
almost sure single-limit point convergence is demonstrated 
and the corresponding convergence rate is derived. 
The obtained results significantly extend the existing results on 
the convergence and convergence rate of stochastic search. 
They also cover several practically important classes of stochastic gradient algorithms
to which the existing results cannot be applied. 
We apply the new results to the asymptotic analysis of 
online algorithms for supervised and temporal-difference learning, 
principal component analysis,  
maximum likelihood estimation
and 
simulation-based optimization of (controlled) Markov chains. 
We also use them to study the asymptotic properties of 
recursive identification methods based on the prediction error and 
maximum likelihood principles. 

The paper is organized as follows. 
In Section \ref{section1}, 
stochastic gradient algorithms with additive noise are considered 
and the main results of the paper are presented. 
In Section \ref{section2}, the main results are applied to 
the asymptotic analysis of stochastic gradient algorithms with Markovian dynamics. 
Sections \ref{section4} -- \ref{section9} contain examples 
of the results presented in Sections \ref{section1} and \ref{section2}. 
In Sections \ref{section4} -- \ref{section8}, online algorithms for 
supervised and temporal-difference learning,  
principal component analysis and  
maximum likelihood estimation  
are studied. 
Recursive identification methods are considered in Sections \ref{section6} and \ref{section10}, 
while simulation-based optimization of (controlled) Markov chains is the subject of 
Section \ref{section9}. 
Section \ref{section0*} contains a detailed outline of the proof of the main results, 
while the proof itself is presented in Section \ref{section1*}. 
Sections \ref{section2*} -- \ref{section6*} contain the proof 
of the results presented in Sections \ref{section2} -- \ref{section6}. 

\section{Main Results} \label{section1}

In this section, the convergence and convergence rate of the following algorithm is analyzed: 
\begin{align} \label{1.1}
	\theta_{n+1}
	=
	\theta_{n}
	-
	\alpha_{n} (\nabla f(\theta_{n} ) + \xi_{n} ), 
	\;\;\; 
	n\geq 0. 
\end{align}
Here, $f:\mathbb{R}^{d_{\theta} } \rightarrow \mathbb{R}$ is a differentiable function, 
while 
$\{\alpha_{n} \}_{n\geq 0}$ is a sequence of positive real numbers. 
$\theta_{0}$ is an $\mathbb{R}^{d_{\theta } }$-valued random variable 
defined on a probability space 
$(\Omega, {\cal F}, P )$, 
while 
$\{\xi_{n} \}_{n\geq 0}$ is an $\mathbb{R}^{d_{\theta} }$-valued 
stochastic process defined on 
the same probability space. 
To allow more generality, we assume for each $n \geq 0$ that 
$\xi_{n}$ is a random function of $\theta_{0},\dots,\theta_{n}$. 
In the area of stochastic optimization, 
recursion (\ref{1.1}) is known as a stochastic gradient search
(or stochastic gradient algorithm),
while function $f(\cdot )$ is referred to as an objective function. 
For further details see \cite{pflug}, \cite{spall} and 
references given therein. 

Throughout the paper, unless otherwise stated, 
the following notation is used. 
The Euclidean norm is denoted by $\|\cdot \|$, 
while $d(\cdot, \cdot )$ stands for the distance induced by the Euclidean norm. 
$S$ is the set of stationary points of $f(\cdot )$, i.e., 
\begin{align*}
	S=\{\theta \in \mathbb{R}^{d_{\theta } }:\nabla f(\theta ) = 0 \}. 
\end{align*}	
Sequence $\{\gamma_{n} \}_{n\geq 0}$
is defined by 
$\gamma_{0}=0$
and 
\begin{align*}
	\gamma_{n} = \sum_{i=0}^{n-1} \alpha_{i}
\end{align*}
for $n\geq 1$. 	
For $t \in (0, \infty )$ and $n\geq 0$, 
$a(n,t)$ is the integer defined as 
\begin{align*}
	a(n,t)=\max\left\{k\geq n: \gamma_{k} - \gamma_{n} \leq t \right\}.  
\end{align*} 

Algorithm (\ref{1.1}) is analyzed under the following assumptions:  

\begin{assumption} \label{a1.1}
$\lim_{n\rightarrow \infty } \alpha_{n} = 0$ 
and 
$\sum_{n=0}^{\infty } \alpha_{n} = \infty$. 
\end{assumption}

\begin{assumption} \label{a1.2}
There exists a real number 
$r \in (1, \infty )$ such that 
\begin{align}\label{a1.2.1}
	\xi 
	=
	\limsup_{n\rightarrow \infty } 
	\max_{n\leq k < a(n,t) }
	\left\|
	\sum_{i=n}^{k} \alpha_{i} \gamma_{i}^{r} \xi_{i} 
	\right\|
	< 
	\infty 
\end{align}
w.p.1 on 
$\{\sup_{n\geq 0} \|\theta_{n} \| < \infty \}$ for all $t\in(0,\infty )$. 
\end{assumption} 

\begin{assumption} \label{a1.3}
For any compact set $Q \subset \mathbb{R}^{d_{\theta } }$ and any 
$a\in f(Q)$, there exist real numbers 
$\delta_{Q,a} \in (0,1]$, 
$\mu_{Q,a} \in (1,2]$, 
$M_{Q,a} \in[1,\infty )$
such that 
\begin{align} \label{a1.3.1}
	|f(\theta ) - a |
	\leq 
	M_{Q,a} \|\nabla f(\theta ) \|^{\mu_{Q,a} } 
\end{align}
for all $\theta \in Q$ satisfying 
$|f(\theta ) - a | \leq \delta_{Q,a}$. 
\end{assumption}

\begin{remark}\label{remark1.1}
As an immediate consequence of Assumption \ref{a1.3}, 
we have that for each $\theta \in R^{d_{\theta } }$, 
there exist real numbers 
$\delta_{\theta } \in (0,1]$, 
$\mu_{\theta } \in (1,2]$, 
$M_{\theta } \in [1,\infty )$ such that 
\begin{align}\label{r1.1.1}
	|f(\theta' ) - f(\theta ) |
	\leq 
	M_{\theta } \|\nabla f(\theta' ) \|^{\mu_{\theta } }
\end{align}
for all $\theta' \in \mathbb{R}^{d_{\theta } }$ 
satisfying $\|\theta' - \theta \| \leq \delta_{\theta }$. 
If $\theta \in S$, 
$\mu_{\theta }$ and $M_{\theta }$ can be selected as 
\begin{align*}
	\mu_{\theta } 
	=
	(1-\varepsilon)
	\liminf_{\theta' \rightarrow \theta } 
	\frac{\log|f(\theta' ) - f(\theta ) | }{\log\|\nabla f(\theta' ) \| }, 
	\;\;\;\;\; 
	M_{\theta } 
	=
	(1+\varepsilon)
	\limsup_{\theta' \rightarrow \theta } 
	\frac{|f(\theta' ) - f(\theta ) |}{\|\nabla f(\theta' ) \|^{\mu_{\theta } } }
\end{align*}
where $\varepsilon$ is a small positive constant
(since $\{\theta_{n} \}_{n\geq 0}$ converges to $S$, 
the values of $\mu_{\theta }$, $M_{\theta }$ 
for $\theta\not\in S$
are not 
relevant to the problems studied in the paper). 
Moreover, if 
$Q\subseteq \{\theta' \in R^{d_{\theta } }: 
\|\theta' - \theta \| \leq \delta_{\theta } \}$
and $a=f(\theta )$ for some $\theta \in R^{d_{\theta } }$, 
$\mu_{Q,a}$ and $M_{Q,a}$ can be assigned the values $\mu_{Q,a} = \mu_{\theta }$, 
$M_{Q,a} = M_{\theta }$. 
\end{remark}

\begin{remark}\label{remark1.2}
In order for Assumption \ref{a1.3} to be true,
it is sufficient that 
the assumption holds locally in an open vicinity of $S$, 
i.e., 
that there exists  
an open set $V \supset S$
with the following property: 
For any compact set $Q \subset V$ and any $a\in f(Q)$, 
there exist real numbers 
$\delta_{Q,a} \in (0,1]$, 
$\mu_{Q,a} \in (1,2]$, 
$M_{Q,a} \in [1,\infty )$ 
such that (\ref{a1.3.1}) holds 
for all $\theta \in Q$ satisfying 
$|f(\theta ) - a | \leq \delta_{Q,a}$
(see Appendix \ref{appendix1} for details). 
\end{remark}

Assumption \ref{a1.1} corresponds to the sequence 
$\{\alpha_{n} \}_{n\geq 0}$ and is widely used in the asymptotic 
analysis of stochastic gradient and stochastic approximation 
algorithms. 
It is fulfilled when 
$\alpha_{n}=n^{-a}$ for $n\geq 1$ and some constant $a\in (0,1]$. 

Assumption \ref{a1.2} is a noise condition. 
In this or a similar form, it is involved in most of 
the results on the convergence and convergence rate of stochastic gradient search
and stochastic approximation. 
It holds for algorithms with Markovian dynamics 
(see the next section). 
It is also satisfied when  
$\{\xi_{n} \}_{n\geq 0}$ is  
a martingale-difference or mixingale sequence. 

Assumption \ref{a1.3} is related to the stability 
of the gradient flow $d\theta/dt = - \nabla f(\theta )$ 
and the geometry of the set of stationary points $S$. 
In differential geometry, 
relations (\ref{a1.3.1}) and (\ref{r1.1.1}) 
are known as the Lojasiewicz gradient inequality
(see \cite{lojasiewicz1} and \cite{lojasiewicz2} for details). 
They hold if $f(\cdot )$ is analytic or subanalytic in 
an open vicinity of $S$
(see \cite{lojasiewicz2} for the proof; 
for the version of Lojasiewicz inequality appearing in 
Assumption \ref{a1.3} and (\ref{a1.3.1}), see \cite[Theorem \L I, page 775]{kurdyka}; 
for the definition and properties of analytic and 
subanalytic functions, see 
\cite{bierstone&milman}, \cite{krantz&parks}). 
In addition to this, 
Assumption \ref{a1.3} and relations (\ref{a1.3.1}), (\ref{r1.1.1}) include as a special case 
all stability conditions adopted by the existing results on the convergence rate 
of $\{\theta_{n} \}_{n\geq 0}$. 
These results are based on 
the following two conditions:  
(i) 
$f(\cdot )$ has a unique minimum $\theta_{*}$,  
and 
(ii) 
there exist a real number $\nu\in [0,\infty )$ and a  
positive definite matrix $A\in R^{d_{\theta } \times d_{\theta } }$ 
such that 
\begin{align}\label{1.101} 
	\nabla f(\theta )
	=
	A(\theta - \theta_{*} ) \|\theta - \theta_{*} \|^{\nu } 
	+
	o(\|\theta - \theta_{*} \|^{\nu+1} )
\end{align}
in an open vicinity of $\theta_{*}$
(see \cite{benveniste}, \cite{chen}, \cite{kushner&yin} 
and references cited therein).\footnote{
Due to (\ref{1.101}), $f(\cdot )$ is strongly convex in 
a vicinity of $\theta_{*}$. 
When $\nu=0$, (\ref{1.101}) is equivalent to 
the positive definiteness of $\nabla^{2} f(\theta_{*} )$. }
Using elementary calculus, 
it is straightforward to show that (i) and (ii) imply Assumption \ref{a1.3}.\footnote{
As a result of (i) and (ii), we also get
\begin{align*}
	0
	\leq 
	f(\theta ) - f(\theta_{*} ) 
	= &
	\int_{0}^{1} (\nabla f(\theta_{*} + t(\theta - \theta_{*} ) ) )^{T} 
	(\theta - \theta_{*} ) dt 
	\\
	\leq &
	(\theta - \theta_{*} )^{T} A (\theta - \theta_{*} ) 
	\|\theta - \theta_{*} \|^{\nu } 
	+ 
	o(\|\theta - \theta_{*} \|^{\nu+2} ) 
	\leq
	2\lambda_{max} \|\theta - \theta_{*} \|^{\nu+2} 
	\\
	\|\nabla f(\theta ) \|
	\geq & 
	\|A(\theta - \theta_{*} ) \| \|\theta - \theta_{*} \|^{\nu } 
	- 
	o(\|\theta - \theta_{*} \|^{\nu+1} )
	\geq 
	\frac{\lambda_{min} }{2} \|\theta - \theta_{*} \|^{\nu+1} 
\end{align*}
in a sufficiently small open vicinity of $\theta_{*}$, 
where $\lambda_{min}$ and $\lambda_{max}$ are the smallest and largest eigenvalue of 
$A$ (respectively). 
Consequently, 
\begin{align*}
	0
	\leq 
	f(\theta ) - f(\theta_{*} ) 
	\leq 
	M\|\nabla f(\theta ) \|^{\mu}
\end{align*}
in an open vicinity of $\theta_{*}$, 
where $\mu=(\nu+2)/(\nu+1)$ 
and $M=2^{\mu+1} \lambda_{max}/\lambda_{min}^{\mu}$. 
Hence, according to Remark \ref{remark1.2}, 
Assumption \ref{a1.3} is satisfied when (i) and (ii) hold. 
}

Although tightly connected to analyticity and subanalyticity 
(which are rather restrictive conditions), 
Assumption \ref{a1.3} covers  
many stochastic gradient 
algorithms routinely used in signal processing, automatic control, communications, statistics,
machine learning and operations research. 
In this paper, 
we show analyticity for 
the objective functions corresponding to online algorithms for 
supervised and temporal-difference learning, 
maximum likelihood estimation and  
principal component analysis (Sections \ref{section4} -- \ref{section8}). 
We also demonstrate the analyticity for 
simulation-based optimization of (controlled) Markov chains (Section \ref{section9}),
as well as for the recursive identification methods based on 
the prediction error and maximum likelihood principles (Sections \ref{section6} and \ref{section10}). 
It is also worth mentioning that 
the objective functions corresponding to many adaptive signal processing algorithms
are usually polynomial or rational, and hence, analytic, too
(see e.g., \cite{cichocki&amari} and references cited therein). 
In order to state the main results of this section, we need further notation. 
For $\theta \in \mathbb{R}^{d_{\theta } }$, 
$C_{\theta } \in [1,\infty )$ 
stands for an upper bound of 
$\|\nabla f(\cdot )\|$ 
on $\{\theta'\in \mathbb{R}^{d_{\theta } }: \|\theta' - \theta \| \leq \delta_{\theta }\}$ 
and for a Lipschitz constant of $\nabla f(\cdot )$ on the same set. 
Moreover, $p_{\theta }$, $q_{\theta }$ and $r_{\theta }$ are real numbers defines as 
\begin{align}\label{1.25}
	r_{\theta }
	= 
	\begin{cases} 
	1/(2 - \mu_{\theta } ), 
	&\text{if } \mu_{\theta } < 2
	\\
	\infty, 
	&\text{if } \mu_{\theta } = 2
	\end{cases}, 
	\;\;\;\;\; 
	p_{\theta } = \mu_{\theta } \min\{r,r_{\theta } \}, 
	\;\;\;\;\;
	q_{\theta } = \min\{r,r_{\theta } \} - 1  
\end{align}
($\delta_{\theta}$, $\mu_{\theta }$ are specified in Remark \ref{remark1.1}).

Our main results on the convergence and convergence rate of the recursion (\ref{1.1}) 
are contained in the next two theorems.  

\begin{theorem}[Convergence] \label{theorem1.1}
Let Assumptions \ref{a1.1} -- \ref{a1.3} hold. 
Then,  
$\hat{\theta } = \lim_{n\rightarrow \infty } \theta_{n}$
exists 
and satisfies 
$\nabla f(\hat{\theta } ) = 0$
w.p.1 on $\{\sup_{n\geq 0} \|\theta_{n} \| < \infty \}$. 
\end{theorem}

\begin{theorem}[Convergence Rate] \label{theorem1.2}
Let Assumptions \ref{a1.1} -- \ref{a1.3} hold. 
Then, 
there exists a random variable 
$\hat{K}$ (which is a deterministic function of $\hat{p}$, $C_{\hat{\theta } }$, $M_{\hat{\theta } }$)
such that $0< \hat{K} < \infty$ everywhere and 
such that the following is true: 
\begin{align}
	& \label{t1.1.1*}
	\limsup_{n\rightarrow \infty } 
	\gamma_{n}^{\hat{p} }
	\|\nabla f(\theta_{n} ) \|^{2} 
	\leq 
	\hat{K}
	\big(
	\varphi(\xi)
	\big)^{\hat{\mu} }, 
	\\
	& \label{t1.1.3*}
	\limsup_{n\rightarrow \infty } 
	\gamma_{n}^{\hat{p} }
	|f(\theta_{n} ) - f(\hat{\theta } ) |
	\leq 
	\hat{K} 
	\big(
	\varphi(\xi)
	\big)^{\hat{\mu} }, 
	\\
	& \label{t1.1.5*}
	\limsup_{n\rightarrow \infty } 
	\gamma_{n}^{\hat{q} } \|\theta_{n} - \hat{\theta } \| 
	\leq 
	\hat{K} \varphi(\xi ) 
\end{align}
w.p.1 on 
$\{\sup_{n\geq 0} \|\theta_{n} \| < \infty \}$, 
where 
$\hat{\mu} = \mu_{\hat{\theta } }$, 
$\hat{p} = p_{\hat{\theta } }$, $\hat{q} = q_{\hat{\theta } }$, 
$\hat{r} = r_{\hat{\theta } }$ and 
\begin{align*}
	\varphi(\xi)
	=
	\begin{cases}
	\xi, 
	&\text{if } r < \hat{r} \\
	1+\xi, 
	&\text{if } r = \hat{r} \\
	1, 
	&\text{if } r> \hat{r}  
	\end{cases}.
\end{align*}
\end{theorem}

A proof of Theorems \ref{theorem1.1} and \ref{theorem1.2} is provided in Section \ref{section1*}, 
while its outline is presented in Section \ref{section0*}. 
As an immediate consequence of the previous theorems, we get the next corollary: 

\begin{corollary} \label{corollary1.1}
Let Assumptions \ref{a1.1} -- \ref{a1.3} hold. 
Then, the following is true: 
\begin{compactenum}[(i)]
\item
$\|\nabla f(\theta_{n} ) \|^{2}	=
o\big(\gamma_{n}^{-\hat{p} } \big)$,
$|f(\theta_{n} ) - f(\hat{\theta } ) | =
o\big(\gamma_{n}^{-\hat{p} } \big)$ 
and
$\|\theta_{n} - \hat{\theta } \| =
o\big(\gamma_{n}^{-\hat{q} } \big)$
w.p.1 on 
$\{\sup_{n\geq 0} \|\theta_{n} \| < \infty \} \cap
\{\xi = 0, \hat{r} > r \}$.  
\item
$\|\nabla f(\theta_{n} ) \|^{2}	=
O\big(\gamma_{n}^{-\hat{p} } \big)$,
$|f(\theta_{n} ) - f(\hat{\theta } ) | =
O\big(\gamma_{n}^{-\hat{p} } \big)$ 
and
$\|\theta_{n} - \hat{\theta } \| =
O\big(\gamma_{n}^{-\hat{q} } \big)$
w.p.1 on 
$\{\sup_{n\geq 0} \|\theta_{n} \| < \infty \} \cap 
\{\xi = 0, \hat{r} > r \}^{c}$. 
\item
$\|\nabla f(\theta_{n} ) \|^{2} = o(\gamma_{n}^{-p} )$
and  
$|f(\theta_{n} ) - f(\hat{\theta } ) | = o(\gamma_{n}^{-p} )$
w.p.1 on 
$\{\sup_{n\geq 0} \|\theta_{n} \| < \infty \}$, 
where $p=\min\{1,r\}$. 
\end{compactenum}
\end{corollary}

\begin{remark}
The estimate of the convergence rate provided in Part (iii) of Corollary \ref{corollary1.1} 
is deterministic and independent of $\theta_{*}$. 
\end{remark}

In the literature on stochastic and deterministic optimization, 
the asymptotic behavior of gradient search is usually 
characterized by
the convergence of sequences 
$\{\nabla f(\theta_{n} ) \}_{n\geq 0}$,  
$\{f(\theta_{n} ) \}_{n\geq 0}$ and   
$\{\theta_{n} \}_{n\geq 0}$
(see e.g., \cite{bertsekas}, \cite{bertsekas&tsitsiklis2}, 
\cite{polyak&tsypkin}, \cite{polyak}  
and references cited therein). 
Similarly, the convergence rate can be described 
by the rates at which 
$\{\nabla f(\theta_{n} ) \}_{n\geq 0}$,  
$\{f(\theta_{n} ) \}_{n\geq 0}$  
and 
$\{\theta_{n} \}_{n\geq 0}$
tend to their limit points. 
In the case of algorithm (\ref{1.1}), 
this kind of information is provided by 
Theorems \ref{theorem1.1}, \ref{theorem1.2} and 
Corollary \ref{corollary1.1}. 
Theorem \ref{theorem1.1} claims that 
algorithm (\ref{1.1}) almost surely converges to a single-limit point  
(and does not exhibit limit cycles). 
Theorem \ref{theorem1.2} and Corollary \ref{corollary1.1} 
provide almost sure upper bounds on 
the convergence rate of 
$\{\nabla f(\theta_{n} ) \}_{n\geq 0}$, 
$\{f(\theta_{n} ) \}_{n\geq 0}$
and 
$\{\theta_{n} \}_{n\geq 0}$. 
The bounds are tightly connected to  
the convergence rate of 
gradient flow
$d\theta/dt = - \nabla f(\theta )$
and of noise average 
$\sum_{i=n}^{k} \alpha_{i} \xi_{i}$. 
Basically, Theorem \ref{theorem1.2} and Corollary \ref{corollary1.1} 
claim that 
the convergence rate of 
$\{\|\nabla f(\theta_{n} ) \|^{2} \}_{n\geq 0}$ 
and 
$\{f(\theta_{n} ) \}_{n\geq 0}$ is the slower of the rates 
$O(\gamma_{n}^{- \hat{r} \hat{\mu} } )$
(the rate of the gradient flow 
$d\theta/dt = - \nabla f(\theta )$ sampled at time-instants 
$\{\gamma_{n} \}_{n\geq 0}$) 
and 
$O(\gamma_{n}^{- r \hat{\mu} } )$
(the rate of the noise average 
$\max_{n\leq k < a(n,1) } \|\sum_{i=n}^{k} \alpha_{i} \xi_{i} \|^{\hat{\mu} }$). 
These estimates of the convergence rate of 
$\{f(\theta_{n} ) \}_{n\geq 0}$ and $\{\nabla f(\theta_{n} ) \}_{n\geq 0}$
seem to be rather tight. 
This is indicated by the arguments the proof of Theorem \ref{theorem1.2} 
is based on (see Section \ref{section0*} for an outline), 
as well as by the following two special cases:

\begin{vcase}{$\xi_{n} = 0 \text{ \rm for each } n\geq 0$}  
{Due to Assumption \ref{a1.3}, we have 
\begin{align*}
	\frac{d(f(\theta(t) ) - f(\hat{\theta} ) ) }{dt}
	=
	-\|\nabla f(\theta(t) ) \|^{2} 
	\leq 
	- 
	\left(
	\frac{f(\theta(t)) - f(\hat{\theta} ) }{\hat{M} } 
	\right)^{2/\hat{\mu} }
\end{align*}
for a solution $\theta(\cdot )$ of 
$d\theta/dt = - \nabla f(\theta )$
satisfying 
$\lim_{t\rightarrow \infty } \theta(t) = \hat{\theta}$ and 
$\|\theta(t) - \hat{\theta} \|\leq \delta_{\hat{\theta} }$ for all $t\in [0,\infty )$
($\delta_{\theta}$ is specified in Remark \ref{remark1.1}). 
Consequently, the Bellman-Gronwall inequality yields 
\begin{align*}
	f(\theta(t) ) - f(\hat{\theta} )
	= 
	O(t^{-\hat{\mu}/(2 - \hat{\mu} ) } ) 
	=
	O(t^{-\hat{\mu} \hat{r} }). 
\end{align*}
As $\{\theta_{n} \}_{n\geq 0}$ is asymptotically equivalent to 
$\theta(\cdot )$ sampled at time-instances  
$\{\gamma_{n} \}_{n\geq 0}$, 
we get 
$f(\theta_{n} ) - f(\hat{\theta} ) = O(\gamma_{n}^{-\hat{\mu} \hat{r} } )$. 
The same result is implied by Theorem \ref{theorem1.1} 
and Corollary \ref{corollary1.1}. 
}
\end{vcase}

\begin{vcase}{$f(\theta ) = \theta^{T} B \theta \text{ \rm for some positive definite
matrix } B$}
{In this case, 
recursion (\ref{1.1}) reduces to a linear stochastic approximation 
algorithm. 
For such an algorithm, the tightest bound on  
the convergence rate of 
$\{f(\theta_{n} ) \}_{n\geq 0}$ 
and $\{\|\nabla f(\theta_{n} ) \|^{2} \}_{n\geq 0}$
is 
$O(\gamma_{n}^{-2r} )$ if $\xi>0$ 
and 
$o(\gamma_{n}^{-2r} )$ if $\xi=0$
(see \cite{tadic1}). 
The same rate is predicted by Theorem \ref{theorem1.2} and 
Corollary \ref{corollary1.1}. 
}
\end{vcase}

The results of Theorems \ref{theorem1.1}, \ref{theorem1.2} and Corollary \ref{corollary1.1} 
are of a local nature. 
They hold only on the event where 
algorithm (\ref{1.1}) is stable
(i.e., where sequence 
$\{\theta_{n} \}_{n\geq 0}$ is bounded). 
Stating results on the convergence and convergence rate 
in such a form is quite sensible due to the following 
reasons. 
The stability of stochastic gradient search is based on 
well-understood arguments which are rather different from 
the arguments used here to analyze convergence and convergence rate. 
Moreover and more importantly, 
it is straightforward to get a global version of 
Theorems \ref{theorem1.1}, \ref{theorem1.2} and Corollary \ref{corollary1.1} 
by combining them with 
the methods for verifying or ensuring stability
(e.g., with the results of 
\cite{borkar&meyn} and \cite{chen}; 
see Appendix \ref{appendix2} for a global version of the results presented in this section). 

In the literature on deterministic optimization, 
a significant attention has recently been given to analytic and subanalytic functions, 
their properties and the methods for their minimization
(see e.g., \cite{absil1}, \cite{absil2}, \cite{bolte}). 
Crucially relying on Lojasiewicz gradient inequality 
and on the fact that $\{f(\theta_{n} ) \}_{n\geq 0}$ is decreasing, 
it has been demonstrated in \cite{absil1} 
that the deterministic gradient search converges to 
a single limit-point when the objective function is analytic. 
Theorems \ref{theorem1.1}, \ref{theorem1.2} and Corollary \ref{corollary1.1} 
provide a generalization of \cite{absil1} to stochastic gradient search. 
Since $\{f(\theta_{n} ) \}_{n\geq 0}$ is not decreasing in the case of stochastic gradient search
(due to noise $\{\xi_{n} \}_{n\geq 0}$), 
the arguments behind the results of \cite{absil1} cannot be applied to the asymptotic analysis
carried out here
(even the classical version of the Lojasiewicz inequality 
(\ref{r1.1.1}) cannot be used, but its generalization (\ref{a1.3.1})). 
Instead, a different and much more sophisticated techniques are needed. 
These techniques are based on a `singular' Lyapunov function 
(function $v(\cdot )$ introduced in (\ref{1.101*})). 

The single limit-point convergence and convergence rate of stochastic gradient search 
(and stochastic approximation) 
have been the subject of a number of papers and books
(see \cite{benveniste}, \cite{kushner&yin}, \cite{ljung2}, 
\cite{polyak}, \cite{spall} 
and references cited therein). 
Although the existing results 
provide a good insight into 
the asymptotic behavior and efficiency 
of stochastic gradient algorithms, 
they are based on restrictive conditions.  
The existing results on the single limit-point convergence of (\ref{1.1}) 
require (explicitly or implicitly) $f(\cdot )$ to have an isolated minimum $\theta_{*}$ 
such that $f(\cdot )$ is strongly convex in on open vicinity of $\theta_{*}$
and such that $\{\theta_{n} \}_{n\geq 0}$ almost surely 
visits the attraction domain of $\theta_{*}$
infinitely often. 
In addition to this, 
the existing results on the convergence rate of (\ref{1.1}) 
require $\nabla f(\cdot )$ to admit the representation (\ref{1.101}) 
in an open vicinity of $\theta_{*}$. 
These conditions are hard to verify for complex stochastic gradient algorithms. 
For such algorithms, it is very difficult even to show the existence of an isolated 
minimum, let alone to verify the representation (\ref{1.101}) or to check if 
$\{\theta_{n} \}_{n\geq 0}$ infinitely often enters the attraction domain of 
$\theta_{*}$. 
Furthermore, the conditions the existing results rely on are unlikely to hold 
for complex stochastic gradient algorithms as the corresponding objective functions 
are prone to non-isolated minima (each of which is a potential limit point of (\ref{1.1})). 
Several practically relevant examples of such a situation 
are provided Sections \ref{section4} -- \ref{section10}. 

Relying on the Lojasiewicz gradient inequality, 
Theorems \ref{theorem1.1}, \ref{theorem1.2} and Corollary \ref{corollary1.1} 
overcome the described difficulties.  
Both theorems and their corollary allow the objective function 
$f(\cdot )$
to be non-convex and have multiple and non-isolated minima. 
They also do not require $\nabla f(\cdot )$ to admit any particular representation
(notice that (\ref{1.101}) cannot hold if $\theta_{*}$ is a non-isolated minimum) 
and $\{\theta_{n} \}_{n\geq 0}$ to exhibit (a priori) any particular behavior 
(i.e., to visit infinitely often the attraction domain of an isolated minimum).  
Furthermore, they cover several practically important classes of stochastic gradient algorithms
which do not fit into the framework of the existing results
(for details, see Sections \ref{section4} -- \ref{section10}). 
To the best or our knowledge, 
these are the only results on the convergence and convergence rate 
of stochastic search which enjoy such features. 

\section{Stochastic Gradient Algorithms with Markovian Dynamics} \label{section2}

In order to illustrate the results of Section \ref{section1} and 
to set up a framework for the analysis carried out in Sections 
\ref{section4} -- \ref{section9}, 
we apply Theorems \ref{theorem1.1}, \ref{theorem1.2} and 
Corollary \ref{corollary1.1} 
to stochastic gradient algorithms with Markovian dynamics. 
These algorithms are defined by the following difference equation: 
\begin{align} \label{2.1}
	\theta_{n+1}
	=
	\theta_{n} 
	-
	\alpha_{n} F(\theta_{n}, Z_{n+1} ), 
	\;\;\; 
	n\geq 0. 
\end{align}
In this recursion, 
$F: \mathbb{R}^{d_{\theta} } \times \mathbb{R}^{d_{z} } \rightarrow \mathbb{R}^{d_{\theta } }$ is 
a measurable function, 
while $\{\alpha_{n} \}_{n\geq 0}$ is a sequence of positive real numbers. 
$\theta_{0}\in \mathbb{R}^{d_{\theta } }$ is an arbitrary vector, 
while 
$\{Z_{n} \}_{n\geq 0}$ is an $\mathbb{R}^{d_{z} }$-valued stochastic process 
defined on a probability space $(\Omega, {\cal F}, P)$. 
$\{Z_{n} \}_{n\geq 0}$ is a Markov process controlled by 
$\{\theta_{n} \}_{n\geq 0}$, i.e., 
there exists a family of transition kernels 
$\{\Pi_{\theta }(\cdot,\cdot ) \}_{\theta \in \mathbb{R}^{d_{\theta } } }$
(defined on $\mathbb{R}^{d_{z} }$)
such that 
\begin{align}\label{2.3} 
	P(Z_{n+1} \in B|\theta_{0},Z_{0},\dots,\theta_{n},Z_{n} ) 
	=
	\Pi_{\theta_{n} }(Z_{n}, B )
\end{align}
w.p.1 for $n\geq 0$ and any measurable set $B \subseteq \mathbb{R}^{d_{z} }$. 
In the context of stochastic gradient search, 
$F(\theta_{n},Z_{n+1} )$ is regarded to as an estimator of 
$\nabla f(\theta_{n} )$. 

The algorithm (\ref{2.1}) is analyzed under the following assumptions. 

\begin{assumption} \label{a2.1}
$\lim_{n\rightarrow \infty } \alpha_{n} = 0$, 
$\limsup_{n\rightarrow \infty } |\alpha_{n+1}^{-1} - \alpha_{n}^{-1} | < \infty$ 
and 
$\sum_{n=0}^{\infty } \alpha_{n} = \infty$. 
Moreover, there exists a real number $r\in (1,\infty )$
such that 
$\sum_{n=0}^{\infty } \alpha_{n}^{2} \gamma_{n}^{2r} < \infty$. 
\end{assumption}

\begin{assumption} \label{a2.2}
There exist a differentiable function 
$f: \mathbb{R}^{d_{\theta } } \rightarrow \mathbb{R}$ 
and 
a measurable function 
$\tilde{F}: \mathbb{R}^{d_{\theta} } \times \mathbb{R}^{d_{z} } \rightarrow \mathbb{R}^{d_{\theta } }$
such that 
$\nabla f(\cdot )$ is locally Lipschitz continuous and 
such that 
\begin{align*}
	F(\theta, z ) 
	-
	\nabla f(\theta )
	=
	\tilde{F}(\theta, z ) 
	-
	(\Pi\tilde{F} )(\theta, z )
\end{align*}
for each $\theta\in \mathbb{R}^{d_{\theta } }$, 
$z \in \mathbb{R}^{d_{z} }$, 
where 
$(\Pi\tilde{F} )(\theta, z ) 
= \int \tilde{F}(\theta, z' ) \Pi_{\theta }(z, dz' )$. 
\end{assumption}

\begin{assumption} \label{a2.3}
For any compact set $Q \subset \mathbb{R}^{d_{\theta } }$ and $s \in (0,1)$, 
there exists a measurable function 
$\varphi_{Q, s }: \mathbb{R}^{d_{z} } \rightarrow [1,\infty )$ such that 
\begin{align*}
	&
	\max\{
	\|F(\theta,z) \|, \|\tilde{F}(\theta,z) \|, \|(\Pi \tilde{F} )(\theta,z) \|
	\}
	\leq 
	\varphi_{Q, s }(z), 
	\\
	&
	\|(\Pi \tilde{F} )(\theta',z) - (\Pi \tilde{F} )(\theta'',z) \|
	\leq 
	\varphi_{Q, s }(z) \|\theta' - \theta'' \|^{s} 
\end{align*}
for all 
$\theta, \theta', \theta'' \in Q$, $z \in \mathbb{R}^{d_{z} }$.  
In addition to this, 
\begin{align*}
	\sup_{n\geq 0}
	E\left(
	\varphi_{Q, s }^{2}(Z_{n} ) 
	I_{\{\tau_{Q} \geq n \} }
	|\theta_{0}=\theta, Z_{0}=z 
	\right)
	< 
	\infty
\end{align*}
for all $\theta \in \mathbb{R}^{d_{\theta } }$, $z \in \mathbb{R}^{d_{z} }$, 
where 
$\tau_{Q} = 
\inf(\{n\geq 0: \theta_{n} \not\in Q \}\cup\{\infty \} )$. 
\end{assumption}

The main results on the convergence rate of recursion 
(\ref{2.1}) are contained in the next theorem. 

\begin{theorem} \label{theorem2.1}
Let Assumptions \ref{a2.1} -- \ref{a2.3} hold, 
and suppose that $f(\cdot )$ (introduced in Assumption \ref{a2.2}) 
satisfies Assumption \ref{a1.3}. 
Then, the following is true: 
\begin{compactenum}[(i)]
\item
$\hat{\theta} = \lim_{n\rightarrow \infty } \theta_{n}$
exists and satisfies 
$\nabla f(\hat{\theta} ) = 0$
w.p.1 on $\{\sup_{n\geq 0} \|\theta_{n} \| < \infty \}$. 
\item
$\|\nabla f(\theta_{n} ) \|^{2} = 
o\big(\gamma_{n}^{-\hat{p} } \big)$,
$|f(\theta_{n} ) - f(\hat{\theta} ) |	=
o\big(\gamma_{n}^{-\hat{p} } \big)$
and
$\|\theta_{n} - \hat{\theta} \| =
o\big(\gamma_{n}^{-\hat{q} } \big)$ 
w.p.1 on 
$\{\sup_{n\geq 0} \|\theta_{n} \| < \infty \} 
\cap \{\hat{r} > r\}$. 
\item 
$\|\nabla f(\theta_{n} ) \|^{2}	=
O\big(\gamma_{n}^{-\hat{p} } \big)$,
$|f(\theta_{n} ) - f(\hat{\theta} ) | =
O\big(\gamma_{n}^{-\hat{p} } \big)$ 
and 
$\|\theta_{n} - \hat{\theta} \| =
O\big(\gamma_{n}^{-\hat{q} } \big)$
w.p.1 on 
$\{\sup_{n\geq 0} \|\theta_{n} \| < \infty \} 
\cap \{\hat{r} \leq r\}$. 
\item
$\|\nabla f(\theta_{n} ) \|^{2} = o(\gamma_{n}^{-p} )$ 
and 
$|f(\theta_{n} ) - f(\hat{\theta} ) | = o(\gamma_{n}^{-p} )$ 
w.p.1 on $\{\sup_{n\geq 0} \|\theta_{n} \| < \infty \}$. 
\end{compactenum}
\end{theorem}

A proof is provided in Section \ref{section2*}. 
$p$, $\hat{p}$, $\hat{q}$ and $\hat{r}$
are defined in Theorem \ref{theorem1.2} and Corollary \ref{corollary1.1}. 

Assumption \ref{a2.1} is related to the sequence  
$\{\alpha_{n} \}_{n\geq 0}$. 
It holds if 
$\alpha_{n} = 1/n^{a}$ for $n\geq 1$ and some constant $a \in (3/4,1]$ 
(in that case, 
$\gamma_{n} = O(n^{1-a} )$
for $n\rightarrow \infty$, 
while $r$ can be any number satisfying
$0 < r < (a-1/2)/(1-a)$). 
On the other side,  
Assumptions \ref{a2.2} and \ref{a2.3} 
correspond to   
the stochastic process 
$\{Z_{n} \}_{n\geq 0}$
and are standard for the asymptotic analysis of 
stochastic approximation algorithms with Markovian dynamics. 
Assumptions \ref{a2.2} and \ref{a2.3} have been introduced by 
Metivier and Priouret in \cite{metivier&priouret1} 
(see also \cite[Part II]{benveniste}), 
and later generalized by Kushner and Yin 
(see \cite{kushner&yin} and references cited therein). 
However, 
neither the results of Metivier and Priouret, 
nor the results of Kushner and Yin 
provide any information on the single limit-point convergence and convergence rate 
of stochastic gradient search in 
the case of multiple and non-isolated 
minima.  

Regarding Theorem \ref{theorem2.1}, 
the following note is also in order.  
As already mentioned in the beginning of the section, 
the purpose of the theorem 
is illustrating the results of Section \ref{section1}  
and providing a framework for studying the examples 
presented in the next few sections. 
Since these examples perfectly fit into the framework 
developed by Metivier and Priouret, 
more general assumptions and settings of 
\cite{kushner&yin} are not considered here 
in order to keep the exposition as concise as possible. 

\section{Example 1: Supervised Learning} \label{section4} 

In this section, 
online algorithms for supervised learning in  
feedforward neural networks are analyzed using 
Theorems \ref{theorem1.1}, \ref{theorem1.2} and \ref{theorem2.1}. 
To avoid unnecessary technical details and complicated notation, 
only two-layer networks are considered here. 
However, 
the obtained results can be extended to the networks with any number of layers. 

The input-output function of a two-layer perceptron can be defined as 
\begin{align*}
	G_{\theta }(x) 
	=
	\sum_{i=1}^{M} a_{i} 
	\psi\left(
	\sum_{j=1}^{N}
	b_{i,j} x_{j} 
	\right).  
\end{align*}
Here, 
$\psi:\mathbb{R} \rightarrow \mathbb{R}$ is a differentiable function, 
while  
$M,N\geq 1$ are integers. 
$a_{1}, \dots, a_{M}$, $b_{1,1}, \dots, b_{M,N}$ and 
$x_{1}, \dots, x_{N}$ are real numbers, 
while 
$\theta = 
[a_{1} \cdots a_{M} \; b_{1,1} \cdots b_{M,N} ]^{T}$, 
$x = [x_{1} \cdots x_{N} ]^{T}$
and 
$d_{\theta } = M(N+1)$. 
In this context, 
$\psi(\cdot )$ represents the network activation function, 
while $x$ and $G_{\theta }(x)$ are the network input and output (respectively). 
$\theta$ is the vector of the network parameters to be tuned 
through the process of supervised learning. 

Let ${\cal X} \subseteq \mathbb{R}^{N}$, ${\cal Y} \subseteq \mathbb{R}$
be measurable sets, 
while $\{(X_{n}, Y_{n} ) \}_{n\geq 0}$ are ${\cal X}\times {\cal Y}$-valued 
i.i.d. random variables defined on a probability space $(\Omega, {\cal F}, P )$. 
Function $f(\cdot )$ is defined as 
\begin{align*}
	f(\theta ) 
	= 
	\frac{1}{2}
	E(Y_{0}-G_{\theta }(X_{0}) )^{2} 
\end{align*}
for $\theta \in \mathbb{R}^{d_{\theta } }$. 
Then, the mean-square error based supervised learning
in feedforward neural networks can be described as the minimization 
of $f(\cdot )$
in a situation when only a realization of $\{(X_{n}, Y_{n} ) \}_{n\geq 0}$
is available. 
In this context, $\{(X_{n}, Y_{n} ) \}_{n\geq 0}$ is referred to as a training sequence. 
For more details on neural networks and supervised learning, 
see e.g., \cite{hastie&tibshirani&friedman}, 
\cite{haykin} and references cited therein. 

Function $f(\cdot )$ is usually minimized by 
the following stochastic gradient algorithm: 
\begin{align} \label{4.1} 
	\theta_{n+1} 
	=
	\theta_{n} 
	+
	\alpha_{n} 
	(Y_{n} - G_{\theta_{n} }(X_{n} ) ) H_{\theta_{n} }(X_{n} ), 
	\;\;\; 
	n\geq 0. 
\end{align}
In this recursion, 
$\{\alpha_{n} \}_{n\geq 0}$
is a sequence of positive real numbers. 
$\theta_{0}\in \mathbb{R}^{d_{\theta } }$ is an arbitrary vector, 
while 
$H_{\theta }(\cdot ) = \nabla_{\theta } G_{\theta }(\cdot )$. 

\begin{remark}\label{r4.1}
Even for relatively small $M$ and $N$, 
function $f(\cdot )$ is prone to multiple and non-isolated minima. 
To illustrate this, we consider the simplest possible case 
when $\psi(\cdot )$ is identity mapping
(i.e., when $\psi(t)=t$ for each $t\in \mathbb{R}$). 
In this situation, the set of global minima of $f(\cdot )$ admits the representation
\begin{align*}
	S_{*}
	=
	\{\theta = [a^{T} \; \text{\rm vec}^{T}(B) ]^{T}: 
	a\in \mathbb{R}^{M}, B\in \mathbb{R}^{M\times N}, B^{T}a=\phi_{*} \}, 
\end{align*}
where 
$\phi_{*} = \arg\min_{\phi\in \mathbb{R}^{N} } 
\int (y-\phi^{T} x)^{2} \pi(dx,dy)$, 
while $\text{\rm vec}(B)$ is the vector whose components are the entries of $B$
(i.e., $\text{\rm vec}(B) = [b_{1,1} \cdots b_{M,N}]^{T}$,  
where $b_{i,j}$ denotes the $(i,j)$-entry of $B$). 
Obviously, $S_{*}$ has uncountably many elements each of which is non-isolated. 
This clearly indicates that 
function $f(\cdot )$ is very likely to have multiple and non-isolated minima in a general case 
when $\psi(\cdot )$ is nonlinear. 
\end{remark}

The asymptotic behavior of algorithm (\ref{4.1}) is analyzed under the 
following assumptions: 

\begin{assumption} \label{a4.1}
$\psi(\cdot )$ is real-analytic. 
Moreover, 
$\psi(\cdot )$ has 
a (complex-valued) continuation 
$\hat{\psi}(\cdot )$ with the following properties: 
\begin{compactenum}[(i)]
\item
$\hat{\psi}(z)$ maps $z\in \mathbb{C}$ to $\mathbb{C}$
($\mathbb{C}$ denotes the set of complex numbers). 
\item
$\hat{\psi}(x) = \psi(x)$
for all $x\in \mathbb{R}$. 
\item
There exists a real number
$\varepsilon \in (0,1)$ such that 
$\hat{\psi}(\cdot )$ 
is analytic
on $V_{\varepsilon }(\mathbb{R} ) = 
\{z\in \mathbb{C}: d(z, \mathbb{R} ) \leq \varepsilon \}$. 
\end{compactenum}
\end{assumption} 

\begin{assumption} \label{a4.2}
${\cal X}$ and ${\cal Y}$ are compact. 
\end{assumption}

Assumption \ref{a4.1} is related to the network activation function. 
It holds when 
$\psi(\cdot )$ is a logistic function\footnote
{Complex-valued logistic function can be defined as 
$\hat{\psi}(z) = (1 + \exp(-z) )^{-1}$ for $z\in \mathbb{C}$. 
Since 
\begin{align*}
	|1 + \exp(-z) |^{2}
	=
	1 
	+ 
	\exp(-2\text{Re}(z) )
	+
	2 \exp(-\text{Re}(z) ) \cos(\text{Im}(z) )
	\geq 
	1 + \exp(-2\text{Re}(z) )
\end{align*}
when $|\text{Im}(z) | \leq \pi/2$, 
$\hat{\psi}(\cdot )$
is analytical on 
$V_{\pi/2}(\mathbb{R} ) \{z\in \mathbb{C}: d(z,\mathbb{R}) \leq \pi/2 \}$. 
} 
or a standard Gaussian density\footnote
{Complex-valued standard Gaussian density can be defined by  
$\hat{\psi}(z) = (2\pi )^{-1/2} \exp(-z^{2}/2 )$ 
for $z\in \mathbb{C}$. 
It is analytical on entire $\mathbb{C}$. 
}, 
which are the most common activation functions for 
feedforward neural networks. 
Assumption \ref{a4.2} corresponds to the 
training sequence 
$\{(X_{n}, Y_{n} ) \}_{n\geq 0}$ 
and practically always holds in real-world applications
(as only bounded signals can be generated by real-world systems). 

Our main results on the properties of 
objective function $f(\cdot )$
and algorithm (\ref{4.1})
are contained in the next two theorems. 

\begin{theorem} \label{theorem4.1}
Let Assumptions \ref{a4.1} and \ref{a4.2} hold.
Then, $f(\cdot )$ is analytic on entire $\mathbb{R}^{d_{\theta } }$. 
\end{theorem}

\begin{theorem} \label{theorem4.2}
Let Assumptions \ref{a2.1}, \ref{a4.1} and \ref{a4.2} hold. 
Then, all conclusions of Theorem \ref{theorem2.1} are true for 
$\{\theta_{n} \}_{n\geq 0}$ defined in this section. 
\end{theorem}

A proof of Theorem \ref{theorem4.1} and \ref{theorem4.2} is provided in Section \ref{section4*}. 

The asymptotic properties of online algorithms for supervised learning have been 
studied in a large number of papers and books 
(see \cite{bertsekas&tsitsiklis1}, \cite{hastie&tibshirani&friedman}, \cite{haykin} 
and references cited therein). 
To the best of out knowledge, the available literature does not provide 
any information on the single limit-point convergence and convergence rate which can be 
verified for feedforward neural networks with nonlinear 
activation functions. 
The reason comes out of the fact that 
the existing asymptotic results for stochastic gradient search 
hold under very restrictive conditions which fail to hold for such networks
(as explained in Remark \ref{r4.1} and Section \ref{section1}). 

\section{Example 2: Principal Component Analysis} \label{section5} 

To illustrate the results of Sections \ref{section1} and \ref{section2}, 
we apply them to the asymptotic analysis of online algorithms for principal component analysis. 

To state the problem of principal component analysis and to define the corresponding online algorithms, 
we use the following notation. 
$M$ and $N$ are integers satisfying $N\geq M > 1$. 
$\{X_{n} \}_{n\geq 0}$ is an $\mathbb{R}^{N}$-valued i.i.d. stochastic process
defined on a probability space $(\Omega, {\cal F}, P )$, 
while $R = E(X_{0} X_{0}^{T} )$. 
Then, the principal component analysis can be stated as the computation of 
the $M$ leading eigenvectors of $R$ 
(i.e., the eigenvectors corresponding to the $M$ largest eigenvalues)
given a realization of $\{X_{n} \}_{n\geq 0}$. 
Online algorithms for principal component analysis are based on the minimization of 
\begin{align*}
	f(\Theta )
	=
	E\|X_{0} - \Theta\Theta^{T} X_{0} \|^{2}
\end{align*}
with respect to $\Theta\in\mathbb{R}^{N\times M}$
(see e.g., \cite{delmas1}, \cite{delmas2}, \cite{yang} and references cited therein). 
Since 
\begin{align*}
	\nabla f(\Theta )
	=
	-
	\left(
	R(2{\boldsymbol I} - \Theta\Theta^{T} ) 
	-
	\Theta\Theta^{T} R
	\right)\Theta
\end{align*}
(
here, ${\boldsymbol I}$ denotes $N\times N$ unit matrix, while 
$\nabla f(\Theta )$ is the $N\times M$ matrix defined by 
$[\nabla f(\Theta ) ]_{i,j} = \partial f/\partial [\Theta]_{i,j}$ 
for $1\leq i \leq N$, $1\leq j\leq M$). 
the minimization can be performed by the following stochastic gradient search: 
\begin{align}\label{5.1}
	\Theta_{n+1} 
	=
	\Theta_{n} 
	+
	\alpha_{n} 
	\left(
	X_{n} X_{n}^{T} (2{\boldsymbol I} - \Theta_{n}\Theta_{n}^{T} ) 
	-
	\Theta_{n}\Theta_{n}^{T} X_{n} X_{n}^{T} 
	\right)\Theta_{n}, 
	\;\;\; n\geq 0. 
\end{align}
In this recursion, $\{\alpha_{n} \}_{n\geq 0}$ is a sequence of positive reals, 
while $\Theta_{0}\in \mathbb{R}^{N\times M}$ is an arbitrary matrix. 
Since 
$\lim_{n\rightarrow\infty } \Theta_{n}^{T}\Theta_{n} = {\boldsymbol I}$
(see \cite{delmas1}, \cite{yang}), 
algorithm (\ref{5.1}) can be simplified to 
\begin{align}\label{5.3}
	\Theta_{n+1} 
	=
	\Theta_{n} 
	+
	\alpha_{n} 
	\left({\boldsymbol I} - \Theta_{n}\Theta_{n}^{T} \right) X_{n} X_{n}^{T} \Theta_{n}, 
	\;\;\; n\geq 0. 
\end{align}
In the literature on principal component analysis, 
recursions (\ref{5.1}) and (\ref{5.3}) are known as the Yang and Oja algorithm (respectively). 
As opposed to (\ref{5.1}), algorithm (\ref{5.3}) is not a stochastic gradient search. 
Despite this, (\ref{5.3}) can still be analyzed using the results of Sections \ref{section1} and \ref{section2}. 
Since such this analysis involves some technical difficulties
(such as bringing (\ref{5.3}) to a form similar to (\ref{5.1}) and analyzing the associated quantities), 
the focus of this section is on recursion (\ref{5.1}). 

\begin{remark}
Let $\lambda_{1},\dots,\lambda_{N}$ be 
eigenvalues of $R$
satisfying $\lambda_{1}\geq\cdots\geq\lambda_{N}$, 
while 
$e_{i}\in \mathbb{R}^{N}$ is an eigenvector corresponding to $\lambda_{i}$. 
Moreover, let 
\begin{align*}
	&
	S_{*}
	=
	\left\{
	\Theta=[e_{1} \cdots e_{M} ]Q: 
	Q\in \mathbb{R}^{M\times M}
	\right\}, 
	\\
	&
	S
	=
	\left\{
	\Theta=[e_{i_{1} } \cdots e_{i_{M} } ]Q: 
	Q\in \mathbb{R}^{M\times M}, 
	1\leq i_{1} < \cdots < i_{M} \leq M
	\right\}.
\end{align*}
Then, if $\lambda_{M} > \lambda_{M+1}$, 
$S_{*}$ and $S$ are the sets of global minima and stationary points 
of $f(\cdot )$, respectively 
(see \cite{delmas1}, \cite{yang}). 
Obviously, both $S_{*}$ and $S$ have uncountably many elements each of which is non-isolated. 
\end{remark}

Algorithm (\ref{5.1}) is analyzed under the following assumption. 

\begin{assumption}\label{a5.1}
$E\|X_{0}\|^{4} < \infty$. 
\end{assumption}

The main results on the properties of $f(\cdot )$ and algorithm (\ref{5.1}) 
are provided in the next two theorems. 

\begin{theorem} \label{theorem5.1}
Let Assumption \ref{a5.1} hold.
Then, $f(\cdot )$ is analytic on entire $\mathbb{R}^{N\times M}$. 
\end{theorem}

\begin{theorem} \label{theorem5.2}
Let Assumptions \ref{a2.1} and \ref{a5.1} hold. 
Then, all conclusion of Theorem \ref{theorem2.1} are true for 
$\{\Theta_{n} \}_{n\geq 0}$ 
(i.e., for $\{\theta_{n} \}_{n\geq 0}$ defined by 
$\theta_{n} = [\vartheta_{n}^{1,1} \cdots \vartheta_{n}^{N,M} ]^{T}$, 
where $\vartheta_{n}^{i,j}$ is the $(i,j)$-entry of $\Theta_{n}$). 
\end{theorem}


\begin{remark}
Theorem \ref{theorem5.1} is an immediate consequence of the fact that 
$f(\Theta )$ is polynomial in $\Theta$. 
On the other hand, Assumptions \ref{a2.2} and \ref{a2.3} hold for algorithm (\ref{5.1}), 
since $\{X_{n} \}_{n\geq 0}$ can be interpreted as a controlled Markov chain 
whose transition kernel $\Pi_{\Theta}(x,\cdot )$ does not depend on 
$(\Theta,x)$. 
As a result of this, Theorem \ref{theorem5.2} directly follows from Theorem \ref{theorem2.1}. 
\end{remark}

The asymptotic behavior of online algorithms for principal component analysis has been 
studied in a number of papers 
(see \cite[Section 10.5]{borkar} and \cite{delmas2} for a recent review). 
Although the existing results provide a good insight into the properties of these algorithms, 
they are mainly concerned with the behavior of 
$\{\Theta_{n} \Theta_{n}^{T} \}_{n\geq 0}$ 
and do not provide any information about the single limit-point convergence and convergence rate of 
$\{\Theta_{n} \}_{n\geq 0}$
(for the difficulties associated with the asymptotic analysis of 
$\{\Theta_{n} \}_{n\geq 0}$, see \cite[Section III]{delmas1}). 
The aim of Theorems \ref{theorem5.1} and \ref{theorem5.2} is to fill this gap in 
the literature on principal component analysis. 

\section{Example 3: Maximum Likelihood Estimation} \label{section7} 

In this section, 
Theorems \ref{theorem1.1}, \ref{theorem1.2} and \ref{theorem2.1} 
are used to analyze the asymptotic behavior of online algorithms 
for maximum likelihood estimation in i.i.d. data. 

To state the problem of maximum likelihood estimation and to define the corresponding 
online algorithm, 
we use the following notation. 
$d_{\theta},N\geq 1$ are integers.  
$\Theta\subseteq\mathbb{R}^{d_{\theta } }$ is an open set,  
while 
${\cal X}\subseteq\mathbb{R}^{N}$ is a measurable sets. 
$\lambda(\cdot )$ is a measure on $\mathbb{R}^{N}$. 
For each $\theta\in\Theta$, 
$p_{\theta }(\cdot )$ is a (parameterized)  
probability density with respect to $\lambda(\cdot )$
(i.e., $p_{\theta}(x)$ is a measurable function mapping 
$(\theta,x) \in \Theta\times\mathbb{R}^{N}$ 
to $[0,\infty )$ and satisfying $\int p_{\theta }(x) \lambda(dx) = 1$ 
for all $\theta\in\Theta$). 
$\{X_{n} \}_{n\geq 0}$ are ${\cal X}$-valued i.i.d. random variables 
which are 
defined on a probability space $(\Omega, {\cal F}, P )$ and 
admit a probability density $p(\cdot )$ with respect to $\lambda(\cdot )$
($p(\cdot )$ is not necessarily an element of 
$\{p_{\theta }(\cdot ) \}_{\theta\in\Theta }$). 

The problem of parameter estimation for i.i.d. data can be stated as follows: 
Given a realization of $\{X_{n} \}_{n\geq 0}$, 
estimate the values of $\theta$ for which $p_{\theta }(\cdot )$ 
provides the best approximation to $p(\cdot )$. 
If the estimation is based on the maximum likelihood principle, 
the estimation reduces to 
the minimization of the negative log-likelihood 
\begin{align*}
	f(\theta )
	=
	-
	\int\log\left(p_{\theta }(x) \right) p(x) \lambda(dx) 
\end{align*}
with respect to $\theta\in\Theta$. 
In online settings, 
$f(\cdot )$ is usually minimized by stochastic gradient (or stochastic Newton) algorithm. 
Such an algorithm is defined by the following recursion: 
\begin{align}\label{7.1}
	\theta_{n+1} 
	=
	\theta_{n} 
	-
	\alpha_{n} F(\theta_{n}, X_{n} ), 
	\;\;\; n\geq 0. 
\end{align}
Here, $\{\alpha_{n} \}_{n\geq 0}$ is a sequence of positive real numbers. 
$\theta_{0}\in\Theta$ is an arbitrary vector, 
while $F(\theta,x) = - \nabla_{\theta } p_{\theta }(x)/ p_{\theta }(x)$ 
for $\theta\in\Theta$, $x\in {\cal X}$. 
In the literature on statistical inference and system identification, 
algorithm (\ref{7.1}) is commonly referred to as the recursive maximum 
likelihood method. 

\begin{remark}\label{r7.1} 
In the case of multivariate parameters, 
negative log-likelihood $f(\cdot )$ is prone to multiple and non-isolated minima. 
This inevitably happens whenever 
$\{p_{\theta }(\cdot ) \}_{\theta\in\Theta }$ is over-parameterized for 
$p(\cdot )$. 
To illustrate this, we consider the situation when $p(\cdot )$ and 
$p_{\theta }(\cdot )$ are finite mixtures of probability densities 
from the same parametric family. 
More specifically, we assume 
\begin{align*}
	&
	p(x)
	=
	\sum_{i=1}^{M} w_{i}^{*} q_{\phi_{i}^{*} }(x), 
	\;\;\;\;\;  
	p_{\theta }(x)
	=
	\sum_{i=1}^{M+1} w_{i} q_{\phi_{i} }(x). 
\end{align*}
Here, $\{q_{\phi }(\cdot ) \}_{\phi\in\Phi }$ are (parameterized)    
probability densities with respect to $\lambda(\cdot )$, 
while $\Phi\subseteq\mathbb{R}^{L}$ is an open set and $L,M\geq 1$ are integers. 
$w_{1}^{*},\dots,w_{M}^{*},w_{1},\dots,w_{M+1} \in (0,1)$ are
real numbers satisfying 
$\sum_{i=1}^{M} w_{i}^{*} = \sum_{i=1}^{M+1} w_{i} = 1$, 
while $\phi_{1}^{*},\dots,\phi_{M}^{*}$, $\phi_{1},\dots,\phi_{M+1}$ $\in \Phi$ 
and $\theta = [w_{1} \cdots w_{M+1} \; \phi_{1}^{T} \cdots \phi_{M+1}^{T} ]^{T}$. 
On the other side, let 
\begin{align*}
	S_{*}^{i}
	=
	\Big\{&
	\theta
	=
	[w_{1} \cdots w_{M+1} \; \phi_{1}^{T} \cdots \phi_{M+1}^{T} ]^{T}
	\in 
	(0,1)^{M+1}\times \Phi^{M+1} 
	\\
	&
	:
	\phi_{i}=\phi_{M+1}=\phi_{i}^{*}, w_{i}+w_{M+1}=w_{i}^{*}, 
	w_{j}=w_{j}^{*}, \phi_{j}=\phi_{j}^{*}, 
	1\leq j \leq M, j\neq i
	\Big\}
\end{align*}
for $1\leq i \leq M$, 
while $S_{*} = \bigcup_{i=1}^{M} S_{*}^{i}$. 
Then, it is straightforward to show that each element of $S_{*}$ 
is a non-isolated global minimum of $f(\cdot )$. 
This strongly suggests that in a general case, when $p(\cdot )$ is not included in 
$\{p_{\theta }(\cdot ) \}_{\theta\in\Theta }$, 
negative log-likelihood $f(\cdot )$ is very likely to be multi-modal and 
has non-isolated minima. 
\end{remark}

Algorithm (\ref{7.1}) is analyzed under the following assumptions. 

\begin{assumption}\label{a7.1}
${\cal X}$ is compact and $\inf_{x\in {\cal X} } p(x) > 0$. 
\end{assumption}

\begin{assumption}\label{a7.2}
$p_{\theta }(x) > 0$ for all $\theta\in\Theta$, $x\in {\cal X}$. 
\end{assumption}

\begin{assumption}\label{a7.3}
For each $x\in {\cal X}$, 
$p_{\theta }(x)$ is real-analytic in $\theta$ on entire $\Theta$. 
Moreover, 
$p_{\theta }(x)$ has 
a (complex-valued) continuation 
$\hat{p}_{\eta}(x)$ with the following properties: 
\begin{compactenum}[(i)]
\item
$\hat{p}_{\eta}(x)$ maps $(\eta,x)\in \mathbb{C}^{d_{\theta } } \times {\cal X}$ 
to $\mathbb{C}$. 
\item
$\hat{p}_{\theta}(x) = p_{\theta}(x)$
for all $\theta\in\Theta$, $x\in {\cal X}$. 
\item
For any $\theta\in\Theta$, 
there exists a real number
$\delta_{\theta } \in (0,1)$ such that 
$\hat{p}_{\eta}(x)$ is 
analytic in $\eta$ 
and 
continuous in $(\eta,x)$
for any 
$\eta\in\mathbb{C}^{d_{\theta } }$, $x\in {\cal X}$ 
satisfying $\|\eta-\theta\|\leq\delta_{\theta }$. 
\end{compactenum}
\end{assumption}

Assumption \ref{a7.1} corresponds to the statistical properties of data 
$\{X_{n} \}_{n\geq 0}$ and covers many practically important applications and situations. 
Assumptions \ref{a7.2} and \ref{a7.3} are related to the 
parameterized family 
$\{p_{\theta }(\cdot ) \}_{\theta\in\Theta }$. 
They hold for many practically relevant statistical models. 
E.g., Assumptions \ref{a7.2} and \ref{a7.3} are satisfied 
if $p_{\theta }(\cdot )$ is a mixture of exponential, gamma, logistic, normal, log-normal, 
Pareto, uniform and Weinbull distributions, 
and if the mixture is parameterized by the mixture weights 
and by the `natural parameters' of the ingredient distributions. 

Let $\Lambda$ be the event defined by 
\begin{align}\label{7.3}
	\Lambda
	=
	\left\{
	\sup_{n\geq 0} \|\theta_{n} \|<\infty, 
	\inf_{n\geq 0} d(\theta_{n}, \Theta^{c} ) > 0 
	\right\}. 
\end{align}
With this notation, the main results on the properties of $f(\cdot )$
and the asymptotic behavior of (\ref{7.1}) read as follows: 

\begin{theorem}\label{theorem7.1}
Let Assumptions \ref{a7.1} -- \ref{a7.3} hold.
Then, $f(\cdot )$ is analytic on entire $\Theta$.  
\end{theorem}

\begin{theorem}\label{theorem7.2}
Let Assumptions \ref{a2.1} and \ref{a7.1} -- \ref{a7.3} hold. 
Then, the following is true: 
\begin{compactenum}[(i)]
\item
$\hat{\theta} = \lim_{n\rightarrow \infty } \theta_{n}$
exists and satisfies 
$\nabla f(\hat{\theta } ) = 0$
w.p.1 on $\Lambda$. 
\item
$\|\nabla f(\theta_{n} ) \|^{2}	=
o\big(\gamma_{n}^{-\hat{p} } \big)$,
$|f(\theta_{n} ) - f(\hat{\theta } ) | =
o\big(\gamma_{n}^{-\hat{p} } \big)$
and 
$\|\theta_{n} - \hat{\theta} \| =
o\big(\gamma_{n}^{-\hat{q} } \big)$
w.p.1 
on 
$\Lambda 
\cap \{\hat{r} > r \}$. 
\item
$	\|\nabla f(\theta_{n} ) \|^{2} =
O\big(\gamma_{n}^{-\hat{p} } \big)$,
$|f(\theta_{n} ) f(\hat{\theta} ) | =
O\big(\gamma_{n}^{-\hat{p} } \big)$
and  
$\|\theta_{n} - \hat{\theta} \| =
O\big(\gamma_{n}^{-\hat{q} } \big)$
w.p.1 on 
$\Lambda 
\cap \{\hat{r} \leq r \}$. 
\item
$\|\nabla f(\theta_{n} ) \|^{2} = o(\gamma_{n}^{-p} )$
and 
$|f(\theta_{n} ) - f(\hat{\theta} ) | = o(\gamma_{n}^{-p} )$
w.p.1 on 
$\Lambda$. 
\end{compactenum}
\end{theorem}

A proof of Theorems \ref{theorem7.1} and \ref{theorem7.2} is provided in Section \ref{section7*}. 
$p$, $\hat{p}$, $\hat{q}$ and $\hat{r}$
are defined in Theorem \ref{theorem1.2} and Corollary \ref{corollary1.1}. 

\begin{remark}
Algorithm (\ref{7.1}) usually involves a projection (or truncation) device which ensures
that estimates $\{\theta_{n} \}_{n\geq 0}$ remain in $\Theta$
(see e.g., \cite[Section 3.44]{ljung2}). 
However, in order to avoid unnecessary technical details and to 
keep the exposition as concise as possible, 
this aspect of algorithm (\ref{7.1}) is not discussed here. 
Instead, similarly as in \cite{benveniste}, \cite{ljung1},  \cite{ljung2}, 
we state our asymptotic results 
in a local form. 
\end{remark}

The minimization of the negative log-likelihood 
using stochastic gradient search has a long tradition in statistical inference, 
system identification and signal and image processing, 
while the asymptotic properties of the corresponding algorithms have studied in 
a number of papers
(see e.g., \cite{benveniste}, \cite{gu}, \cite{ljung2}, \cite{nevelson}, \cite{younes} and references 
cited therein). 
Although the available literature provides a good insight into 
the asymptotic behavior of the recursive maximum likelihood method, 
the existing results on the convergence and convergence rate 
(of algorithm (\ref{7.1})) rely on very restrictive conditions: 
These results 
require the negative log-likelihood $f(\cdot )$ 
to have an isolated minimum $\theta_{*}$ and its gradient $\nabla f(\cdot )$ 
to admit representation (\ref{1.101}). 
As such, the existing results do not cover the case when the negative log-likelihood 
$f(\cdot )$ has multiple and non-isolated minima, 
which, as explained in Remark \ref{r7.1}, often happens in practice. 
The aim of Therems \ref{theorem7.1} and \ref{theorem7.2} is to fill this gap 
in the literature on maximum likelihood estimation. 



\section{Example 4: Temporal-Difference Learning} \label{section8} 

In this section, the asymptotic behavior of online algorithms for 
temporal-difference learning is analyzed using 
Theorems \ref{theorem1.1}, \ref{theorem1.2} and 
\ref{theorem2.1}. 

In order to explain temporal-difference learning and 
to define the corresponding algorithm, 
we use the following notation. 
$N \geq 1$ is an integer, 
while 
${\cal X} \subseteq \mathbb{R}^{N}$ is a measurable set. 
$\{X_{n} \}_{n\geq 0}$
is an ${\cal X}$-valued Markov chain defined on a 
probability space 
$(\Omega, {\cal F}, P )$, while $P(\cdot,\cdot )$ is its transition kernel. 
$c:\mathbb{R}^{N} \rightarrow \mathbb{R}$ is a locally Lipschitz continuous function. 
$\beta \in (0,1)$ is a constant, while function $g(x)$ is defined as 
\begin{align*}
	g(x)
	=
	E\left(\left. 
	\sum_{n=0}^{\infty } \beta^{n} c(X_{n} ) 
	\right| X_{0}=x
	\right)
\end{align*}
for $x\in {\cal X}$. 
$d_{\theta }\geq 1$ is an integer, 
while $G_{\theta }(x)$ is a real-valued measurable function of 
$(\theta,x) \in \mathbb{R}^{d_{\theta } } \times {\cal X}$. 
$f(\cdot )$ is the function defined by 
\begin{align}\label{8.101} 
	f(\theta )
	=
	\frac{1}{2}
	\lim_{n\rightarrow \infty } 
	E(g(X_{n} ) - G_{\theta }(X_{n} ) )^{2}
\end{align}
for $\theta\in \mathbb{R}^{d_{\theta } }$. 
With this notation, the problem of temporal-difference learning can be 
posed as the minimization of 
$f(\cdot )$. 
In this context, 
$c(x)$ is considered as a cost of visiting state
$x$, 
while 
$g(x)$ is regarded to as the total discounted cost incurred by 
$\{X_{n} \}_{n\geq 0}$ 
when $\{X_{n} \}_{n\geq 0}$ starts from state $x$. 
$G_{\theta }(\cdot )$
is a parameterized approximation of $g(\cdot )$, 
while $\theta$ is the parameter to be tuned 
through the process of temporal-difference learning. 
For more details on temporal-difference learning, 
see e.g., 
\cite{bertsekas&tsitsiklis1}, \cite{powell} 
and references cited therein. 

Function $f(\cdot )$ can be minimized by the following 
algorithm: 
\begin{align}
	& \label{8.1}
	Y_{n+1}
	=
	\beta Y_{n} 
	+
	H_{\theta_{n} }(X_{n} ), 
	\\
	& \label{8.3}
	\theta_{n+1} 
	=
	\theta_{n} 
	+
	\alpha_{n} 
	(c(X_{n} ) + \beta G_{\theta_{n} }(X_{n+1} ) - G_{\theta_{n} }(X_{n} ) ) Y_{n+1}, 
	\;\;\; 
	n\geq 0. 
\end{align}
In this recursion, 
$\{\alpha_{n} \}_{n\geq 0}$ is a sequence of positive reals. 
$\theta_{0}\in \mathbb{R}^{d_{\theta } }$ is an arbitrary vector, 
while 
$H_{\theta }(\cdot ) = \nabla_{\theta } G_{\theta }(\cdot )$. 
In the literature on reinforcement learning, 
recursion (\ref{8.1}), (\ref{8.3}) is known as 
$TD(1)$ temporal-difference learning algorithm
with a nonlinear function approximation, 
while $G_{\theta }(\cdot )$ is referred to as 
a function approximation (or just as an `approximator'). 

We analyze algorithm 
(\ref{8.1}), (\ref{8.3}) under the following assumptions: 

\begin{assumption} \label{a8.1}
${\cal X}$ is compact. 
\end{assumption} 

\begin{assumption} \label{a8.2} 
$\{X_{n} \}_{n\geq 0}$ has a unique invariant probability measure $\pi(\cdot )$. 
Moreover, there exist real numbers $\rho\in (0,1)$, $C\in [1,\infty )$ such that 
\begin{align*}
	|P^{n}(x,B) - \pi(B) |
	\leq 
	C\rho^{n} 
\end{align*}
for all $x\in {\cal X}$, $n\geq 0$ and any measurable set $B\subseteq{\cal X}$  
(here, $P^{n}(\cdot,\cdot)$ denotes the $n$-th transition probability of 
$\{X_{n} \}_{n\geq 0}$).  
\end{assumption}

\begin{assumption}\label{a8.3}
For each $x\in {\cal X}$, 
$G_{\theta }(x)$ is real-analytic in $\theta$ on entire $\mathbb{R}^{d_{\theta } }$. 
Moreover, 
$G_{\theta }(x)$ has 
a (complex-valued) continuation 
$\hat{G}_{\eta}(x)$ with the following properties: 
\begin{compactenum}[(i)]
\item
$\hat{G}_{\eta}(x)$ maps $(\eta,x)\in \mathbb{C}^{d_{\theta } } \times {\cal X}$ 
to $\mathbb{C}$. 
\item
$\hat{G}_{\theta}(x) = G_{\theta}(x)$
for all $\theta\in\mathbb{R}^{d_{\theta } }$, $x\in {\cal X}$. 
\item
For any $\theta\in\mathbb{R}^{d_{\theta } }$, 
there exist a real number
$\delta_{\theta } \in (0,1)$ such that 
$\hat{G}_{\eta}(x)$ is 
analytic in $\eta$ 
and
continuous in $(\eta,x)$
for any 
$\eta\in\mathbb{C}^{d_{\theta } }$, $x\in {\cal X}$ 
satisfying $\|\eta-\theta\|\leq\delta_{\theta }$. 
\end{compactenum}
\end{assumption}

Our main results on the properties of 
$f(\cdot )$
and asymptotic behavior of the algorithm 
(\ref{8.1}), (\ref{8.3}) are presented in the next two theorems. 

\begin{theorem} \label{theorem8.1}
Let Assumptions \ref{a8.1} -- \ref{a8.3} hold.
Then, $f(\cdot )$ is analytic on entire $\mathbb{R}^{d_{\theta } }$.  
\end{theorem}

\begin{theorem} \label{theorem8.2}
Let Assumptions \ref{a2.1} and \ref{a8.1} -- \ref{a8.3} hold. 
Then, all conclusions of Theorem \ref{theorem2.1} are true for 
$\{\theta_{n} \}_{n\geq 0}$ defined in this section. 
\end{theorem}

A proof of Theorems \ref{theorem8.1} and \ref{theorem8.2} is provided in Section \ref{section8*}. 

Assumptions \ref{a8.1} and \ref{a8.2} correspond to the stability of Markov chain 
$\{X_{n} \}_{n\geq 0}$. 
In this or similar form, they are involved in any result on the 
asymptotic behavior of temporal-difference learning. 
On the other side, 
Assumption \ref{a8.3} is related to the properties of  
$G_{\theta }(\cdot )$. 
It covers some of the most popular function approximations 
used in reinforcement learning 
(e.g., 
feedforward neural networks with analytic activation functions; 
for details see \cite{bertsekas&tsitsiklis1}, \cite{powell}). 

Asymptotic properties of temporal-difference learning 
have been the subject of a number of papers
(see \cite{bertsekas&tsitsiklis1}, \cite{powell} and references cited therein). 
However, the available literature on reinforcement learning 
does not offer any information 
on the single limit-point convergence and convergence rate 
which can be verified for temporal-difference learning algorithms with non-linear 
function approximation
(i.e., for $G_{\theta }(\cdot )$ being nonlinear in $\theta$). 
Similarly as in the case of supervised learning, 
the reason comes out of the fact that 
the existing asymptotic results for stochastic gradient search 
hold under very restrictive conditions which are hard (if possible at all) 
to demonstrate for such algorithms. 
The aim of Theorems \ref{theorem8.1} and \ref{theorem8.2} is to fill this gap 
in the literature on reinforcement learning. 

\section{Example 5: Identification of Hidden Markov Models}\label{section10} 

In this section, Theorems \ref{theorem1.1}, \ref{theorem1.2} and \ref{theorem2.1}
are applied to the asymptotic analysis of 
recursive maximum split-likelihood algorithm. 
Recursive maximum split-likelihood algorithm is a method for the identification of hidden Markov models. 

In order to define hidden Markov models and to state the problem of their identification, 
we need the following notation. 
$N_{x}>1$ and $N_{y}>1$ are integers, 
while ${\cal X} = \{1,\dots,N_{x} \}$ and ${\cal Y} = \{1,\dots,N_{y} \}$. 
$p(x'|x)$ and $q(y|x)$ are non-negative functions of 
(respectively)
$(x,x') \in {\cal X} \times {\cal X}$, $(x,y) \in {\cal X} \times {\cal Y}$
which satisfy 
\begin{align*}
	\sum_{x' \in {\cal X} } p(x'|x) = 1, 
	\;\;\;\;\; 
	\sum_{y \in {\cal Y} } q(y|x) = 1
\end{align*}
for each $x\in {\cal X}$. 
$\{ (X_{n}, Y_{n} ) \}_{n\geq 0}$ is an ${\cal X} \times {\cal Y}$-valued 
Markov chain 
which is defined on a (canonical) probability space 
$(\Omega, {\cal F}, P )$ and admits 
\begin{align*}
	P(X_{n+1} = x', Y_{n+1} = y' |X_{n} = x, Y_{n} = y )
	=
	q(y'|x') p(x'|x)
\end{align*}
for all $x,x' \in {\cal X}$, $y,y' \in {\cal Y}$. 
On the other side, 
$d_{\theta } \geq 1$ is an integer, while $\Theta \subseteq R^{d_{\theta } }$ is an open set. 
$\pi_{\theta }(x)$,
$p_{\theta }(x'|x)$ and $q_{\theta }(y|x)$ are non-negative functions of
(respectively)
$(\theta, x ) \in \Theta \times {\cal X}$, 
$(\theta, x, x' ) \in \Theta \times {\cal X} \times {\cal X}$, 
$(\theta, x, y ) \in \Theta \times {\cal X} \times {\cal Y}$
with the following properties: 
They are differentiable in $\theta$ 
for all $\theta\in \Theta$, $x,x' \in {\cal X}$, $y \in {\cal Y}$
and satisfy 
\begin{align*}
	\sum_{x'\in {\cal X} } \pi_{\theta }(x') = 1, 
	\;\;\;\;\; 
	\sum_{x' \in {\cal X} } p_{\theta }(x'|x) = 1,  
	\;\;\;\;\;
	\sum_{y \in {\cal Y} } q_{\theta }(y|x) = 1
\end{align*}
for each $\theta\in\Theta$, $x\in {\cal X}$. 
For $\theta\in\Theta$, 
$\{ (X_{n}^{\theta }, Y_{n}^{\theta } ) \}_{n\geq 0}$
is an ${\cal X} \times {\cal Y}$-valued Markov chain 
which is defined on a (canonical) probability space $(\Omega, {\cal F}, P_{\theta } )$
and satisfies 
\begin{align*}
	&
	P_{\theta }(X_{0}^{\theta } = x, Y_{0}^{\theta } = y ) 
	=
	q_{\theta }(y|x) \pi_{\theta }(x), 
	\\
	&
	P_{\theta }(X_{n+1}^{\theta } = x', Y_{n+1}^{\theta } = y'
	|X_{n}^{\theta } = x, Y_{n}^{\theta } = y )
	= 
	q_{\theta }(y'|x') p_{\theta }(x'|x)
\end{align*}
for all $x,x' \in {\cal X}$, $y,y' \in {\cal Y}$, $n\geq 0$. 
For $\theta\in\Theta$, $y_{1:N} = (y_{1}, \dots, y_{N} ) \in {\cal Y}^{N}$ 
and $N\geq 1$, 
function $\phi_{N,\theta}(y_{1:N} )$ 
is defined as 
\begin{align*}
	\phi_{N, \theta} (y_{1:N} )
	= &
	-
	\frac{1}{N}
	\log\left(
	\sum_{x_{0},\dots,x_{N} \in {\cal X} } 
	\left(
	\prod_{i=1}^{N}
	\left(
	q_{\theta}(y_{i}|x_{i} ) p_{\theta }(x_{i} |x_{i-1} ) 
	\right)
	\right)
	\pi_{\theta }(x_{0} ) 
	\right), 
\end{align*}
while 
\begin{align*}
	&
	f_{N}(\theta ) 
	= 
	\lim_{n\rightarrow \infty } 
	E(\phi_{N,\theta}(Y_{nN+1:(n+1)N } ) ), 
	\;\;\;\;\; 
	f_{\infty }(\theta ) 
	=
	\lim_{N\rightarrow\infty } 
	E(\phi_{N,\theta}(Y_{1:N} ) )
\end{align*}
(here, $Y_{nN+1:(n+1)N}$ stands for $(Y_{nN+1}, \dots, Y_{(n+1)N} )$). 

In the statistics and engineering literature, 
$\{ (X_{n}, Y_{n} ) \}_{n\geq 0}$  
(as well as 
$\{ (X_{n}^{\theta }, Y_{n}^{\theta } ) \}_{n\geq 0}$)
is known as a hidden Markov model, 
while 
$X_{n}$ and $Y_{n}$ are 
its (unobservable) state and (observable) output at discrete-time $n$. 
The identification of $\{ (X_{n}, Y_{n} ) \}_{n\geq 0}$ can be stated as follows: 
Given a realization of the output sequence $\{Y_{n} \}_{n\geq 0}$, 
estimate $\{p(x'|x) \}_{x,x' \in {\cal X} }$, $\{q(y|x) \}_{x \in {\cal X}, y \in {\cal Y} }$. 
If the identification is based on the maximum likelihood principle
and the parameterized model 
$\{p_{\theta }(x'|x) \}_{x,x' \in {\cal X} }$, $\{q_{\theta }(y|x) \}_{x \in {\cal X}, y \in {\cal Y} }$, 
the estimation reduces to the minimization of the negative log-likelihood 
$f_{\infty }(\cdot )$ over $\Theta$. 
In this context, 
$\{ (X_{n}^{\theta }, Y_{n}^{\theta } ) \}_{n\geq 0}$
is considered as a candidate model of 
$\{ (X_{n}, Y_{n} ) \}_{n\geq 0}$. 
For more details on hidden Markov models and their identification see 
\cite[Part II]{cappe&moulines&ryden} and references cited therein. 

As the negative log-likelihood $f(\cdot )$ and its gradient are rarely 
available analytically, 
$f_{\infty }(\cdot )$ is usually minimized by stochastic gradient search. 
The consistent estimation of $\nabla f_{\infty }(\cdot )$
is computationally expensive (even for moderately large $N_{x}$, $N_{y}$), 
since it is based on the optimal filter and the filter derivatives 
(see e.g., \cite{cappe&moulines&ryden}). 
To reduce the computational complexity, 
a number of approaches based on approximate maximum likelihood (also known 
as pseudo-likelihood) has been proposed. 
Among them, 
the maximum split-likelihood method 
\cite{ryden1}, \cite{ryden2} has attracted a considerable attention in 
the literature. 
This approach is based on the following fact. 
If $\{X_{n} \}_{n\geq 0}$ is geometrically ergodic and 
if the optimal filter for the candidate model 
$\{ (X_{n}^{\theta }, Y_{n}^{\theta } ) \}_{n\geq 0}$ is stable, 
then 
\begin{align*}
	\nabla f_{\infty }(\theta ) 
	=
	\lim_{N\rightarrow\infty } \nabla f_{N}(\theta )
	=
	\lim_{N\rightarrow \infty } \lim_{n\rightarrow \infty } 
	E\left(
	\nabla_{\theta } \phi_{N,\theta}(Y_{nN+1:(n+1)N} )
	\right)
\end{align*}
(see Theorem \ref{theorem10.1}, below). 
Hence, 
$\nabla_{\theta }\phi_{N,\theta}(Y_{nN+1:(n+1)N} )$
is a reasonably good estimator of $\nabla f_{\infty }(\theta )$ when $n,N\gg 1$. 
Combining this estimator 
with stochastic gradient search, 
we get the recursive maximum split-likelihood algorithm: 
\begin{align}\label{10.1}	
	\theta_{n+1}
	=
	\theta_{n} 
	-
	\alpha_{n} 
	\psi_{N,\theta_{n} }(Y_{nN+1:(n+1)N} ), 
	\;\;\; n\geq 0. 
\end{align}
Here, $\{\alpha_{n} \}_{n\geq 0}$ is a sequence of positive real numbers, 
$N\geq 1$ is a fixed integer, and 
$\psi_{N,\theta }(\cdot ) = \nabla_{\theta } \phi_{N,\theta }(\cdot )$. 

To analyze algorithm (\ref{10.1}), we need the following assumptions: 

\begin{assumption}\label{a10.1} 
$\{X_{n} \}_{n\geq 0}$ is geometrically ergodic. 
\end{assumption}

\begin{assumption}\label{a10.2} 
$p_{\theta }(x'|x)>0$, 
$q_{\theta }(y|x)>0$ and $\pi_{\theta }(x)>0$
for all $\theta\in\Theta$, $x,x'\in {\cal X}$, $y\in {\cal Y}$. 
\end{assumption}

\begin{assumption}\label{a10.3} 
For each $x,x'\in {\cal X}$, $y\in {\cal Y}$, 
$p_{\theta }(x'|x)$, 
$q_{\theta }(y|x)$ and $\pi_{\theta }(x)$
are real-analytic in $\theta$ on entire $\Theta$. 
\end{assumption}

Assumption \ref{a10.1} is related to the stability 
of the system being identified $\{ (X_{n}, Y_{n} ) \}_{n\geq 0}$. 
In this or similar form, it is often involved in the asymptotic analysis 
of the identification methods for hidden Markov models
(see e.g. \cite{cappe&moulines&ryden} and references cited therein). 
Assumptions \ref{a10.2} and \ref{a10.3} correspond to the stability and 
parameterization of the 
candidate model $\{ (X_{n}^{\theta }, Y_{n}^{\theta } ) \}_{n\geq 0}$. 
It is satisfied for many commonly used parameterizations 
(e.g., natural, trigonometric and exponential).  

Our main results on the properties of $f_{N}(\cdot )$, $f_{\infty }(\cdot )$ 
and the asymptotic behavior of 
algorithm (\ref{10.1}) are provided in the next two theorems. 

\begin{theorem} \label{theorem10.1}
Let Assumptions \ref{a10.1} and \ref{a10.2} hold. 
Then, the following is true: 
\begin{compactenum}[(i)]
\item
$f_{N}(\cdot )$ and $f_{\infty }(\cdot )$ are analytic on entire $\Theta$.  
\item
For each $\theta\in\Theta$, there exists a real number $L_{\theta } \in (0,\infty )$
(not depending on $N$)\footnote
{$L_{\theta }$ depends only on the mixing rate of $\{X_{n} \}_{n\geq 0}$ 
and on the upper bounds and Lipschitz constants of 
$p_{\theta}(x'|x)$, $q_{\theta }(y|x)$, $\pi_{\theta }(x)$ and their derivatives. 
For further details, see the proof of Theorem \ref{theorem10.1}. }
such that 
\begin{align}\label{t10.1.1*}
	\max\{|f_{N}(\theta ) - f_{\infty } (\theta ) |, 
	\|\nabla f_{N}(\theta ) - \nabla f_{\infty } (\theta ) \| \}
	\leq 
	L_{\theta }/N. 
\end{align}
\end{compactenum}
\end{theorem}

\begin{theorem} \label{theorem10.2}
Let Assumptions \ref{a2.1}, \ref{a10.1} and \ref{a10.2} hold. 
Then, the following is true: 
\begin{compactenum}[(i)]
\item
$\hat{\theta} = \lim_{n\rightarrow \infty } \theta_{n}$
exists and satisfies 
$\nabla f_{N}(\hat{\theta } ) = 0$
w.p.1 on $\Lambda$. 
\item
$\|\nabla f(\theta_{n} ) \|^{2}	=
o\big(\gamma_{n}^{-\hat{p} } \big)$,
$|f(\theta_{n} ) - f(\hat{\theta } ) | =
o\big(\gamma_{n}^{-\hat{p} } \big)$
and 
$\|\theta_{n} - \hat{\theta} \| =
o\big(\gamma_{n}^{-\hat{q} } \big)$
w.p.1 
on 
$\Lambda 
\cap \{\hat{r} > r \}$. 
\item
$	\|\nabla f(\theta_{n} ) \|^{2} =
O\big(\gamma_{n}^{-\hat{p} } \big)$,
$|f(\theta_{n} ) f(\hat{\theta} ) | =
O\big(\gamma_{n}^{-\hat{p} } \big)$
and  
$\|\theta_{n} - \hat{\theta} \| =
O\big(\gamma_{n}^{-\hat{q} } \big)$
w.p.1 on 
$\Lambda 
\cap \{\hat{r} \leq r \}$. 
\item
$\|\nabla f(\theta_{n} ) \|^{2} = o(\gamma_{n}^{-p} )$
and 
$|f(\theta_{n} ) - f(\hat{\theta} ) | = o(\gamma_{n}^{-p} )$
w.p.1 on 
$\Lambda$. 
\end{compactenum}
\end{theorem}

A proof of Theorems \ref{theorem10.1} and \ref{theorem10.2} is  provided in Section \ref{section10*}. 
$p$, $\hat{p}$, $\hat{q}$ and $\hat{r}$
are defined in Theorem \ref{theorem1.2} and Corollary \ref{corollary1.1},\footnote
{In this case, $f(\cdot )$ should be replaced by $f_{N}(\cdot )$.} 
while $\Lambda$ is specified in (\ref{7.3}). 

\begin{remark}\label{r10.2}
Similarly as (\ref{7.1}), 
algorithm (\ref{10.1}) usually involves a projection (or truncation) device which ensures
that estimates $\{\theta_{n} \}_{n\geq 0}$ remain in $\Theta$
(see e.g., \cite[Section 3.44]{ljung2}). 
However, in order to avoid unnecessary technical details and to 
keep the exposition as concise as possible, 
this aspect of algorithm (\ref{7.1}) is not discussed here. 
Instead, similarly as in \cite{benveniste}, \cite{ljung1},  \cite{ljung2}, 
we state our asymptotic results 
in a local form. 
\end{remark}

\begin{remark}
As in the case of the maximum likelihood estimation for i.i.d. data, 
$f_{N}(\cdot )$, $f_{\infty }(\cdot )$ are usually multimodal 
(notice that $f_{N}(\cdot )$, $f_{\infty }(\cdot )$ 
are invariant in the order of their arguments). 
In addition to this, $f_{N}(\cdot )$, $f_{\infty }(\cdot )$ 
are likely to have non-isolated minima 
(which inevitably happens whenever 
model $\{(X_{n}^{\theta }, Y_{n}^{\theta } ) \}_{n\geq 0}$ is over-parameterized for 
$\{(X_{n}, Y_{n} )\}_{n\geq 0}$). 
\end{remark}

The asymptotic properties of the maximum split-likelihood have throughly been studies 
in \cite{ryden1}, \cite{ryden2} (see also \cite{cappe&moulines&ryden}. 
Although the results of \cite{ryden1}, \cite{ryden2} offer a good insight into the method, 
they hold under restrictive conditions: 
These results guarantee the single limit-point convergence of 
$\{\theta_{n} \}_{n\geq 0}$ and provide the convergence rate only if 
$f_{N}(\cdot )$ has a unique minimum at which $\nabla^{2} f_{N}(\cdot )$ is positive definite. 
Hence, \cite{ryden1}, \cite{ryden2} do not cover 
the case where $f_{N}(\cdot )$ has multiple non-isolated minima, 
which, as explained in Remark \ref{r10.2} often happen in practice. 
The purpose of Theorems \ref{theorem10.1} and \ref{theorem10.2} is to overcome 
these difficulties. 

\section{Example 6: Identification of Linear Stochastic Systems} \label{section6}

To illustrate the general results of Sections \ref{section1} and \ref{section2}, 
we apply them to the asymptotic analysis of the recursive prediction error 
method for identification of linear stochastic systems. 
To avoid unnecessary technical details and complicated notation, 
only the identification of univariate ARMA models is considered here. 
However, it is straightforward to generalize the obtained results to 
any linear stochastic system. 

To define the recursive prediction error methods for 
ARMA models, 
we use the following notation. 
$M,N\geq 1$ are integers, while $d_{\theta } = M+N$. 
$A_{\theta }(\cdot )$ and $B_{\theta }(\cdot )$ are the polynomials defined by 
\begin{align*}
	A_{\theta }(z) 
	=
	1 
	-
	\sum_{k=1}^{M} a_{k} z^{-k}, 
	\;\;\;\;\; 
	B_{\theta }(z) 
	=
	1 
	+
	\sum_{k=1}^{N} b_{k} z^{-k}
\end{align*}
for $z\in \mathbb{C}$,  
$a_{1}, \dots, a_{M}, b_{1}, \dots, b_{N} \in \mathbb{R}$ and 
$\theta = [a_{1} \cdots a_{M} \; b_{1} \cdots b_{N} ]^{T}$
($\mathbb{C}$ denotes the set of complex numbers). 
${\cal Y}\subseteq\mathbb{R}$ is a measurable set, while 
\begin{align*}
	\Theta_{a}
	=
	\{\theta \in \mathbb{R}^{d_{\theta} }: A_{\theta}(z) = 0 \Rightarrow |z| < 1 \}, 
	\;\;\;\;\;  
	\Theta_{b}
	=
	\{\theta \in \mathbb{R}^{d_{\theta} }: B_{\theta}(z) = 0 \Rightarrow |z| < 1 \}  
\end{align*}
and $\Theta=\Theta_{a} \bigcap \Theta_{b}$. 
$\{Y_{n} \}_{n\geq 0}$ is a ${\cal Y}$-valued stochastic process 
which represents 
the signal generated by the system being identified. 
For $\theta \in \Theta$, 
$\{Y_{n}^{\theta} \}_{n\geq 0}$ is the output of the ARMA model 
\begin{align} \label{6.1'}
	A_{\theta }(q) Y_{n}^{\theta } 
	=
	B_{\theta }(q) U_{n}, 
	\;\;\; n \geq 0, 
\end{align}
where 
$\{U_{n} \}_{\geq 0}$ is a real-valued white noise and 
$q^{-1}$ is the (backward) time-shift operator.  
For the same $\theta$, 
$\{\varepsilon_{n}^{\theta } \}_{n\geq 0}$ is the stochastic process 
generated by the recursion 
\begin{align} \label{6.1''}
	B_{\theta }(q) \varepsilon_{n}^{\theta } 
	=
	A_{\theta }(q) Y_{n}, 
	\;\;\; n\geq 0.    
\end{align}
In that case, 
$\hat{Y}_{n}^{\theta } = Y_{n} - \varepsilon_{n}^{\theta }$ 
is the mean-square optimal prediction  
of $Y_{n}$ given 
$Y_{0},\dots,Y_{n-1}$ and model (\ref{6.1'}) 
(for details see e.g., \cite{ljung2}, \cite{ljung3}).   
On the other side, $\varepsilon_{n}^{\theta }$ can be interpreted 
as the prediction error. 

The parametric identification of ARMA models can be stated as follows: 
Given a realization of $\{Y_{n} \}_{n\geq 0}$, 
estimate the values of $\theta$ for which model (\ref{6.1'})
provides the best approximation to signal 
$\{Y_{n} \}_{n\geq 0}$. 
If the identification is based on the prediction error principle, 
this estimation problem reduces to the minimization of 
the mean-square prediction error 
\begin{align*}
	f(\theta )
	=
	\frac{1}{2} 
	\lim_{n\rightarrow \infty } 
	E\left((\varepsilon_{n}^{\theta } )^{2} \right) 
\end{align*}
with respect to $\theta\in\Theta$. 
In online settings, 
$f(\cdot )$ is usually minimized by stochastic gradient (or stochastic Newton) 
algorithm. 
Such an algorithm is defined by the following recursion: 
\begin{align}
	& \label{6.1}
	\phi_{n}
	=
	[Y_{n} \cdots Y_{n-M+1 } \; \varepsilon_{n} \cdots \varepsilon_{n-N+1 } ]^{T}, 
	\\
	& \label{6.3}
	\varepsilon_{n+1} 
	=
	Y_{n+1} - \phi_{n}^{T} \theta_{n}, 
	\\
	& \label{6.5}
	\psi_{n+1} 
	=
	\phi_{n}
	-
	[\psi_{n} \cdots \psi_{n-N+1 } ] 
	D \: \theta_{n}, 
	\\
	& \label{6.7}
	\theta_{n+1} 
	=
	\theta_{n} 
	+
	\alpha_{n} \psi_{n+1} \varepsilon_{n+1}, 
	\;\;\;\;\; n\geq 0.  
\end{align}
In this recursion, 
$\{\alpha_{n} \}_{n\geq 0}$ denotes a sequence of positive reals.  
$D$ is the $N\times (M+N)$ block-matrix 
defined by $D = [{\boldsymbol 0} \; {\boldsymbol I} ]$, 
where ${\boldsymbol I}$ and ${\boldsymbol 0}$ denote $N\times N$ unit matrix 
and $N\times M$ zero matrix (respectively). 
$\{Y_{n} \}_{n\geq -M }$ is a real-valued stochastic process 
defined on a probability space
$(\Omega, {\cal F}, P)$.  
$\theta_{0} \in \Theta$, $\psi_{0},\dots,\psi_{-N+1} \in \mathbb{R}^{d_{\theta } }$ 
are arbitrary vectors, while 
$\varepsilon_{0},\dots,\varepsilon_{-N+1} \in \mathbb{R}$ 
are arbitrary numbers. 
$\theta_{0}, \varepsilon_{0},\dots,\varepsilon_{-N+1}, 
\psi_{0},\dots,\psi_{-N+1}$ 
represent the initial conditions of the algorithm
(\ref{6.1}) -- (\ref{6.7}). 
In the literature on system identification, 
recursion (\ref{6.1}) -- (\ref{6.7}) is known as 
the recursive prediction error algorithm for ARMA models. 
$\varepsilon_{n}$ is referred to as the prediction error, 
while $\psi_{n}$ is the negative gradient of $\varepsilon_{n}$
with respect to $\theta$ 
(for more details see 
\cite{ljung2}, \cite{ljung3} 
and references cited therein). 

We study the asymptotic behavior of algorithm (\ref{6.1}) -- (\ref{6.7}) for the case 
where $\{Y_{n} \}_{n\geq 0}$ is an output of a Markovian system. 
More specifically, we assume that there exist an integer $L\geq 1$, 
a measurable set ${\cal X} \subseteq \mathbb{R}^{L}$ and 
an ${\cal X}$-valued stochastic process $\{X_{n} \}_{n\geq 0}$ 
defined on $(\Omega, {\cal F}, P )$ 
such that $\{(X_{n}, Y_{n} ) \}_{n\geq 0}$ is a Markov chain. 
In this context, $\{X_{n} \}_{n\geq 0}$ can be interpreted as
unobservable states of the system being identified. 

Let ${\cal W} = {\cal X} \times {\cal Y}$, 
while $\{W_{n} \}_{n\geq 0}$ is the stochastic process defined by 
$W_{n} = [X_{n}^{T} \; Y_{n} ]^{T}$ for $n\geq 0$. 
To analyze algorithm (\ref{6.1}) -- (\ref{6.7}), we rely on the following 
assumptions: 

\begin{assumption} \label{a6.1}
${\cal W}$ is compact. 
\end{assumption}

\begin{assumption} \label{a6.2} 
$\{W_{n} \}_{n\geq 0}$ has a unique invariant probability measure 
$\pi(\cdot )$. 
Moreover, there exist real numbers $\rho\in (0,1)$, $C\in [1,\infty )$
such that 
\begin{align*}
	|P^{n}(w,B) - \pi(B) |
	\leq 
	C\rho^{n}
\end{align*}
for all $w\in {\cal W}$, $n\geq 0$ and any measurable set $B\subseteq {\cal W}$ 
(here, $P^{n}(\cdot,\cdot)$ denotes the $n$-th step transition probability of 
$\{W_{n} \}_{n\geq 0}$). 
\end{assumption}

\begin{assumption} \label{a6.3}
For any compact set $Q \subset \Theta$, 
\begin{align} \label{6.9}
	\sup_{n\geq 0} 
	E\left(
	(\varepsilon_{n}^{4} + \|\psi_{n} \|^{4} )
	I_{ \{\tau_{Q} \geq n \} }
	\right)
	< \infty,  
\end{align}
where $\tau_{Q} = \inf\{n\geq 0: \theta_{n} \notin Q \}$. 
\end{assumption}

Our main results on the properties of $f(\cdot )$ and the asymptotic behavior of 
algorithm (\ref{6.1}) -- (\ref{6.7}) are provided in the next two theorems. 

\begin{theorem} \label{theorem6.1}
Let Assumptions \ref{a6.1} -- \ref{a6.3} hold. 
Then, $f(\cdot )$ is analytic on entire $\Theta$.  
\end{theorem}

\begin{theorem} \label{theorem6.2}
Let Assumptions \ref{a2.1}, \ref{a6.1} and \ref{a6.2} hold. 
Then, the following is true: 
\begin{compactenum}[(i)]
\item
$\hat{\theta} = \lim_{n\rightarrow \infty } \theta_{n}$
exists and satisfies 
$\nabla f(\hat{\theta } ) = 0$
w.p.1 on $\Lambda$. 
\item
$\|\nabla f(\theta_{n} ) \|^{2}	=
o\big(\gamma_{n}^{-\hat{p} } \big)$,
$|f(\theta_{n} ) - f(\hat{\theta } ) | =
o\big(\gamma_{n}^{-\hat{p} } \big)$
and 
$\|\theta_{n} - \hat{\theta} \| =
o\big(\gamma_{n}^{-\hat{q} } \big)$
w.p.1 
on 
$\Lambda 
\cap \{\hat{r} > r \}$. 
\item
$	\|\nabla f(\theta_{n} ) \|^{2} =
O\big(\gamma_{n}^{-\hat{p} } \big)$,
$|f(\theta_{n} ) f(\hat{\theta} ) | =
O\big(\gamma_{n}^{-\hat{p} } \big)$
and  
$\|\theta_{n} - \hat{\theta} \| =
O\big(\gamma_{n}^{-\hat{q} } \big)$
w.p.1 on 
$\Lambda 
\cap \{\hat{r} \leq r \}$. 
\item
$\|\nabla f(\theta_{n} ) \|^{2} = o(\gamma_{n}^{-p} )$
and 
$|f(\theta_{n} ) - f(\hat{\theta} ) | = o(\gamma_{n}^{-p} )$
w.p.1 on 
$\Lambda$. 
\end{compactenum}
\end{theorem}

A proof of Theorems \ref{theorem6.1} and \ref{theorem6.2} is  provided in Section \ref{section6*}. 
$p$, $\hat{p}$, $\hat{q}$ and $\hat{r}$
are defined in Theorem \ref{theorem1.2} and Corollary \ref{corollary1.1}, 
while $\Lambda$ is specified in (\ref{7.3}). 

\begin{remark}
Similarly as (\ref{7.1}), 
algorithm (\ref{6.1}) -- (\ref{6.7}) 
involves a projection (or truncation) device which prevents 
$\{\theta_{n} \}_{n\geq 0}$ from leaving $\Theta$  
(see \cite[Section 3.44]{ljung2}), 
i.e., which ensures the stability of the parameterized model 
$\{Y_{n}^{\theta } \}_{n\geq 0}$ 
(condition $\theta_{n}\in\Theta_{a}$)
and the stability of the prediction error 
$\{\varepsilon_{n}^{\theta } \}_{n\geq 0}$ 
and subrecursion (\ref{6.1}) -- (\ref{6.5}) 
(condition $\theta_{n}\in\Theta_{b}$). 
However, in order to avoid unnecessary technical details and to 
keep the exposition as concise as possible, 
this aspect of algorithm (\ref{6.1}) -- (\ref{6.7}) is not studied here. 
Instead, similarly as in \cite{benveniste}, \cite{ljung1}, \cite{ljung2}, 
we state our asymptotic results 
in a local form. 
Since the stability of algorithm (\ref{6.1}) -- (\ref{6.7}) is not affected 
by the stability of $\{Y_{n}^{\theta } \}_{n\geq 0}$, 
Theorems \ref{theorem6.1} and \ref{theorem6.2} remain valid 
if $\Theta$ is defined by $\Theta=\Theta_{b}$. 
\end{remark} 

\begin{remark}\label{r6.2}
As well-documented in the literature on system identification 
(see e.g., \cite[Section 3.7]{stoica}), 
the mean-square prediction error $f(\cdot )$ is multimodal for 
ARMA models. 
In addition to this, $f(\cdot )$ is likely to have non-isolated minima 
and stationary points (which inevitably happens whenever 
model (\ref{6.1'}) is over-parameterized for 
$\{Y_{n} \}_{n\geq 0}$). 
\end{remark}

Assumptions \ref{a6.1} and \ref{a6.2} correspond to the system being identified. 
They hold whenever the system is a geometrically ergodic hidden Markov model 
(in that case, $\{X_{n} \}_{n\geq 0}$ is the hidden Markov chain).  
They also cover a number of linear and nonlinear
stochastic systems encountered in real-world applications
(including ARMA models driven by bounded i.i.d. or Markovian noise). 
In addition to this, 
Assumptions \ref{a6.1} and \ref{a6.2} allow for the possibility that  
$\{Y_{n} \}_{n\geq 0}$ is not a member of the parametric family of ARMA models 
(\ref{6.1'}) 
(which is rather important from the practical point of view, 
as such models cannot provide an exact representation of a real-world system, 
but only an accurate approximation).   
Unfortunately, Assumption \ref{a6.1} requires states $\{X_{n} \}_{n\geq 0}$
and outputs $\{Y_{n} \}_{n\geq 0}$ to be compactly supported 
(i.e., almost surely bounded). 
Although this may seem restrictive from theoretical point of view, 
it is always satisfied in practice 
(as systems met in real-world applications generate only bounded signals). 
Anyway, relying on the concept of $V$-uniform ergodicity 
(see e.g., \cite[Chapter 16]{meyn&tweedie}), it is relatively straightforward to extend the results of this section 
to Markovian systems with non-compactly supported states and outputs. 

Assumption \ref{a6.3} is related to the stability of 
subrecursion (\ref{6.1}) -- (\ref{6.5}) 
and of sequences  
$\{\varepsilon_{n} \}_{\geq 0}$, $\{\psi_{n} \}_{n\geq 0}$. 
In this or a similar form, 
Assumption \ref{a6.3} is involved in practically all 
asymptotic results for the recursive prediction error 
identification methods. 
E.g., 
\cite[Theorems 4.1 -- 4.3]{ljung2} 
(probably the most general result of this kind)  
require $\{(\varepsilon_{n}, \psi_{n} ) \}_{n\geq 0}$
to visit a fixed compact set infinitely often w.p.1 
on event $\Lambda$. 
When $\{Y_{n} \}_{n\geq 0}$ is generated by a  
Markovian system, such a requirement is practically equivalent to (\ref{6.9}). 

Various aspects of the recursive prediction error identification in 
linear stochastic systems have been the subject 
of numerous papers and books
(see \cite{ljung2}, \cite{ljung3} 
and references cited therein). 
Although the available literature offers a good insight into 
the asymptotic behavior of the recursive prediction error method, 
the existing results on the convergence and convergence rate 
(of algorithm (\ref{6.3}) -- (\ref{6.7})) hold under very restrictive 
conditions: 
These results require the mean-square prediction error $f(\cdot )$ to 
have an isolated minimum $\theta_{*}$ at which 
$\nabla^{2} f(\cdot )$ is positive definite 
(see \cite{ljung2}, probably the strongest result of this type). 
As such, the existing results cannot cover the case when $f(\cdot )$
has multiple and non-isolated minima, 
which, as explained in Remark \ref{r6.2}, often happens in practice. 
The aim of Theorems \ref{theorem6.1} and \ref{theorem6.2} is to fill 
this gap in the literature on system identification. 


\section{Example 7: Simulation-Based Optimization of Controlled 
\\
Markov Chains}\label{section9}

In this section, we explain how 
Theorems \ref{theorem1.1}, \ref{theorem1.2} and \ref{theorem2.1}
can be used to analyze 
the actor-critic algorithms proposed 
by Tsitsiklis and Konda in \cite{konda&tsitsiklis}. 
These algorithms fall into the category of  
reinforcement learning. 
They can be considered as simulation-based methods
for average-cost Markov decision problems, too. 

To state average-cost Markov decision problems and to define 
the actor-critic algorithms of Tsitsiklis and Konda, 
we use the following notation. 
$d_{\theta }\geq 1$ and $M, N > 1$ are integers, 
while ${\cal X} = \{1,\dots,N\}$ and
${\cal Y} = \{1,\dots,M\}$. 
$c(x,y)$, $p(x'|x,y)$ and $q_{\theta }(y|x)$ 
are functions mapping 
$\theta\in\mathbb{R}^{d_{\theta } }$, $x,x' \in {\cal X}$, $y\in {\cal Y}$
to $[0,\infty )$. 
For each $\theta\in\mathbb{R}^{d_{\theta } }$, $x\in {\cal X}$, $y\in {\cal Y}$, 
$p(\cdot |x,y)$ and $q_{\theta }(\cdot |x)$ are 
probability mass functions on ${\cal X}$, ${\cal Y}$ (respectively). 
For each $\theta\in\mathbb{R}^{d_{\theta } }$, 
$\{ (X_{n}^{\theta }, Y_{n}^{\theta } ) \}_{n\geq 0}$ 
is an ${\cal X}\times {\cal Y}$-valued Markov chain
which is defined on a (canonical) probability space 
$(\Omega, {\cal F}, P )$ and satisfies 
\begin{align}\label{9.1'}
	P(X_{n+1}^{\theta }=x, Y_{n+1}^{\theta } = y 
	|X_{n}^{\theta }, Y_{n}^{\theta } )
	=
	q_{\theta }(y|x) p(x|X_{n}^{\theta }, Y_{n}^{\theta } )
\end{align}
for all $x\in {\cal X}$, $y\in {\cal Y}$, 
$n\geq 0$. 
As an immediate consequence of (\ref{9.1'}), 
$\{X_{n}^{\theta } \}_{n\geq 0}$ is also a Markov chain
whose transition kernel
$p_{\theta }(x'|x)$ is defined by  
\begin{align*}
	p_{\theta }(x'|x)
	=
	\sum_{y\in{\cal Y} } 
	p(x'|x,y) q_{\theta }(y|x) 
\end{align*}
for $x,x'\in{\cal X}$. 

An average-cost Markov decision problem with 
parameterized randomized policy can be posed as 
the minimization of 
\begin{align*}
	f(\theta )
	=
	\lim_{n\rightarrow\infty } 
	E\left(
	\frac{1}{n} \sum_{k=1}^{n} 
	c(X_{k}^{\theta }, Y_{k}^{\theta } )
	\right) 
\end{align*}
with respect to $\theta\in\mathbb{R}^{d_{\theta } }$. 
In this context, 
$\{(X_{n}^{\theta }, Y_{n}^{\theta } ) \}_{n\geq 0}$ is referred to as a controlled 
Markov chain with parameterized randomized policy. 
$\{X_{n}^{\theta } \}_{n\geq 0}$ represent the chain states, 
while $\{Y_{n}^{\theta } \}_{n\geq 0}$ are the control actions. 
$p(x'|x,y)$ is the state transition kernel, 
while $q_{\theta }(y|x)$ is the action likelihood. 
$c(x,y)$ is the cost of state-action pair $(x,y)$. 
For further details on controlled Markov chains 
and Markov decision problems, see \cite{bertsekas&tsitsiklis1}, \cite{powell}
and references cited therein. 

In \cite{konda&tsitsiklis}, 
Tsitsiklis and Konda have proposed a class of actor-critic algorithms 
for the minimization of $f(\cdot )$. 
These algorithms are based on Markov chain regeneration
and 
can be defined by the following difference equations: 
\begin{align}
	V_{n+1}
	=& \label{9.1} 
	c(X_{n},Y_{n} ) - \eta_{2,n} 
	+
	(s_{\theta_{n} }(X_{n+1}, Y_{n+1} ) - s_{\theta_{n} }(X_{n}, Y_{n} ) )^{T} \eta_{1,n}, 
	\\
	W_{n+1}  
	=& \label{9.3} 
	W_{n} I_{\{X_{n+1} \neq x_{*} \} }
	+
	s_{\theta_{n} }(X_{n+1}, Y_{n+1} ), 
	\\
	\theta_{n+1}
	=& \label{9.5} 
	\theta_{n}
	-
	\alpha_{n} 
	s_{\theta_{n} }(X_{n+1}, Y_{n+1} ) s_{\theta_{n} }^{T}(X_{n+1}, Y_{n+1} )  
	\eta_{1,n}, 
	\\ 
	\eta_{1,n+1}
	=& \label{9.7} 
	\eta_{1,n} 
	+
	\beta_{n} 
	W_{n+1} V_{n+1}, 
	\\
	\eta_{2,n+1}
	=&
	\eta_{2,n} 
	+
	\beta_{n} (c(X_{n+1}, Y_{n+1} ) - \eta_{2,n} ), 
	\;\;\;\;\; n\geq 0. 
\end{align}
$\{\alpha_{n} \}_{n\geq 0}$ and 
$\{\beta_{n} \}_{n\geq 0}$ are sequences of positive real numbers. 
$\theta_{0}, \eta_{1,0}, W_{0} \in \mathbb{R}^{d_{\theta } }$ are 
arbitrary vectors, while $\eta_{2,0}\in \mathbb{R}$ is an arbitrary number. 
$s_{\theta }(x,y)$ is defined by 
$s_{\theta }(x,y) = \nabla_{\theta } q_{\theta }(y|x) / q_{\theta }(y|x)$
for $\theta\in\mathbb{R}^{d_{\theta } }$, $x\in{\cal X}$, 
$y\in{\cal Y}$. 
$x_{*}$ is a fixed element of ${\cal X}$. 
$\{X_{n} \}_{n\geq 0}$ and $\{Y_{n} \}_{n\geq 0}$ are 
stochastic processes generated through the following Monte Carlo 
simulations: 
For each $n\geq 0$, 
$X_{n+1}$ is simulated from $p(\cdot|X_{n}, Y_{n} )$
(independently of 
$\{\theta_{i},\eta_{1,i},\eta_{2,i} \}_{1\leq i\leq n}$
and 
$\{X_{j},Y_{j} \}_{1\leq j < n }$), 
while $Y_{n+1}$ is simulated from 
$q_{\theta_{n} }(\cdot|X_{n+1} )$
(independently from 
$\{\theta_{i} \}_{1\leq i < n}$
and 
$\{\eta_{1,j},\eta_{2,j}, X_{j},Y_{j} \}_{1\leq j \leq n}$). 
Hence, $\{ (X_{n}, Y_{n} ) \}_{n\geq 0}$ satisfies 
\begin{align*}
	&
	P(X_{n+1}=x,Y_{n+1}=y
	|\theta_{0},\eta_{1,0},\eta_{2,0},X_{0},Y_{0},
	\dots, 
	\theta_{n},\eta_{1,n},\eta_{2,n},X_{n},Y_{n} )
	\\
	&
	=
	q_{\theta_{n} }(y|x) p(x|X_{n}, Y_{n} )
\end{align*}
w.p.1 for $n\geq 0$. 

Algorithm (\ref{9.1}) -- (\ref{9.7}) is analyzed under the following assumptions: 

\begin{assumption}\label{a9.1}
$\lim_{n\rightarrow\infty } \alpha_{n} = 
\lim_{n\rightarrow\infty } \beta_{n} = 
\lim_{n\rightarrow\infty } \alpha_{n}\beta_{n}^{-1} = 0$, 
$\limsup_{n\rightarrow\infty } |\alpha_{n+1}^{-1} - \alpha_{n}^{-1} | <
\infty$, 
$\limsup_{n\rightarrow\infty } |\beta_{n+1}^{-1} - \beta_{n}^{-1} | <
\infty$
and 
$\sum_{n=0}^{\infty } \alpha_{n} = \infty$. 
Moreover, there exists a real number $r\in [1,\infty )$ 
such that 
$\sum_{n=0}^{\infty } \beta_{n}^{2} \gamma_{n}^{2r} < \infty$. 
\end{assumption}

\begin{assumption}\label{a9.2}
For each $\theta\in\mathbb{R}^{d_{\theta } }$, 
$p_{\theta }(x'|x)$ is an irreducible and aperiodic transition kernel. 
\end{assumption}

\begin{assumption}\label{a9.3}
For any compact set $Q\subset\mathbb{R}^{d_{\theta } }$, 
there exists an integer $n_{Q}\geq 1$ and a real number $\varepsilon_{Q} \in (0,1)$
such that 
\begin{align*}
	\sum_{n=1}^{n_{Q} } \:
	\sum_{x_{1},\dots,x_{n} \in {\cal X} }
	p_{\vartheta_{n} }(x_{*}|x_{n} ) \cdots p_{\vartheta_{0} }(x_{1}|x )
	\geq 
	\varepsilon_{Q} 
\end{align*}
for all $x\in{\cal X}$ and any sequence 
$\{\vartheta_{n} \}_{0\leq n \leq n_{Q} }$ in $Q$. 
\end{assumption} 

\begin{assumption}\label{a9.4}
For any compact set $Q\subset\mathbb{R}^{d_{\theta } }$, 
there exists a real number $K_{Q}\in [1,\infty )$ such that 
\begin{align*}
 &
 \|\nabla_{\theta } q_{\theta }(y|x) \| 
 \leq 
 K_{Q} q_{\theta }(y|x), 
 \\
 &
 \|s_{\theta' }(x,y) - s_{\theta''}(x,y) \| 
 \leq 
 K_{Q} \|\theta'-\theta'' \|
\end{align*}
for all $\theta,\theta',\theta'' \in Q$, $x\in{\cal X}$, $y\in{\cal Y}$. 
\end{assumption}

\begin{assumption}\label{a9.5}
For each $x\in{\cal X}$, $y\in{\cal Y}$, 
$q_{\theta}(y|x)$ is real-analytic in $\theta$ on entire 
$\mathbb{R}^{d_{\theta } }$. 
\end{assumption} 

\begin{assumption}\label{a9.6}
For each $\theta\in\mathbb{R}^{d_{\theta } }$, 
\begin{align*}
	\sum_{x\in{\cal X}, y\in{\cal Y} } 
	s_{\theta }(x,y) s_{\theta }^{T}(x,y) 
	\pi_{\theta }(x)
\end{align*}
is positive definite, 
where $\pi_{\theta }(x)$ is the invariant probability mass function of 
$\{X_{n}^{\theta } \}_{n\geq 0}$
(i.e., 
$\pi_{\theta }(x) = \lim_{n\rightarrow\infty } P(X_{n}^{\theta } = x )$). 
\end{assumption}

To the best of our knowledge, the strongest result on 
the asymptotic behavior of 
algorithm (\ref{9.1}) -- (\ref{9.7}) have been provided by 
Tsitsiklis and Konda in \cite{konda&tsitsiklis}. 
They have analyzed algorithm (\ref{9.1}) -- (\ref{9.7}) 
under conditions slightly weaker than 
Assumptions \ref{a9.1} -- \ref{a9.6}.\footnote{ 
The only difference between the conditions adopted in \cite{konda&tsitsiklis} and here 
is that the results of \cite{konda&tsitsiklis} 
hold whenever $q_{\theta}(y|x)$ is twice differentiable in $\theta$, 
while 
Assumption \ref{a9.5} requires $q_{\theta}(y|x)$ to be analytical in $\theta$. 
However, Assumption \ref{a9.5} covers a number of parameterizations of the action likelihood 
$q_{\theta}(y|x)$ such as `natural,' trigonometric of logistic. 
}
As a main result of their analysis, Tsitsiklis and Konda have demonstrated 
that $\liminf_{n\rightarrow\infty } \nabla f(\theta_{n} ) = 0$ w.p.1. 
Using the arguments of Theorems \ref{theorem1.1}, \ref{theorem1.2} and \ref{theorem2.1}, 
much stronger asymptotic results are possible. 
These results are presented in the next two theorems. 

\begin{theorem} \label{theorem9.1}
Let Assumptions \ref{a9.2} and \ref{a9.5} hold.
Then, $f(\cdot )$ is analytic on entire $\mathbb{R}^{d_{\theta } }$.  
\end{theorem}

\begin{theorem} \label{theorem9.2}
Let Assumptions \ref{a9.1} -- \ref{a9.6} hold. 
Then, all conclusions of Theorem \ref{theorem2.1} are true for 
$\{\theta_{n} \}_{n\geq 0}$ defined in this section. 
\end{theorem}

Algorithm (\ref{9.1}) -- (\ref{9.7}) falls into the category of 
two time-scale stochastic approximation (see e.g., \cite{borkar}) 
and does not fit exactly 
into the framework studied in Sections \ref{section1} and \ref{section2}. 
Fortunately, the algorithm is asymptotically equivalent to recursion 
(\ref{1.1}) and (\ref{2.1}), and hence, 
with some modifications, Theorems \ref{theorem1.1}, \ref{theorem1.2} and
\ref{theorem2.1} can be applied to its asymptotic analysis. 
Although intuitively straightforward, these modifications involve 
a number of technical details. 
Therefore, complete proof of Theorems \ref{theorem9.1} and \ref{theorem9.2} 
are provided in separate paper \cite{tadic10}. 
Here, in Section \ref{section9*}, only an outline of the proof is presented. 

\section{Outline of the Proof of 
Theorems \ref{theorem1.1} and \ref{theorem1.2}} \label{section0*}

Theorems \ref{theorem1.1} and \ref{theorem1.2} are proved in several steps. 
These steps can be summarized as follows: 

\begin{step}
\rm The asymptotic properties of 
$\{\theta_{n} \}_{n\geq 0}$, 
$\{f(\theta_{n} ) \}_{n\geq 0}$ and 
$\{\nabla f(\theta_{n} ) \}_{n\geq 0}$  
are analyzed 
(Lemmas \ref{lemma1.1}, \ref{lemma1.2}). 
The analysis is based on Taylor formula and Bellman-Gronwall inequality. 
The obtained results are a prerequisite for Steps 2, 3. 
\end{step}

\begin{step}
\rm $\lim_{n\rightarrow\infty } \nabla f(\theta_{n} ) = 0$
and the convergence of $f(\theta_{n} )$ are demonstrated 
(Lemmas \ref{lemma1.3}, \ref{lemma1.3'}).  
The proof is based on Lojasiewicz inequaltiy (\ref{1.3*})
(which is a consequence of Assumption \ref{a1.3}),  
Lemma \ref{lemma1.2} (relations (\ref{l1.2.3*}), (\ref{l1.2.5*}))
and standard stochastic approximation arguments. 
This result is used later at Steps 3, 4. 
\end{step}

\begin{step}
\rm The asymptotic behavior of $\{u(\theta_{n} ) \}_{n\geq 0}$, 
$\{v(\theta_{n} ) \}_{n\geq 0}$ is studied (Lemma \ref{lemma1.4}; 
$u(\cdot )$, $v(\cdot )$ are defined in (\ref{1.101*})).  
The obtained results crucially rely on 
Lojasiewicz inequality (\ref{1.3*}) and Steps 1, 2 
(Lemmas \ref{lemma1.2}, \ref{lemma1.3'}). 
The results are a corner-stone of the analysis 
carried out at Steps 4, 5, 6. 
\end{step}

\begin{step}
\rm $\liminf_{n\rightarrow\infty } \gamma_{n}^{\hat{p} } (f(\theta_{n} ) - \hat{f} ) > -\infty$
is demonstrated (Lemma \ref{lemma1.5}; 
$\hat{p}$ is defined in (\ref{1.105*})).  
The idea of the proof can be described as follows. 
If the previous relation is not true, 
then there exists 
a sufficiently large integer $n_{0} > 0$ such that 
$u(\theta_{n_{0} } ) < 0$ and 
\begin{align}\label{0.1}
	\hat{M} 
	\max_{n \leq k < a(n,1) }
	\left\|
	\sum_{i=n}^{k-1} \alpha_{i} \xi_{i} 
	\right\|^{\hat{\mu}}
	\leq 
	|u(\theta_{n_{0} } ) |
\end{align}
for $n\geq n_{0}$ 
(notice that 
$\max_{n\leq k < a(n,1) } 
\!\big\|\!\sum_{i=n}^{k} \alpha_{i} \xi_{i} \big\|^{\hat{\mu} }$ 
$= O(\gamma_{n} ^{-\hat{\mu} r } ) = O(\gamma_{n}^{-\hat{p} } )$ 
follows from 
Lemma \ref{lemma1.1}; 
$\hat{\mu}$, $\hat{M}$, $u(\cdot )$ are defined in 
(\ref{1.103*}), (\ref{1.101*})). 
Define sequence $\{n_{k} \}_{k\geq 0}$ recursively by 
$n_{k+1} = a(n_{k}, 1 )$ for $k\geq 0$.   
Let us show by induction that 
$u(\theta_{n_{k} } ) \leq u(\theta_{n_{0} } )$
for each $k\geq 0$. 
Obviously, this is true for $k=0$. 
Assume that  
$u(\theta_{n_{k} } ) \leq  u(\theta_{n_{0} } )$ for some $k\geq 0$. 
As $|u(\theta_{n_{k} } ) |\geq |u(\theta_{n_{0} } ) |$
(due to $u(\theta_{n_{0} } ) < 0$), 
the Lojasiewicz inequality (\ref{1.3*}) and (\ref{0.1}) imply 
$\big\|\sum_{i=n_{k} }^{n_{k+1} - 1} \alpha_{i} \xi_{i} \big\|
\leq \|\nabla f(\theta_{n_{k} } ) \|$. 
On the other side, Taylor formula yields 
\begin{align*}
	u(\theta_{n_{k+1} } )
	\approx &
	u(\theta_{n_{k} } ) 
	-
	(\nabla f(\theta_{n_{k} } ) )^{T} 
	\sum_{i=n_{k} }^{n_{k+1} - 1} \alpha_{i} (\nabla f(\theta_{i} ) + \xi_{i} )
	\\
	\approx &
	u(\theta_{n_{k} } )
	-
	(\nabla f(\theta_{n_{k} } ) )^{T}
	\left( 
	(\gamma_{n_{k+1} } - \gamma_{n_{k} } ) \nabla f(\theta_{n_{k} } )
	+
	\sum_{i=n_{k} }^{n_{k+1}-1} \alpha_{i} \xi_{i} 
	\right)
	\\
	\leq &
	u(\theta_{n_{k} } )
	-
	\|\nabla f(\theta_{n_{k} } ) \|
	\left( 
	\|\nabla f(\theta_{n_{k} } ) \|
	-
	\left\|
	\sum_{i=n_{k} }^{n_{k+1}-1} \alpha_{i} \xi_{i} 
	\right\| 
	\right)
	\\
	\leq &
	u(\theta_{n_{k} } )
\end{align*}
(notice that $\gamma_{n_{k+1} } - \gamma_{n_{k} } \approx 1$). 
Hence, $u(\theta_{n_{k+1} } ) \leq u(\theta_{n_{k} } )$. 
Then, by mathematical induction, we conclude 
$u(\theta_{n_{k+1} } ) \leq u(\theta_{n_{0} } )$ for any $k\geq 0$. 
However, this is not possible as $u(\theta_{n_{0} } ) < 0$
and 
$\lim_{n\rightarrow \infty } u(\theta_{n} ) = 0$
(due to Lemma \ref{lemma1.3'}). 
\end{step}

\begin{step}
\rm $\liminf_{n\rightarrow\infty } \gamma_{n}^{\hat{p} } (f(\theta_{n} ) - \hat{f} ) < \infty$
is proved (Lemma \ref{lemma1.6}). 
The idea of the proof can be summarized as follows. 
If the previous relation is not satisfied, 
then 
$\lim_{n\rightarrow \infty } \gamma_{n}^{-1} v(\theta_{n} ) = 0$
and there exists a sufficiently large integer $n_{0} > 0$
such that $u(\theta_{n} ) > 0$
and 
\begin{align*}
	2^{\hat{\mu } } \hat{M} 
	\max_{n\leq k < a(n,1) }
	\left\|
	\sum_{i=n}^{k-1} \alpha_{i} \xi_{i} 
	\right\|^{\hat{\mu} } 
	\leq 
	u(\theta_{n} )
\end{align*}
for $n\geq n_{0}$
(notice again that due to Lemma \ref{1.1}, 
$\max_{n\leq k < a(n,1) } 
\!\big\|\!\sum_{i=n}^{k} \alpha_{i} \xi_{i} \big\|^{\hat{\mu} } 
= O(\gamma_{n}^{-\hat{p} } )$; 
$u(\cdot )$, $v(\cdot )$ are defined in (\ref{1.101*})). 
Let $\{n_{k} \}_{k\geq 0}$ be defined in the same way as in Step 4. 
Then, the Lojasiewicz inequality (\ref{1.3*}) yields 
$\big\|\sum_{i=n_{k} }^{n_{k+1}-1} \alpha_{i} \xi_{i} \big\|
\leq \|\nabla f(\theta_{n_{k} } ) \|/2$ 
for $k\geq 0$. 
The same inequality also implies 
\begin{align*}
	\|\nabla f(\theta_{n_{k} } ) \|^{2} 
	\geq 
	\hat{M}^{-2/\hat{\mu} } 
	(\hat{f} - f(\theta_{n_{k} } ) )^{2/\hat{\mu} }
	\geq 
	2\hat{p} \hat{L} 
	(u(\theta_{n_{k} } ) )^{1+1/\hat{p} }
\end{align*} 
for $k\geq 0$, where $\hat{L} = 2^{-1} \hat{p}^{-1} \hat{M}^{-2/\hat{\mu} }$
(notice that 
$u(\theta_{n_{k} } ) \approx 0$ and 
$2/\hat{\mu} = 1+ 1/(\hat{\mu}\hat{r} ) \leq 1 + 1/\hat{p}$; 
$\hat{r}$ is defined in (\ref{1.105*})). 
Then, owing to Taylor formula, 
we have 
\begin{align*}
	v(\theta_{n_{k+1} } ) 
	\approx &
	v(\theta_{n_{k} } ) 
	+
	\frac{(\nabla f(\theta_{n_{k} } ) )^{T} }
	{\hat{p} (u(\theta_{n_{k} } ) )^{1+1/\hat{p} } }
	\left(
	(\gamma_{n_{k+1} } - \gamma_{n_{k} } ) \nabla f(\theta_{n_{k} } ) 
	+
	\sum_{i=n_{k} }^{n_{k+1}-1} 
	\alpha_{i} \xi_{i} 
	\right)
	\\
	\geq &
	v(\theta_{n_{k} } ) 
	+
	\frac{(\gamma_{n_{k+1} } - \gamma_{n_{k} } ) \|\nabla f(\theta_{n_{k} } ) \|^{2} }
	{2\hat{p} (u(\theta_{n_{k} } ) )^{1+1/\hat{p} } }
	\\
	&
	+
	\frac{\|\nabla f(\theta_{n_{k} } ) \| }
	{\hat{p} (u(\theta_{n_{k} } ) )^{1+1/\hat{p} } }
	\left(
	\frac{\|\nabla f(\theta_{n_{k} } ) \| }{2}
	-
	\left\|
	\sum_{i=n_{k} }^{n_{k+1}-1} 
	\alpha_{i} \xi_{i} 
	\right\|
	\right)
	\\
	\geq &
	v(\theta_{n_{k} } ) 
	+
	\hat{L} (\gamma_{n_{k+1} } - \gamma_{n_{k} } )
\end{align*}
for $k\geq 0$. 
Therefore, 
$\liminf_{k\rightarrow \infty } \gamma_{n_{k} }^{-1} v(\theta_{n_{k} } ) \geq \hat{L} > 0$. 
However, this is not possible due to  
\linebreak 
$\lim_{n\rightarrow\infty } \gamma_{n}^{-1} v(\theta_{n} ) = 0$. 
\end{step} 

\begin{step}
\rm  
$\limsup_{n\rightarrow\infty } \gamma_{n}^{\hat{p} } (f(\theta_{n} ) - \hat{f} ) < \infty$
is proved (Lemma \ref{lemma1.7}). 
The idea of the proof can be described as follows. 
Let $\hat{L}$ have the same meaning as in Step 5. 
If the previous relation is not satisfied, 
then, owing to the results of Step 5, 
there exist sufficiently large integer $m_{0}$ 
and sufficiently small real number $t\in (0,1)$ 
with the following properties: 
$(1/\hat{L} )^{\hat{p} } 
< \gamma_{m_{0} }^{\hat{p} } u(\theta_{m_{0} } ) 
\leq \gamma_{a(m_{0},t) }^{\hat{p} } u(\theta_{a(m_{0}, t ) } )$
and 
\begin{align*}
	(2/t)^{\hat{\mu } } \hat{M} 
	\left\|
	\sum_{i=m_{0} }^{a(m_{0}, t ) -1} \alpha_{i} \xi_{i} 
	\right\|^{\hat{\mu} } 
	\leq 
	u(\theta_{m_{0} } ) 
\end{align*}
(notice again that due to Lemma \ref{1.1}, 
$\max_{n\leq k < a(n,1) } 
\!\big\|\!\sum_{i=n}^{k} \alpha_{i} \xi_{i} \big\|^{\hat{\mu} } 
= O(\gamma_{n}^{-\hat{p} } )$). 
Let $n_{0} = a(m_{0}, t )$. 
Consequently, 
$\gamma_{n_{0} }^{-1} v(\theta_{n_{0} } ) 
\leq \gamma_{m_{0} }^{-1} v(\theta_{m_{0} } ) 
< \hat{L}$, 
while 
the Lojasiewicz inequality (\ref{1.3*}) implies  
$\big\|\sum_{i=m_{0} }^{n_{0}-1} \alpha_{i} \xi_{i} \big\|
\leq (t/2) \|\nabla f(\theta_{m_{0} } ) \|$ 
and 
\begin{align*}
	\|\nabla f(\theta_{m_{0} } ) \|^{2} 
	\geq 
	\hat{M}^{-2/\hat{\mu} } 
	(\hat{f} - f(\theta_{m_{0} } ) )^{2/\hat{\mu} }
	\geq 
	2\hat{p} \hat{L} 
	(u(\theta_{m_{0} } ) )^{1+1/\hat{p} }. 
\end{align*} 
Combining this with Taylor formula, we get 
\begin{align*}
	v(\theta_{n_{0} } ) 
	\approx &
	v(\theta_{m_{0} } ) 
	+
	\frac{(\nabla f(\theta_{m_{0} } ) )^{T} }
	{\hat{p} (u(\theta_{m_{0} } ) )^{1+1/\hat{p} } }
	\left(
	(\gamma_{n_{0} } - \gamma_{m_{0} } ) \nabla f(\theta_{m_{0} } ) 
	+
	\sum_{i=m_{0} }^{n_{0}-1} 
	\alpha_{i} \xi_{i} 
	\right)
	\\
	\geq &
	v(\theta_{n_{k} } ) 
	+
	\frac{(\gamma_{n_{0} } - \gamma_{m_{0} } ) \|\nabla f(\theta_{m_{0} } ) \|^{2} }
	{2\hat{p} (u(\theta_{m_{0} } ) )^{1+1/\hat{p} } }
	\\
	&
	+
	\frac{\|\nabla f(\theta_{m_{0} } ) \| }
	{\hat{p} (u(\theta_{m_{0} } ) )^{1+1/\hat{p} } }
	\left(
	\frac{t\|\nabla f(\theta_{m_{0} } ) \| }{2}
	-
	\left\|
	\sum_{i=m_{0} }^{n_{0}-1} 
	\alpha_{i} \xi_{i} 
	\right\|
	\right)
	\\
	\geq &
	v(\theta_{n_{k} } ) 
	+
	\hat{L} (\gamma_{n_{0} } - \gamma_{m_{0} } )
\end{align*}
(notice that $\gamma_{n_{0} } - \gamma_{m_{0} } \approx t$). 
Therefore, 
\begin{align*}
	\gamma_{n_{0} }^{-1} v(\theta_{n_{0} } ) 
	\geq 
	\gamma_{m_{0} }^{-1} v(\theta_{m_{0} } ) 
	+
	(1 - \gamma_{m_{0} }/ \gamma_{n_{0} } )
	(\hat{L} - \gamma_{m_{0} }^{-1} v(\theta_{m_{0} } ) )
	> 
	\gamma_{m_{0} }^{-1} v(\theta_{m_{0} } ).  
\end{align*}
However, this is impossible as 
$\gamma_{n_{0} }^{-1} v(\theta_{n_{0} } ) 
\leq 
\gamma_{m_{0} }^{-1} v(\theta_{m_{0} } )$. 
\end{step} 

\begin{step}
\rm $\|\nabla f(\theta_{n} ) \|^{2} = O(\gamma_{n}^{-\hat{p} } )$ is demonstrated 
(Lemma \ref{lemma1.5}).  
The proof is based on the following idea. 
Due to Taylor formula, we have 
\begin{align*}
	\|\nabla f(\theta_{n} ) \|^{2}
	\approx &
	\frac{u(\theta_{n} ) - u(\theta_{a(n,1) } ) }{\gamma_{a(n,1) } - \gamma_{n} }
	-
	\frac{(\nabla f(\theta_{n} ) )^{T} }{\gamma_{a(n,1) } - \gamma_{n} }
	\sum_{i=n}^{a(n,1)-1} \alpha_{i} \xi_{i}
	\\
	\leq &
	|u(\theta_{a(n,1) } ) |
	+
	|u(\theta_{n} ) |
	+
	\frac{\|\nabla f(\theta_{n} ) \|^{2} }{2}
	+
	\frac{1}{2} 
	\left\|
	\sum_{i=n}^{a(n,1)-1} \alpha_{i} \xi_{i} 
	\right\|^{2} 
\end{align*}
for all sufficiently large $n$
(notice that $\gamma_{a(n,1)} - \gamma_{n} \approx 1$). 
Consequently, 
\begin{align*}
	\|\nabla f(\theta_{n} ) \|^{2}	
	\leq 
	2|u(\theta_{a(n,1) } ) |
	+
	2|u(\theta_{n} ) |	
	+
	\left\|
	\sum_{i=n}^{a(n,1)-1} \alpha_{i} \xi_{i} 
	\right\|^{2} 
\end{align*}
for the same $n$. 
Then, $\|\nabla f(\theta_{n} ) \|^{2} = O(\gamma_{n}^{-\hat{p} } )$
directly follows from the results of Steps 4 and 6
(also notice that 
$\max_{n\leq k < a(n,1) } \big\|\sum_{i=n}^{k} \alpha_{i} \xi_{i} \big\|^{2} 
= O(\gamma_{n} ^{-2r} ) 
= O(\gamma_{n}^{-\hat{p} } )$ follows from 
Lemma \ref{lemma1.1}).  
\end{step} 

\begin{step} 
\rm $\max_{k\geq n} \|\theta_{k} - \theta_{n} \| = O(\gamma_{n}^{-\hat{q} } )$ 
is proved (Lemmas \ref{lemma1.21}, \ref{lemma1.23}; 
$\hat{q}$ is defined in (\ref{1.103*})). 
The idea of the proof can be summarized as follows. 
Let $\{n_{k} \}_{k\geq 0}$ be the sequence recursively defined by 
$n_{0} = 0$ and $n_{k+1} = a(n_{k}, 1 )$ for $k\geq 0$. 
Owing to Taylor formula, we have 
\begin{align}\label{s8.3} 
	u(\theta_{k} ) 
	-
	u(\theta_{n} )
	\approx &
	- 
	(\gamma_{k} - \gamma_{n} ) \|\nabla f(\theta_{n_{k} } ) \|^{2} 
	-
	(\nabla f(\theta_{n_{k} } ) )^{T}
	\sum_{i=n}^{k-1} \alpha_{i} \xi_{i} 	
\end{align}
for $n\leq k \leq a(n,1)$ and all sufficiently large $n$. 
We also have 
\begin{align}\label{s8.1}
	\|\theta_{k} - \theta_{n} \| 
	\approx &
	\left\|
	(\gamma_{k} - \gamma_{n} ) \nabla f(\theta_{n} ) 
	+
	\sum_{i=n}^{k-1} \alpha_{i} \xi_{i} 
	\right\|
	\leq 
	(\gamma_{k} - \gamma_{n} ) 
	\|\nabla f(\theta_{n} ) \| 
	+
	\left\|
	\sum_{i=n}^{k-1} \alpha_{i} \xi_{i} 
	\right\|	
\end{align}
for the same $n,k$. 
Combining (\ref{s8.3}), (\ref{s8.1}), we get 
\begin{align}
	&\label{s8.7}
	\begin{aligned}[b]
	\|\theta_{k} - \theta_{n} \|
	\leq & 
	\frac{1}{\|\nabla f(\theta_{n} ) \| }
	\left(
	u(\theta_{n} )
	-
	u(\theta_{k} ) 
	-
	(\nabla f(\theta_{n} ) )^{T} 
	\sum_{i=n}^{k-1} \alpha_{i} \xi_{i} 
	\right)
	+ 
	\left\|
	\sum_{i=n}^{k-1} \alpha_{i} \xi_{i} 
	\right\|
	\\
	\leq &
	\frac{u(\theta_{n} ) - u(\theta_{k} ) }
	{\|\nabla f(\theta_{n} ) \| }
	+
	2
	\left\|
	\sum_{i=n}^{k-1} \alpha_{i} \xi_{i} 
	\right\|
	\end{aligned}
\end{align}
for $n\leq k \leq a(n,1)$ and all sufficiently large $n$. 
Similarly, using the results of Step 7 and (\ref{s8.1}), we obtain    
\begin{align}\label{s8.101}
	\max_{n \leq k \leq a(n,1) }
	\|\theta_{k} - \theta_{n} \| 
	=
	O(\gamma_{n}^{-\hat{p}/2 } )
	=
	o(\gamma_{n}^{-(\hat{q}+1) } )
\end{align}
(notice that 
$\hat{q} < \hat{p}/2$, $\hat{q} + 1 \leq r$, 
$\gamma_{a(n,1)} - \gamma_{n} \approx 1$
and that 
$\max_{n\leq k < a(n,1)} \big\|\sum_{i=n}^{k} \alpha_{i} \xi_{i} \big\| 
= O(\gamma_{n} ^{-r} )$ 
follows from 
Lemma \ref{lemma1.1}). 
On the other side, if $\|\nabla f(\theta_{n} ) \| \geq \gamma_{n}^{-(\hat{q}+1)}$, 
(\ref{s8.7}) yields 
\begin{align}\label{s8.9} 
	\|\theta_{k} - \theta_{n} \| 
	\leq &
	\gamma_{n}^{\hat{q}+1}
	(u(\theta_{n} ) - u(\theta_{k} ) )
	+
	2 
	\left\|
	\sum_{i=n}^{k-1} \alpha_{i} \xi_{i} 
	\right\|
	\nonumber\\
	\leq &
	\hat{L}_{1} 
	\left(
	\gamma_{n}^{\hat{q}+1}
	(u(\theta_{n} ) - u(\theta_{k} ) )
	+
	\gamma_{n}^{-(\hat{q}+1)}
	\right)
\end{align}	
for $n\leq k \leq a(n,1)$, all sufficiently large $n$ and some $\hat{L}_{1} \in [1,\infty )$. 
If $\|\nabla f(\theta_{n} ) \| \leq \gamma_{n}^{-(\hat{q}+1)}$, 
a similar relation results from 
(\ref{s8.3}), (\ref{s8.1}): 
\begin{align}\label{s8.11}
	\|\theta_{k} - \theta_{n} \|
	\leq &
	\|\nabla f(\theta_{n} ) \|
	+
	\left\|\sum_{i=n}^{k-1} \alpha_{i} \xi_{i} \right\| 
	+
	\gamma_{n}^{\hat{q}+1 }
	(u(\theta_{n} ) - u(\theta_{k} ) )
	+ 
	\gamma_{n}^{\hat{q}+1 }
	|u(\theta_{n} ) - u(\theta_{k} ) |
	\nonumber\\
	\leq &
	\hat{L}_{2} 
	\left(
	\gamma_{n}^{\hat{q}+1 }
	(u(\theta_{n} ) - u(\theta_{k} ) )	
	+
	\gamma_{n}^{-(\hat{q}+1) }
	\right)
\end{align}
for the same $n,k$ and some $\hat{L}_{2} \in [1,\infty )$. 
Combining (\ref{s8.9}), (\ref{s8.11}), we get 
\begin{align}\label{s8.103}
	\|\theta_{n_{j} } - \theta_{n_{k} } \| 
	\leq 
	\sum_{i=k}^{j-1} \|\theta_{n_{i+1} } - \theta_{n_{i} } \| 
	\leq &
	\hat{L} 
	\sum_{i=k}^{\infty } \gamma_{n_{i} }^{-(\hat{q}+1) } 
	+
	\hat{L}
	\sum_{i=k+1}^{\infty } (\gamma_{n_{i} }^{\hat{q}+1 } - \gamma_{n_{i-1} }^{\hat{q}+1 } )
	|u(\theta_{n_{i} } ) |
	\nonumber\\
	&
	+
	\hat{L} \gamma_{n_{k} }^{\hat{q} + 1} |u(\theta_{n_{k} } ) |
	+
	\hat{L} \gamma_{n_{j} }^{\hat{q} + 1} |u(\theta_{n_{j} } ) |
\end{align}
for $j\geq k$ and all sufficiently large $k$, 
where $\hat{L} = \max\{\hat{L}_{1}, \hat{L}_{2} \}$. 
As $u(\theta_{n} ) = O(\gamma_{n}^{-\hat{p} } )$ 
(due to the results of Steps 4, 6) and 
\begin{align*}
	\sum_{i=k}^{\infty } \gamma_{n_{i} }^{-(\hat{q}+1) }
	=
	O(\gamma_{n_{k} }^{-\hat{q} } ), 
	\;\;\;\;\; 
	\sum_{i=k+1}^{\infty } 
	\gamma_{n_{i} }^{-\hat{p} } 
	(\gamma_{n_{i} }^{\hat{q}+1} - \gamma_{n_{i-1} }^{\hat{q}+1} ) 
	=
	O(\gamma_{n_{k} }^{-\hat{q} } )
\end{align*}
(see (\ref{l1.23.201}), (\ref{l1.23.701})), 
we conclude from (\ref{s8.101}), (\ref{s8.103}) 
that $\max_{k\geq n} \|\theta_{k} - \theta_{n} \| = O(\gamma_{n}^{-\hat{q} } )$. 
\end{step} 

\begin{step} 
\rm 
Theorems \ref{theorem1.1} and \ref{theorem1.2} are proved.  
The convergence and convergence rate of $\{\theta_{n} \}_{n\geq 0}$
directly follow from the results of Step 8, 
while the convergence rates of 
$\{f(\theta_{n} ) \}_{n\geq 0}$, 
$\{\nabla f(\theta_{n} ) \}_{n\geq 0}$ 
are immediate consequences of Steps 4 -- 7. 
\end{step} 

\section{Proof of Theorems \ref{theorem1.1} and \ref{theorem1.2}} \label{section1*}

In this section, the following notation is used. 
$\Lambda$ is the event defined as 
\begin{align*}
	\Lambda = \left\{\sup_{n\geq 0} \|\theta_{n} \| < \infty \right\}. 
\end{align*}	
For $k > n \geq 1$, let 
$\zeta_{n,n}=\zeta'_{n,n}=\zeta''_{n,n}=0$ and 
\begin{align*}	
	\zeta'_{n,k}
	=
	\sum_{i=n}^{k-1} \alpha_{i} \xi_{i}, 
	\;\;\; 
	\zeta''_{n,k}
	=
	\sum_{i=n}^{k-1} 
	\alpha_{i} (\nabla f(\theta_{i} ) - \nabla f(\theta_{n} ) ), 
\end{align*}
while 
$
	\zeta_{n,k} 
	= 
	\zeta'_{n,k} + \zeta''_{n,k}.  
$
For the same $k,n$, let $\phi_{n,n}=\phi'_{n,n}=\phi''_{n,n}=0$ and 
\begin{align*}
	\phi'_{n,k}
	=
	(\nabla f(\theta_{n} ) )^{T} 
	\zeta_{n,k}, 
	\;\;\; 
	\phi''_{n,k} 
	=
	-
	\int_{0}^{1} 
	(\nabla f(\theta_{n} + s (\theta_{k} - \theta_{n} ) ) - \nabla f(\theta_{n} ) )^{T}
	(\theta_{k} - \theta_{n} ) ds, 
\end{align*}
while 
$
	\phi_{n,k} = \phi'_{n,k} + \phi''_{n,k}.
$
Then, it is straightforward to show 
\begin{align} 
	&\label{1.1'*}
	\begin{aligned}[b]
	\theta_{k} - 	\theta_{n} 
	=&
	-
	\sum_{i=n}^{k-1} \alpha_{i} \nabla f(\theta_{i} ) 
	-
	\zeta'_{n,k}
	=
	-
	(\gamma_{k} - \gamma_{n} ) \nabla f(\theta_{n} ) 
	-
	\zeta_{n,k}, 
	\end{aligned} 
	\\ 
	&\label{1.1*}
	f(\theta_{k} ) - f(\theta_{n} )
	= 
	-
	(\gamma_{k} - \gamma_{n} ) 
	\|\nabla f(\theta_{n} ) \|^{2} 
	-
	\phi_{n,k}			
\end{align}
for $0 \leq n \leq k$. 

In this section, 
we also rely on the following notation. 
For a compact set $Q \subset \mathbb{R}^{d_{\theta } }$, 
$C_{Q}$ stands for an upper bound of 
$\|\nabla f(\cdot )\|$ on $Q$ 
and for a Lipschitz constant of $\nabla f(\cdot )$ on the same set. 
$\hat{A}$ is the set of accumulation points of 
$\{\theta_{n} \}_{n\geq 0}$, 
while 
\begin{align*}
	\hat{f} 
	=
	\liminf_{n\rightarrow \infty } 
	f(\theta_{n} ). 
\end{align*}
$\hat{B}$ and $\hat{Q}$ are random sets defined by 
\begin{align*}
	\hat{B}
	=
	\bigcup_{\theta \in \hat{A} } 
	\left\{ 
	\theta' \in \mathbb{R}^{d_{\theta } }: 
	\|\theta' - \theta \| \leq \delta_{\theta }/2
	\right\}, 
	\;\;\;\;\; 
	\hat{Q} 
	=
	\text{cl}(\hat{B} )
\end{align*}
on event $\Lambda$, 
and by 
\begin{align*}
	\hat{B} 
	=
	\hat{A}, 
	\;\;\;\;\;
	\hat{Q} 
	=
	\hat{A}
\end{align*}
outside $\Lambda$
($\delta_{\theta }$ is specified in Remark \ref{remark1.1}). 
Overriding the definition of $\hat{\mu}$, $\hat{p}$, $\hat{q}$, $\hat{r}$, 
in Theorem \ref{theorem1.2}, 
we define random quantities 
$\hat{\delta}$, $\hat{\mu}$, $\hat{p}$, $\hat{q}$, $\hat{r}$, $\hat{C}$, $\hat{M}$ 
as 
\begin{align}
	&\label{1.103*}
	\hat{\delta} 
	=
	\delta_{\hat{Q}, \hat{f} }, 
	\;\;\;\;\;  
	\hat{\mu} 
	=
	\mu_{\hat{Q}, \hat{f} }, 
	\;\;\;\;\;  
	\hat{C} 
	=
	C_{\hat{Q} }, 
	\;\;\;\;\;  
	\hat{M}
	=
	\hat{M}_{\hat{Q}, \hat{f} }, 
\end{align}
\begin{align}
	&\label{1.105*}
	\hat{r} 
	=
	\begin{cases}
	1/(2-\hat{\mu} ), 
	&\text{ if } \hat{\mu} < 2
	\\
	\infty, 
	&\text{ if } \hat{\mu} = 2
	\end{cases}, 
	\;\;\;
	\hat{p}
	=
	\hat{\mu} \min\{r, \hat{r} \}, 
	\;\;\; 
	\hat{q}
	=
	\min\{r,\hat{r} \} - 1
\end{align}
on $\Lambda$
($\delta_{Q,a}$, $\mu_{Q,a}$, $M_{Q,a}$ are specified in Assumption \ref{a1.3}), and as 
\begin{align*}
	\hat{\delta} = 1, 
	\;\;\; 
	\hat{\mu} = 2, 
	\;\;\; 
	\hat{C} = 1, 
	\;\;\; 
	\hat{M} =1, 
	\;\;\; 
	\hat{r} = \infty, 
	\;\;\; 
	\hat{p} = 2r, 
	\;\;\; 
	\hat{q}
	=
	r - 1
\end{align*}
outside $\Lambda$
(later, when Theorem \ref{theorem1.1} is proved, 
it will be clear that $\hat{\mu}$, $\hat{p}$, $\hat{r}$
specified here coincide with 
$\hat{\mu}$, $\hat{p}$, $\hat{r}$ defined in Theorem \ref{theorem1.2}).  
$u(\cdot )$, $v(\cdot )$ are functions defined by 
\begin{align}\label{1.101*}
	u(\theta ) 
	=
	f(\theta ) - \hat{f}, 
	\;\;\;\;\;  
	v(\theta ) 
	=
	\begin{cases} 
	(f(\theta ) - \hat{f} )^{-1/\hat{p} }, 
	&\text{if } f(\theta ) > \hat{f}
	\\
	0, 
	&\text{otherwise} 
	\end{cases}  
\end{align}
for 
$\theta \in \mathbb{R}^{d_{\theta } }$. 
For $\varepsilon \in (0, \infty )$, 
$\varphi_{\varepsilon }(\xi)$ and $\phi_{\varepsilon }(\xi )$ are random quantities defined as 
\begin{align}\label{1.201*}
	\varphi_{\varepsilon }(\xi) = \varphi(\xi) + \varepsilon, 
	\;\;\;\;\; 
	\phi_{\varepsilon }(\xi )
	=
	\begin{cases}
	\varphi_{\varepsilon }(\xi ), 
	&\text{ if } r\leq \hat{r}
	\\
	\left(\varphi_{\varepsilon }(\xi ) \right)^{\hat{\mu} - 1}, 
	&\text{ if } r>\hat{r}
	\end{cases} 
\end{align}
($\xi$ is specified in Assumption \ref{a1.2}, 
while $\varphi(\xi)$ is defined in the statement of Theorem \ref{theorem1.2}). 

\begin{remark}
On event $\Lambda$, 
$\hat{Q}$ is compact and satisfies 
$\hat{A} \subset \text{\rm int} \hat{Q}$. 
Thus, 
$\hat{\delta}$, $\hat{p}$, $\hat{r}$, $\hat{C}$, $\hat{M}$, 
$v(\cdot )$
are well-defined on $\Lambda$
(what happens with these quantities outside $\Lambda$
does not affect the results presented in this section). 
Then, Assumption \ref{a1.3} implies 
\begin{align}\label{1.3*}
	|f(\theta ) - \hat{f} |
	\leq 
	\hat{M} \|\nabla f(\theta ) \|^{\hat{\mu} }
\end{align}
on $\Lambda$ for all $\theta \in \hat{Q}$ satisfying 
$|f(\theta ) - \hat{f} | \leq \hat{\delta }$. 
\end{remark}

\begin{remark}
Regarding the notation, the following note is also in order: 
Diacritic $\:\tilde{}\:$ is used for a locally defined quantity, i.e., 
for a quantity whose definition holds only in the proof where such a quantity 
appears. 
\end{remark}

\begin{lemma} \label{lemma1.1}
Let Assumptions \ref{a1.1} and \ref{a1.2} hold. 
Then, there exists an event $N_{0} \in {\cal F}$
such that 
$P(N_{0} ) = 0$ and 
\begin{align*} 
	\limsup_{n\rightarrow \infty } \gamma_{n}^{r} 
	\max_{n\leq k \leq a(n,1) }
	\|\zeta'_{n,k} \|
	\leq 
	\xi
	< 
	\infty 
\end{align*}
on $\Lambda\setminus N_{0}$. 
\end{lemma}

\begin{sproof}
It is straightforward to verify 
\begin{align*}
	\zeta'_{n,k}
	=
	\sum_{i=n}^{k-1} 
	(\gamma_{i}^{-r} - \gamma_{i+1}^{-r} ) 
	\left(
	\sum_{j=n}^{i} \alpha_{j} \gamma_{j}^{r} \xi_{j} 
	\right)
	+
	\gamma_{k}^{-r} \sum_{i=n}^{k-1} \alpha_{i} \gamma_{i}^{r} \xi_{i} 
\end{align*}
for $0\leq n < k$. 
Consequently, 
\begin{align*}
	\|\zeta'_{n,k} \|
	\leq &
	\left(
	\gamma_{k}^{-r} 
	+ 
	\sum_{i=n}^{k-1} (\gamma_{i}^{-r} - \gamma_{i+1}^{-r} )
	\right)
	\max_{n\leq j < a(n,1) } 
	\left\|
	\sum_{i=n}^{j} \alpha_{i} \gamma_{i}^{r} \xi_{i} 
	\right\|
	= 
	\gamma_{n}^{-r} 
	\max_{n\leq j < a(n,1) } 
	\left\|
	\sum_{i=n}^{j} \alpha_{i} \gamma_{i}^{r} \xi_{i} 
	\right\|
\end{align*}
for $0 \leq n \leq k \leq a(n,1)$. 
Thus, 
\begin{align*}
	\gamma_{n}^{r} 
	\|\zeta'_{n,k} \|
	\leq &
	\max_{n\leq j < a(n,1) } 
	\left\|
	\sum_{i=n}^{j} \alpha_{i} \gamma_{i}^{r} \xi_{i} 
	\right\|
\end{align*}
for $0 \leq n \leq k \leq a(n,1)$. 
Then, the lemma's assertion directly follows from Assumption \ref{a1.2}. 
\end{sproof} 

\begin{lemma} \label{lemma1.2}
Suppose that Assumptions \ref{a1.1} -- \ref{a1.3} hold. 
Then, 
there exist random quantities  
$\hat{C}_{1}$, $\hat{t}$
(which are deterministic functions of 
$\hat{C}$) 
and for any real number $\varepsilon \in (0,\infty )$, 
there exists  
a non-negative integer-valued 
random quantity  
$\tau_{1,\varepsilon}$ such that the following is true:
$1\leq \hat{C}_{1} < \infty$, 
$0 < \hat{t} < 1$, 
$0\leq \tau_{1,\varepsilon } < \infty$
everywhere and 
\begin{align}
	& \label{l1.2.1*} 
	\max_{n\leq k \leq a(n,\hat{t} ) }
	\|\theta_{k} - \theta_{n} \|
	\leq 
	\hat{C}_{1} 
	\left(
	\|\nabla f(\theta_{n} ) \|
	+
	\gamma_{n}^{-r} (\xi + \varepsilon ) 
	\right), 
	\\
	& \label{l1.2.3*} 	
	\begin{aligned}[b]
		\max_{n\leq k \leq a(n,\hat{t} ) }
		(f(\theta_{k} ) - f(\theta_{n} ) )
		\leq &	
		\hat{C}_{1} 
		\left( 
		\gamma_{n}^{-r} 
		\|\nabla f(\theta_{n} ) \| 
		(\xi + \varepsilon )  	
		+ 
		\gamma_{n}^{-2r} 
		(\xi + \varepsilon )^{2}
		\right), 
	\end{aligned}
	\\
	& \label{l1.2.5*}
	\begin{aligned}[b]
		f(\theta_{a(n,\hat{t} ) } ) - f(\theta_{n} ) 
		+
		\hat{t} \|\nabla f(\theta_{n} ) \|^{2}/2
		\leq 
		\hat{C}_{1} 
		\left( 
		\gamma_{n}^{-r} 
		\|\nabla f(\theta_{n} ) \| (\xi + \varepsilon )
		+  
		\gamma_{n}^{-2r} 
		(\xi + \varepsilon )^{2} 
		\right) 
	\end{aligned} 	
	\\
	& \label{l1.2.7*}
	\begin{aligned}[b]
		&
		2\left( f(\theta_{a(n,\hat{t} ) } ) - f(\theta_{n} ) \right) 
		+
		\hat{t} \|\nabla f(\theta_{n} ) \|^{2}/2
		+
		\|\nabla f(\theta_{n} ) \| \|\theta_{a(n,\hat{t} ) } - \theta_{n} \|
		\\
		&
		\;\;\;\;\; 
		\leq 
		\hat{C}_{1} 
		\left( 
		\gamma_{n}^{-r} 
		\|\nabla f(\theta_{n} ) \| (\xi + \varepsilon )
		+  
		\gamma_{n}^{-2r} 
		(\xi + \varepsilon )^{2} 
		\right) 
	\end{aligned} 	
\end{align}
on $\Lambda\setminus N_{0}$ 
for  
$n>\tau_{1,\varepsilon }$. 
\end{lemma}

\begin{sproof}
Let 
$\tilde{C}_{1} = 2 \hat{C} \exp(\hat{C} )$, 
$\tilde{C}_{2} = 2 \hat{C} \tilde{C}_{1}$, 
$\tilde{C}_{3} = 2 \hat{C} \tilde{C}_{1}^{2} + \hat{C}_{2}$, 
$\tilde{C}_{4} = \tilde{C}_{2} + 2 \tilde{C}_{3}$, 
while 
$\hat{C}_{1} = \tilde{C}_{4}$, 
$\hat{t} = 1/(4 \tilde{C}_{4} )$. 
Moreover, let $\varepsilon \in (0,\infty )$ be an arbitrary real number.  
Then, owing to Lemma \ref{lemma1.1} and the fact that 
$\gamma_{a(n,\hat{t} ) } - \gamma_{n} = 
\hat{t} + O(\alpha_{a(n,\hat{t} ) } )$ for $n\rightarrow \infty$, 
it is possible to construct a non-negative integer-valued random quantity 
$\tau_{1,\varepsilon}$ such that 
$0 \leq \tau_{1,\varepsilon} < \infty$ everywhere
and such that 
$\theta_{n} \in \hat{Q}$, 
\begin{align} 
	& \label{l1.2.3}
	\gamma_{a(n,\hat{t} )} - \gamma_{n} 
	\geq 
	2\hat{t}/3, 
	\\
	& \label{l1.2.1}
	\max_{n\leq k \leq a(n,1) } 
	\|\zeta'_{n,k} \|
	\leq 
	\gamma_{n}^{-r} (\xi + \varepsilon ) 
\end{align}
on $\Lambda \setminus N_{0}$ for $n>\tau_{1,\varepsilon }$. 

Let $\omega$ be an arbitrary sample from $\Lambda\setminus N_{0}$
(notice that all formulas which follow in the proof correspond to this sample).
Since $\theta_{n} \in \hat{Q}$ for $n>\tau_{1,\varepsilon}$, 
(\ref{1.1'*}), (\ref{l1.2.1}) yield
\begin{align}\label{l1.2.1501}
	\|\nabla f(\theta_{k} ) \| 
	\leq &
	\|\nabla f(\theta_{n} ) \| 
	+
	\|\nabla f(\theta_{k} ) - \nabla f(\theta_{n} ) \|
	\nonumber\\
	\leq &
	\|\nabla f(\theta_{n} ) \| 
	+
	\hat{C} \|\theta_{k} - \theta_{n} \|
	\nonumber\\
	\leq &
	\|\nabla f(\theta_{n} ) \| 
	+
	\hat{C} 
	\sum_{i=n}^{k-1} \alpha_{i} \|\nabla f(\theta_{i} ) \| 
	+
	\hat{C} \|\zeta'_{n,k} \|
	\nonumber\\
	\leq &
	\|\nabla f(\theta_{n} ) \| 
	+
	\hat{C} \gamma_{n}^{-r} (\xi + \varepsilon) 
	+
	\hat{C} 
	\sum_{i=n}^{k-1} \alpha_{i} \|\nabla f(\theta_{i} ) \| 
\end{align}
for $\tau_{1,\varepsilon } < n \leq k \leq a(n,1)$.  
Then, Bellman-Gronwall inequality implies 
\begin{align}\label{l1.2.1505}
	\|\nabla f(\theta_{k} ) \|	
	\leq &
	\left(
	\|\nabla f(\theta_{n} ) \|
	+
	\hat{C} \gamma_{n}^{-r} (\xi + \varepsilon ) 
	\right)
	\exp\left(
	\hat{C} (\gamma_{k} - \gamma_{n} )
	\right) 
	\leq 
	\hat{C} \exp(\hat{C} )
	\left(
	\|\nabla f(\theta_{n} ) \|
	+
	\gamma_{n}^{-r} (\xi + \varepsilon ) 
	\right)
\end{align}
for $\tau_{1,\varepsilon } < n \leq k \leq a(n,1)$
(notice that $\gamma_{k} - \gamma_{n} \leq \gamma_{a(n,1)} - \gamma_{n} \leq 1$
when $n \leq k \leq a(n,1)$). 
Consequently, (\ref{l1.2.1}) gives 
\begin{align}\label{l1.2.5}
	\|\theta_{k} - \theta_{n} \|
	\leq &
	\sum_{i=n}^{k-1} 
	\alpha_{i} \|\nabla f(\theta_{i} ) \| 
	+
	\|\zeta'_{n,k} \|
	\nonumber\\
	\leq &
	\hat{C} \exp(\hat{C} )
	\left(
	\|\nabla f(\theta_{n} ) \|
	+
	\gamma_{n}^{-r} (\xi + \varepsilon ) 
	\right)
	(\gamma_{k} - \gamma_{n} )
	+ 
	\gamma_{n}^{-r} (\xi + \varepsilon ) 
	\nonumber\\
	\leq &
	\tilde{C}_{1} 
	\left(
	(\gamma_{k} - \gamma_{n} )	\|\nabla f(\theta_{n} ) \|
	+ 
	\gamma_{n}^{-r} (\xi + \varepsilon ) 
	\right)
\end{align}
for $\tau_{1,\varepsilon } < n \leq k \leq a(n,1)$.
Therefore, (\ref{l1.2.1}) yields 
\begin{align}	\label{l1.2.7}
	\|\zeta_{n,k} \|
	\leq &
	\|\zeta'_{n,k} \|
	+
	\hat{C} 
	\sum_{i=n}^{k-1} 
	\alpha_{i} \|\theta_{i} - \theta_{n} \| 
	\nonumber \\
	\leq &
	\gamma_{n}^{-r} (\xi + \varepsilon )
	+
	\hat{C} \tilde{C}_{1} 
	\left(
	(\gamma_{k} - \gamma_{n} ) \|\nabla f(\theta_{n} ) \| 
	+
	\gamma_{n}^{-r} (\xi + \varepsilon ) 
	\right)
	(\gamma_{k} - \gamma_{n} ) 
	\nonumber \\
	\leq &
	\tilde{C}_{2} 
	\left(
	(\gamma_{k} - \gamma_{n} )^{2} \|\nabla f(\theta_{n} ) \| 
	+
	\gamma_{n}^{-r} (\xi + \varepsilon ) 
	\right)
\end{align}
for $\tau_{1,\varepsilon } < n \leq k \leq a(n,1)$. 
Thus, 
\begin{align} 
	\label{l1.2.9}
	|\phi_{n,k} |
	\leq &
	\|\nabla f(\theta_{n} ) \| \|\zeta_{n,k} \| 
	+
	\hat{C} \|\theta_{k} - \theta_{n} \|^{2} 
	\nonumber \\
	\leq &
	\tilde{C}_{2} 
	\left(
	(\gamma_{k} - \gamma_{n} )^{2} 
	\|\nabla f(\theta_{n} ) \|^{2} 
	+
	\gamma_{n}^{-r} 
	\|\nabla f(\theta_{n} ) \| 
	(\xi + \varepsilon ) 
	\right)
	\nonumber\\
	&+
	\hat{C} \tilde{C}_{1}^{2} 
	\left( 
	(\gamma_{k} - \gamma_{n} ) \|\nabla f(\theta_{n} ) \| 
	+
	\gamma_{n}^{-r} (\xi + \varepsilon ) 
	\right)^{2} 
	\nonumber\\
	\leq &
	\tilde{C}_{3} 
	\left(
	(\gamma_{k} - \gamma_{n} )^{2} 
	\|\nabla f(\theta_{n} ) \|^{2} 
	+
	\gamma_{n}^{-r} 
	\|\nabla f(\theta_{n} ) \| 
	(\xi + \varepsilon ) 
	+
	\gamma_{n}^{-2r} (\xi + \varepsilon )^{2} 
	\right)
\end{align}
for $\tau_{1,\varepsilon } < n \leq k \leq a(n,1)$. 

Owing to (\ref{1.1*}), (\ref{l1.2.9}), we have 
\begin{align}\label{l1.2.21} 
	f(\theta_{k} )
	-
	f(\theta_{n} ) 
	\leq &
	- (\gamma_{k} - \gamma_{n} ) \|\nabla f(\theta_{n} ) \|^{2} 
	+
	|\phi_{n,k} | 
	\nonumber\\
	\leq &
	-
	\left(
	1 - \tilde{C}_{3} (\gamma_{k} - \gamma_{n} ) 
	\right)
	(\gamma_{k} - \gamma_{n} )
	\|\nabla f(\theta_{n} ) \|^{2} 
	\nonumber\\
	&+
	\tilde{C}_{3}
	\left(
	\gamma_{n}^{-r} \|\nabla f(\theta_{n} ) \| (\xi + \varepsilon ) 
	+
	\gamma_{n}^{-2r} (\xi + \varepsilon )^{2} 
	\right)
\end{align}
for $\tau_{1,\varepsilon } < n \leq k \leq a(n,1)$. 
Since 
\begin{align}\label{l1.2.23}
	\tilde{C}_{3} (\gamma_{k}-\gamma_{n} )	
	\leq 
	\tilde{C}_{4} (\gamma_{k}-\gamma_{n} )	
	\leq 
	\tilde{C}_{4} (\gamma_{a(n,\hat{t} ) }-\gamma_{n} )	
	\leq 
	\tilde{C}_{4} \hat{t} \leq 1/4
\end{align}
for $0\leq n \leq k \leq a(n,\hat{t} )$, 
(\ref{l1.2.21}) yields 
\begin{align} \label{l1.2.25}
	f(\theta_{k} ) - f(\theta_{n} ) 
	\leq &
	-
	3(\gamma_{k} - \gamma_{n} )
	\|\nabla f(\theta_{n} ) \|^{2}/4 
	+
	\tilde{C}_{3}
	\left(
	\gamma_{n}^{-r} 
	\|\nabla f(\theta_{n} ) \| (\xi + \varepsilon )  	
	+ 
	\gamma_{n}^{-2r} (\xi + \varepsilon )^{2} 
	\right)
\end{align}
for 
$\tau_{1,\varepsilon } < n \leq k \leq a(n,\hat{t} )$. 
As an immediate consequence of (\ref{l1.2.3}), (\ref{l1.2.5}), (\ref{l1.2.25})
we get that 
(\ref{l1.2.1*}) - (\ref{l1.2.5*})  
hold for $n > \tau_{1,\varepsilon }$
(notice that $\gamma_{k} - \gamma_{n} \leq 1$ 
for $n\leq k \leq a(n,1)$).  

Due to (\ref{1.1'*}), we have 
\begin{align*}
	(\gamma_{k} - \gamma_{n} ) \|\nabla f(\theta_{n} ) \|^{2} 
	= &
	\|\nabla f(\theta_{n} ) \| 
	\|(\gamma_{k} - \gamma_{n} ) \nabla f(\theta_{n} ) \| 
	= 
	\|\nabla f(\theta_{n} ) \| 
	\|\theta_{k} - \theta_{n} + \zeta_{n,k} \|
\end{align*}
for $0\leq n \leq k$. 
Combining this with 
(\ref{1.1*}), (\ref{l1.2.9}) and the first part of 
(\ref{l1.2.7}),  we get 
\begin{align*}
	2\left( f(\theta_{k} ) - f(\theta_{n} ) \right) 
	= &
	-
	\|\nabla f(\theta_{n} ) \| 
	\|\theta_{k} - \theta_{n} + \zeta_{n,k} \| 
	-
	(\gamma_{k} - \gamma_{n} ) \|\nabla f(\theta_{n} ) \|^{2} 
	-
	2\phi_{n,k} 
	\\
	\leq &
	-
	\|\nabla f(\theta_{n} ) \| 
	\|\theta_{k} - \theta_{n} \| 
	-
	(\gamma_{k} - \gamma_{n} ) \|\nabla f(\theta_{n} ) \|^{2} 
	\\
	&+
	\|\nabla f(\theta_{n} ) \| \|\zeta_{n,k} \| 
	+
	2|\phi_{n,k} |
	\\
	\leq &
	-
	\|\nabla f(\theta_{n} ) \| 
	\|\theta_{k} - \theta_{n} \| 
	-
	(\gamma_{k} - \gamma_{n} ) \|\nabla f(\theta_{n} ) \|^{2} 
	+
	\tilde{C}_{4} 
	(\gamma_{k} - \gamma_{n} )^{2} 
	\|\nabla f(\theta_{n} ) \|^{2} 
	\\
	&+
	\tilde{C}_{4} 
	\left(
	\gamma_{n}^{-r} \|\nabla f(\theta_{n} ) \| (\xi + \varepsilon ) 
	+
	\gamma_{n}^{-2r} (\xi + \varepsilon )^{2} 
	\right)
	\\
	= & 
	-
	\|\nabla f(\theta_{n} ) \| 
	\|\theta_{k} - \theta_{n} \| 
	-
	\left(
	1 
	-
	\tilde{C}_{4} (\gamma_{k} - \gamma_{n} ) 
	\right) 
	(\gamma_{k} - \gamma_{n} )
	\|\nabla f(\theta_{n} ) \|^{2} 
	\\
	&+
	\tilde{C}_{4} 
	\left(
	\gamma_{n}^{-r} \|\nabla f(\theta_{n} ) \| (\xi + \varepsilon ) 
	+
	\gamma_{n}^{-2r} (\xi + \varepsilon )^{2} 
	\right)
\end{align*}
for $\tau_{1,\varepsilon} < n \leq k \leq a(n,1)$. 
Consequently, (\ref{l1.2.23}) yields 
\begin{align*}
	2\left( f(\theta_{k} ) - f(\theta_{n} ) \right) 
	\leq &
	-
	\|\nabla f(\theta_{n} ) \| 
	\|\theta_{k} - \theta_{n} \| 
	- 3(\gamma_{k} - \gamma_{n} ) \|\nabla f(\theta_{n} ) \|^{2}/4
	\\
	&+
	\tilde{C}_{4} 
	\left(
	\gamma_{n}^{-r} \|\nabla f(\theta_{n} ) \| (\xi + \varepsilon ) 
	+
	\gamma_{n}^{-2r} (\xi + \varepsilon )^{2} 
	\right)
\end{align*}
for $\tau_{1,\varepsilon} < n \leq k \leq a(n,\hat{t} )$. 
Then, (\ref{l1.2.3}) implies that 
(\ref{l1.2.7*}) is true for $n > \tau_{1,\varepsilon }$. 
\end{sproof}

\begin{lemma} \label{lemma1.3}
Suppose that Assumptions \ref{a1.1} -- \ref{a1.3} hold. 
Then, 
$\lim_{n\rightarrow \infty } \nabla f(\theta_{n} ) = 0$ 
on $\Lambda \setminus N_{0}$. 
\end{lemma}

\begin{sproof}
The lemma's assertion is proved by contradiction. 
We assume that 
$\limsup_{n\rightarrow \infty } \|\nabla f(\theta_{n} ) \| > 0$
for some sample 
$\omega \in \Lambda\setminus N_{0}$ 
(notice that all formulas which follow in the proof correspond to 
this sample). 
Then, there exists $a\in (0,\infty )$ and an 
increasing sequence $\{l_{k} \}_{k\geq 0}$
(both depending on $\omega$)
such that 
$\liminf_{k\rightarrow \infty } \|\nabla f(\theta_{l_{k} } ) \| > a$. 
Since $\liminf_{k\rightarrow \infty } f(\theta_{a(l_{k},\hat{t} ) } ) \geq \hat{f}$, 
Lemma \ref{lemma1.2} (inequality (\ref{l1.2.5*})) gives 
\begin{align*}
	\hat{f} 
	-
	\liminf_{k\rightarrow \infty } f(\theta_{l_{k} } ) 
	\leq &
	\limsup_{k\rightarrow \infty } 
	(f(\theta_{a(l_{k},\hat{t} ) } ) - f(\theta_{l_{k} } ) )
	\leq 
	-
	(\hat{t}/2 ) 
	\liminf_{k\rightarrow \infty } 
	\|\nabla f(\theta_{l_{k} } ) \|^{2}
	\leq 
	- a^{2} \hat{t}/2.   
\end{align*}
Therefore, 
$\liminf_{k\rightarrow \infty } f(\theta_{l_{k} } ) \geq \hat{f} + a \hat{t}^{2}/2$. 
Consequently, 
there exist $b,c \in \mathbb{R}$ 
(depending on $\omega$)
such that 
$\hat{f} < b < c < \hat{f} + a \hat{t}^{2}/2$, 
$b < \hat{f} + \hat{\delta}$
and $\limsup_{n\rightarrow \infty } f(\theta_{n} ) > c$. 
Thus, there exist sequences 
$\{m_{k} \}_{k\geq 0}$, $\{n_{k} \}_{k\geq 0}$
(depending on $\omega$)
with the following properties: 
$m_{k} < n_{k} < m_{k+1}$, 
$f(\theta_{m_{k} } ) < b$, 
$f(\theta_{n_{k} } ) > c$ and 
\begin{align} \label{l1.3.5}
	\max_{m_{k} < n \leq n_{k} } f(\theta_{n} ) 
	\geq 
	b
\end{align}
for $k\geq 0$. 
Then, Lemma \ref{lemma1.2} (inequality (\ref{l1.2.3*})) implies 
\begin{align}
	&\label{l1.3.7}
	\limsup_{k\rightarrow \infty } 
	(f(\theta_{m_{k} + 1 } ) - f(\theta_{m_{k} } ) )
	\leq 
	0, 
	\\
	&\label{l1.3.9} 
	\limsup_{k\rightarrow \infty } 
	\max_{m_{k} \leq n \leq a(m_{k}, \hat{t} ) }
	(f(\theta_{n} ) - f(\theta_{m_{k} } ) )
	\leq
	0. 
\end{align}
Since 
\begin{align*}
	b
	>
	f(\theta_{m_{k} } ) 
	=
	f(\theta_{m_{k} + 1 } ) 
	-
	(f(\theta_{m_{k} + 1 } ) - f(\theta_{m_{k} } ) )
	\geq 
	b 
	-
	(f(\theta_{m_{k} + 1 } ) - f(\theta_{m_{k} } ) )
\end{align*}
for $k\geq 0$, 
(\ref{l1.3.7}) yields 
$\lim_{k\rightarrow \infty } f(\theta_{m_{k} } ) = b$. 
As
$f(\theta_{n_{k} } ) - f(\theta_{m_{k} } ) > c-b$
for $k\geq 0$, 
(\ref{l1.3.9}) implies 
$a(m_{k},\hat{t} ) < n_{k}$ for all, 
but infinitely many $k$
(otherwise, 
$\liminf_{k\rightarrow \infty } (f(\theta_{n_{k} } ) - f(\theta_{m_{k} } ) ) \leq 0$
would follow from (\ref{l1.3.9})). 
Consequently, 
$\liminf_{k\rightarrow \infty } f(\theta_{a(m_{k},\hat{t} ) } ) \geq b$
(due to (\ref{l1.3.5})), 
while Lemma \ref{lemma1.2} 
(inequality (\ref{l1.2.5*})) gives 
\begin{align*}	
	0
	\leq
	\limsup_{k\rightarrow \infty } 
	f(\theta_{a(m_{k},\hat{t} )} ) 
	- 
	b 
	= &
	\limsup_{k\rightarrow \infty } 
	(f(\theta_{a(m_{k},\hat{t} )} ) - f(\theta_{m_{k} } ) )
	\leq 
	- 
	(\hat{t}/2)
	\liminf_{k\rightarrow \infty } 
	\|\nabla f(\theta_{m_{k} } ) \|^{2}.   
\end{align*}
Therefore, 
$\lim_{k\rightarrow \infty } \|\nabla f(\theta_{m_{k} } ) \| = 0$. 
Moreover, 
there exists 
$k_{0} \geq 0$ 
(depending on $\omega$)
such that 
$\theta_{m_{k} } \in \hat{Q}$
and 
$f(\theta_{m_{k} } ) \geq (\hat{f} + b ) / 2$
for $k\geq k_{0}$
(notice that 
$\lim_{k\rightarrow \infty } f(\theta_{m_{k} } ) = b > (\hat{f} + b )/2$). 
Consequently, 
$\theta_{m_{k} } \in \hat{Q}$
and 
$0 < (b - \hat{f} )/2 \leq f(\theta_{m_{k} } ) - \hat{f} \leq \hat{\delta}$
for $k\geq k_{0}$
(notice that $f(\theta_{m_{k} } ) < b < \hat{f} + \hat{\delta}$
for $k\geq 0$). 
Then, owing to (\ref{1.3*}) (i.e., to Assumption \ref{a2.3}), 
we have 
\begin{align*} 
	0
	<
	(b - \hat{f} )/2
	\leq 
	f(\theta_{m_{k} } ) - \hat{f} 
	\leq 
	\hat{M} 
	\|\nabla f(\theta_{m_{k} } ) \|^{\hat{\mu} }
\end{align*}
for $k\geq k_{0}$. 
However, this directly contradicts the fact
$\lim_{k\rightarrow \infty } \|\nabla f(\theta_{m_{k} } ) \| = 0$. 
Hence, 
$\lim_{n\rightarrow \infty } \nabla f(\theta_{n} ) = 0$
on $\Lambda \setminus N_{0}$. 
\end{sproof}

\begin{lemma} \label{lemma1.3'}
Suppose that Assumptions \ref{a1.1} -- \ref{a1.3} hold. 
Then, $\lim_{n\rightarrow \infty } f(\theta_{n} ) = \hat{f}$
on $\Lambda\setminus N_{0}$. 
\end{lemma}

\begin{sproof}
We use contradiction to prove the lemma's assertion: 
Suppose that 
$\hat{f} < \limsup_{n\rightarrow \infty } f(\theta_{n} )$
for some sample $\omega \in \Lambda\setminus N_{0}$
(notice that all formulas which follow in the proof correspond to this sample). 
Then, there exists 
$a\in \mathbb{R}$ 
(depending on $\omega$)
such that 
$\hat{f} < a < \hat{f} + \hat{\delta}$
and 
$\limsup_{n\rightarrow \infty } f (\theta_{n} ) > a$. 
Thus, there exists an increasing sequence 
$\{n_{k} \}_{k\geq 0}$ 
(depending on $\omega$)
such that 
$f(\theta_{n_{k} } ) < a$ and 
$f(\theta_{n_{k}+1} ) \geq a$ for $k\geq 0$. 
On the other side, 
Lemma \ref{lemma1.2} (inequality (\ref{l1.2.3*})) implies 
\begin{align}\label{l1.3'.1}
	\limsup_{k\rightarrow \infty } 
	(f(\theta_{n_{k}+1 } ) - f(\theta_{n_{k} } ) )
	\leq 
	0. 
\end{align}
Since
\begin{align*}
	a
	>
	f(\theta_{n_{k} } )
	=
	f(\theta_{n_{k}+1 } )
	-
	(f(\theta_{n_{k}+1 } ) - f(\theta_{n_{k} } ) )
	\geq 
	a 
	-
	(f(\theta_{n_{k}+1 } ) - f(\theta_{n_{k} } ) )
\end{align*}
for $k\geq 0$, 
(\ref{l1.3'.1}) yields 
$\lim_{k\rightarrow \infty } f(\theta_{n_{k} } ) = a$. 
Moreover, 
there exists 
$k_{0} \geq 0$
(depending on $\omega$)
such that 
$\theta_{n_{k} } \in \hat{Q}$
and 
$f(\theta_{n_{k} } ) \geq (\hat{f} + a ) /2$
for $k\geq k_{0}$
(notice that 
$\lim_{k\rightarrow \infty } f(\theta_{n_{k} } ) = a > (\hat{f} + a )/2$). 
Thus, 
$\theta_{n_{k} } \in \hat{Q}$
and 
$0 < (a - \hat{f} )/ 2 \leq f(\theta_{n_{k} } ) - \hat{f} \leq \hat{\delta}$
for $k\geq k_{0}$
(notice that $f(\theta_{n_{k} } ) < a < \hat{f} + \hat{\delta}$ for $k\geq 0$). 
Then, due to (\ref{1.3*})
(i.e., to Assumption \ref{a1.3}), we have 
\begin{align*}
	0
	<
	(a - \hat{f} )/2
	\leq 
	f(\theta_{n_{k} } ) - \hat{f} 
	\leq 
	\hat{M} \|\nabla f(\theta_{n_{k} } ) \|^{\hat{\mu} }
\end{align*}
for $k\geq k_{0}$. 
However, this directly contradicts the fact
$\lim_{n\rightarrow \infty } \nabla f(\theta_{n} ) = 0$. 
Hence, $\lim_{n\rightarrow \infty } f(\theta_{n} ) = \hat{f}$
on $\Lambda \setminus N_{0}$.
\end{sproof}

\begin{lemma} \label{lemma1.4}
Suppose that Assumptions \ref{a1.1} -- \ref{a1.3} hold. 
Then, there exist random quantities  
$\hat{C}_{2}$, $\hat{C}_{3}$ 
(which are deterministic functions of $\hat{p}$, $\hat{C}$, $\hat{M}$)
and 
for any real number $\varepsilon \in (0,\infty )$, there exists 
a non-negative integer-valued random quantity  
$\tau_{2,\varepsilon}$ such that the following is true: 
$1\leq \hat{C}_{2}, \hat{C}_{3} < \infty$, 
$0\leq \tau_{2,\varepsilon } < \infty$
everywhere and  
\begin{align}
	&\label{l1.4.1*}
	\left(
	u(\theta_{a(n,\hat{t} ) } ) 
	- 
	u(\theta_{n} ) 
	+
	\hat{t} \|\nabla f(\theta_{n} ) \|^{2}/4
	\right)
	I_{A_{n,\varepsilon } }
	\leq 
	0, 
	\\
	&\label{l1.4.3*}
	\left(
	u(\theta_{a(n,\hat{t} ) } ) 
	- 
	u(\theta_{n} ) 
	+
	(\hat{t}/\hat{C}_{3} ) \: u(\theta_{n} ) 
	\right)
	I_{B_{n,\varepsilon } }
	\leq 
	0, 
	\\
	&\label{l1.4.5*}
	\left(
	v(\theta_{a(n,\hat{t} ) } ) 
	- 
	v(\theta_{n} ) 
	- 
	(\hat{t}/\hat{C}_{3} )
	(\varphi_{\varepsilon }(\xi) )^{-\hat{\mu}/\hat{p} }
	\right)
	I_{C_{n,\varepsilon } }
	\geq 
	0
\end{align}
on $\Lambda\setminus N_{0}$ for $n\geq \tau_{2,\varepsilon}$, 
where 
\begin{align*}
	&
	\begin{aligned}[b]
		A_{n,\varepsilon}
		= &
		\left\{
		\gamma_{n}^{\hat{p} } 
		|u(\theta_{n} ) |
		\geq 
		\hat{C}_{2} 
		(\varphi_{\varepsilon }(\xi) )^{\hat{\mu} }
		\right\}
		\cup
		\Big\{
		\gamma_{n}^{\hat{p} } 
		\|\nabla f(\theta_{n} ) \|^{2}
		\geq 
		\hat{C}_{2} 
		(\varphi_{\varepsilon }(\xi) )^{\hat{\mu} } 
		\Big\}, 
	\end{aligned}
	\\
	&
	B_{n,\varepsilon}
	= 
	\left\{
	\gamma_{n}^{\hat{p} } 
	u(\theta_{n} ) 
	\geq 
	\hat{C}_{2} 
	(\varphi_{\varepsilon }(\xi) )^{\hat{\mu} }
	\right\}
	\cap
	\{\hat{\mu} = 2 \}, 
	\\
	&
	C_{n,\varepsilon }
	=
	\left\{
	\gamma_{n}^{\hat{p} } u(\theta_{n} ) 
	\geq 
	\hat{C}_{2} 
	(\varphi_{\varepsilon }(\xi) )^{\hat{\mu} }
	\right\}
	\cap
	\left\{
	u(\theta_{a(n,\hat{t} ) } ) > 0
	\right\}
	\cap
	\left\{
	\hat{\mu} < 2
	\right\}. 
\end{align*} 
\end{lemma}

\begin{remark}
Inequalities (\ref{l1.4.1*}) -- (\ref{l1.4.5*}) can be represented in 
the following equivalent form: 
Relations 
\begin{align}
	&\label{l1.4.21*}
	\left(
	\gamma_{n}^{\hat{p} } 
	|u(\theta_{n} ) |
	\geq 
	\hat{C}_{2} 
	(\varphi_{\varepsilon }(\xi) )^{\hat{\mu} }
	\;\vee\;
	\gamma_{n}^{\hat{p} } 
	\|\nabla f(\theta_{n} ) \|^{2}
	\geq 
	\hat{C}_{2} 
	(\varphi_{\varepsilon }(\xi) )^{\hat{\mu} } 
	\right)
	\; \wedge \; 
	n > \tau_{2,\varepsilon}
	\nonumber\\
	&
	\Longrightarrow
	u(\theta_{a(n,\hat{t} ) } ) 
	\leq
	u(\theta_{n} ) 
	-
	\hat{t} \|\nabla f(\theta_{n} ) \|^{2}/4, 
	\\
	&\label{l1.4.23*}
	\gamma_{n}^{\hat{p} } 
	u(\theta_{n} ) 
	\geq 
	\hat{C}_{2} 
	(\varphi_{\varepsilon }(\xi) )^{\hat{\mu} }
	\;\wedge\;
	\hat{\mu}=2
	\; \wedge \; 
	n > \tau_{2,\varepsilon}
	\nonumber\\
	&
	\Longrightarrow
	u(\theta_{a(n,\hat{t} ) } ) 
	\leq 
	\left(
	1
	-
	\hat{t} /\hat{C}_{3} 
	\right)
	u(\theta_{n} ), 
	\\
	&\label{l1.4.25*}
	\gamma_{n}^{\hat{p} } u(\theta_{n} ) 
	\geq 
	\hat{C}_{2} 
	(\varphi_{\varepsilon }(\xi) )^{\hat{\mu} }
	\;\wedge\;
	u(\theta_{a(n,\hat{t} ) } ) > 0
	\;\wedge\;
	\hat{\mu} < 2
	\; \wedge \; 
	n > \tau_{2,\varepsilon}
	\nonumber\\
	&
	\Longrightarrow
	v(\theta_{a(n,\hat{t} ) } ) 
	\geq
	v(\theta_{n} ) 
	+
	(\hat{t}/\hat{C}_{3} )
	(\varphi_{\varepsilon }(\xi) )^{-\hat{\mu}/\hat{p} }
\end{align}
are true on $\Lambda\setminus N_{0}$. 
\end{remark}

\begin{sproof}
Let 
$\tilde{C} = 8\hat{C}_{1}/\hat{t}$,   
$\hat{C}_{2} = \tilde{C}^{2} \hat{M}$
and 
$\hat{C}_{3} = 4 \hat{p} \hat{M}^{2}$. 
Moreover, let $\varepsilon \in (0,\infty )$ be an arbitrary real number.  
Then, owing to Lemma \ref{lemma1.1} and \ref{lemma1.3'}, 
it is possible to construct a non-negative inter-valued random quantity 
$\tau_{2,\varepsilon}$ 
such that $\tau_{1,\varepsilon} \leq \tau_{2,\varepsilon} < \infty$ everywhere
and such that 
$\theta_{n} \in \hat{Q}$, 
$|u(\theta_{n} ) | \leq \hat{\delta}$, 
\begin{align}
	&\label{l1.4.1'}
	\gamma_{n}^{-\hat{p}/2} 
	(\varphi_{\varepsilon }(\xi) )^{\hat{\mu}/2 }
	\geq 
	\gamma_{n}^{-r} (\xi + \varepsilon), 
	\\
	&\label{l1.4.1''}
	\gamma_{n}^{-\hat{p}/\hat{\mu} } 
	\varphi_{\varepsilon }(\xi) 
	\geq  
	\gamma_{n}^{-r} (\xi + \varepsilon )
\end{align}
on $\Lambda \setminus N_{0}$ for 
$n>\tau_{2,\varepsilon}$.\footnote{
To conclude that 
(\ref{l1.4.1'}) holds on $\Lambda\setminus N_{0}$
for all but finitely many $n$, 
notice that 
$\hat{p}/2 < \min\{r,\hat{r} \} \leq r$
when $\hat{\mu}<2$ 
and that the left and right hand sides of the inequality in (\ref{l1.4.1'})
are equal when $\hat{\mu}=2$. 
In order to deduce that 
(\ref{l1.4.1''}) is true on $\Lambda\setminus N_{0}$ for all but finitely many $n$, 
notice that 
$\hat{p}/\hat{\mu} = r$, $\varphi_{\varepsilon}(\xi ) \geq \xi + \varepsilon$
when $r\leq \hat{r}$
and that $\hat{p}/\hat{\mu} = \hat{r} < r$
when $r>\hat{r}$.}  
Since 
$\tau_{2,\varepsilon } \geq \tau_{1,\varepsilon }$
on $\Lambda \setminus N_{0}$, 
Lemma \ref{lemma1.2} (inequality (\ref{l1.2.5*}))
yields
\begin{align}\label{l1.4.3'}
		u(\theta_{a(n,\hat{t} ) } ) - u(\theta_{n} ) 
		\leq &
		-
		\hat{t} \|\nabla f(\theta_{n} ) \|^{2}/2
		+
		\hat{C}_{1} 
		\left(
		\gamma_{n}^{-r} 
		\|\nabla f(\theta_{n} ) \| (\xi + \varepsilon) 
		+ 
		\gamma_{n}^{-2r} 
		(\xi + \varepsilon )^{2} 
		\right)
\end{align}
on $\Lambda\setminus N_{0}$ for 
$n>\tau_{2,\varepsilon}$. 
As $\theta_{n} \in \hat{Q}$
and $|u(\theta_{n} ) | \leq \hat{\delta}$
on $\Lambda\setminus N_{0}$ for 
$n>\tau_{2,\varepsilon }$, 
(\ref{1.3*})
(i.e., Assumption \ref{a1.3}) implies
\begin{align}\label{l1.4.5'}
	|u(\theta_{n} ) |
	\leq 
	\hat{M} \|\nabla f(\theta_{n} ) \|^{\hat{\mu} }
\end{align}
on $\Lambda\setminus N_{0}$ for 
$n>\tau_{2,\varepsilon }$.

Let $\omega$ be an arbitrary sample from 
$\Lambda\setminus N_{0}$
(notice that all formulas which follow in the proof correspond to this 
sample). 
First, we show (\ref{l1.4.1*}). 
We proceed by contradiction: 
Suppose that 
(\ref{l1.4.1*}) is violated for some 
$n>\tau_{2,\varepsilon}$. 
Therefore, 
\begin{align} \label{l1.4.1} 
	u(\theta_{a(n,\hat{t} ) } ) - u(\theta_{n} )
	>
	-\hat{t} \|\nabla f(\theta_{n} ) \|^{2}/4 
\end{align}
and at least one of the following two inequalities is true: 
\begin{align}
		&\label{l1.4.3}
		|u(\theta_{n} ) |
		\geq 
		\hat{C}_{2} 
		\gamma_{n}^{-\hat{p} } 
		(\varphi_{\varepsilon }(\xi) )^{\hat{\mu} },
		\\
		&\label{l1.4.5}
		\|\nabla f(\theta_{n} ) \|^{2}
		\geq 
		\hat{C}_{2} 
		\gamma_{n}^{-\hat{p} } 
		(\varphi_{\varepsilon }(\xi) )^{\hat{\mu} }.  
\end{align}
If (\ref{l1.4.3}) holds, 
then (\ref{l1.4.1''}), (\ref{l1.4.5'}) imply  
\begin{align*}
	\|\nabla f(\theta_{n} ) \|
	\geq 
	(|u(\theta_{n} ) |/\hat{M} )^{1/\hat{\mu} }
	\geq 
	(\hat{C}_{2}/\hat{M} )^{1/\hat{\mu} } \gamma_{n}^{-\hat{p}/\hat{\mu} } 
	\varphi_{\varepsilon }(\xi)
	\geq 
	\tilde{C} \gamma_{n}^{-r} (\xi + \varepsilon ) 
\end{align*}
(notice that 
$(\hat{C}_{2}/\hat{M} )^{1/\hat{\mu} } = \tilde{C}^{2/\hat{\mu} } \geq \tilde{C}$ 
owing to $\hat{\mu}\leq 2$). 
On the other side, if (\ref{l1.4.5}) is satisfied, 
then (\ref{l1.4.1'}) yields 
\begin{align*}
	\|\nabla f(\theta_{n} ) \|
	\geq 
	\hat{C}_{2}^{1/2} 
	\gamma_{n}^{-\hat{p}/2} 
	(\varphi_{\varepsilon }(\xi) )^{\hat{\mu}/2 }
	\geq 
	\tilde{C} 
	\gamma_{n}^{-r} (\xi + \varepsilon ). 
\end{align*}
Thus, as a result of one of (\ref{l1.4.3}), (\ref{l1.4.5}), 
we get 
\begin{align*}
	\|\nabla f(\theta_{n} ) \|
	\geq 
	\tilde{C} 
	\gamma_{n}^{-r} (\xi + \varepsilon ). 
\end{align*}
Consequently, 
\begin{align*}
	&
	\hat{t} \|\nabla f(\theta_{n} ) \|^{2}/8
	\geq 
	(\tilde{C} \hat{t} /8 ) \gamma_{n}^{-r} 
	\|\nabla f(\theta_{n} ) \|
	(\xi + \varepsilon ) 
	= 
	\hat{C}_{1} 
	\gamma_{n}^{-r} 
	\|\nabla f(\theta_{n} ) \|
	(\xi + \varepsilon ), 
	\\
	&
	\hat{t} \|\nabla f(\theta_{n} ) \|^{2}/8
	\geq 
	(\tilde{C}^{2} \hat{t} /8 ) \gamma_{n}^{-2r}
	(\xi + \varepsilon )^{2}
	\geq 
	\hat{C}_{1} \gamma_{n}^{-2r} (\xi + \varepsilon )^{2}
\end{align*}
(notice that 
$\tilde{C} \hat{t} /8 = \hat{C}_{1}$, 
$\tilde{C}^{2} \hat{t} /8 \geq \tilde{C}\hat{t}/8 = \hat{C}_{1}$). 
Combining this with (\ref{l1.4.3'}), 
we get 
\begin{align}\label{l1.4.7}
	u(\theta_{a(n,\hat{t} ) } ) - u(\theta_{n} )
	\leq 
	-\hat{t} \|\nabla f(\theta_{n} ) \|^{2}/4,  
\end{align}
which directly contradicts (\ref{l1.4.1}). 
Hence, (\ref{l1.4.1*}) is true for $n > \tau_{2,\varepsilon }$. 
Then, as a result of (\ref{l1.4.5'}) and 
the fact that $B_{n,\varepsilon } \subseteq A_{n,\varepsilon }$
for $n\geq 0$, 
we get
\begin{align*}
	\left(
	u(\theta_{a(n,\hat{t} ) } ) 
	- 
	u(\theta_{n} ) 
	+
	(\hat{t}/\hat{C}_{3} ) \: u(\theta_{n} ) 
	\right)
	I_{B_{n,\varepsilon } }
	&
	\leq
	\left(
	u(\theta_{a(n,\hat{t} ) } ) 
	- 
	u(\theta_{n} ) 
	+
	(\hat{M} \hat{t} /\hat{C}_{3} ) \: \|\nabla f(\theta_{n} ) \|^{2}
	\right)
	I_{B_{n,\varepsilon } }
	\\
	&
	\leq
	\left(
	u(\theta_{a(n,\hat{t} ) } ) 
	- 
	u(\theta_{n} ) 
	+
	\hat{t} \|\nabla f(\theta_{n} ) \|^{2}/4
	\right)
	I_{B_{n,\varepsilon } }
	\leq 
	0
\end{align*}
for $n> \tau_{2,\varepsilon}$
(notice that $u(\theta_{n} ) > 0$ on $B_{n,\varepsilon}$ for 
each $n\geq 0$; also notice that $\hat{C}_{3} \geq 4 \hat{M}$). 
Thus, (\ref{l1.4.3*}) is true for $n> \tau_{2,\varepsilon}$. 

Now, let us prove (\ref{l1.4.5*}). 
To do so, we again use contradiction: 
Suppose that (\ref{l1.4.3*}) does not hold for some $n>\tau_{2,\varepsilon }$. 
Consequently, we have 
$\hat{\mu} < 2$, $u(\theta_{a(n,\hat{t} ) } ) > 0$ and 
\begin{align}
	&\label{l1.4.9}
	\gamma_{n}^{\hat{p} } \: u(\theta_{n} ) 
	\geq 
	\hat{C}_{2} 
	(\varphi_{\varepsilon }(\xi) )^{\hat{\mu} }
	>
	0, 
	\\
	&\label{l1.4.11}
	v(\theta_{a(n,\hat{t} ) } ) 
	- 
	v(\theta_{n} ) 
	<  
	(\hat{t} / \hat{C}_{3} ) 
	(\varphi_{\varepsilon }(\xi) )^{-\hat{\mu}/\hat{p} }.  
\end{align}
Combining (\ref{l1.4.9}) with
(already proved) (\ref{l1.4.1*}), 
we get (\ref{l1.4.7}), 
while $\hat{\mu} < 2$ implies 
 \begin{align} \label{l1.4.11'}
	2/\hat{\mu} 
	= 
	1 + 1/(\hat{\mu} \hat{r} ) 
	\leq 
	1 + 1/\hat{p}
\end{align}
(notice that 
$\hat{r} = 1/(2 - \hat{\mu} )$
owing to $\hat{\mu} < 2$; 
also notice that 
$\hat{p} = \hat{\mu} \min\{r,\hat{r} \} \leq \hat{\mu} \hat{r}$). 
As $0 < u(\theta_{n} ) \leq \hat{\delta} \leq 1$
(due to (\ref{l1.4.9}) and the definition of $\tau_{2,\varepsilon }$), 
inequalities (\ref{l1.4.5'}), (\ref{l1.4.11'})  
yield
\begin{align}\label{l1.4.15}
	\|\nabla f(\theta_{n} ) \|^{2}
	\geq
	\left(
	u(\theta_{n} )/\hat{M} 
	\right)^{2/\hat{\mu} }
	\geq 
	\left(
	u(\theta_{n} )
	\right)^{1+1/\hat{p} }
	/\hat{M}^{2}  
\end{align}
(notice that $\hat{M}^{2/\hat{\mu} } \leq \hat{M}^{2}$ due to 
$\hat{\mu} < 2$, $\hat{M}\geq 1$). 
Since 
$\|\nabla f(\theta_{n} ) \| > 0$
and 
$0 < u(\theta_{a(n,\hat{t} ) } ) < u(\theta_{n} )$
(due to (\ref{l1.4.5'}), (\ref{l1.4.7})), 
inequalities (\ref{l1.4.7}), (\ref{l1.4.15}) give
\begin{align*}
	\frac{\hat{t} }{4}
	\leq 
	\frac{u(\theta_{n} ) - u(\theta_{a(n,\hat{t} ) } ) }
	{\|\nabla f(\theta_{n} ) \|^{2} }
	\leq &
	\hat{M}^{2}
	\frac{u(\theta_{n} ) - u(\theta_{a(n,\hat{t} ) } ) }
	{\left(u(\theta_{n} ) \right)^{1+1/\hat{p} } }
	\\
	= &
	\hat{M}^{2}
	\int^{u(\theta_{n} ) }_{u(\theta_{a(n,\hat{t} ) } ) }
	\frac{du}{\left(u(\theta_{n} ) \right)^{1+1/\hat{p} } }
	\\
	\leq &
	\hat{M}^{2}
	\int^{u(\theta_{n} ) }_{u(\theta_{a(n,\hat{t} ) } ) }
	\frac{du}{u^{1+1/\hat{p} } }
	\\
	= &
	\hat{p}
	\hat{M}^{2}
	\left(
	v(\theta_{a(n,\hat{t} ) } )
	-
	v(\theta_{n} )
	\right).
\end{align*}
Therefore, 
\begin{align*}
	v(\theta_{a(n,\hat{t} ) } )
	-
	v(\theta_{n} )
	\geq 
	\hat{t}/(4\hat{p} \hat{M}^{2} )
	=
	(\hat{t}/\hat{C}_{3} ),  
\end{align*}
which directly contradicts (\ref{l1.4.11}). 
Thus, (\ref{l1.4.5*}) is satisfied for $n>\tau_{2,\varepsilon}$. 
\end{sproof}

\begin{lemma} \label{lemma1.21}
Suppose that Assumptions \ref{a1.1} -- \ref{a1.3} hold. 
Then, there exists a random quantity $\hat{C}_{4}$
(which is a deterministic function of $\hat{C}$)
and for any $\varepsilon\in (0,\infty )$ 
there exists a non-negative integer-valued random quantity 
$\tau_{3,\varepsilon }$
such that the following is true: 
$1\leq \hat{C}_{4} < \infty$, $0\leq \tau_{3,\varepsilon } < \infty$
everywhere 
and 
\begin{align}\label{l1.21.1*}
	\|\theta_{a(n,\hat{t}) } - \theta_{n} \|
	\leq 
	-\gamma_{n}^{\hat{q}+1} 
	\left(
	u(\theta_{a(n,\hat{t} ) } ) - u(\theta_{n} ) 
	\right)
	(\phi_{\varepsilon }(\xi ) )^{-1} 
	+ 
	\hat{C}_{4} \gamma_{n}^{-(\hat{q}+1)} \phi_{\varepsilon}(\xi ) 
\end{align}
on $\Lambda\setminus N_{0}$ for $n > \tau_{3,\varepsilon}$
and any $\varepsilon \in (0,\infty )$. 
\end{lemma} 

\begin{sproof}
Let $\varepsilon \in (0,\infty )$ be an arbitrary real number, 
while 
$\hat{C}_{4} = 10\hat{C}_{1}^{2}/\hat{t}$. 
Then, it is possible to construct a non-negative integer-valued random quantity 
$\tau_{3,\varepsilon }$ such that 
$\tau_{1,\varepsilon } \leq \tau_{3,\varepsilon } < \infty$ everywhere 
and such that 
\begin{align}\label{l1.21.301}
	\gamma_{n}^{-(\hat{q}+1) } \phi_{\varepsilon }(\xi )
	\geq 
	\gamma_{n}^{-r} (\xi+\varepsilon )
\end{align}
on $\Lambda\setminus N_{0}$ for $n>\tau_{3,\varepsilon }$.\footnote{
To deduce that (\ref{l1.21.301}) holds on $\Lambda\setminus N_{0}$ 
for all but finitely many $n$, notice that 
$\hat{q}+1=\hat{r}<r$ when $r>\hat{r}$
and that 
$\hat{q}+1=r$, 
$\phi_{\varepsilon }(\xi ) = \varphi_{\varepsilon }(\xi ) \geq \xi + \varepsilon$
when $r\leq \hat{r}$.} 

Let $\omega$ be an arbitrary sample from $\Lambda\setminus N_{0}$
(notice that all formulas which follow in the proof correspond to this sample), 
while $n>\tau_{3,\varepsilon}$ is an arbitrary integer. 
To prove (\ref{l1.21.1*}), 
we consider separately the cases
$\|\nabla f(\theta_{n} ) \| \geq 
(4\hat{C}_{1}/\hat{t} ) \gamma_{n}^{-(\hat{q}+1 ) } \phi_{\varepsilon}(\xi )$
and 
$\|\nabla f(\theta_{n} ) \| \leq 
(4\hat{C}_{1}/\hat{t} ) \gamma_{n}^{-(\hat{q} + 1 ) } \phi_{\varepsilon}(\xi )$. 

{\em Case 
$\|\nabla f(\theta_{n} ) \| \geq 
(4\hat{C}_{1}/\hat{t} ) \gamma_{n}^{-(\hat{q} + 1 ) } \phi_{\varepsilon}(\xi )$:} 
Owing to (\ref{l1.21.301}), 
we have 
\begin{align*}
	\|\nabla f(\theta_{n} ) \| 
	\geq 
	(4\hat{C}_{1}/\hat{t} )
	\gamma_{n}^{-r} (\xi + \varepsilon ). 
\end{align*}
Therefore, 
\begin{align*}
	&
	(\hat{t}/4 ) \|\nabla f(\theta_{n} ) \|^{2} 
	\geq 
	\hat{C}_{1} \gamma_{n}^{-r} \|\nabla f(\theta_{n} ) \| (\xi + \varepsilon ), 
	\\
	&
	(\hat{t}/4 ) \|\nabla f(\theta_{n} ) \|^{2} 
	\geq 
	(4\hat{C}_{1}^{2}/\hat{t} )
	\gamma_{n}^{-2r} (\xi + \varepsilon )^{2} 
	\geq 
	\hat{C}_{1} 
	\gamma_{n}^{-2r} (\xi + \varepsilon )^{2}.
\end{align*}
Then, Lemma \ref{lemma1.2} (inequality (\ref{l1.2.7*})) yields 
\begin{align*}
	\|\nabla f(\theta_{n} ) \| 
	\|\theta_{a(n,\hat{t} ) } - \theta_{n} \| 
	\leq &
	-
	2\left(
	u(\theta_{a(n,\hat{t} ) } ) - u(\theta_{n} ) 
	\right)
	-
	\hat{t} \|\nabla f(\theta_{n} ) \|^{2}/2 
	\\
	&+
	\hat{C}_{1} 
	\left(
	\gamma_{n}^{-r} \|\nabla f(\theta_{n} ) \| 
	(\xi + \varepsilon ) 
	+
	\gamma_{n}^{-2r} (\xi + \varepsilon )^{2} 
	\right)
	\\
	\leq &
	-
	2\left(
	u(\theta_{a(n,\hat{t} ) } ) - u(\theta_{n} ) 
	\right). 
\end{align*}
Consequently, 
\begin{align*}
	\|\theta_{a(n,\hat{t} ) } - \theta_{n} \| 
	\leq &
	-
	2 \|\nabla f(\theta_{n} ) \|^{-1} 
	\left(
	u(\theta_{a(n,\hat{t} ) } ) - u(\theta_{n} ) 
	\right)
	\\
	\leq &
	- 
	(2\hat{C}_{1}/\hat{t} )^{-1} \gamma_{n}^{\hat{q} + 1 } 
	\left(
	u(\theta_{a(n,\hat{t} ) } ) - u(\theta_{n} ) 
	\right)
	(\phi_{\varepsilon}(\xi ) )^{-1}
	\\
	\leq & 
	-\gamma_{n}^{\hat{q} + 1 } 
	\left(
	u(\theta_{a(n,\hat{t} ) } ) - u(\theta_{n} ) 
	\right)
	(\phi_{\varepsilon}(\xi ) )^{-1} 
	+ 
	\hat{C}_{4} \gamma_{n}^{-(\hat{q} + 1 ) } \phi_{\varepsilon}(\xi ).  
\end{align*}
Hence, (\ref{l1.21.1*}) is true when 
$\|\nabla f(\theta_{n} ) \| \geq 
(4\hat{C}_{1}/\hat{t} ) \gamma_{n}^{-(\hat{q} + 1 ) } \phi_{\varepsilon}(\xi )$. 

{\em Case 
$\|\nabla f(\theta_{n} ) \| \leq 
(4\hat{C}_{1}/\hat{t} ) \gamma_{n}^{-(\hat{q} + 1 ) } \phi_{\varepsilon}(\xi )$:} 
Due to
Lemma \ref{lemma1.2} (inequalities (\ref{l1.2.1*}), (\ref{l1.2.3*}))
and (\ref{l1.21.301}), we have  
\begin{align}
	&\label{l1.21.3}
	\begin{aligned}[b]
	\|\theta_{a(n,\hat{t} ) } - \theta_{n} \|
	\leq &
	\hat{C}_{1}
	\left(
	\|\nabla f(\theta_{n} ) \| 
	+
	\gamma_{n}^{-r} (\xi + \varepsilon ) 
	\right)
	\leq 
	(\hat{C}_{4}/2 ) \gamma_{n}^{-(\hat{q} + 1 ) } \phi_{\varepsilon}(\xi ),  
	\end{aligned}
	\\
	&
	\begin{aligned}[b]
	u(\theta_{a(n,\hat{t} ) } ) - u(\theta_{n} )
	\leq &
	\hat{C}_{1}
	\left(
	\gamma_{n}^{-r} \|\nabla f(\theta_{n} ) \| (\xi + \varepsilon ) 
	+
	\gamma_{n}^{-2r} (\xi + \varepsilon )^{2} 
	\right)
	\nonumber
	\leq 
	(\hat{C}_{4}/2 ) \gamma_{n}^{-2(\hat{q} + 1 ) } 
	(\phi_{\varepsilon}(\xi ) )^{2}
	\end{aligned}
\end{align}
Hence, 
\begin{align*}
	\gamma_{n}^{\hat{q} + 1 } 
	\left(u(\theta_{a(n,\hat{t} ) } ) - u(\theta_{n} ) \right) 
	(\phi_{\varepsilon}(\xi) )^{-1} 
	\leq 
	(\hat{C}_{4}/2 ) 
	\gamma_{n}^{-(\hat{q} + 1 ) } \phi_{\varepsilon}(\xi). 
\end{align*}
Combining this with (\ref{l1.21.3}), 
we get 
\begin{align*}
	\|\theta_{a(n,\hat{t} ) } - \theta_{n} \|
	\leq &
	(\hat{C}_{4}/2 )
	\gamma_{n}^{-(\hat{q} + 1 ) } \phi_{\varepsilon}(\xi) 
	-
	\gamma_{n}^{\hat{q} + 1 } 
	\left(
	u(\theta_{a(n,\hat{t} ) } ) - u(\theta_{n} )
	\right)
	(\phi_{\varepsilon}(\xi) )^{-1}
	\\
	&+
	\gamma_{n}^{\hat{q} + 1 } 
	\left(
	u(\theta_{a(n,\hat{t} ) } ) - u(\theta_{n} )
	\right)
	(\phi_{\varepsilon}(\xi) )^{-1}
	\\
	\leq &
	-
	\gamma_{n}^{\hat{q} + 1 } 
	\left(
	u(\theta_{a(n,\hat{t} ) } ) - u(\theta_{n} )
	\right)
	(\phi_{\varepsilon}(\xi) )^{-1}
	+
	\hat{C}_{4} 
	\gamma_{n}^{-(\hat{q} + 1 ) } \phi_{\varepsilon}(\xi).  
\end{align*}
Thus, (\ref{l1.21.1*}) holds when 
$\|\nabla f(\theta_{n} ) \| \leq 
(4\hat{C}_{1}/\hat{t} ) \gamma_{n}^{-(\hat{q} + 1 ) } \phi_{\varepsilon}(\xi )$. 
\end{sproof}

\begin{lemma} \label{lemma1.5}
Suppose that Assumptions \ref{a1.1} -- \ref{a1.3} hold. 
Then, 
\begin{align}	\label{l1.5.1*} 
	u(\theta_{n} ) 
	\geq 
	-
	\hat{C}_{2} \gamma_{n}^{-\hat{p} } 
	(\varphi_{\varepsilon }(\xi) )^{\hat{\mu} }
\end{align}
on $\Lambda\setminus N_{0}$
for $n > \tau_{2,\varepsilon }$ 
and any $\varepsilon \in (0,\infty )$. 
Furthermore, 
there exists a random quantity  
$\hat{C}_{5} \in [1,\infty )$ 
(which is a deterministic function of $\hat{p}$, $\hat{C}$, $\hat{M}$)
such that the following is true: 
$1\leq \hat{C}_{5} < \infty$ everywhere 
and 
\begin{align} \label{l1.5.3*}
	\|\nabla f(\theta_{n} ) \|^{2}
	\leq 
	\hat{C}_{5} 
	\left(
	\psi(u(\theta_{n} ) )
	+ 
	\gamma_{n}^{-\hat{p} } (\varphi_{\varepsilon }(\xi) )^{\hat{\mu} } 
	\right)
\end{align}
on $\Lambda\setminus N_{0}$ 
for $n> \tau_{2,\varepsilon }$
and any $\varepsilon \in (0,\infty )$, 
where
function $\psi(\cdot )$
is defined by 
$\psi(x) = x \:{\rm I}_{(0,\infty )}(x)$, $x\in \mathbb{R}$. 
\end{lemma}

\begin{sproof}
Let 
$\hat{C}_{5} = 4 \hat{C}_{2} / \hat{t}$, 
while $\varepsilon \in (0,\infty )$ is an arbitrary real number. 
Moreover, 
$\omega$ is an arbitrary sample from $\Lambda\setminus N_{0}$
(notice that all formulas which follow in the proof correspond to 
this sample). 

First, 
we prove (\ref{l1.5.1*}). 
To do so, we use contradiction:  
Assume that (\ref{l1.5.1*}) is not satisfied for 
some $n > \tau_{2,\varepsilon }$. 
Define $\{n_{k} \}_{k\geq 0}$ recursively by
$n_{0}=n$
and  
$n_{k} = a(n_{k-1},\hat{t} )$ 
for $k\geq 1$.  
Let us show by induction that 
$\{u(\theta_{n_{k} } ) \}_{k\geq 0}$ is non-increasing: 
Suppose that 
$u(\theta_{n_{l} } ) \leq u(\theta_{n_{l-1} } )$ for 
$0\leq l \leq k$. 
Consequently, 
\begin{align*}
	u(\theta_{n_{k} } ) 
	\leq 
	u(\theta_{n_{0} } ) 
	\leq 
	-
	\hat{C}_{2} \gamma_{n_{0} }^{-\hat{p} }
	(\varphi_{\varepsilon }(\xi) )^{\hat{\mu} }
	\leq 
	-
	\hat{C}_{2} \gamma_{n_{k} }^{-\hat{p} }
	(\varphi_{\varepsilon }(\xi) )^{\hat{\mu} } 
\end{align*}
(notice that $\{\gamma_{n} \}_{n\geq 0}$ is increasing). 
Then, Lemma \ref{lemma1.4} (relations (\ref{l1.4.1*}), (\ref{l1.4.21*})) yields 
\begin{align*}
	u(\theta_{n_{k+1} } ) - u(\theta_{n_{k} } ) 
	\leq 
	-
	\hat{t} \|\nabla f(\theta_{n_{k} } ) \|^{2}/4 
	\leq 
	0,  
\end{align*}
i.e., $u(\theta_{n_{k+1} } ) \leq u(\theta_{n_{k} } )$. 
Thus, $\{u(\theta_{n_{k} } ) \}_{k\geq 0}$ is non-increasing. 
Therefore,
\begin{align*}
	\limsup_{n\rightarrow \infty } 
	u(\theta_{n_{k} } ) 
	\leq 
	u(\theta_{n_{0} } ) 
	< 
	0. 
\end{align*}
However, this is not possible, as 
$\lim_{n\rightarrow \infty } u(\theta_{n} ) = 0$
(due to Lemma \ref{lemma1.3'}). 
Hence, (\ref{l1.5.1*}) indeed holds 
for $n> \tau_{2,\varepsilon }$. 

Now, (\ref{l1.5.3*}) is demonstrated. 
Again, we proceed by contradiction: 
Suppose that (\ref{l1.5.3*}) is violated for some 
$n> \tau_{2,\varepsilon }$. 
Consequently, 
\begin{align*}
	\|\nabla f(\theta_{n} ) \|^{2}
	\geq 
	\hat{C}_{5} \gamma_{n}^{-\hat{p} } (\varphi_{\varepsilon }(\xi) )^{\hat{\mu} }
	\geq 
	\hat{C}_{2} \gamma_{n}^{-\hat{p} } (\varphi_{\varepsilon }(\xi) )^{\hat{\mu} } 
\end{align*}
(notice that $\hat{C}_{5} \geq \hat{C}_{2}$),
which, together with Lemma \ref{lemma1.4}
(relations (\ref{l1.4.1*}), (\ref{l1.4.21*})), yields 
\begin{align*}
	u(\theta_{a(n,\hat{t} ) } ) - u(\theta_{n} )
	\leq 
	-\hat{t} \|\nabla f(\theta_{n} ) \|^{2}/4. 
\end{align*}
Then, (\ref{l1.5.1*}) implies 
\begin{align*}
	\|\nabla f(\theta_{n} ) \|^{2}
	\leq &
	(4/\hat{t} )
	\left(
	u(\theta_{n} ) 
	-
	u(\theta_{a(n,\hat{t} )} ) 
	\right)
	\\
	\leq &
	(4/\hat{t} )
	\left(
	\psi(u(\theta_{n} ) )
	+
	\hat{C}_{2} \gamma_{a(n,\hat{t} ) }^{-\hat{p} } 
	(\varphi_{\varepsilon }(\xi) )^{\hat{\mu} }
	\right)
	\\
	\leq &
	\hat{C}_{5} 
	\left(
	\psi(u(\theta_{n} ) )
	+ 
	\gamma_{n}^{-\hat{p} } 
	(\varphi_{\varepsilon }(\xi) )^{\hat{\mu} }
	\right). 
\end{align*}
However, this directly contradicts our assumption 
that $n$ violates (\ref{l1.5.3*}). 
Thus, (\ref{l1.5.3*}) is indeed satisfied for 
$n>\tau_{2,\varepsilon }$. 
\end{sproof}

\begin{lemma} \label{lemma1.6}
Suppose that Assumptions \ref{a1.1} -- \ref{a1.3} hold. 
Then, 
there exists a random quantity 
$\hat{C}_{6}$ 
(which is a deterministic function of $\hat{p}$, $\hat{C}$, $\hat{M}$)
such that the following is true: 
$1\leq \hat{C}_{6} < \infty$ everywhere
and 
\begin{align} \label{l1.6.1*}
	\liminf_{n\rightarrow \infty } 
	\gamma_{n}^{\hat{p} } \:
	u(\theta_{n} ) 
	\leq 
	\hat{C}_{6}
	(\varphi_{\varepsilon}(\xi) )^{\hat{\mu} }
\end{align}
on $\Lambda\setminus N_{0}$
for any $\varepsilon \in (0,\infty )$. 
\end{lemma}

\begin{sproof}
Let 
$\hat{C}_{6} = \hat{C}_{2}+\hat{C}_{3}^{\hat{p} }$. 
We prove (\ref{l1.6.1*}) by contradiction: 
Assume that (\ref{l1.6.1*}) is violated for some 
sample $\omega$ from $\Lambda\setminus N_{0}$
(notice that the formulas which follow in the proof correspond to 
this sample)
and some real number $\varepsilon \in (0,\infty )$. 
Consequently, there exists  
$n_{0} > \tau_{2,\varepsilon }$ 
(depending on $\omega$, $\varepsilon$)
such that 
\begin{align} \label{l1.6.1}
	u(\theta_{n} ) 
	\geq 
	\hat{C}_{6} \gamma_{n}^{-\hat{p} } 
	(\varphi_{\varepsilon }(\xi) )^{\hat{\mu} }
\end{align}
for $n\geq n_{0}$. 
Let $\{n_{k} \}_{k\geq 0}$ be defined recursively by 
$n_{k} = a(n_{k-1},\hat{t} )$ for $k\geq 1$. 	
In what follows in the proof, we consider separately 
the cases $\hat{\mu} < 2$
and $\hat{\mu} = 2$. 

{\em Case $\hat{\mu} < 2$:}
Due to (\ref{l1.6.1}), we have 
\begin{align*}
		v(\theta_{n_{k} } ) 
		\leq &
		\hat{C}_{6}^{-1/\hat{p} } 
		\gamma_{n_{k} }
		(\varphi_{\varepsilon }(\xi) )^{-\hat{\mu} /\hat{p} }. 
\end{align*}
On the other side, Lemma \ref{lemma1.4} 
(relations (\ref{l1.4.5*}), (\ref{l1.4.25*})) 
and (\ref{l1.6.1}) yield
\begin{align*}
	&
	\begin{aligned}[b]
		v(\theta_{n_{k+1} } ) - v(\theta_{n_{k} } ) 
		\geq 
		(\hat{t}/\hat{C}_{3} )
		(\varphi_{\varepsilon }(\xi) )^{-\hat{\mu}/\hat{p} }
		\geq 
		(1/\hat{C}_{3} ) (\gamma_{n_{k+1} } - \gamma_{n_{k} } )
		(\varphi_{\varepsilon }(\xi) )^{-\hat{\mu}/\hat{p} }
	\end{aligned}
\end{align*}
for $k\geq 0$
(notice that $\hat{t} \geq \gamma_{n_{k+1} } - \gamma_{n_{k} }$).  
Therefore, 
\begin{align*}
	(1/\hat{C}_{3} ) (\gamma_{n_{k} }  - \gamma_{n_{0} } )
	(\varphi_{\varepsilon }(\xi) )^{-\hat{\mu}/\hat{p} }
	\leq &
	\sum_{i=0}^{k-1} 
	(v(\theta_{n_{i+1} } ) - v(\theta_{n_{i} } ) )
	= 
	v(\theta_{n_{k} } ) - v(\theta_{n_{0} } ) 
	\leq 
	\hat{C}_{6}^{-1/\hat{p} } \gamma_{n_{k} }
	(\varphi_{\varepsilon }(\xi) )^{-\hat{\mu}/\hat{p} }
\end{align*}
for $k\geq 1$. 
Thus, 
\begin{align*}
	(1 - \gamma_{n_{0} }/\gamma_{n_{k} } )
	\leq 
	\hat{C}_{3} \hat{C}_{6}^{-1/\hat{p} }
\end{align*}
for $k\geq 1$. 
However, this is impossible, since the limit process 
$k\rightarrow \infty$ (applied to the previous relation)  
yields 
$\hat{C}_{3} \geq \hat{C}_{6}^{1/\hat{p} }$
(notice that $\hat{C}_{6} > \hat{C}_{3}^{\hat{p} }$). 
Hence, (\ref{l1.6.1*}) holds when $\hat{\mu} < 2$. 

{\em Case $\hat{\mu} = 2$:} 
As a result of 
Lemma \ref{lemma1.4} (relations (\ref{l1.4.3*}), (\ref{l1.4.23*})) and (\ref{l1.6.1}), we get 
\begin{align*}
	u(\theta_{n_{k+1} } )
	\leq
	(1 - \hat{t}/\hat{C}_{3} ) u(\theta_{n_{k} } )
	\leq
	\left(
	1
	-
	(\gamma_{n_{k+1} } - \gamma_{n_{k} } )/ \hat{C}_{3}
	\right)
	u(\theta_{n_{k} } )
\end{align*}
for $k\geq 0$. 
Consequently, 
\begin{align*}
	u(\theta_{n_{k} } )
	\leq 
	u(\theta_{n_{0} } )
	\prod_{i=1}^{k} 
	\left(
	1 
	-  
	(\gamma_{n_{i} } - \gamma_{n_{i-1} } )/\hat{C}_{3} 	
	\right)
	\leq &
	u(\theta_{n_{0} } )	
	\exp\left(
	-
	(1/\hat{C}_{3} )
	\sum_{i=1}^{k} (\gamma_{n_{i} } - \gamma_{n_{i-1} } )	
	\right)
	\\
	= & 
	u(\theta_{n_{0} } )	
	\exp\left(
	- 
	(\gamma_{n_{k} } - \gamma_{n_{0} } )/\hat{C}_{3}	
	\right)
\end{align*}
for $k\geq 0$. 
Then, (\ref{l1.6.1}) yields 
\begin{align*}
	\hat{C}_{6} 
	(\varphi_{\varepsilon}(\xi) )^{\hat{\mu} }
	\leq 
	u(\theta_{n_{0} } )
	\gamma_{n_{k} }^{\hat{p} } 
	\exp\left(
	-(\gamma_{n_{k} } - \gamma_{n_{0} } )/\hat{C}_{3} 
	\right)
\end{align*}
for $k\geq 0$. 
However, this is not possible, as the limit 
process $k\rightarrow \infty$
(applied to the previous relation) 
implies 
$\hat{C}_{6} (\varphi_{\varepsilon}(\xi) )^{\hat{\mu} }
\leq 0$. 
Thus, (\ref{l1.6.1*}) holds also when $\hat{\mu} = 2$. 
\end{sproof}

\begin{lemma} \label{lemma1.7}
Suppose that Assumptions \ref{a1.1} -- \ref{a1.3} hold. 
Then, 
there exists a random quantity  
$\hat{C}_{7}$ (which is a deterministic function of $\hat{p}$, $\hat{C}$, $\hat{M}$)
such that 
the following is true: 
$1\leq \hat{C}_{7} < \infty$ everywhere 
and 
\begin{align} \label{l1.7.1*}
	\limsup_{n\rightarrow \infty } 
	\gamma_{n}^{\hat{p} } \:
	u(\theta_{n} ) 
	\leq 
	\hat{C}_{7} 
	(\varphi_{\varepsilon}(\xi) )^{\hat{\mu} }	
\end{align}
on $\Lambda\setminus N_{0}$
for any $\varepsilon \in (0,\infty )$. 
\end{lemma}

\begin{sproof}
Let 
$\tilde{C}_{1} = 3 \hat{C}_{1}  \hat{C}_{5}$, 
$\tilde{C}_{2} = 6 \tilde{C}_{1} \hat{C}_{2} + \hat{C}_{3}^{\hat{p} } + \hat{C}_{6}$
and 
$\hat{C}_{7} = 2 (\tilde{C}_{1} + \tilde{C}_{2} )^{2}$. 
We use contradiction to show (\ref{l1.7.1*}): 
Suppose that (\ref{l1.7.1*}) is violated for some sample 
$\omega$ from $\Lambda\setminus N_{0}$
(notice that the formulas which appear in the proof correspond to 
this sample)
and some real number $\varepsilon \in (0,\infty )$.
Then, it can be deduced from Lemma \ref{lemma1.6} that 
there exist 
$n_{0} > m_{0} > \tau_{2,\varepsilon }$ 
(depending on $\omega$, $\varepsilon$)
such that 
\begin{align}
	& \label{l1.7.1}
	\gamma_{m_{0} }^{\hat{p} } u(\theta_{m_{0} } ) 
	\leq 
	\tilde{C}_{2} 
	(\varphi_{\varepsilon }(\xi) )^{\hat{\mu} }, 
	\\
	& \label{l1.7.3}
	\gamma_{n_{0} }^{\hat{p} } u(\theta_{n_{0} } ) 
	\geq 
	\hat{C}_{7} 
	(\varphi_{\varepsilon }(\xi) )^{\hat{\mu} }, 
	\\
	& \label{l1.7.5}
	\min_{m_{0}<n\leq n_{0} }
	\gamma_{n}^{\hat{p} } \: u(\theta_{n} ) 
	>
	\tilde{C}_{2}  
	(\varphi_{\varepsilon }(\xi) )^{\hat{\mu} },   
	\\
	& \label{l1.7.5'}
	\max_{m_{0}\leq n < n_{0} }
	\gamma_{n}^{\hat{p} } \: u(\theta_{n} ) 
	<
	\hat{C}_{7} 
	(\varphi_{\varepsilon }(\xi) )^{\hat{\mu} }  
\end{align}
(notice that $\tilde{C}_{2} > \hat{C}_{6}$) 
and such that 
\begin{align} 
	& \label{l1.7.7'}
	(
	\gamma_{a(m_{0},\hat{t} ) }/\gamma_{m_{0} }
	)^{\hat{p} }
	\leq 
	\min\{2, (1 - \hat{t}/\hat{C}_{3} )^{-1} \}, 
	\\
	&	\label{l1.7.7}
	\gamma_{m_{0} }^{-2r} 
	(\xi + \varepsilon )^{2}
	\leq 
	\gamma_{m_{0} }^{-\hat{p} } 
	(\varphi_{\varepsilon }(\xi) )^{\hat{\mu} }
\end{align}
(to see that (\ref{l1.7.7'}) holds for all, but finitely many
$m_{0}$, 
notice that 
$\lim_{n\rightarrow \infty } \gamma_{a(n,\hat{t} ) }/\gamma_{n} = 1$; 
to conclude that (\ref{l1.7.7}) is true for all, but finitely many $m_{0}$, 
notice that 
$\hat{p} < 2\min\{r,\hat{r} \} \leq 2r$ if $\hat{\mu} < 2$
and
that the left and right-hand sides of (\ref{l1.7.7}) are 
equal when $\hat{\mu}=2$). 

Let $l_{0} = a(m_{0},\hat{t} )$. 
As a direct consequence of Lemmas \ref{lemma1.2}, \ref{lemma1.5} 
(relations (\ref{l1.2.3*}), (\ref{l1.5.3*})) and 
(\ref{l1.7.7}), we get
\begin{align} \label{l1.7.9}
	u(\theta_{n} )
	-
	u(\theta_{m_{0} } )
	\leq &
	\hat{C}_{1}\left(
	\gamma_{m_{0} }^{-r} 
	\|\nabla f(\theta_{m_{0} } ) \| 
	(\xi + \varepsilon ) 
	+
	\gamma_{m_{0} }^{-2r} 
	(\xi + \varepsilon )^{2} 
	\right)
	\nonumber\\
	\leq &
	\hat{C}_{1}\left( 
	\|\nabla f(\theta_{m_{0} } ) \|^{2}/2 
	+
	3\gamma_{m_{0} }^{-2r} 
	(\xi + \varepsilon )^{2}/2 
	\right)
	\nonumber \\
	\leq &
	\hat{C}_{1} \hat{C}_{5} \: \psi(u(\theta_{m_{0} } ) )
	+
	(2\hat{C}_{1} + \hat{C}_{1} \hat{C}_{5} )
	\gamma_{m_{0} }^{-\hat{p} }
	(\varphi_{\varepsilon }(\xi) )^{\hat{\mu} }
	\nonumber \\
	\leq &
	\tilde{C}_{1} 
	\left(
	\psi(u(\theta_{m_{0} } ) )
	+
	\gamma_{m_{0} }^{-\hat{p} }
	(\varphi_{\varepsilon }(\xi) )^{\hat{\mu} }
	\right) 
\end{align}
for $m_{0} \leq n \leq l_{0}$. 
Then, 
(\ref{l1.7.5}), (\ref{l1.7.7'}), (\ref{l1.7.9}) yield 
\begin{align} \label{l1.7.21'}
	u(\theta_{m_{0} } ) 
	+
	\tilde{C}_{1} 
	\psi(u(\theta_{m_{0} } ) )
	\geq &
	u(\theta_{m_{0} + 1 } ) 
	-
	\tilde{C}_{1} 
	\gamma_{m_{0} }^{-\hat{p} }
	(\varphi_{\varepsilon }(\xi) )^{\hat{\mu} }
	\nonumber \\
	\geq &
	(\tilde{C}_{2} 
	\gamma_{m_{0} + 1 }^{-\hat{p} }
	-
	\tilde{C}_{1} 
	\gamma_{m_{0} }^{-\hat{p} })
	(\varphi_{\varepsilon }(\xi) )^{\hat{\mu} }
	\nonumber \\
	= &
	\left(
	\tilde{C}_{2} 
	(\gamma_{m_{0} + 1 }/\gamma_{m_{0} } )^{-\hat{p} } 
	-
	\tilde{C}_{1} 
	\right)
	\gamma_{m_{0} }^{-\hat{p} }
	(\varphi_{\varepsilon }(\xi) )^{\hat{\mu} }
	\nonumber \\
	\geq &
	(\tilde{C}_{2}/2 - \tilde{C}_{1} )  
	\gamma_{m_{0} }^{-\hat{p} }
	(\varphi_{\varepsilon }(\xi) )^{\hat{\mu} } 
	> 0 
\end{align}
(notice that 
$(\gamma_{m_{0} + 1 } /\gamma_{m_{0} } )^{\hat{p} } \leq 
(\gamma_{l_{0} } /\gamma_{m_{0} } )^{\hat{p} } \leq 2$; 
also notice that $\tilde{C}_{2}/2 \geq 3 \tilde{C}_{1}$), 
while (\ref{l1.7.1}), (\ref{l1.7.7'}), (\ref{l1.7.9}) imply 
\begin{align} \label{l1.7.21}
	u(\theta_{n} )  
	\leq 
	(1 + \tilde{C}_{1} )
	u(\theta_{m_{0} } ) 
	+
	\tilde{C}_{1} 
	\gamma_{m_{0} }^{-\hat{p} }
	(\varphi_{\varepsilon }(\xi) )^{\hat{\mu} }
	\nonumber 
	\leq &
	(\tilde{C}_{1} + \tilde{C}_{2} + \tilde{C}_{1} \tilde{C}_{2} )
	\gamma_{m_{0} }^{-\hat{p} } (\varphi_{\varepsilon }(\xi) )^{\hat{\mu} } 
	\nonumber \\
	< &
	(\hat{C}_{7}/2) 
	(\gamma_{n}/\gamma_{m_{0} } )^{\hat{p} }
	\gamma_{n}^{-\hat{p} } (\varphi_{\varepsilon }(\xi) )^{\hat{\mu} } 
	\nonumber \\
	\leq  &
	\hat{C}_{7} 
	\gamma_{n}^{-\hat{p} } (\varphi_{\varepsilon }(\xi) )^{\hat{\mu} } 
\end{align}
for $m_{0} \leq n \leq l_{0}$
(notice that 
$(\gamma_{n}/\gamma_{m_{0} } )^{\hat{p} } \leq
(\gamma_{l_{0} } /\gamma_{m_{0} } )^{\hat{p} } \leq 2$
for $m_{0}\leq n\leq l_{0}$; 
also notice that $\hat{C}_{7}/2 = (\tilde{C}_{1} + \tilde{C}_{2} )^{2} >
\tilde{C}_{1} + \tilde{C}_{2} + \tilde{C}_{1} \tilde{C}_{2}$).  
Due to  
(\ref{l1.7.3}), (\ref{l1.7.5'}), (\ref{l1.7.21}), we have 
$l_{0} < n_{0}$.
On the other side, 
since $x + \tilde{C}_{1} \psi(x) \geq 0$ only if $x\geq 0$
and since $x + \tilde{C}_{1} \psi(x) = (1 + \tilde{C}_{1} ) x$ for $x\geq 0$, 
inequality (\ref{l1.7.21'}) implies   
\begin{align} \label{l1.7.23} 
	u(\theta_{m_{0} } )
	\geq &
	(1 + \tilde{C}_{1} )^{-1} (\tilde{C}_{2}/2 - \tilde{C}_{1} )
	\gamma_{m_{0} }^{-\hat{p} }
	(\varphi_{\varepsilon }(\xi) )^{\hat{\mu} }
	\geq 
	\hat{C}_{2} 
	\gamma_{m_{0} }^{-\hat{p} }
	(\varphi_{\varepsilon }(\xi) )^{\hat{\mu} } 
\end{align}
(notice that 
$\tilde{C}_{2}/2 - \tilde{C}_{1}
\geq \tilde{C}_{1} (3 \hat{C}_{2} - 1 )
\geq 2 \tilde{C}_{1} \hat{C}_{2} 
\geq (1+\tilde{C}_{1} ) \hat{C}_{2}$).

In what follows in the proof, we consider separately the cases 
$\hat{\mu} < 2$ and $\hat{\mu} = 2$. 

{\em Case $\hat{\mu} < 2$:}
Owing to Lemma \ref{lemma1.4} (relations (\ref{l1.4.5*}), (\ref{l1.4.25*})) and 
(\ref{l1.7.1}), 
(\ref{l1.7.23}), we have 
\begin{align*}
	v(\theta_{l_{0} } ) 
	\geq 
	v(\theta_{m_{0} } ) 
	+
	(\hat{t}/\hat{C}_{3} ) 
	(\varphi_{\varepsilon }(\xi) )^{-\hat{\mu}/\hat{p} }
	\geq &
	\left(
	\tilde{C}_{2}^{-1/\hat{p} } 
	\gamma_{m_{0} } 
	+
	\hat{C}_{3}^{-1} 
	(\gamma_{l_{0} } - \gamma_{m_{0} } )
	\right)
	(\varphi_{\varepsilon }(\xi) )^{-\hat{\mu}/\hat{p} }
	\\
	> &
	\min\{\tilde{C}_{2}^{-1/\hat{p} }, \hat{C}_{3}^{-1} \}
	\gamma_{l_{0} } 
	(\varphi_{\varepsilon}(\xi) )^{-\hat{\mu}/\hat{p} }
	\\
	= &
	\tilde{C}_{2}^{-1/\hat{p} } 
	\gamma_{l_{0} } 	
	(\varphi_{\varepsilon }(\xi) )^{-\hat{\mu}/\hat{p} } 
\end{align*}
(notice that 
$\hat{t} \geq \gamma_{l_{0} } - \gamma_{m_{0} }$; 
also notice 
$\tilde{C}_{2}^{-1/\hat{p} } < \hat{C}_{3}^{-1}$).
Consequently, 
\begin{align*}
	u(\theta_{l_{0} } ) 
	=
	\left(
	v(\theta_{l_{0} } ) 
	\right)^{-\hat{p} }
	< 
	\tilde{C}_{2} 
	\gamma_{l_{0} }^{-\hat{p} }
	(\varphi_{\varepsilon }(\xi) )^{\hat{\mu} }. 
\end{align*}
However, this directly contradicts 
(\ref{l1.7.5}) and the fact that 
$l_{0} < n_{0}$. 
Thus, (\ref{l1.7.1*}) holds 
when $\hat{\mu} < 2$. 

{\em Case $\hat{\mu} = 2$:}
Using Lemma \ref{lemma1.4} (relations (\ref{l1.4.3*}), (\ref{l1.4.23*})) and (\ref{l1.7.23}), 
we get
\begin{align*}
	u(\theta_{l_{0} } )
	\leq 
	\left(
	1
	- 
	\hat{t}/\hat{C}_{3} 
	\right)
	u(\theta_{{m}_{0} } ). 
\end{align*}
Then, (\ref{l1.7.1}), (\ref{l1.7.7'}) yield 
\begin{align*}
	u(\theta_{l_{0} } )
	\leq &
	\tilde{C}_{2} 
	(1 - \hat{t}/\hat{C}_{3} )
	(\gamma_{l_{0} }/\gamma_{m_{0} } )^{\hat{p} }
	\gamma_{l_{0} }^{-\hat{p} }
	(\varphi_{\varepsilon}(\xi) )^{\hat{\mu} }
	\leq 
	\tilde{C}_{2} 
	\gamma_{l_{0} }^{-\hat{p} }
	(\varphi_{\varepsilon}(\xi) )^{\hat{\mu} }. 
\end{align*}
However, this is impossible due to (\ref{l1.7.5})
and the fact that $l_{0} < n_{0}$. 
Hence, 
(\ref{l1.7.1*}) 
also in the case $\hat{\mu} = 2$. 
\end{sproof}

\begin{lemma} \label{lemma1.23}
Suppose that Assumptions \ref{a1.1} -- \ref{a1.3} hold. 
Then, 
there exists a random quantity 
$\hat{C}_{8}$ 
(which is a deterministic function of $\hat{p}$, $\hat{C}$, $\hat{M}$)
such that the following is true: 
$1\leq \hat{C}_{8} < \infty$ everywhere 
and 
\begin{align} \label{l1.23.1*}
	\limsup_{n\rightarrow \infty } 
	\gamma_{n}^{\hat{q} }
	\sup_{k\geq n}
	\|\theta_{k} - \theta_{n} \|
	\leq 
	\hat{C}_{8}
	\varphi_{\varepsilon}(\xi ) 
\end{align}
on $\Lambda\setminus N_{0}$. 
\end{lemma}

\begin{sproof}
Let 
$\varepsilon \in (0,\infty )$
be an arbitrary real number, 
while 
$\tilde{C}_{1} = 2 (\hat{C}_{2} + \hat{C}_{7} )$, 
$\tilde{C}_{2} = 2(\hat{q} + 1 ) \tilde{C}_{1} + \hat{C}_{4}$, 
$\tilde{C}_{3} = 2\tilde{C}_{1} + 3\hat{q}^{-1} \hat{t}^{-1} \tilde{C}_{2}$, 
$\hat{C}_{8} = 2\tilde{C}_{1} + \tilde{C}_{3}$. 
Moreover, let 
$\omega$ is an arbitrary sample from 
$\Lambda\setminus N_{0}$
(notice that all formulas which follow in the proof correspond to 
this sample). 

Owing to Lemmas \ref{lemma1.5} and \ref{lemma1.7}, we have 
\begin{align}
	&\label{l1.23.501}
	\limsup_{n\rightarrow\infty } 
	\gamma_{n}^{\hat{p} } |u(\theta_{n} ) |
	\leq 
	\max\{\hat{C}_{2}, \hat{C}_{7} \} 
	(\varphi_{\varepsilon}(\xi) )^{\hat{\mu} },  
	\\
	&\label{l1.23.503}
	\begin{aligned}[b]
	\limsup_{n\rightarrow\infty } 
	\gamma_{n}^{\hat{p} } \|\nabla f(\theta_{n} ) \|^{2} 
	\leq & 
	\hat{C}_{5} 
	\limsup_{n\rightarrow\infty } 
	\gamma_{n}^{\hat{p} } \psi(u(\theta_{n} ) ) 
	+
	\hat{C}_{5} (\varphi_{\varepsilon}(\xi) )^{\hat{\mu} }
	\\
	\leq & 
	2\hat{C}_{5} 
	\max\{\hat{C}_{2}, \hat{C}_{7} \} 
	(\varphi_{\varepsilon}(\xi) )^{\hat{\mu} }.   
	\end{aligned}
\end{align}
We also conclude that $\hat{q}<r$, $\hat{p}/2>\hat{q}$ 
and that 
\begin{align}\label{l1.23.505}
	\gamma_{n}^{-\hat{p} } (\varphi_{\varepsilon}(\xi) )^{\hat{\mu} }
	\leq 
	\gamma_{n}^{-(2\hat{q} + 1 ) } 
	\varphi_{\varepsilon}(\xi) \phi_{\varepsilon}(\xi) 
\end{align}
for all but finitely many $n$.\footnote
{To conclude that $\hat{p}/2>\hat{q}$ and that (\ref{l1.23.505}) holds for all but finitely many $n$, 
notice the following: 
\begin{compactenum}[(i)]
\item
If $\hat{\mu}=2$, then $\hat{r}=\infty$, $\hat{p}=2r$, $\hat{q}=r-1$, 
$\varphi_{\varepsilon}(\xi) = \phi_{\varepsilon}(\xi)$, 
and thus, 
$\hat{p}=2\hat{q}+2$, 
$(\varphi_{\varepsilon}(\xi) )^{\hat{\mu} } 
= \varphi_{\varepsilon}(\xi) \phi_{\varepsilon}(\xi)$. 
Consequently, $\hat{\mu}=2$ implies that $\hat{p}/2>\hat{q}$
and that (\ref{l1.23.505}) is true for each $n$. 
\item
If $\hat{\mu}<2$, $r\geq \hat{r}$, 
then $\hat{r}=1/(2-\hat{\mu} )$, 
$\hat{p}=\hat{\mu}\hat{r}$, 
$\hat{q}=\hat{r}-1$, 
$\varphi_{\varepsilon}(\xi)>1$, 
and hence, 
$\hat{p}=2\hat{q}+1$, 
$(\varphi_{\varepsilon}(\xi) )^{\hat{\mu} } 
\leq \varphi_{\varepsilon}(\xi) \phi_{\varepsilon}(\xi)$. 
Therefore, $\hat{\mu}<2$, $r\geq \hat{r}$ yields that $\hat{p}/2>\hat{q}$
and that (\ref{l1.23.505}) is satisfied for any $n$. 
\item
If $\hat{\mu}<2$, $r<\hat{r}$, 
then $\hat{r}=1/(2-\hat{\mu} )$, 
$\hat{p}=\hat{\mu}r$, 
$\hat{q}=r-1$, 
and thus, 
$\hat{p}=2r - r/\hat{r}>2r-1=2\hat{q}+1$. 
Consequently, when $\hat{\mu}<2$, $r<\hat{r}$, 
we have that $\hat{p}/2>\hat{q}$
and that (\ref{l1.23.505}) holds for all but finitely many $n$. 
\end{compactenum}
}
Consequently, 
Lemma \ref{lemma1.2} (inequality (\ref{l1.2.1*})) 
and (\ref{l1.23.503}) imply 
\begin{align}\label{l1.23.507}
	\limsup_{n\rightarrow\infty } 
	\gamma_{n}^{\hat{q} } 
	\max_{n\leq k \leq a(n,\hat{t} ) } 
	\|\theta_{k} - \theta_{n} \| 
	\leq &
	2\hat{C}_{1} \hat{C}_{5} 
	\max\{\hat{C}_{2}, \hat{C}_{7} \} 
	(\varphi_{\varepsilon}(\xi) )^{\hat{\mu}/2 } 
	\lim_{n\rightarrow\infty } \gamma_{n}^{\hat{q}-\hat{p}/2 }  
	\nonumber\\
	&
	+
	\hat{C}_{1} (\xi+\varepsilon ) 
	\lim_{n\rightarrow\infty } \gamma_{n}^{\hat{q}-r } 
	=
	0. 
\end{align}
On the other side, 
it is straightforward to show  
$\gamma_{a(n,\hat{t} ) } - \gamma_{n} = 
\hat{t} + O(\alpha_{a(n,\hat{t} ) } )$
and 
\begin{align} \label{l1.23.201}
	\gamma_{a(n,\hat{t} ) }^{\hat{q}+1 } - \gamma_{n}^{\hat{q}+1 } 
	= &
	\gamma_{a(n,\hat{t} ) }^{\hat{q}+1 }
	\left(
	1
	-
	\left(
	1 - (\gamma_{a(n,\hat{t} ) } - \gamma_{n} )/\gamma_{a(n,\hat{t} ) } 
	\right)^{\hat{q}+1 }
	\right)
	= 
	\gamma_{a(n,\hat{t} ) }^{\hat{q}+1 }
	\left(
	(\hat{q}+1 ) \hat{t} \gamma_{a(n,\hat{t} ) }^{-1}
	+
	O(\gamma_{a(n,\hat{t} ) }^{-2} )
	\right) 
\end{align}
for $n\rightarrow \infty$. 
Combining this with (\ref{l1.23.501}), (\ref{l1.23.505}), (\ref{l1.23.507}), 
we deduce that there exist 
$n_{0}>0$ (depending on $\omega$, $\varepsilon$) 
such that $n_{0}>\tau_{3,\varepsilon}$ and such that 
\begin{align}
	&\label{l1.23.1}
	\gamma_{a(n,\hat{t} ) } - \gamma_{n} 
	\geq 
	\hat{t}/2, 
	\\
	&\label{l1.23.3}
	\gamma_{a(n,\hat{t} ) }^{\hat{q}+1 } - \gamma_{n}^{\hat{q}+1 } 
	\leq 
	2(\hat{q}+1 ) \gamma_{a(n,\hat{t} ) }^{\hat{q} }, 
	\\
	&\label{l1.23.5}
	|u(\theta_{n} ) |
	\leq 
	\tilde{C}_{1} \gamma_{n}^{-(2\hat{q}+1 ) }
	\varphi_{\varepsilon}(\xi) \phi_{\varepsilon}(\xi), 
	\\
	&\label{l1.23.103}
	\max_{n\leq k \leq a(n,\hat{t} ) } \|\theta_{k} - \theta_{n} \|
	\leq 
	\tilde{C}_{1} 
	\gamma_{n}^{-\hat{q} } 
	\varphi_{\varepsilon}(\xi) 
\end{align}
on $\Lambda\setminus N_{0}$
for $n> n_{0}$. 

Let $\{n_{k} \}_{k\geq 0}$ be recursively defined by 
$n_{k+1} = a(n_{k}, \hat{t} )$ for $k\geq 0$. 
Then, due to Lemma \ref{lemma1.21}, we have 
\begin{align*}
	\|\theta_{n_{l} } - \theta_{n_{k} } \|
	\leq 
	\sum_{i=k}^{l-1} \|\theta_{n_{i+1} } - \theta_{n_{i} } \|
	\leq &
	\sum_{i=k}^{l-1} 
	\gamma_{n_{i} }^{\hat{q}+1 } 
	\left(
	u(\theta_{n_{i} } ) - u(\theta_{n_{i+1} } )
	\right)
	(\phi_{\varepsilon}(\xi) )^{-1}
	+
	\hat{C}_{4} 
	\sum_{i=k}^{l-1} 
	\gamma_{n_{i} }^{-(\hat{q}+1 ) } 
	\phi_{\varepsilon}(\xi) 
	\\
	\leq &
	\sum_{i=k+1}^{l} 
	(\gamma_{n_{i} }^{\hat{q}+1 } - \gamma_{n_{i-1} }^{\hat{q}+1 } )
	|u(\theta_{n_{i} } ) |
	(\phi_{\varepsilon}(\xi) )^{-1}
	+
	\hat{C}_{4} 
	\sum_{i=k}^{l-1} 
	\gamma_{n_{i} }^{-(\hat{q}+1 )} 
	\phi_{\varepsilon}(\xi)  
	\\
	&+
	\gamma_{n_{l} }^{\hat{q}+1 } 
	|u(\theta_{n_{l} } ) |
	(\phi_{\varepsilon}(\xi) )^{-1}
	+
	\gamma_{n_{k} }^{\hat{q}+1 } 
	|u(\theta_{n_{k} } ) |
	(\phi_{\varepsilon}(\xi) )^{-1}
\end{align*}
for $0\leq k \leq l$. 
As $\phi_{\varepsilon}(\xi) \leq \varphi_{\varepsilon}(\xi)$, 
(\ref{l1.23.3}), (\ref{l1.23.5}) yield 
\begin{align} \label{l1.23.105}
	\|\theta_{n_{l} } - \theta_{n_{k} } \|
	\leq &
	2\tilde{C}_{1} (\hat{q}+1 ) 
	\varphi_{\varepsilon}(\xi) 
	\sum_{i=k+1}^{l} 
	\gamma_{n_{i} }^{-(\hat{q}+1 ) } 
	+
	\hat{C}_{4} \varphi_{\varepsilon}(\xi ) 
	\sum_{i=k}^{l-1} \gamma_{n_{i} }^{-(\hat{q}+1 ) } 
	+
	\tilde{C}_{1} 
	(\gamma_{n_{k} }^{-\hat{q} } + \gamma_{n_{l} }^{-\hat{q} } )
	\varphi_{\varepsilon}(\xi ) 
	\nonumber \\
	\leq &
	\tilde{C}_{2} \varphi_{\varepsilon}(\xi) 
	\sum_{i=k+1}^{\infty } 
	\gamma_{n_{i} }^{-(\hat{q}+1 ) } 
	+
	2\tilde{C}_{1} \gamma_{n_{k} }^{-\hat{q} } 
	\varphi_{\varepsilon}(\xi) 
\end{align}
for $0\leq k \leq l$. 
Since 
\begin{align*}
	\gamma_{n_{l} } 
	=
	\gamma_{n_{k} } 
	+
	\sum_{i=k}^{l-1} 
	(\gamma_{n_{i+1} } - \gamma_{n_{i} } )
	\geq 
	\gamma_{n_{k} } 
	+
	2^{-1} \hat{t} (l-k)
\end{align*}
for $0\leq k \leq l$
(owing to (\ref{l1.23.1})), we get
\begin{align}\label{l1.23.701}
	\sum_{i=k}^{\infty } \gamma_{n_{i} }^{-(\hat{q}+1 ) }
	\leq & 
	\sum_{i=0}^{\infty } 
	(\gamma_{n_{k} } + i\hat{t}/2 )^{-(\hat{q}+1 ) } 
	\leq 
	\gamma_{n_{k} }^{-(\hat{q}+1 ) }
	+
	\int_{0}^{\infty } 
	(\gamma_{n_{k} } + u\hat{t}/2 )^{-(\hat{q}+1 ) } du 
	\leq 
	3 \hat{q}^{-1} \hat{t}^{-1}
	\gamma_{n_{k} }^{-\hat{q} }
\end{align}
for $k\geq 0$. 
Then, (\ref{l1.23.105}) implies
\begin{align} \label{l1.23.107}
	\|\theta_{n_{l} } - \theta_{n_{k} } \| 
	\leq &
	\tilde{C}_{3} 
	\gamma_{n_{k} }^{-\hat{q} }
	\varphi_{\varepsilon}(\xi ) 
\end{align}
for $0\leq k \leq l$. 
Combining this with (\ref{l1.23.103}), we obtain
\begin{align*}
	\|\theta_{k} - \theta_{n} \|
	\leq &
	\|\theta_{k} - \theta_{n_{j} } \| 
	+
	\|\theta_{n_{j} } - \theta_{n_{i} } \| 
	+ 
	\|\theta_{n_{i} } - \theta_{n} \| 
	\\
	\leq &
	\tilde{C}_{3} 
	\gamma_{n_{i} }^{-\hat{q} }
	\varphi_{\varepsilon}(\xi ) 
	+
	\tilde{C}_{1}  
	(\gamma_{n}^{-\hat{q} } + \gamma_{n_{j} }^{-\hat{q} } )
	\varphi_{\varepsilon}(\xi ) 
	\\
	\leq & 
	\hat{C}_{8} 
	\gamma_{n}^{-\hat{q} }
	\varphi_{\varepsilon}(\xi ) 
\end{align*}
for $n_{0} < n \leq k$, $1\leq i \leq j$
satisfying 
$n_{i-1} \leq n < n_{i}$, $n_{j} \leq k < n_{j+1}$. 
Then, it is obvious that (\ref{l1.23.1*}) is true. 
\end{sproof}

\begin{vproof}{Theorems \ref{theorem1.1} and \ref{theorem1.2}}
{Owing to Lemmas \ref{lemma1.3} and \ref{lemma1.23}, 
$\hat{\theta } = \lim_{n\rightarrow \infty } \theta_{n}$
exists and satisfies $\nabla f(\hat{\theta } ) = 0$
on $\Lambda\setminus N_{0}$. 
Thus, Theorem \ref{theorem1.1} holds. 
In addition, we have 
$\hat{Q} \subseteq \{\theta \in \mathbb{R}^{d_{\theta } }: 
\|\theta - \hat{\theta } \| \leq \delta_{\hat{\theta } } \}$
on $\Lambda\setminus N_{0}$
($\delta_{\theta }$ is specified in Remark \ref{remark1.1}). 
Therefore, on $\Lambda\setminus N_{0}$, random quantities 
$\hat{\mu}$, $\hat{p}$, $\hat{r}$ defined in the beginning of this section 
coincide with 
$\hat{\mu}$, $\hat{p}$, $\hat{r}$ specified in Theorem 
\ref{theorem1.2}
(see Remark \ref{remark1.1}). 
Similarly, on $\Lambda\setminus N_{0}$, 
$\hat{C}$, $\hat{M}$ introduced in this section 
are identical to $C_{\hat{\theta } }$, $M_{\hat{\theta } }$
(specified in Section \ref{section1}). 

Let
$\hat{K} = 
2\hat{C}_{5} (\hat{C}_{2} + \hat{C}_{7} ) + \hat{C}_{8}$.  
Then, Lemmas \ref{lemma1.5}, \ref{lemma1.7} 
and the limit process $\varepsilon\rightarrow 0$
imply 
\begin{align*} 
	\limsup_{n\rightarrow \infty } 
	\gamma_{n}^{\hat{p} } |u(\theta_{n} ) |
	\leq 
	\max\{\hat{C}_{2}, \hat{C}_{7} \} 
	(\varphi(\xi) )^{\hat{\mu} }
	\leq 
	\hat{K} (\varphi(\xi) )^{\hat{\mu} }
\end{align*}
on $\Lambda\setminus N_{0}$. 
Consequently, 
Lemma \ref{lemma1.5} yields 
\begin{align*} 
	\limsup_{n\rightarrow \infty } 
	\gamma_{n}^{\hat{p} } \|\nabla f(\theta_{n} ) \|^{2}
	\leq &
	\hat{C}_{5} 
	(\varphi(\xi) )^{\hat{\mu} }
	+
	\hat{C}_{5}
	\limsup_{n\rightarrow \infty } 
	\gamma_{n}^{\hat{p} } 
	\psi(u(\theta_{n} ) )
	\leq 
	\hat{K}  
	(\varphi(\xi) )^{\hat{\mu} }
\end{align*}
on $\Lambda\setminus N_{0}$. 
On the other side, 
using Lemma \ref{lemma1.23}, we 
get 
\begin{align*}
	\limsup_{n\rightarrow \infty } 
	\gamma_{n}^{\hat{q} }
	\|\theta_{n} - \hat{\theta } \| 
	\leq 
	\hat{C}_{8} 
	\varphi(\xi ) 
	\leq 
	\hat{K}\varphi(\xi ) 
\end{align*}
on $\Lambda\setminus N_{0}$. 
Hence, Theorem \ref{theorem1.2} holds, too.
}
\end{vproof}

\section{Proof of Theorem \ref{theorem2.1}} \label{section2*} 

The following notation is used in this section. 
For $\theta \in \mathbb{R}^{d_{\theta } }$, $z \in \mathbb{R}^{d_{z} }$, 
$E_{\theta, z }(\cdot )$ denotes the conditional expectation given $\theta_{0}=\theta$, $Z_{0}=z$. 
For $n\geq 1$, $\xi_{n}$ is the random variable defined as
\begin{align*}
	\xi_{n} 
	=
	F(\theta_{n}, Z_{n+1} ) - \nabla f(\theta_{n} ), 
\end{align*}
while $\xi_{1,n}$, $\xi_{2,n}$, $\xi_{3,n}$ are the random variables defined by 
\begin{align*} 
	&
	\xi_{1,n} 
	=
	\tilde{F}(\theta_{n}, Z_{n+1} ) - (\Pi\tilde{F} )(\theta_{n}, Z_{n} ), 
	\;\;\; 
	\xi_{2,n}
	=
	(\Pi\tilde{F} )(\theta_{n}, Z_{n} ) - (\Pi\tilde{F} )(\theta_{n-1}, Z_{n} ), 
	\;\;\; 
	\xi_{3,n}
	=
	-(\Pi\tilde{F} )(\theta_{n}, Z_{n+1} ).  
\end{align*}
Then, it is straightforward to show that 
algorithm (\ref{2.1}) admits the form (\ref{1.1}). 
On the other side, Assumption \ref{a2.2} yields  
\begin{align} \label{2.3*}
	\sum_{i=n}^{k} 
	\alpha_{i} \gamma_{i}^{r} \xi_{i} 
	= &
	\sum_{i=n}^{k} 
	\alpha_{i} \gamma_{i}^{r} \xi_{1,i} 
	+
	\sum_{i=n}^{k} 
	\alpha_{i} \gamma_{i}^{r} \xi_{2,i} 
	-
	\sum_{i=n}^{k} 
	(\alpha_{i} \gamma_{i}^{r} - \alpha_{i+1} \gamma_{i+1}^{r} ) \xi_{3,i} 
	\nonumber \\
	&
	-
	\alpha_{k+1} \gamma_{k+1}^{r} \xi_{3,k} 
	+
	\alpha_{n} \gamma_{n}^{r} \xi_{3,n-1} 
\end{align}
for $1\leq n \leq k$. 

\begin{lemma} \label{lemma11.1}
Let Assumption \ref{a2.1} hold. 
Then, there exists a real number $s \in (0,1)$ such that 
\linebreak
$\sum_{n=0}^{\infty } \alpha_{n}^{1+s} \gamma_{n}^{r} < \infty$. 
\end{lemma} 

\begin{sproof}
Let $p=(2+2r)/(2+r)$, 
$q=(2+2r)/r$, 
$s=(2+r)/(2+2r)$. 
Then, using the H\"{o}lder inequality, we get
\begin{align*}
	\sum_{n=0}^{\infty } 
	\alpha_{n}^{1+s} \gamma_{n}^{r} 
	=
	\sum_{n=1}^{\infty } 
	(\alpha_{n}^{2} \gamma_{n}^{2r} )^{1/p}
	\left(
	\frac{\alpha_{n} }{\gamma_{n}^{2} } 
	\right)^{1/q} 
	\leq
	\left(
	\sum_{n=1}^{\infty }
	\alpha_{n}^{2} \gamma_{n}^{2r} 
	\right)^{1/p}
	\left(
	\sum_{n=1}^{\infty } 
	\frac{\alpha_{n} }{\gamma_{n}^{2} } 
	\right)^{1/q}. 
\end{align*}
Since 
$\gamma_{n+1}/\gamma_{n} = 1 + \alpha_{n}/\gamma_{n} = O(1)$ 
for $n\rightarrow \infty$ and 
\begin{align*}
	\sum_{n=1}^{\infty } 
	\frac{\alpha_{n} }{\gamma_{n}^{2} }
	=
	\sum_{n=1}^{\infty } 
	\frac{\gamma_{n+1} - \gamma_{n} }{\gamma_{n}^{2} }
	\leq 
	\sum_{n=1}^{\infty } 
	\left(\frac{\gamma_{n+1} }{\gamma_{n} } \right)^{2}
	\int_{\gamma_{n} }^{\gamma_{n+1} } 
	\frac{dt}{t^{2} }
	\leq 
	\frac{1}{\gamma_{1} }
	\max_{n\geq 0} \left(\frac{\gamma_{n+1} }{\gamma_{n} } \right)^{2}, 
\end{align*}
it is obvious that 
$\sum_{n=0}^{\infty } \alpha_{n}^{1+s} \gamma_{n}^{r}$ converges. 
\end{sproof} 

\begin{vproof}{Theorem \ref{theorem2.1}}
{Let $Q \subset \mathbb{R}^{d_{\theta } }$ be an arbitrary compact set, 
while $\Lambda_{Q}=\bigcap_{n=0}^{\infty } \{\theta_{n}\in Q \}$. 
Moreover, let $s\in (0,1)$ be a real number such that 
$\sum_{n=0}^{\infty} \alpha_{n}^{1+s} \gamma_{n}^{r} < \infty$.  
To prove the theorem's assertion, it is sufficient to show that  
$\sum_{n=0}^{\infty } \alpha_{n} \gamma_{n}^{r} \xi_{n}$ 
converges w.p.1 on $\bigcap_{n=0}^{\infty} \{\theta_{n} \in Q \}$. 

Owing to
Assumption \ref{a2.1}, we have 
\begin{align*}
	&
	\alpha_{n-1}^{s} \alpha_{n} \gamma_{n}^{r} 
	=
	\alpha_{n}^{1+s} \gamma_{n}^{r}
	\left(
	1 + \alpha_{n-1} (\alpha_{n}^{-1} - \alpha_{n-1}^{-1} )
	\right)^{s}
	=
	O(\alpha_{n}^{1+s} \gamma_{n}^{r} ) 
\end{align*}
as $n\rightarrow\infty$. 
The same assumption also yields
\begin{align*}
	&
	(\alpha_{n} - \alpha_{n+1} ) \gamma_{n+1}^{r}
	=
	\alpha_{n}^{2}\gamma_{n}^{r}
	(\alpha_{n+1}^{-1} - \alpha_{n}^{-1} )
	\left(1 + \alpha_{n+1} (\alpha_{n}^{-1} - \alpha_{n+1}^{-1} ) \right)
	(1+\alpha_{n}/\gamma_{n} )^{r}
	=
	O(\alpha_{n}^{2} \gamma_{n}^{r} ), 	
	\\
	&
	\alpha_{n} (\gamma_{n+1}^{r} - \gamma_{n}^{r} )
	=
	\alpha_{n} \gamma_{n}^{r} 
	\left(
	(1 + \alpha_{n}/\gamma_{n} )^{r} - 1 
	\right)
	=
	\alpha_{n} \gamma_{n}^{r} 
	\left(
	r\alpha_{n}/\gamma_{n} + o(\alpha_{n}/\gamma_{n} ) 
	\right)
	=
	o(\alpha_{n}^{2} \gamma_{n}^{r} )
\end{align*}
as $n\rightarrow \infty$. 
Hence, $\alpha_{n}\gamma_{n}^{r}-\alpha_{n+1}\gamma_{n+1}^{r} = O(\alpha_{n}^{2}\gamma_{n}^{r})$ 
as $n\rightarrow \infty$. 
Consequently, 
\begin{align}\label{t2.1.701}
	& 
	\sum_{n=0}^{\infty } 
	\alpha_{n}^{s} \alpha_{n+1} \gamma_{n+1}^{r} 	
	< 
	\infty, 
	\;\;\;\;\; 
	\sum_{n=0}^{\infty } 
	|\alpha_{n} \gamma_{n}^{r} - \alpha_{n+1} \gamma_{n+1}^{r} |
	< 
	\infty. 
\end{align}

Let ${\cal F}_{n}=\sigma\{\theta_{0},Z_{0},\dots,\theta_{n},Z_{n} \}$ for $n\geq 0$. 
Since $\{\tau_{Q}>n\}$ is measurable with respect to ${\cal F}_{n}$, 
Assumption \ref{a2.2} implies 
\begin{align*}
	E_{\theta,z}
	\left(
	\xi_{1,n} 
	I_{ \{\tau_{Q}>n \} }
	|
	{\cal F}_{n} 
	\right)
	=
	\left(
	E_{\theta,z}
	(
	\tilde{F}(\theta_{n}, Z_{n+1} )
	|
	{\cal F}_{n} 
	)
	-
	(\Pi\tilde{F} )(\theta_{n},Z_{n} )
	\right) 
	I_{ \{\tau_{Q}>n \} }
	=
	0
\end{align*}
w.p.1 for each $\theta\in\mathbb{R}^{d_{\theta } }$, $z\in\mathbb{R}^{d_{z} }$, $n\geq 0$. 
On the other side, Assumption \ref{a2.3} yields 
\begin{align*}
	\|\xi_{1,n} \| I_{ \{\tau_{Q}>n \} }
	\leq 
	\varphi_{Q}(Z_{n} ) I_{ \{\tau_{Q}>n-1 \} } + \varphi_{Q}(Z_{n+1} ) I_{ \{\tau_{Q}>n \} }  
\end{align*}
for $n\geq 0$. 
Combining this with Assumptions \ref{a2.1}, \ref{a2.3}, we get 
\begin{align*}
	E_{\theta,z}\left(
	\sum_{n=0}^{\infty } \alpha_{n}^{2}\gamma_{n}^{2r} \|\xi_{1,n} \|^{2} I_{ \{\tau_{Q}>n \} }
	\right)
	\leq
	2E_{\theta,z}\left(
	\sum_{n=0}^{\infty } (\alpha_{n}^{2}\gamma_{n}^{2r} + \alpha_{n+1}^{2}\gamma_{n+1}^{2r} ) 
	\varphi_{Q}^{2}(Z_{n+1} ) I_{ \{\tau_{Q}>n \} }
	\right)
	<\infty
\end{align*}
for all $\theta\in\mathbb{R}^{d_{\theta } }$, $z\in\mathbb{R}^{d_{z} }$. 
Then, using Doob theorem, we conclude that 
$\sum_{n=0}^{\infty } \alpha_{n}\gamma_{n}^{r} \xi_{1,n} I_{ \{\tau_{Q}>n \} }$
converges w.p.1. 
Since $\{\tau_{Q}>n\}\subseteq\Lambda_{Q}$ for $n\geq 0$, 
$\sum_{n=0}^{\infty } \alpha_{n}\gamma_{n}^{r}\xi_{1,n}$ converges w.p.1 on 
$\Lambda_{Q}$. 

As a result of Assumption \ref{a2.3} and (\ref{t2.1.701}), we get
\begin{align*}
	E_{\theta,z}
	\left(
	\sum_{n=0}^{\infty } 
	\alpha_{n}\alpha_{n+1}\gamma_{n+1}^{r} \varphi_{Q}^{2}(Z_{n+1} ) I_{ \{\tau_{Q}>n \} }
	\right)
	<\infty
\end{align*}
for all $\theta\in\mathbb{R}^{d_{\theta } }$, $z\in\mathbb{R}^{d_{z} }$. 
On the other side, 
Assumption \ref{a2.3} implies 
\begin{align*}
	\|\xi_{2,n} \| I_{\Lambda_{Q} }
	\leq &
	\varphi_{Q}(Z_{n} ) \|\theta_{n} - \theta_{n-1} \| I_{\Lambda_{Q} } 
	\\
	\leq &
	\alpha_{n-1} \varphi(Z_{n} ) 
	\|F(\theta_{n-1},Z_{n} ) \| I_{\Lambda_{Q} }
	\\
	\leq &
	\alpha_{n-1} \varphi_{Q}^{2}(Z_{n} ) I_{\Lambda_{Q} }
\end{align*}
for $n\geq 1$.  
Thus, 
\begin{align*}
	&
	\sum_{n=1}^{j} \alpha_{n}\gamma_{n}^{r} \|\xi_{2,n} \| I_{\Lambda_{Q} }
	\leq 
	\sum_{n=0}^{\infty} \alpha_{n} \alpha_{n+1}\gamma_{n+1}^{r} 
	\varphi_{Q}^{2}(Z_{n+1} ) I_{ \{\tau_{Q}>n \} }.  
\end{align*}
Therefore, 
$\sum_{n=0}^{\infty } \alpha_{n}\gamma_{n}^{r}\xi_{2,n}$ converges w.p.1 on 
$\Lambda_{Q}$. 

Due to Assumptions \ref{a2.1}, \ref{a2.3} and (\ref{t2.1.701}), we have 
\begin{align*}
	E_{\theta,z}
	\left(
	\sum_{n=0}^{\infty } 
	\alpha_{n}^{2}\gamma_{n}^{2r} \varphi_{Q}^{2}(Z_{n+1} )
	\right)
	<\infty, 
	\;\;\;\;\;
	E_{\theta,z}
	\left(
	\sum_{n=0}^{\infty } |\alpha_{n}\gamma_{n}^{r} - \alpha_{n+1}\gamma_{n+1}^{r} | 
	\varphi_{Q}^{2}(Z_{n+1} ) 
	\right)
	<\infty
\end{align*}
for all $\theta\in\mathbb{R}^{d_{\theta } }$, $z\in\mathbb{R}^{d_{z} }$. 
On the other side, 
owing to Assumption \ref{a2.3}, we have 
\begin{align*}
	&
	\|\xi_{3,n} \| I_{\Lambda_{Q} } 
	\leq 
	\varphi_{Q}(Z_{n+1} ) I_{\Lambda_{Q} }
	\leq 
	\varphi_{Q}^{2}(Z_{n+1} ) I_{\Lambda_{Q} }
\end{align*}
for $n\geq 0$. 
Hence, 
\begin{align*}
	&
	\sum_{n=0}^{\infty} |\alpha_{i}\gamma_{i}^{r} - \alpha_{i+1}\gamma_{i+1}^{r} | \: 
	\|\xi_{3,i} \| I_{\Lambda_{Q} }
	\leq 
	\sum_{n=0}^{\infty} |\alpha_{i}\gamma_{i}^{r} - \alpha_{i+1}\gamma_{i+1}^{r} | 
	\varphi_{Q}^{2}(Z_{i+1} ) I_{ \{\tau_{Q}>i \} }. 
\end{align*}
Consequently, $\sum_{n=0}^{\infty } (\alpha_{n}\gamma_{n}^{r} - \alpha_{n+1}\gamma_{n+1}^{r} ) \xi_{3,n}$
converges w.p.1 on $\Lambda_{Q}$. 
We also get 
\begin{align}\label{ta2.1.23}
	\lim_{n\rightarrow\infty } 
	\alpha_{n+1}\gamma_{n+1}^{r} \|\xi_{3,n} \| I_{\Lambda_{Q} }
	=
	0
\end{align}
w.p.1. 

As $\sum_{n=0}^{\infty } \alpha_{n}\gamma_{n}^{r}\xi_{1,n}$, 
$\sum_{n=1}^{\infty } \alpha_{n}\gamma_{n}^{r}\xi_{2,n}$, 
$\sum_{n=0}^{\infty } (\alpha_{n}\gamma_{n}^{r}-\alpha_{n+1}\gamma_{n+1}^{r} )\xi_{3,n}$ 
are convergent w.p.1 on $\Lambda_{Q}$, 
(\ref{2.3*}), (\ref{ta2.1.23}) imply that 
$\sum_{n=0}^{\infty } \alpha_{n}\gamma_{n}^{r}\xi_{n}$ converges w.p.1 on $\Lambda_{Q}$, too. 
}
\end{vproof}

\section{Proof of Theorems \ref{theorem4.1} and \ref{theorem4.2}} \label{section4*} 

In this section, we use the following notation. 
For 
$\theta \in \mathbb{R}^{d_{\theta } }$, $x \in \mathbb{R}^{N}$, $y \in \mathbb{R}$
and 
$z = [x^{T} \; y]^{T}$, 
let 
\begin{align*}
	F(\theta, z ) = -(y - G_{\theta }(x) ) H_{\theta }(x), 
\end{align*}
while 
$Z_{n+1} = [X_{n}^{T} \; Y_{n} ]^{T}$ for $n\geq 0$. 
With this notation, it is obvious that algorithm (\ref{4.1})
admits the form of (\ref{2.1}). 

\begin{vproof}{Theorem \ref{theorem4.1}}
{Owing to Assumption \ref{a4.2}, 
there exists a real number $K\in [1,\infty )$ such that 
$\max\{\|x\|,|y| \} \leq K$ for any $x\in {\cal X}$, $y\in {\cal Y}$. 

Let $\delta = \varepsilon/(2KN)$, 
while 
\begin{align*}
	\hat{G}_{\eta}(x)
	=
	\sum_{i=1}^{M} c_{i} 
	\hat{\psi}\left(\sum_{j=1}^{N} d_{i,j} x_{j} \right), 
	\;\;\; 
	\hat{H}_{\eta}(x,y) 
	= 
	\frac{1}{2} (y - \hat{G}_{\eta}(x) )^{2}
\end{align*}
and 
$\hat{f}(\eta ) = E(\hat{H}_{\eta}(X_{0},Y_{0} ) )$
for 
$\eta = [c_{1} \cdots c_{M} \; d_{1,1} \cdots d_{M,N} ]^{T} 
\in \mathbb{C}^{d_{\theta } }$, 
$x = [x_{1} \cdots x_{N} ]^{T} \in {\cal X}$, $y\in {\cal Y}$. 
On the other side, let 
$\theta = [a_{1} \cdots a_{M} \; b_{1,1} \cdots b_{M,N} ]^{T}
\in \mathbb{R}^{d_{\theta } }$
be an arbitrary vector. 
Obviously, it is sufficient to show that $\hat{f}(\cdot )$ is analytic on $V_{\delta }(\theta )$
(here, $V_{\delta }(\theta )$ denotes 
$V_{\delta }(\theta ) = \{\eta\in\mathbb{C}^{d_{\theta } }: \|\eta - \theta \| \leq \delta \}$). 

We have 
\begin{align*}
	\left|
	\sum_{j=1}^{N} d_{i,j} x_{j} 
	- 
	\sum_{j=1}^{N} b_{i,j} x_{j} 
	\right|
	\leq 
	K\sum_{j=1}^{N} 
	|d_{i,j} - b_{i,j} |
	\leq  
	\varepsilon/2
\end{align*}
for each   
$\eta = [c_{1} \cdots c_{M} \; d_{1,1} \cdots d_{M,N} ]^{T} \in V_{\delta }(\theta )$, 
$x = [x_{1} \cdots x_{N} ]^{T}\in {\cal X}$, 
$1\leq i \leq M$. 
Hence, 
$\sum_{j=1}^{N} d_{i,j} x_{j} \in V_{\varepsilon/2}(\mathbb{R} )$
whenever $\eta = [c_{1} \cdots c_{M} \; d_{1,1} \cdots d_{M,N} ]^{T} \in V_{\delta }(\theta )$, 
$x = [x_{1} \cdots x_{N} ]^{T}\in {\cal X}$, 
$1\leq i \leq M$. 
Consequently, Assumption \ref{a4.1} implies that
$\hat{G}_{\eta }(x)$ is analytical in $\eta$ and 
continuous in $(\eta,x)$ for all 
$\eta\in V_{\delta}(\theta )$, $x\in {\cal X}$. 
Therefore, 
$\hat{H}_{\eta}(x,y)$ is analytical in $\eta$ and 
continuous in $(\eta,x,y)$ for each 
$\eta\in V_{\delta}(\theta )$, $x\in {\cal X}$, $y\in {\cal Y}$. 
Since $V_{\delta}(\theta ) \times {\cal X} \times {\cal Y}$ is a compact set, 
there exists a real number $L_{1,\theta } \in [1,\infty )$
such that $|\hat{H}_{\eta}(x,y) | \leq L_{1,\theta }$ 
for any $\eta\in V_{\delta}(\theta )$, $x\in {\cal X}$, $y\in {\cal Y}$. 
Then, Cauchy inequality for complex analytic functions 
(see e.g., \cite[Proposition 2.1.3]{taylor}) 
implies that there exists another real number $L_{2,\theta } \in [1,\infty )$
such that $\|\nabla_{\eta } \hat{H}_{\eta}(x,y) \| \leq L_{2,\theta }$ 
for all $\eta\in V_{\delta}(\theta )$, $x\in {\cal X}$, $y\in {\cal Y}$. 
As a result of this and the dominated convergence theorem, 
$\hat{f}(\eta )$ is differentiable for all $\eta\in V_{\delta}(\theta )$. 
Hence, $\hat{f}(\cdot )$ is analytic on $V_{\delta }(\theta )$. 
}
\end{vproof}

\begin{vproof}{Theorem \ref{theorem4.2}}
{As $\{Z_{n} \}_{n\geq 0}$
can be interpreted as a controlled Markov chain whose 
transition kernel $\Pi_{\theta }(z,\cdot )$ does not depend 
on $(\theta,z)$, 
it is straightforward to show that Assumptions \ref{a2.2} 
and \ref{a2.3} hold. 
Then, the theorem's assertion follows directly from 
Theorem \ref{theorem2.1}. 
}
\end{vproof}

\section{Proof of Theorems \ref{theorem7.1} and \ref{theorem7.2}} \label{section7*}

\begin{vproof}{Theorem \ref{theorem7.1}}
{Let 
\begin{align*}	
	\hat{G}_{\eta}(x) = \log\hat{p}_{\eta}(x),
	\;\;\;\;\; 
	\hat{f}(\eta ) = \int \hat{G}_{\eta}(x) p(x) \lambda(dx)
\end{align*} 
for $\eta\in\mathbb{C}^{d_{\theta } }$, $x\in {\cal X}$, 
while   
$\theta\in\Theta$ is an arbitrary vector.  
Obviously, it sufficient to show that $\hat{f}(\cdot )$ is analytic 
in an open vicinity of $\theta$. 

Since $V_{\delta_{\theta } }(\theta ) \times {\cal X}$ is a compact set 
(here, $V_{\delta_{\theta } }(\theta )$ denotes 
$V_{\delta_{\theta } }(\theta ) = \{\eta\in\mathbb{C}^{d_{\theta } }: 
\|\eta-\theta\|\leq \delta_{\theta } \}$, 
while $\delta_{\theta }$ is specified in Assumption \ref{a7.3}), 
Assumptions \ref{a7.2}, \ref{a7.3} imply that there exist 
real numbers $\varepsilon_{\theta} \in (0,\delta_{\theta } )$, 
$L_{1,\theta } \in [1,\infty )$ such that 
$L_{1,\theta}^{-1} \leq |\hat{p}_{\eta}(x) | \leq L_{1,\theta}$
for all $\eta\in V_{\varepsilon_{\theta } }(\theta )$, $x\in {\cal X}$. 
Therefore, $\hat{G}_{\eta}(x)$ is analytic in $\eta$ for all 
$\eta\in V_{\varepsilon_{\theta } }(\theta )$, $x\in {\cal X}$. 
Moreover, $|\hat{G}_{\eta}(x) | \leq \log L_{1,\theta }$ for all 
$\eta\in V_{\varepsilon_{\theta } }(\theta )$, $x\in {\cal X}$. 
Then, using Cauchy inequality for complex analytic functions, 
we deduce that there exists a real number $L_{2,\theta } \in [1,\infty )$
such that $\|\nabla_{\eta } \hat{G}_{\eta}(x) \| \leq L_{2,\theta }$ for all 
$\eta\in V_{\varepsilon_{\theta } }(\theta )$, $x\in {\cal X}$. 
Consequently, the dominated convergence theorem implies that $\hat{f}(\eta )$
is differentiable for all $\eta\in V_{\varepsilon_{\theta} }(\theta )$. 
Hence, $\hat{f}(\cdot )$ is analytic on $V_{\varepsilon_{\theta} }(\theta )$. 
}
\end{vproof}

\begin{vproof}{Theorem \ref{theorem7.2}}
{Similarly as in the proof of Theorem \ref{theorem4.2}, $\{X_{n} \}_{n\geq 0}$
can be interpreted as a controlled Markov chain whose 
transition kernel $\Pi_{\theta }(x,\cdot )$ does not depend 
on $(\theta,x)$. 
Therefore,  
Assumptions \ref{a2.2} and \ref{a2.3} 
are satisfied for algorithm (\ref{7.1}). 
Hence, the theorem's assertion is a straightforward consequence of  
Theorem \ref{theorem2.1}. 
}
\end{vproof}

\section{Proof of Theorems \ref{theorem8.1} and \ref{theorem8.2}} \label{section8*}

In this section, we rely on the following notation. 
Let $d_{w}=2N$, $d_{z}=d_{\theta } + d_{w}$, 
while $W_{n} = [X_{n}^{T} \; X_{n-1}^{T} ]^{T}$, 
$Z_{n} = [Y_{n}^{T} \; X_{n}^{T} \; X_{n-1}^{T} ]^{T}$ for $n\geq 1$. 
Moreover, let 
\begin{align*}
	&
	\tilde{G}_{\theta }(x,x') 
	=
	c(x') + \beta G_{\theta }(x) - G_{\theta }(x'), 
	\;\;\;  
	F(\theta, z ) 
	= 
	-\tilde{G}_{\theta }(x,x') y, 
\end{align*}
for $\theta, y \in \mathbb{R}^{d_{\theta } }$, $x,x' \in {\cal X}$, 
$z = [y^{T} x^{T} (x')^{T} ]^{T}$, 
while 
\begin{align*}
	\\
	&
	\Pi_{\theta }(z,B ) 
	=
	\int I_{B}(\beta y + H_{\theta }(x), x'' , x ) P(x,dx'')
\end{align*}
for the same $\theta,y,x,x',z$ and a measurable set 
$B\subseteq \mathbb{R}^{d_{\theta} }\times {\cal X}\times {\cal X}$. 
Then, it is straightforward to verify that algorithm  
(\ref{8.1}), (\ref{8.3}) admits the form of the recursion 
studied in Section \ref{section2}
(i.e., $\{\theta_{n} \}_{n\geq 0}$, $\{Z_{n} \}_{n\geq 0}$, $\Pi_{\theta}(z,B)$, $F(\theta,z)$
defined here and in Section \ref{section8} satisfy (\ref{2.1}), (\ref{2.3})). 

The following notation is also used in this section.  
Function $B_{\theta }(w)$ is defined by $B_{\theta }(w) = H_{\theta }(x')$
for $\theta\in\mathbb{R}^{d_{\theta } }$, $x,x' \in\mathbb{R}^{N}$, 
$w=[x^{T} \; (x')^{T} ]^{T}$. 
Stochastic processes $\{V_{n}^{\theta } \}_{n\geq 0}$, $\{Z_{n}^{\theta } \}_{n\geq 0}$ 
are recursively defined by 
\begin{align*}
	V_{n+1}^{\theta } 
	=
	\beta V_{n}^{\theta } + B_{\theta }(W_{n+1} )
\end{align*}
and 
$Z_{n}^{\theta } = [(V_{n}^{\theta } )^{T} \; W_{n}^{T} ]^{T}$ 
for $\theta\in\mathbb{R}^{d_{\theta } }$, $n\geq 0$, 
where $V_{0}^{\theta }\in\mathbb{R}^{d_{z} } $ is an arbitrary vector. 
Then, it is straightforward to show 
that $B_{\theta }(w)$ is locally Lipschitz continuous in $(\theta,w)$ 
and 
that $\Pi_{\theta}(\cdot,\cdot )$ 
is a transition kernel of $\{Z_{n}^{\theta } \}_{n\geq 0}$.  

\begin{lemma}\label{lemma8.1}
Let Assumptions \ref{a8.1} -- \ref{a8.3} hold. 
Then, 
\begin{align}\label{l8.1.3*} 
	\lim_{n\rightarrow\infty } (\Pi^{n} F)(\theta,z) 
	=
	\nabla f(\theta )
\end{align}
for all $\theta\in\mathbb{R}^{d_{\theta } }$, 
$z\in\mathbb{R}^{d_{z} }$. 
Moreover, 
for any compact set $Q\subset\mathbb{R}^{d_{\theta } }$, 
there exists a real number $L_{Q}\in [1,\infty )$ such that 
\begin{align}\label{l8.1.1*}
	\|Y_{n} \| 
	I_{ \{\tau_{Q}\geq n \} }
	\leq 
	L_{Q} (1 + \|Y_{0} \| )
\end{align}
for $n\geq 0$
($\tau_{Q}$ is specified in Assumption \ref{a2.3}). 
\end{lemma}

\begin{sproof}
Let $Q\subset\mathbb{R}^{d_{\theta } }$ be an arbitrary compact set. 
Then, owing to Assumption \ref{a8.3}, 
there exists a real number $M_{Q}\in [1,\infty )$ such that 
$\max\{|c(x) |, |G_{\theta }(x) |, \|H_{\theta }(x) \| \} \leq M_{Q}$ 
for all $\theta\in Q$, $x\in{\cal X}$. 
Since 
\begin{align*}
	Y_{n+1} 
	=
	\beta^{n+1} Y_{0} 
	+
	\sum_{k=0}^{n} 
	\beta^{n-k} H_{\theta_{k} }(X_{k} )
\end{align*}
for $n\geq 0$, 
we get 
\begin{align*}
	\|Y_{n+1} \| I_{ \{\tau_{Q}\geq n \} }
	\leq 
	\|Y_{0} \| 
	+
	M_{Q} 
	\sum_{k=0}^{n} \beta^{n-k} 
	\leq 
	\|Y_{0} \| 
	+
	M_{Q} (1-\beta )^{-1}
\end{align*}
for the same $n$. 
Consequently, there exists a real number $L_{Q}\in [1,\infty )$
such that (\ref{l8.1.1*}) is true for $n\geq 0$. 
On the other side, 
it is straightforward to verify 
\begin{align*}
	(\Pi^{n} F )(\theta,z) 
	&
	=
	E(F(\theta, Z_{n+1}^{\theta } )|Z_{1}^{\theta } = z )
	\\
	&
	=
	\begin{aligned}[t]
	-
	E\left(
	\tilde{G}_{\theta }(X_{n+1}, X_{n} )
	\left. 
	\left( 
	\beta^{n} y 
	+
	\sum_{k=0}^{n-1} \beta^{k} H_{\theta }(X_{n-k} ) 
	\right)
	\right|
	X_{1}=x
	\right)
	\end{aligned}
	\\
	&
	=
	-
	\sum_{k=0}^{n-1}
	\beta^{k}  
	\int \tilde{G}_{k,\theta }(x'') H_{\theta }(x'') P^{n-k-1}(x,dx'') 
	+ 
	\beta^{n} \tilde{G}_{n-1,\theta }(x) y
\end{align*}
for all $\theta,y\in\mathbb{R}^{d_{\theta } }$, $x,x'\in{\cal X}$, 
$z=[y^{T} x^{T} (x')^{T} ]^{T}$, $n\geq 1$, 
where 
\begin{align*}
	\tilde{G}_{k,\theta }(x) 
	=
	(P^{k} c )(x) 
	+
	\beta (P^{k+1} G )_{\theta }(x) 
	-
	(P^{k} G )_{\theta }(x). 
\end{align*}
It is also easy to show 
\begin{align*}
	\nabla f(\theta )
	=
	-
	\int (g(x) - G_{\theta }(x) ) H_{\theta }(x) \pi(dx)
	=
	-
	\sum_{k=0}^{\infty } \beta^{k} 
	\int \tilde{G}_{k,\theta }(x) H_{\theta }(x) \pi(dx)
\end{align*}
for each $\theta\in\mathbb{R}^{d_{\theta } }$. 
As $\|\tilde{G}_{k,\theta }(x) H_{\theta }(x) \| \leq 3M_{Q}^{2}$
for any $\theta\in Q$, $x\in{\cal X}$, $k\geq 0$, 
Assumption \ref{a8.2} implies 
\begin{align*}
	\left\|
	\int \tilde{G}_{k,\theta }(x') H_{\theta }(x') (P^{l}-\pi)(x,dx') 
	\right\|
	\leq 
	3CM_{Q}^{2} \rho^{l}
\end{align*}
for all $\theta\in Q$, $x\in{\cal X}$, $k,l\geq 0$. 
Consequently, 
\begin{align*}
	\|(\Pi^{n} F)(\theta,z) - \nabla f(\theta ) \|
	\leq &
	\sum_{k=0}^{n-1} 
	\beta^{k} 
	\left\|
	\int \tilde{G}_{k,\theta }(x'') H_{\theta }(x'') (P^{n-k-1}-\pi)(x,dx'') 
	\right\|
	\\
	&
	+ \!
	\sum_{k=n}^{\infty} 
	\beta^{k} 
	\left\|
	\! \int \! \tilde{G}_{k,\theta }(x'') H_{\theta }(x'') \pi(x,dx'') 
	\right\|
	\!+
	\beta^{n} \|\tilde{G}_{n-\!1,\theta}(x) \| \|y\|
	\\
	\leq &
	3CM_{Q}^{2} 
	\sum_{k=0}^{n-1} \beta^{k} \rho^{n-k-1} 
	+
	3M_{Q}^{2} \beta^{n} 
	(\|y\| + (1-\beta )^{-1} )
\end{align*}
for each $\theta\in Q$, $y\in\mathbb{R}^{d_{\theta } }$, $x,x'\in{\cal X}$, 
$z=[y^{T} x^{T} (x')^{T} ]^{T}$, $n\geq 1$. 
Hence, (\ref{l8.1.3*}) holds for all $\theta\in\mathbb{R}^{d_{\theta } }$, 
$z\in\mathbb{R}^{d_{z} }$. 
\end{sproof}

\begin{vproof}{Theorem \ref{theorem8.1}}
{Let 
\begin{align*}
	\hat{H}_{\eta}(x) = 2^{-1} (g(x) - \hat{G}_{\eta}(x) )^{2},
	\;\;\;\;\; 
	\hat{f}(\eta ) = \int \hat{H}_{\eta}(x) \pi(dx)
\end{align*}
for $\eta\in\mathbb{C}^{d_{\theta } }$, $x\in {\cal X}$, 
while   
$\theta\in\Theta$ is an arbitrary vector.  
Obviously, it sufficient to show that $\hat{f}(\cdot )$ is analytic 
on 
$V_{\delta_{\theta } }(\theta )= \{\eta\in\mathbb{C}^{d_{\theta } }: 
\|\eta-\theta\|\leq \delta_{\theta } \}$
($\delta_{\theta }$ is specified in Assumption \ref{a8.3}). 

Owing to Assumption \ref{a8.3}, 
$\hat{H}_{\eta}(x)$ is analytic in $\eta$ 
for all $\eta\in V_{\delta_{\theta } }(\theta )$. 
Due to the same assumption, there exists a real number 
$L_{1,\theta } \in [1,\infty )$ such that 
$|\hat{H}_{\eta}(x) | \leq L_{1,\theta}$
for all $\eta\in V_{\delta_{\theta } }(\theta )$, $x\in {\cal X}$. 
Combining this with Cauchy inequality for complex analytic functions, 
we deduce that there exists a real number $L_{2,\theta } \in [1,\infty )$
such that $\|\nabla_{\eta } \hat{H}_{\eta}(x) \| \leq L_{2,\theta }$ for all 
$\eta\in V_{\delta_{\theta } }(\theta )$, $x\in {\cal X}$. 
Consequently, the dominated convergence theorem implies that
$\hat{f}(\eta )$
is differentiable for all $\eta\in V_{\delta_{\theta} }(\theta )$. 
Thus, $\hat{f}(\cdot )$ is analytic on $V_{\delta_{\theta} }(\theta )$. 
}
\end{vproof}

\begin{vproof}{Theorem \ref{theorem8.2}}
{Owing to Assumptions \ref{a8.1} -- \ref{a8.3}, 
$\{V_{n}^{\theta } \}_{n\geq 0}$ and $\{Z_{n}^{\theta } \}_{n\geq 0}$ defined here
satisfy all conditions of Theorem \ref{theorema4.1} (Appendix \ref{appendix4}). 
Moreover, for any compact set $Q\subset\mathbb{R}^{d_{\theta } }$, 
there exists a real number $K_{Q} \in [1,\infty )$ 
such that (\ref{ta4.1.1*}) -- (\ref{ta4.1.5*}) are satisfied for $p=1$, $K_{2,Q}=K_{Q}$ and all 
$\theta,\theta',\theta''\in Q$, 
$z,z',z'' \in \mathbb{R}^{d_{\theta } }\times {\cal X}\times {\cal X}$. 
Consequently, Theorem \ref{theorema4.1} and Lemma \ref{lemma8.1} imply 
that Assumptions \ref{a2.2} and \ref{a2.3} hold. 
Then, the theorem's assertion directly follows from 
Theorem \ref{theorem2.1}. 
}
\end{vproof}

\section{Proof of Theorems \ref{theorem10.1} and \ref{theorem10.2}}\label{section10*}

In this section, we rely on the following notation. 
For $\theta\in\Theta$, $z=(x_{1:N},y_{1:N} ) \in {\cal X}^{N}\times {\cal Y}^{N}$, 
let 
$F(\theta, z ) = \psi_{N,\theta }(y_{1:N} )$, 
while $f(\theta )=f_{N}(\theta )$. 
For $n\geq 0$, let 
$Z_{n} = (X_{nN+1:(n+1)N}, Y_{nN+1:(n+1)N} )$. 
Obviously, 
$\{Z_{n} \}_{n\geq 0}$ is a Markov chain. 
Let $\Pi(z,z')$ be the transition kernel of $\{Z_{n} \}_{n\geq 0}$. 
Then, it is easy to show that algorithm (\ref{10.1}) admits the form of the recursion 
studied in Section \ref{section2}
(i.e., $\{\theta_{n} \}_{n\geq 0}$, $\{Z_{n} \}_{n\geq 0}$, 
$\Pi(z,z')$, $F(\theta, z)$ defined here and in Section \ref{section10} satisfy 
(\ref{2.1}), (\ref{2.3})). 

\begin{vproof}{Theorem \ref{theorem10.1}}
{Using \cite[Theorem 1, Proposition 1]{tadic5}, 
it can easily be demonstrated that $f_{\infty }(\cdot )$ is real-analytic on entire $\Theta$
(notice that all conditions of \cite[Theorem 1, Proposition 1]{tadic5} 
hold when Assumptions \ref{a10.1}, \ref{a10.2} are satisfied). 

Owing to Assumption \ref{a10.1}, 
$\{Z_{n} \}_{n\geq 0}$ is geometrically ergodic. 
Let $\pi(\cdot )$, $\nu(\cdot )$ be the invariant probabilities of 
$\{X_{n} \}_{n\geq 0}$, 
$\{Z_{n} \}_{n\geq 0}$ (respectively). 
Then, there exist real numbers $\rho\in(0,1)$, $C\in[1,\infty )$
such that 
\begin{align*}
	|p^{n}(x'|x) - \pi(x') |
	\leq 
	C\rho^{n}, 
	\;\;\;\;\; 
	|\Pi^{n}(z,z') - \nu(z') |
	\leq 
	C\rho^{n}
\end{align*}
for each $x,x'\in{\cal X}$, $z,z'\in{\cal X}^{N}\times {\cal Y}^{N}$, $n\geq 0$. 
Therefore, 
\begin{align*}
	f_{N}(\theta )
	=
	\sum_{\stackrel{\scriptstyle x_{1:N}\in {\cal X}^{N} }{y_{1:N}\in {\cal Y}^{N} } } 
	\phi_{N,\theta }(y_{1:N} ) 
	\nu(x_{1:N}, y_{1:N} )
\end{align*}
for all $\theta\in\Theta$. 
On the other side, 
Assumptions \ref{a10.2}, \ref{a10.3} imply that for each $y_{1:N}\in {\cal Y}^{N}$, 
$\phi_{N,\theta }(y_{1:N} )$ is real-analytic in $\theta$ on entire $\Theta$. 
Consequently, $f_{N}(\cdot )$ is real-analytic on entire $\Theta$, too. 

Let $Q\subset\mathbb{R}^{d_{\theta } }$ be any compact set. 
Moreover, let ${\cal P}^{N_{x} }$ be the set of $N_{x}$-dimensional probability vectors, 
while $e=[1\cdots 1]^{T}\in \mathbb{R}^{N_{x} }$. 
For $\theta\in\Theta$, $x,x'\in{\cal X}$, $y\in{\cal Y}$, let 
\begin{align*}
	r_{\theta }(y,x'|x)
	=
	q_{\theta }(y|x') p_{\theta }(x'|x), 
\end{align*}
while $R_{\theta }(y)$ is the $N_{x}\times N_{x}$ matrix 
whose $(i,j)$ entry is $r_{\theta }(y,i|j)$. 
For $\theta\in\Theta$, $y\in{\cal Y}$, $u\in{\cal P}^{N_{x} }$, 
$V\in\mathbb{R}^{N_{x}\times N_{x} }$, let 
\begin{align*}
	\Phi_{\theta }(y,u)
	=
	\log(e^{T} R_{\theta }(y) u ), 
	\;\;\;\;\; 
	\Psi_{\theta }(y,u,V)
	=
	\nabla_{\theta }\Phi_{\theta }(y,u) 
	+
	V\: \nabla_{u}\Phi_{\theta }(y,u). 
\end{align*}
Then, owing to Assumption \ref{a10.2}, there exists a real number $\delta_{Q}\in(0,1)$ such that 
\begin{align}\label{t10.1.71}
	r_{\theta}(y,x'|x)
	\geq
	\delta_{Q}
\end{align}
for all $\theta\in Q$, $x,x'\in{\cal X}$, $y\in{\cal Y}$. 
Combining this with Assumption \ref{a10.3}, 
we conclude that there exists a real number $\tilde{C}_{1,Q}\in[1,\infty )$ such that 
\begin{align}
	&\label{t10.1.73}
	\|\Psi_{\theta}(y,u,V) \|
	\leq 
	\tilde{C}_{1,Q} (1 + \|V\| ), 
	\\
	&\label{t10.1.75}
	|\Phi_{\theta}(y,u' ) - \Phi_{\theta}(y,u'' ) |
	\leq 
	\tilde{C}_{1,Q} \|u' - u'' \|, 
	\\
	&\label{t10.1.77}
	\|\Psi_{\theta}(y,u',V') - \Psi_{\theta}(y,u'',V'') \|
	\leq 
	\tilde{C}_{1,Q} 
	(\|u'-u''\| + \|V'-V''\| ) (1 + \|V'\| + \|V''\| )
\end{align}
for all $\theta\in Q$, $x,x'\in{\cal X}$, $y\in{\cal Y}$, 
$u, u', u'' \in{\cal P}^{N_{x} }$, $V, V', V'' \in\mathbb{R}^{N_{x}\times N_{x} }$.

For $\theta\in\Theta$, $y_{1:n}\in{\cal Y}^{n}$, $n\geq 1$, 
let $u_{0,\theta }$, $u_{n,\theta }(y_{1:n} )$ be the $N_{x}$-dimensional vectors 
whose $i$-th components are 
\begin{align*}
	u_{0,i,\theta }
	=
	\pi_{\theta }(i), 
	\;\;\;\;\; 
	u_{n,i,\theta }(y_{1:n} )
	=
	P_{\theta }(X_{n}^{\theta } = i|Y_{1:n}^{\theta } = y_{1:n} )
\end{align*}
(notice that $\{ u_{n,\theta }(y_{1:n} ) \}_{n\geq 1}$ is the optimal filter for the model 
$\{(X_{n}^{\theta }, Y_{n}^{\theta } ) \}_{n\geq 0})$. 
For the same $\theta$, $y_{1:n}$, $n$, let 
\begin{align*}
	V_{0,\theta } = \nabla_{\theta } u_{0,\theta }, 
	\;\;\;\;\; 
	V_{n,\theta }(y_{1:n} ) = \nabla_{\theta } u_{n,\theta }(y_{1:n} ). 
\end{align*}
Then, it is straightforward to verify 
\begin{align*}
	\log P_{\theta }(Y_{1}^{\theta } = y )
	=
	\Phi_{\theta }(y,u_{0,\theta } ), 
	\;\;\;\;\; 
	\log\left(
	\frac{P_{\theta }(Y_{1:n+1}^{\theta } = y_{1:n+1} ) }{P_{\theta }(Y_{1:n}^{\theta } = y_{1:n} ) }
	\right)
	=
	\Phi_{\theta }(y_{n+1}, u_{n,\theta }(y_{1:n} ) )
\end{align*}
for $\theta\in\Theta$, $y\in{\cal Y}$, $y_{1:n+1}=(y_{1},\dots,y_{n+1} )\in{\cal Y}^{n+1}$, $n\geq 0$. 
As 
\begin{align*}
	\phi_{n,\theta }(y_{1:n} ) 
	=
	-
	\frac{1}{n}
	\left(
	\log P_{\theta }(Y_{1}^{\theta } = y_{1} )
	+
	\sum_{i=1}^{n-1}
	\log\left(
	\frac{P_{\theta }(Y_{1:i+1}^{\theta } = y_{1:i+1} ) }{P_{\theta }(Y_{1:i}^{\theta } = y_{1:i} ) }
	\right)
	\right)
\end{align*}
for $\theta\in\Theta$, $y_{1:n}=(y_{1},\dots,y_{n} )\in{\cal Y}^{n}$, $n\geq 1$, 
we conclude 
\begin{align}\label{t10.1.3'}
	\phi_{n,\theta }(y_{1:n} )
	=
	-
	\frac{1}{n}
	\sum_{i=0}^{n-1} 
	\Phi_{\theta }(y_{i+1}, u_{i,\theta}(y_{1:i} ) )
\end{align} 
for the same $\theta$, $y_{1:n}$, $n$. 
Differentiating (\ref{t10.1.3'}) (in $\theta$), we get
\begin{align}\label{t10.1.3''}
	\psi_{n,\theta }(y_{1:n} )
	=
	-
	\frac{1}{n}
	\sum_{i=0}^{n-1} 
	\Psi_{\theta }(y_{i+1}, u_{i,\theta}(y_{1:i} ), V_{i,\theta}(y_{1:i} ) ) 
\end{align}
for $\theta\in\Theta$, $y_{1:n} = (y_{1},\dots,y_{n} ) \in {\cal Y}^{n}$, 
$n\geq 0$. 

Let $U_{0}^{\theta } = u_{0,\theta }$, $U_{n}^{\theta } = u_{n,\theta }(Y_{1:n} )$ and  
$V_{0}^{\theta } = V_{0,\theta }$, 
$V_{n}^{\theta } = V_{n,\theta }(Y_{1:n} )$ for $\theta\in\Theta$, $n\geq 1$. 
Then, using \cite[Theorems 4.1, 4.2]{tadic&doucet} and (\ref{t10.1.71}) -- (\ref{t10.1.77}), 
we conclude that 
$\{(X_{n+1}, Y_{n+1}, U_{n}^{\theta }, V_{n}^{\theta } ) \}_{n\geq 0}$
is geometrically ergodic for each $\theta\in\Theta$. 
We also deduce that there exist functions $g(\cdot )$, $h(\cdot )$
and 
real numbers $\varepsilon_{Q}\in(0,1)$, 
$\tilde{C}_{2,Q}\in[1,\infty )$ (depending on $\rho$, $\delta_{Q}$, $C$, $\tilde{C}_{1,Q}$) 
such that 
\begin{align}\label{t10.1.5}
	\max\{
	|E(\Phi_{\theta }(Y_{n+1}, U_{n}^{\theta } ) ) - g(\theta ) |, 
	\|E(\Psi_{\theta }(Y_{n+1}, U_{n}^{\theta }, V_{n}^{\theta } ) ) - h(\theta ) \|
	\}
	\leq 
	\tilde{C}_{2,Q} \varepsilon_{Q}^{n}
\end{align}
for all $\theta\in Q$, $n\geq 0$. 
As a result of (\ref{t10.1.3'}) -- (\ref{t10.1.5}), we get 
\begin{align*}
	g(\theta ) 
	=
	\lim_{n\rightarrow\infty } 
	E(\phi_{n,\theta }(Y_{1:n} ) ), 
	\;\;\;\;\; 
	h(\theta )
	=
	\lim_{n\rightarrow\infty } 
	E(\psi_{n,\theta }(Y_{1:n} ) )
	=
	\lim_{n\rightarrow\infty } 
	\nabla_{\theta } 
	E(\phi_{n,\theta }(Y_{1:n} ) )
\end{align*}
for all $\theta\in\Theta$. 
Therefore, $g(\cdot )=f_{\infty }(\cdot )$, $h(\cdot ) = \nabla f_{\infty }(\cdot )$
(notice that $E(\psi_{n,\theta }(Y_{1:n} ) )$ converges to $h(\theta )$ 
uniformly in $\theta$ on each compact subset of $\Theta$). 

In the rest of the proof, we assume that $\{X_{n} \}_{n\geq 0}$ is in steady-state
(i.e., $X_{0}$ is distributed according to $\pi(\cdot )$). 
Then, we have 
\begin{align*}
	f_{N}(\theta ) 
	=
	E(\phi_{N,\theta }(Y_{1:N} ) ), 
	\;\;\;\;\; 
	\nabla f_{N}(\theta ) 
	=
	E(\psi_{N,\theta }(Y_{1:N} ) )
\end{align*}
for each $\theta\in\Theta$. 
Combining this with (\ref{t10.1.3'}) -- (\ref{t10.1.5}), we get 
\begin{align*}
	&
	|f_{N}(\theta ) - f_{\infty }(\theta ) |
	=
	\left|
	\frac{1}{N} \sum_{i=0}^{N-1} 
	\left(
	E(\Phi_{\theta }(Y_{i+1}, U_{i}^{\theta } ) )
	-
	g(\theta )
	\right)
	\right|
	\leq 
	\frac{\tilde{C}_{2,Q} }{N}
	\sum_{i=1}^{N} \varepsilon_{Q}^{i}
	\leq 
	\frac{\tilde{C}_{2,Q} }{(1 - \varepsilon_{Q} ) N }, 
	\\
	&
	\|\nabla f_{N}(\theta ) - \nabla f_{\infty }(\theta ) \|
	=
	\left\|
	\frac{1}{N} \sum_{i=0}^{N-1} 
	\left(
	E(\Psi_{\theta }(Y_{i+1}, U_{i}^{\theta }, V_{i}^{\theta } ) )
	-
	h(\theta )
	\right)
	\right\|
	\leq 
	\frac{\tilde{C}_{2,Q} }{N}
	\sum_{i=1}^{N} \varepsilon_{Q}^{i}
	\leq 
	\frac{\tilde{C}_{2,Q} }{(1 - \varepsilon_{Q} ) N }
\end{align*}
for all $\theta\in Q$. 
Then, it can easily be deduced that for each $\theta\in\Theta$, 
there exists $L_{\theta }\in (0,\infty )$ such that 
(\ref{t10.1.1*}) holds. 
}
\end{vproof}

\begin{vproof}{Theorem \ref{theorem10.2}}
{Let $\nu(z)$ have the same meaning as in the proof of Theorem \ref{theorem10.1}, 
while $\tilde{\Pi}^{n}(z,z')= \Pi^{n}(z,z') - \nu(z')$ 
for $z,z'\in {\cal X}^{N}\times {\cal Y}^{N}$, $n\geq 0$. 
Since $\{Z_{n} \}_{n\geq 0}$ is geometrically ergodic, 
there exist real numbers $\rho\in(0,1)$, $C\in[1,\infty )$ such that 
$|\tilde{\Pi}^{n}(z,z')| \leq C\rho^{n}$ 
for each $z,z'\in {\cal X}^{N}\times {\cal Y}^{N}$, $n\geq 0$. 

Let $Q\subset\Theta$ be an arbitrary compact set, 
while $C_{Q}\in[1,\infty )$ stands for an upper bound of 
$\|F(\cdot, z)\|$ on $Q$ and 
for a Lipschitz constant of $F(\cdot,z)$ on the same set
(here, $z$ is any element of ${\cal X}^{N}\times{\cal Y}^{N}$). 
For $\theta\in\Theta$, $z\in{\cal X}^{N}\times{\cal Y}^{N}$, $n\geq 0$, 
let 
\begin{align*}
	(\tilde{\Pi}^{n} F)(\theta,z)
	=
	\sum_{z'\in{\cal X}^{N}\times{\cal Y}^{N} }
	F(\theta, z' ) \tilde{\Pi}^{n}(z,z'), 
\end{align*}
while $\tilde{F}(\theta, z ) = \sum_{n=0}^{\infty } (\tilde{\Pi}^{n} F)(\theta,z)$. 
Then, we have 
\begin{align*}
	&
	\|(\tilde{\Pi}^{n} F)(\theta,z) \|
	\leq 
	CC_{Q}\rho^{n}, 
	\\
	&
	\|(\tilde{\Pi}^{n} F)(\theta',z) - (\tilde{\Pi}^{n} F)(\theta'',z) \|
	\leq 
	CC_{Q}\rho^{n} \|\theta' - \theta'' \|
\end{align*}
for all $\theta,\theta',\theta''\in Q$, $z\in{\cal X}^{N}\times{\cal Y}^{N}$, $n\geq 0$. 
Therefore, $\sum_{n=0}^{\infty } \|(\tilde{\Pi}^{n} F)(\theta,z) \| < \infty$
for any $\theta,\in \Theta$, $z\in{\cal X}^{N}\times{\cal Y}^{N}$. 
Consequently, for each $\theta,\in \Theta$, $z\in{\cal X}^{N}\times{\cal Y}^{N}$, 
$\tilde{F}(\theta,z)$ is well-defined and satisfies
$(\Pi\tilde{F} )(\theta, z ) = \sum_{n=1}^{\infty } (\tilde{\Pi}^{n} F)(\theta,z)$. 
Thus, Assumption \ref{a2.2} holds. 
We also have 
\begin{align*}
	&
	\|\tilde{F}(\theta,z) \|
	\leq 
	CC_{Q}(1-\rho )^{-1}, 
	\\
	&
	\|\tilde{F}(\theta',z) - \tilde{F}(\theta'',z) \|
	\leq 
	CC_{Q}(1-\rho )^{-1} \|\theta' - \theta'' \|
\end{align*}
for all $\theta,\theta',\theta''\in Q$, $z\in{\cal X}^{N}\times{\cal Y}^{N}$. 
Hence, Assumption \ref{a2.3} holds, too. 
Then, the theorem's assertion directly follows from Theorem \ref{theorem2.1}. 
}
\end{vproof}

\section{Proof of Theorems \ref{theorem6.1} and \ref{theorem6.2}} \label{section6*}

In this section, we use the following notation. 
$d_{v}$, $d_{w}$, $d_{z}$ are integers defined by 
$d_{v} = (M+N)(N+1)$, 
$d_{w} = L+1$, 
$d_{z} = d_{v}+d_{w}$. 
Stochastic processes $\{\varepsilon_{n}^{\theta } \}_{n\geq 0}$, 
$\{\phi_{n}^{\theta } \}_{n\geq 0}$, 
$\{\psi_{n}^{\theta } \}_{n\geq 0}$
are recursively defined by 
\begin{align*}
	&
	\phi_{n}^{\theta } 
	=
	[Y_{n} \cdots Y_{n-M+1} \; \varepsilon_{n}^{\theta } \cdots \varepsilon_{n-N+1}^{\theta } ]^{T}, 
	\\
	&
	\psi_{n+1}^{\theta } 
	=
	\phi_{n}^{\theta } 
	-
	[\psi_{n}^{\theta } \cdots \psi_{n-N+1}^{\theta } ] D \theta, 
	\\
	& 
	\varepsilon_{n+1}^{\theta } 
	=
	Y_{n+1} - (\phi_{n}^{\theta } )^{T} \theta
\end{align*}
for $n\geq 0$, $\theta\in\Theta$, 
where $\varepsilon_{0}^{\theta }, \dots, \varepsilon_{-N+1}^{\theta } \!\in \mathbb{R}$
are arbitrary numbers 
and $\psi_{0}^{\theta }, \dots, \psi_{-N+1}^{\theta } \!\in \mathbb{R}^{d_{\theta } }$
are arbitrary vectors.  
$\{V_{n}^{\theta } \}_{n\geq 0}$, $\{Z_{n}^{\theta } \}_{n\geq 0}$ are stochastic processes 
defined by 
\begin{align*}
	V_{n}^{\theta } 
	=
	[Y_{n} \cdots Y_{n-M+1} \; \varepsilon_{n}^{\theta } \cdots \varepsilon_{n-N+1}^{\theta } 
	\; (\psi_{n}^{\theta } )^{T} \cdots (\psi_{n-N+1}^{\theta } )^{T} ]^{T}
\end{align*}
and $Z_{n}^{\theta } = [(V_{n}^{\theta } )^{T} \; W_{n}^{T} ]^{T}$
for $n\geq 0$, $\theta\in\Theta$. 
Similarly, stochastic processes 
$\{V_{n} \}_{n\geq 0}$, $\{Z_{n} \}_{n\geq 0}$ are defined as 
\begin{align}\label{6.501*}
	V_{n} 
	=
	[Y_{n} \cdots Y_{n-M+1} \; \varepsilon_{n} \cdots \varepsilon_{n-N+1} 
	\; \psi_{n}^{T} \cdots \psi_{n-N+1}^{T} ]^{T}
\end{align}
and $Z_{n} = [V_{n}^{T} \; W_{n}^{T} ]^{T}$ for $n\geq 0$. 
Then, it can easily be deduced that there exists a matrix $B\in\mathbb{R}^{d_{v}\times d_{w} }$ 
and a function $A_{\theta }$ mapping $\theta\in\Theta$ to $\mathbb{R}^{d_{v}\times d_{v} }$
such that the following holds: 
\begin{compactenum}[(i)]
\item
$A_{\theta }$ is linear in $\theta$. 
\item
The eigenvalues of $A_{\theta }$ lie in $\{z\in\mathbb{C}: |z|<1 \}$
for all $\theta\in\Theta$. 
\item
$V_{n+1} = A_{\theta_{n} } V_{n} + B W_{n+1}$
and 
$V_{n+1}^{\theta } = A_{\theta } V_{n}^{\theta } + B W_{n+1}$
for each $n\geq 0$, $\theta\in\Theta$. 
\end{compactenum}

In this section, besides the notation introduced in the previous paragraph, 
we also rely on the following notation. 
$F(\theta,z)$, $\phi(z)$ are the functions defined by 
\begin{align*}
	F(\theta,z)
	=
	-\tilde{\psi}_{1} \tilde{\varepsilon}_{1}, 
	\;\;\;\;\; 
	\phi(z)
	=
	\frac{1}{2} \tilde{\varepsilon}_{1}^{2}
\end{align*}
for $\theta\in\Theta$, 
$y_{1},\dots,y_{M} \in \mathbb{R}$, $\tilde{\varepsilon}_{1},\dots,\tilde{\varepsilon}_{N}\in \mathbb{R}$, 
$\tilde{\psi}_{1},\dots,\tilde{\psi}_{N} \in \mathbb{R}^{d_{\theta } }$, 
$w\in {\cal W}$, 
$v=[y_{1} \cdots y_{M} \; \tilde{\varepsilon}_{1} \cdots \tilde{\varepsilon}_{N} \; 
\tilde{\psi}_{1}^{T} \cdots \tilde{\psi}_{N}^{T} ]^{T}$, 
$z=[v^{T} \; w^{T} ]^{T}$
(here, $y_{i}$, $\tilde{\varepsilon}_{j}$, $\tilde{\psi}_{k}$
are deterministic variables corresponding to 
$Y_{n-i+1}$, $\varepsilon_{n-j+1}$, $\psi_{n-k+1}$ in (\ref{6.501*})). 
$\Pi_{\theta }(z,B)$ is the transition kernel defined as 
\begin{align*}
	\Pi_{\theta }(z,B)
	=
	\int I_{B}(A_{\theta } v + Bw', w' ) P(w,dw')
\end{align*}
for a measurable set $B\subseteq \mathbb{R}^{d_{v} } \times {\cal W}$ and 
$v\in\mathbb{R}^{d_{v} }$, $w\in {\cal W}$, $z=[v^{T} \; w^{T}]^{T}$. 
Then, it is easy to show that algorithm  
(\ref{6.1}) -- (\ref{6.7}) admits the form of the recursion 
studied in Section \ref{section2}
(i.e., $\{\theta_{n} \}_{n\geq 0}$, $\{Z_{n} \}_{n\geq 0}$, 
$\Pi_{\theta}(z,B)$, $F(\theta,z)$
defined here and in Section \ref{section6} satisfy (\ref{2.1}), (\ref{2.3})). 
It is also straightforward to vefity that 
$\Pi_{\theta}(z,B)$ is a transition kernel of 
$\{Z_{n}^{\theta } \}_{n\geq 0}$ for all $\theta\in\Theta$
and that $B_{\theta }(q) \varepsilon_{n}^{\theta } = A_{\theta }(q) Y_{n}$
for each $\theta\in\Theta$, $n\geq 0$. 
In addition to this, it is easy to demonstrate that if 
$\varepsilon_{0}^{\theta } = \cdots = \varepsilon_{-N+1}^{\theta } = 0$, 
$\psi_{0}^{\theta } = \cdots = \psi_{-N+1}^{\theta } = 0$ for all $\theta\in\Theta$, 
then 
$\psi_{n}^{\theta } = -\nabla_{\theta } \varepsilon_{n}^{\theta }$ 
for each $\theta\in\Theta$, $n\geq 0$. 
Consequently, if $\varepsilon_{0}^{\theta } = \cdots = \varepsilon_{-N+1}^{\theta } = 0$, 
$\psi_{0}^{\theta } = \cdots = \psi_{-N+1}^{\theta } = 0$ for all $\theta\in\Theta$, 
then 
\begin{align} \label{6.1*}
	& 
	(\Pi^{n} \phi)(\theta, 0 ) 
	= 
	\frac{1}{2}
	E\big((\varepsilon_{n}^{\theta } )^{2} \big), 
	\;\;\;\;\; 
	(\Pi^{n} F)(\theta, 0 ) 
	=
	\frac{1}{2}
	\nabla_{\theta } 
	E\big((\varepsilon_{n}^{\theta } )^{2} \big) 
\end{align}
for each $\theta\in\Theta$, $n\geq 0$.

\begin{vproof}{Theorem \ref{theorem6.1}}
{Let 
$r_{k} = r_{-k} = \lim_{n\rightarrow\infty } \text{Cov}(Y_{n},Y_{n+k} )$ 
for $k\geq 0$, 
while $m=\lim_{n\rightarrow\infty } E(Y_{n} )$. 
Moreover, let 
$
	\varphi(\omega )
	=
	\sum_{k=-\infty }^{\infty } r_{k} e^{-i\omega k }
$
for $\omega \in [-\pi,\pi]$. 
Then, Assumptions \ref{a6.1}, \ref{a6.2} imply 
$\sum_{k=0}^{\infty } |r_{k} | < \infty$, 
and consequently, 
$\varphi(\cdot )$ is real-analytic on $[-\pi,\pi ]$. 

Let $\theta \in \Theta$ be an arbitrary vector, 
while $C_{\theta }(z) = A_{\theta }(z)/B_{\theta }(z)$, 
for $z \in \mathbb{C}$. 
Since $\varepsilon_{n}^{\theta } = C_{\theta }(q) Y_{n}$ 
for $n\geq 0$, 
and since $C_{\theta }(\cdot )$ has poles only in 
$\{z\in \mathbb{C}: |z|< 1 \}$, 
the spectral theory for stationary processes 
(see e.g., \cite[Chapter II]{ljung3}) 
yields 
$\lim_{n\rightarrow \infty } E(\varepsilon_{n}^{\theta } ) 
= m C_{\theta }(1)$
and 
\begin{align*}
	&
	\lim_{n\rightarrow \infty } 
	\text{Cov}(\varepsilon_{n}^{\theta }, \varepsilon_{n+k}^{\theta } )
	=
	\frac{1}{2\pi} 
	\int_{-\pi}^{\pi} 
	|C_{\theta }(e^{i\omega } ) |^{2} \varphi(\omega ) e^{i\omega k } d\omega
\end{align*}
for $k\geq 0$.  
Therefore, 
\begin{align*} 
	f(\theta ) 
	=
	\frac{1}{2} 
	\lim_{n\rightarrow\infty } 
	\left(
	\text{Var}(\varepsilon_{n}^{\theta } ) 
	+
	\left(E(\varepsilon_{n}^{\theta } ) \right)^{2} 
	\right)
	=
	\frac{m^{2} |C_{\theta }(1) |^{2} }{2}
	+
	\frac{1}{4\pi} 
	\int_{-\pi}^{\pi} 
	|C_{\theta }(e^{i\omega } ) |^{2} \varphi(\omega ) d\omega.  
\end{align*}

For $\eta =[c_{1} \cdots c_{M} \; d_{1} \cdots d_{N} ]^{T} \in \mathbb{C}^{M+N}$, 
$z\in\mathbb{C}$, let 
\begin{align*}
	\hat{A}_{\eta }(z) 
	=
	1
	-
	\sum_{k=1}^{M} c_{k} z^{-k}, 
	\;\;\; 
	\hat{B}_{\eta }(z) 
	=
	1
	+
	\sum_{k=1}^{N} d_{k} z^{-k}
\end{align*}
and $\hat{C}_{\eta }(z) = \hat{A}_{\eta }(z)/\hat{B}_{\eta }(z)$, 
while 
\begin{align*}
	\hat{f}(\eta )
	=
	\frac{m^{2} |\hat{C}_{\eta }(1) |^{2} }{2}
	+
	\frac{1}{4\pi} 
	\int_{-\pi}^{\pi} 
	|\hat{C}_{\eta }(e^{i\omega } ) |^{2} \varphi(\omega ) d\omega.  
\end{align*}
Then, to prove the theorem's assertion, 
it is sufficient to show that $\hat{f}(\cdot )$
is analytic in an open vicinity of $\theta$. 

Obviously, $\hat{A}_{\eta}(z)$, $\hat{B}_{\eta}(z)$ 
are analytic in $(\eta,z)$ for all $\eta\in\mathbb{C}^{d_{\theta } }$, 
$z\in\mathbb{C}$, 
while $\hat{B}_{\theta }(z) = B_{\theta }(z) \neq 0$ 
for any $z\in\mathbb{C}$ satisfying $|z|=1$. 
Consequently, there exists a real number $\delta_{\theta }\in (0,1)$
such that 
$\hat{C}_{\eta }(e^{i\omega } )$ is analytic in $\eta$
and continuous in $(\eta, \omega )$
for all $\eta\in V_{\delta_{\theta } }(\theta )$, $\omega\in [-\pi,\pi]$
(here, $V_{\delta_{\theta } }(\theta )$ denotes 
$V_{\delta_{\theta } }(\theta ) 
= \{\eta\in\mathbb{C}^{d_{\theta } }: \|\eta-\theta \| \leq \delta_{\theta } \}$). 
Thus, there exists a real number $L_{1,\theta } \in [1,\infty )$
such that 
$|\hat{C}_{\eta }(e^{i\omega } ) | \leq L_{1,\theta }$
for all $\eta\in V_{\delta }(\theta )$, $\omega\in [-\pi,\pi]$. 
Then, Cauchy inequality for complex analytic functions 
(see e.g., \cite[Proposition 2.1.3]{taylor}) implies 
that there exists a real number $L_{2,\theta } \in [1,\infty )$
such that $\|\nabla_{\eta } \hat{C}_{\eta}(e^{i\omega } ) \| \leq L_{2,\theta }$ 
for each  
$\eta\in V_{\delta_{\theta } }(\theta )$, $\omega\in [-\pi,\pi]$. 
As a result of this and the dominated convergence theorem, 
$\int_{-\pi}^{\pi} |\hat{C}_{\eta }(e^{i\omega } ) |^{2} \varphi(\omega ) d\omega$
is differentiable in $\eta$ for any $\eta\in V_{\delta_{\theta } }(\theta )$. 
Hence, $\hat{f}(\cdot )$ is analytic on $V_{\delta_{\theta} }(\theta )$. 
}
\end{vproof}

\begin{vproof}{Theorem \ref{theorem6.2}}
{Owing to Assumptions \ref{a6.1} -- \ref{a6.3}, 
$\{V_{n}^{\theta } \}_{n\geq 0}$ and $\{Z_{n}^{\theta } \}_{n\geq 0}$ defined here
satisfy all conditions of Theorem \ref{theorema4.1} (Appendix \ref{appendix4}). 
Moreover, 
there exists a real number $K \in [1,\infty )$ 
such that (\ref{ta4.1.1*}) -- (\ref{ta4.1.5*}) are satisfied for $p=1$, $K_{2,Q}=K$ and all 
$\theta,\theta',\theta''\in \Theta$, 
$z,z',z'' \in \mathbb{R}^{d_{v} }\times {\cal W}$. 
Thus, all conclusions of Theorem \ref{theorema4.1} are true for 
$F(\theta,z)$, $\Pi_{\theta }(z,B)$ specified here. 
On the other side, (\ref{6.1*}) implies that in the case studied here, 
function $g(\theta )$ introduced in Theorem \ref{theorema4.1} is the gradient of $f(\theta )$. 
Consequently, Assumptions \ref{a2.2} and \ref{a2.3} hold. 
Then, the theorem's assertion directly follows from 
Theorem \ref{theorem2.1}. 
}
\end{vproof}

\section{Outline of the Proof of 
Theorems \ref{theorem9.1} and \ref{theorem9.2}}\label{section9*}

Theorem \ref{theorem9.1} is a direct consequence of Assumptions \ref{a9.2}, \ref{a9.5}. 
Owing to Assumption \ref{a9.2}, 
$\{X_{n}^{\theta } \}_{n\geq 0}$ has a unique invariant probability mass function 
$\pi_{\theta }(x)$ for any $\theta\in\mathbb{R}^{d_{\theta } }$. 
Consequently, 
$\pi_{\theta}(x)$ is a rational function of 
$\{p_{\theta }(x''|x') \}_{x',x'' \in{\cal X} }$. 
As $p_{\theta }(x'|x)$ is a polynomial function of 
$\{p(x''|x',y) \}_{x',x''\in{\cal X}, y\in{\cal Y}  }$ and 
$\{q_{\theta}(x'|y) \}_{x'\in{\cal X}, y\in{\cal Y} }$, 
Assumption \ref{a9.5} implies that for any $x\in{\cal X}$, 
$\pi_{\theta}(x)$ is analytic in $\theta$ on entire $\mathbb{R}^{d_{\theta } }$. 
Since 
\begin{align*}
	f(\theta )
	=
	\sum_{x\in{\cal X}, y\in{\cal Y} } c(x,y) q_{\theta }(y|x) \pi_{\theta }(x)
\end{align*}
for any $\theta\in\mathbb{R}^{d_{\theta } }$, 
$f(\cdot )$ is analytic on entire $\mathbb{R}^{d_{\theta } }$. 

To explain how Theorem \ref{theorem9.2} is proved, 
we use the following notation. 
$d_{\eta } = d_{\theta } + 1$ and 
$d_{\vartheta } = d_{\theta } + d_{\eta }$, 
while 
${\cal Z} = 
{\cal X}\times{\cal Y}\times{\cal X}\times{\cal Y}\times\mathbb{R}^{d_{\theta } }$. 
Stochastic processes  
$\{\eta_{n} \}_{n\geq 0}$, 
$\{\vartheta_{n} \}_{n\geq 0}$, $\{Z_{n} \}_{n\geq 0}$ are defined as 
$\eta_{n} = [\eta_{1,n}^{T} \: \eta_{2,n} ]^{T}$,
$\vartheta_{n} = [\theta_{n}^{T} \: \eta_{n}^{T} ]^{T}$,  
$Z_{n+1} = (X_{n}, Y_{n}, X_{n+1}, Y_{n+1}, W_{n+1} )$ for $n\geq 0$. 
$A_{1,\theta }(z)$, $c_{1,\theta }(z)$, $c_{2,\theta }(z)$ are the functions defined as 
\begin{align*}
	A_{1,\theta }(z) 
	=
	s_{\theta }(x',y') s_{\theta }^{T}(x',y'), 
	\;\;\;\;\; 
	c_{1,\theta }(z) 
	=
	w c(x,y), 
	\;\;\;\;\; 
	c_{2,\theta }(z) 
	=
	c(x',y')
\end{align*}
for $\theta,w\in\mathbb{R}^{d_{\theta } }$, 
$x,x'\in{\cal X}$, $y,y'\in{\cal Y}$, 
$z= (x,y,x',y',w)$, 
while  functions $B_{1,\theta }(z)$, $B_{2,\theta }(z)$
are defined by  
\begin{align*}
	B_{1,\theta }(z)
	=
	w (s_{\theta }(x,y) - s_{\theta }(x',y') )^{T}, 
	\;\;\;\;\; 
	B_{2,\theta }(z)
	=
	w
\end{align*}
for the same $\theta,w,x,x',y,y',z$. 
$A_{\theta }(z)$, $B_{\theta }(z)$, $c_{\theta }(z)$ 
are the functions defined as 
\begin{align*}
	A_{\theta }(z) 
	=
	-
	\left[
	\begin{array}{cc}
	A_{1,\theta }(z) & {\boldsymbol 0} 
	\end{array}
	\right], 
	\;\;\; 
	B_{\theta }(z)
	=
	-
	\left[
	\begin{array}{cc}
	B_{1,\theta }(z) & B_{2,\theta }(z) \\
	{\boldsymbol 0}^{T} & 1
	\end{array}
	\right], 
	\;\;\; 
	c_{\theta }(z)
	=
	\left[
	\begin{array}{c}
	c_{1,\theta }(z) \\
	c_{2,\theta }(z)
	\end{array}
	\right]
\end{align*}
for $\theta\in\mathbb{R}^{d_{\theta } }$, $z\in{\cal Z}$, 
where ${\boldsymbol 0}$ denotes $d_{\theta }$-dimensional zero (column) vector
(notice that $A_{\theta }(z) \in \mathbb{R}^{d_{\theta }\times d_{\eta } }$, 
$B_{\theta }(z) \in \mathbb{R}^{d_{\eta }\times d_{\eta } }$). 
$\Pi_{\theta }(z,B)$ is a transition kernel defined by 
\begin{align*}
	\Pi_{\theta }(z,B)
	=
	\sum_{x''\in{\cal X}, y''\in{\cal Y} }
	I_{B}\left(x',y',x'',y'',w I_{\{x''\neq x_{*} \} } + s_{\theta }(x'',y'') \right)
	q_{\theta }(y''|x'') p(x''|x',y') 
\end{align*}
for a measurable set $B\subseteq{\cal Z}$ and 
$\theta,w\in\mathbb{R}^{d_{\theta } }$, 
$x,x'\in{\cal X}$, $y,y'\in{\cal Y}$, 
$z= (x,y,x',y',w)$. 
For $\theta\in\mathbb{R}^{d_{\theta } }$, 
$\{Z_{n}^{\theta } \}_{n\geq 0}$ is a ${\cal Z}$-valued Markov chain 
whose transition kernel is $\Pi_{\theta }(\cdot,\cdot )$. 
$\bar{A}(\theta )$, $\bar{B}(\theta )$, $\bar{c}(\theta )$ are the functions defined as 
\begin{align*}
	\bar{A}(\theta )
	=
	\lim_{n\rightarrow\infty } E(A_{\theta }(Z_{n}^{\theta } ) ), 
	\;\;\; 
	\bar{B}(\theta )
	=
	\lim_{n\rightarrow\infty } E(B_{\theta }(Z_{n}^{\theta } ) ), 
	\;\;\; 
	\bar{c}(\theta )
	=
	\lim_{n\rightarrow\infty } E(c_{\theta }(Z_{n}^{\theta } ) ) 
\end{align*}
for $\theta\in\mathbb{R}^{d_{\theta } }$, 
while functions 
$r(\theta )$, $S(\theta )$ are defined by 
\begin{align*}
	r(\theta ) 
	=
	\lim_{n\rightarrow\infty } E(B_{2,\theta }(Z_{n}^{\theta } ) ), 
	\;\;\; 
	S(\theta ) 
	=
	\lim_{n\rightarrow\infty } E(A_{1,\theta }(Z_{n}^{\theta } ) )
\end{align*}
for the same $\theta$. 
Under the introduced notation, algorithm (\ref{9.1}) -- (\ref{9.7}) can be re-written as 
\begin{align}
	\label{9.1*}
	\theta_{n+1}
	=&
	\theta_{n}
	+
	\alpha_{n} A_{\theta_{n} }(Z_{n+1} ) \eta_{n}, 
	\\
	\label{9.3*} 
	\eta_{n+1} 
	=&
	\eta_{n} 
	+
	\beta_{n} (B_{\theta_{n} }(Z_{n+1} ) \eta_{n} + c_{\theta_{n} }(Z_{n+1} ) ), 
	\;\;\; n\geq 0. 
\end{align}
It can also be shown that $\{\theta_{n} \}_{n\geq 0}$, $\{Z_{n} \}_{n\geq 0}$,
$\Pi_{\theta }(z,B)$ defined here satisfy (\ref{2.3}). 
Hence, recursion (\ref{9.1*}), (\ref{9.3*}) fits into the framework studied in 
\cite{benveniste}, \cite{konda&tsitsiklis}. 
Then, using the results of \cite[Section 5.1]{konda&tsitsiklis}, 
we conclude that 
$\bar{A}(\theta )$, $\bar{B}(\theta )$, $\bar{c}(\theta )$, 
$r(\theta )$, $S(\theta )$  
are well-defined 
and satisfy 
\begin{align}\label{9.21*}
	\bar{A}(\theta )
	=
	-
	\left[
	\begin{array}{cc}
	S(\theta ) & {\boldsymbol 0} 
	\end{array}
	\right], 
	\;\;
	\bar{B}(\theta ) 
	=
	-
	\left[
	\begin{array}{cc}
	S(\theta ) & r(\theta ) \\ 
	{\boldsymbol 0}^{T} & 1
	\end{array}
	\right], 
	\;\;
	\bar{c}(\theta )
	=
	\left[
	\begin{array}{c}
	\nabla f(\theta ) + r(\theta ) f(\theta ) \\ 
	f(\theta )
	\end{array} 
	\right]
\end{align}
for each $\theta\in\mathbb{R}^{d_{\theta } }$. 
Combining the same results with the arguments 
behind Theorem \ref{theorem2.1}, we deduce  
\begin{align}
	&\label{9.23*}
	\lim_{n\rightarrow\infty } 
	\max_{n\leq k < a(n,1) }
	\left\|
	\sum_{i=n}^{k} 
	\alpha_{i} \gamma_{i}^{r} 
	\left(A_{\theta_{i} }(Z_{i+1} ) - \bar{A}(\theta_{i} ) \right)
	\right\|
	=0, 
	\\
	&\label{9.25*}
	\lim_{n\rightarrow\infty } 
	\max_{n\leq k < a(n,1) }
	\left\|
	\sum_{i=n}^{k} 
	\beta_{i} \gamma_{i}^{r} 
	\left(B_{\theta_{i} }(Z_{i+1} ) - \bar{B}(\theta_{i} ) \right)
	\right\|
	=0, 
	\\
	&\label{9.27*}
	\lim_{n\rightarrow\infty } 
	\max_{n\leq k < a(n,1) }
	\left\|
	\sum_{i=n}^{k} 
	\beta_{i} \gamma_{i}^{r} 
	\left(c_{\theta_{i} }(Z_{i+1} ) - \bar{c}(\theta_{i} ) \right)
	\right\|
	=0 
\end{align}
w.p.1 on $\{\sup_{n\geq 0} \|\theta_{n} \| < \infty \}$. 

Recursion (\ref{9.3*}) can be viewed as linear stochastic approximation in 
$\{\eta_{n} \}_{n\geq 0}$. 
As 
$\sup_{\theta\in Q} \lambda_{\text{max} }(\bar{B}(\theta ) ) < 0$ for any compact set 
$Q\subset\mathbb{R}^{d_{\theta } }$
(due to Assumption \ref{a9.6}; $\lambda_{\text{max} }(\bar{B}(\theta ) )$
stands for the maximal eigenvalue of $\bar{B}(\theta )$), 
standard asymptotic results for linear stochastic approximation 
(see \cite{tadic1} or \cite[Appendix A]{konda&tsitsiklis}) 
imply that $\{\eta_{n} \}_{n\geq 0}$ is bounded whenever 
$\{\theta_{n} \}_{n\geq 0}$ is bounded. 
More specifically, for any compact set $Q\subset\mathbb{R}^{d_{\theta } }$, 
there exists a real number $\rho_{Q}\in[1,\infty )$
such that w.p.1 on event $\Lambda_{Q} = \bigcap_{n=0}^{\infty } \{\theta_{n}\in Q\}$, 
$\|\eta_{n} \|\leq\rho_{Q}$ for all, but finitely many $n$. 
Consequently, (\ref{9.21*}) -- (\ref{9.25*}) imply that 
recursion (\ref{9.1*}), (\ref{9.3*}) (i.e., algorithm (\ref{9.1}) -- (\ref{9.7})) 
admits representation 
\begin{align}\label{9.5*}
	\vartheta_{n+1} 
	=
	\vartheta_{n} 
	+
	\alpha_{n} D_{n} (h(\vartheta_{n} ) + \xi_{n} ), 
	\;\;\; n\geq 0.  
\end{align}
Here, $h(\vartheta )$ is the function defined by 
\begin{align*}
	h(\vartheta ) = 
	\left[
	\begin{array}{c}
	\bar{A}(\theta ) \eta \\ \bar{B}(\theta ) \eta + \bar{c}(\theta ) 
	\end{array}
	\right] 
\end{align*}
for $\theta\in\mathbb{R}^{d_{\theta } }$, $\eta\in\mathbb{R}^{d_{\eta } }$, 
$\vartheta = [\theta^{T} \; \eta^{T} ]^{T}$. 
$\{D_{n} \}_{n\geq 0}$ are diagonal matrices defined as 
$D_{n} = \text{diag}\{{\boldsymbol I}', \alpha_{n}^{-1} \beta_{n} {\boldsymbol I}'' \}$ 
for $n\geq 0$, 
where ${\boldsymbol I}'$ and ${\boldsymbol I}''$
denote $d_{\theta }\times d_{\theta }$ and $d_{\eta }\times d_{\eta }$ unit matrices (respectively). 
$\{\xi_{n} \}_{n\geq 0}$ is an $\mathbb{R}^{d_{\vartheta } }$-valued stochastic process satisfying  
\begin{align*}
	\lim_{n\rightarrow\infty } 
	\max_{n\leq k < a(n,1) } 
	\left\|
	\sum_{i=n}^{k} \alpha_{i}\gamma_{i}^{r}D_{i}\xi_{i} 
	\right\|
	=
	0
\end{align*}
w.p.1 on $\{\sup_{n\geq 0} \|\theta_{n} \|<\infty \}$
($a(n,1)$ is defined in Section \ref{section1}). 

To explain how the asymptotic behavior of (\ref{9.5*}) can be analyzed, we rely on the following notation. 
$K_{1}(\theta,\eta)$, $K_{2}(\theta,\eta)$, $K_{3}(\theta,\eta)$, $K_{4}(\theta,\eta)$ 
are the functions defined by   
\begin{align*}
	K_{1}(\theta,\eta)
	=
	\nabla_{\theta }\left(S(\theta )\eta \right) - \nabla^{2} f(\theta ), 
	\;\;\;\;\;
	K_{2}(\theta,\eta)
	=
	{\boldsymbol I} 
	+
	\nabla_{\theta }\left(S(\theta )\eta \right) - \nabla^{2} f(\theta )
\end{align*}
and 
$K_{3}(\theta,\eta) = -S(\theta )\eta$, 
$K_{4}(\theta,\eta) = S(\theta )r(\theta )$
for $\theta,\eta\in\mathbb{R}^{d_{\theta } }$, 
where ${\boldsymbol I}$ denotes $d_{\theta }\times d_{\theta }$ unit matrix. 
For a compact set $Q\subset\mathbb{R}^{d_{\theta } }$, 
$L_{Q}\in[1,\infty )$ stands for a real number satisfying 
\begin{align}\label{9.31*}
	\lambda_{min}(S(\theta ) )\geq L_{Q}^{-1/2}, 
	\;\;\;\;\; 
	\max_{1\leq i \leq 4 } \|K_{i}(\theta,\eta ) \|
	\leq 
	2^{-1} L_{Q}^{1/4}
\end{align}
for all $\theta\in Q$, $\eta\in\mathbb{R}^{d_{\theta } }$
satisfying $\|\eta \|\leq \rho_{Q}$
($\lambda_{\text{min} } (S(\theta ) )$
denotes the smallest eigenvalue of $S(\theta )$; 
notice that $S(\theta )$ is positive definite and continuous for
each $\theta\in\mathbb{R}^{d_{\theta}}$).  

To study the asymptotic behavior of (\ref{9.5*}), 
for each compact set $Q\subset\mathbb{R}^{d_{\theta } }$, 
we construct the following Lyapunov function: 
\begin{align*}
	v_{Q}(\vartheta )
	=
	f(\theta ) 
	+
	\frac{1}{2} \|S(\theta ) \eta_{1} - \nabla f(\theta ) \|^{2}  
	+
	\frac{L_{Q} }{2} 
	(\eta_{2} - f(\theta ) )^{2},  
\end{align*}
where 
$\theta,\eta_{1}\in\mathbb{R}^{d_{\theta } }$, 
$\eta_{2}\in\mathbb{R}$ and 
$\vartheta = [\theta^{T}\; \eta_{1}^{T}\; \eta_{2} ]^{T}$. 
Then, it is straightforward to verify 
\begin{align*}
	(\nabla v_{Q}(\vartheta ) )^{T} D_{n} h(\vartheta )
	=&
	-
	\|\nabla f(\theta ) \|^{2} 
	-
	L_{Q} \alpha_{n}^{-1}\beta_{n} 
	(\eta_{2} - f(\theta ) )^{2} 
	\\
	&
	-
	(S(\theta ) \eta_{1} \!-\! \nabla f(\theta ) )^{T}
	\left(\alpha_{n}^{-1}\beta_{n} S(\theta ) \!+\! K_{1}(\theta,\eta_{1} ) \right) 
	(S(\theta ) \eta_{1} \!-\! \nabla f(\theta ) ) 
	\\
	&
	-
	(\nabla f(\theta ) )^{T} 
	K_{2}(\theta, \eta_{1} ) 
	(S(\theta ) \eta_{1} - \nabla f(\theta ) )
	\\
	&
	-
	L_{Q} 
	(\nabla f(\theta ) )^{T} 
	K_{3}(\theta, \eta_{1} ) 
	(\eta_{2} - f(\theta ) )
	\\
	&
	-
	\alpha_{n}^{-1}\beta_{n} 
	(S(\theta ) \eta_{1} - \nabla f(\theta ) )^{T} 
	K_{4}(\theta, \eta_{1} ) 
	(\eta_{2} - f(\theta ) )
\end{align*}
for all $\theta,\eta_{1}\in\mathbb{R}^{d_{\theta } }$, 
$\eta_{2}\in\mathbb{R}$ and 
$\vartheta = [\theta^{T}\: \eta_{1}^{T}\: \eta_{2} ]^{T}$.  
Owing to (\ref{9.31*}), we have  
\begin{align*}
	&
	\left|
	(S(\theta ) \eta_{1} - \nabla f(\theta ) )^{T} 
	K_{4}(\theta, \eta_{1} )
	(\eta_{2} - f(\theta ) )
	\right|
	\\
	&
	\leq 
	2^{-1} L_{Q}^{1/4} 
	\|S(\theta ) \eta_{1} - \nabla f(\theta ) \| \:
	|\eta_{2} - f(\theta ) |
	\\
	&
	\leq 
	4^{-1} L_{Q}^{-1/2} 
	\|S(\theta ) \eta_{1} - \nabla f(\theta ) \|^{2} 
	+
	4^{-1} L_{Q} 
	(\eta_{2} - f(\theta ) )^{2} 
\end{align*}
for all $\theta\in Q$, 
$\eta_{1} \in \mathbb{R}^{d_{\theta } }$, 
$\eta_{2}\in\mathbb{R}$
satisfying $\|\eta_{1} \|\leq \rho_{Q}$. 
Similarly, we get  
\begin{align*}
	&
	\begin{aligned}[t] 
	\left|
	(\nabla f(\theta ) )^{T} 
	K_{2}(\theta, \eta_{1} )
	(S(\theta ) \eta_{1} - \nabla f(\theta ) )
	\right| 
	\leq &
	2^{-1} L_{Q}^{1/4}  
	\|\nabla f(\theta ) \| \|S(\theta ) \eta_{1} - \nabla f(\theta ) \| 
	\\
	\leq &
	4^{-1} \|\nabla f(\theta ) \|^{2} 
	+
	4^{-1} L_{Q}^{1/2} \|S(\theta ) \eta_{1} - \nabla f(\theta ) \|^{2},  
	\end{aligned} 
	\\
	&
	\begin{aligned}[t] 
	\left|
	(\nabla f(\theta ) )^{T} 
	K_{3}(\theta, \eta_{1} )
	(\eta_{2} - f(\theta ) )
	\right| 
	\leq &
	2^{-1} L_{Q}^{1/4}  
	\|\nabla f(\theta ) \|\: |\eta_{2} - f(\theta ) | 
	\\
	\leq &
	4^{-1} \|\nabla f(\theta ) \|^{2} 
	+
	4^{-1} L_{Q}^{1/2} |\eta_{2} - f(\theta ) |^{2}  
	\end{aligned} 
\end{align*}
for the same $\theta, \eta_{1}, \eta_{2}$. 
We also have 
\begin{align*}
	&
	(S(\theta ) \eta_{1} - \nabla f(\theta ) )^{T}
	\left(\alpha_{n}^{-1}\beta_{n} S(\theta ) + K_{1}(\theta,\eta_{1} ) \right) 
	(S(\theta ) \eta_{1} - \nabla f(\theta ) ) 
	\\
	&\geq 
	\left(
	\alpha_{n}^{-1}\beta_{n} \lambda_{min}(S(\theta ) ) 
	-
	\|K_{1}(\theta,\eta_{1} ) \|
	\right)
	\|S(\theta ) \eta_{1} - \nabla f(\theta ) \|^{2} 
	\\
	&
	\geq 
	\left(
	L_{Q}^{-1/2} \alpha_{n}^{-1}\beta_{n} 
	-
	2^{-1} L_{Q}^{1/4} 
	\right)
	\|S(\theta ) \eta_{1} - \nabla f(\theta ) \|^{2} 
\end{align*}
for all $\theta\in Q$, 
$\eta_{1} \in \mathbb{R}^{d_{\theta } }$ 
satisfying $\|\eta_{1} \|\leq \rho_{Q}$.
Hence, 
\begin{align*}
	(\nabla v_{Q}(\vartheta ) )^{T} D_{n} h(\vartheta )
	\leq &
	-
	2^{-1} \|\nabla f(\theta ) \|^{2} 
	-
	\left(
	2^{-1} L_{Q}^{-1} \alpha_{n}^{-1}\beta_{n} 
	-
	L_{Q}^{1/2} 
	\right)
	\|S(\theta ) \eta_{1} - \nabla f(\theta ) \|^{2} 
	\\
	&
	-
	\left(
	2^{-1} L_{Q}\alpha_{n}^{-1}\beta_{n}  
	-
	L_{Q}^{1/2} 
	\right)
	(\eta_{2} - f(\theta ) )^{2} 
\end{align*} 
for each $\theta\in Q$, $[\eta_{1}^{T} \; \eta_{2}]^{T} \in V_{\rho_{Q} }$, 
$\vartheta = [\theta^{T} \; \eta_{1}^{T} \; \eta_{2} ]^{T}$, $n\geq 0$, 
where $V_{\rho_{Q} } = \{\eta\in\mathbb{R}^{d_{\eta } }: \|\eta\|\leq \rho_{Q} \}$. 
As $\lim_{n\rightarrow\infty } \alpha_{n}^{-1}\beta_{n} = \infty$, 
we deduce that there exists an integer $m_{Q}\geq 1$ such that 
\begin{align}\label{9.33*}
	(\nabla v_{Q}(\vartheta ) )^{T} D_{n} h(\vartheta ) 
	\leq &
	-2^{-1} 
	\left(
	\|\nabla f(\theta ) \|^{2} 
	+
	\|S(\theta ) \eta_{1} - \nabla f(\theta ) \|^{2}	
	+
	L_{Q}(\eta_{2} - f(\theta ) )^{2}
	\right)
	\leq 0
\end{align}
for all 
$\theta\in Q$, 
$[\eta_{1}^{T} \; \eta_{2}]^{T} \in V_{\rho_{Q} }$, 
$\vartheta = [\theta^{T} \: \eta_{1}^{T} \; \eta_{2} ]^{T}$, 
$n\geq m_{Q}$. 
Combining this with standard stochastic approximation arguments, we conclude 
that $\{v_{Q}(\vartheta_{n} ) \}_{n\geq 0}$ converges w.p.1 on $\Lambda_{Q}$
and that 
\begin{align}\label{9.35*}
	\lim_{n\rightarrow\infty } \|\nabla f(\theta_{n} ) \|
	=
	\lim_{n\rightarrow\infty } 
	\|S(\theta_{n} ) \eta_{1,n} - \nabla f(\theta_{n} ) \| 
	=
	\lim_{n\rightarrow\infty } 
	|\eta_{2,n} - f(\theta_{n} ) | 
	=0
\end{align}
w.p.1 on the same event. 
Hence, w.p.1 on $\Lambda_{Q}$, (\ref{9.1*}) asymptotically behaves as a gradient search 
minimizing $v_{Q}(\cdot )$. 
On the other side, Lojasiewicz inequality (\ref{a1.3.1}) and (\ref{9.33*}) yield 
\begin{align}\label{9.37*}
	(\nabla v_{Q}(\vartheta ) )^{T} D_{n} h(\vartheta ) 
	\leq & 
	-
	\left(
	v_{Q}(\vartheta ) - f(\theta ) + 2^{-1} \|\nabla f(\theta ) \|^{2}
	\right)
	\nonumber\\
	\leq &
	-
	\left(
	v_{Q}(\vartheta ) - f(\theta ) 
	\right)
	-
	2^{-1} M_{Q,a}^{-2} |f(\theta ) - a |^{2/\mu_{Q,a} }
	\nonumber\\
	\leq & 
	-
	L_{Q,a}^{-1} 
	\left( 
	v_{Q}(\vartheta ) - f(\theta ) 
	+
	|f(\theta ) - a |
	\right)^{2/\mu_{Q,a} }
	\nonumber\\
	\leq & 
	-
	L_{Q,a}^{-1} 
	|v_{Q}(\vartheta ) - a |^{2/\mu_{Q,a} }
\end{align}
for all $a\in f(Q)$, $\theta\in Q$, $\eta\in V_{\rho_{Q} }$, 
$\vartheta=[\theta^{T} \; \eta^{T} ]^{T}$
satisfying $|f(\theta ) - a|\leq \delta_{Q,a}$
($\delta_{Q,a}$ is specified in Assumption \ref{a1.3}), 
where $L_{Q,a}\in [1,\infty )$ is a suitably chosen real number.\footnote
{$L_{Q,a}$ can be selected as 
$L_{Q,a} \!=\! \max\{2M_{Q,a}^{2}, \!K_{Q,a} \}$, where 
$
	K_{Q,a} 
	\!=\! 
	\sup\{(v_{Q}(\vartheta ) \!-\! f(\theta ) )^{2/\mu_{Q,a} -\! 1 }\!:\! 
	\vartheta \!=\! [\theta^{T} \; \eta^{T} ]^{T}, \theta\!\in\! Q,\eta\!\in\! V_{\rho_{Q} } \} 
$. 
}
Thus, w.p.1 on $\Lambda_{Q}$, 
\begin{align*}
	(\nabla v_{Q}(\vartheta_{n} ) )^{T} D_{n} h(\vartheta_{n} ) 
	\leq 
	\hat{L}^{-1} |v_{Q}(\vartheta_{n} ) - \hat{v}_{Q} |^{2/\hat{\mu} }
\end{align*}
for all, but finitely many $n$, 
where $\hat{v}_{Q} = \lim_{n\rightarrow\infty } v(\vartheta_{n} )$, 
$\hat{\mu} = \mu_{Q,\hat{v}_{Q} }$, $\hat{L} = L_{Q,\hat{v}_{Q} }$. 
As (\ref{9.37*}) and Lojasiewicz inequality (\ref{a1.3.1})
have very similar forms, 
(\ref{9.37*}) can be considered as a Lojasiewicz-type inequality for 
$v_{Q}(\cdot )$. 

The conclusions drawn about recursion (\ref{9.5*}) 
(asymptotic equivalence with a gradient search minimizing $v_{Q}(\cdot )$)
and Lyapunov function $v_{Q}(\cdot )$
(Lojasiewicz-type inequality (\ref{9.37*})) 
strongly suggest that Theorems \ref{theorem1.1}, \ref{theorem1.2}, \ref{theorem2.1}
can be extended to algorithm (\ref{9.1}) -- (\ref{9.7}) 
and that Theorem \ref{theorem9.2} is true. 
A detailed proof of this assertion is provided in \cite{tadic10}. 

\refstepcounter{appendixcounter}\label{appendix1}
\section*{Appendix \arabic{appendixcounter} }

In this section, we prove the claim stated in Remark \ref{remark1.2}. 
If open set $V$ specified in Remark \ref{remark1.2} exists,
we can define the following quantities for 
any compact set $Q\subset \mathbb{R}^{d_{\theta } }$ 
and any $a\in f(Q)$: 
\begin{align*}
	&
	\tilde{\delta}_{Q,a} 
	=
	\begin{cases}
	\delta_{\tilde{Q},a}, 
	&\text{ if } Q\cap S \neq \emptyset, \; a \in f(S)
	\\
	1, 
	&\text{ if } Q \cap S = \emptyset
	\\
	\min\{1, d(a, f(S) )/2 \}, 
	&\text{ if } a \not\in f(S)
	\end{cases}
	\\
	&
	\tilde{\mu}_{Q,a} 
	=
	\begin{cases}
	\mu_{\tilde{Q},a}, 
	&\text{ if } Q \cap S \neq \emptyset, \; a \in f(S)
	\\
	2, &\text{ otherwise }
	\end{cases}
	\\
	&
	\tilde{M}_{Q,a} 
	=
	1
	+
	\sup\left\{ 
	\frac{|f(\theta ) - a |}{\|\nabla f(\theta ) \|^{\tilde{\mu}_{Q,a} } }: 
	\theta \in Q\setminus S, 
	|f(\theta ) - a | \leq \tilde{\delta}_{Q,a}
	\right\}
\end{align*}
where 
$\tilde{Q} = Q$ if $Q\subset V$
and 
$\tilde{Q} = \{\theta\in Q: 
d(\theta, S ) \leq d(Q\setminus V, S )/2 \}$ otherwise. 
Then, it is straightforward to show
\begin{align*}
	a\not\in f(S) 
	\Longrightarrow &
	\inf\{\|\nabla f(\theta ) \|: 
	\theta \in Q, |f(\theta ) - a | \leq \tilde{\delta}_{Q,a} \}
	> 0, 
	\\
	Q\setminus V \neq \emptyset
	\Longrightarrow &
	\inf\{\|\nabla f(\theta ) \|: 
	\theta \in Q\setminus \tilde{Q} \} 
	>0, 
	\\
	Q\cap S \neq \emptyset
	\Longrightarrow & 
	\sup\left\{
	\frac{|f(\theta ) - a |}{\|\nabla f(\theta ) \|^{\tilde{\mu}_{Q,a} } }: 
	\theta \in \tilde{Q}, 
	|f(\theta ) - a | \leq \tilde{\delta}_{Q,a} 
	\right\}
	\leq 
	M_{\tilde{Q},a}
	< 
	\infty. 
\end{align*}
Consequently, 
$\tilde{\delta}_{Q,a}$, $\tilde{\mu}_{Q,a}$, $\tilde{M}_{Q,a}$
are well-defined
and enjoy the following properties: 
$0 < \tilde{\delta}_{Q,a} \leq 1$, 
$1 < \tilde{\mu}_{Q,a} \leq 2$, 
$1\leq \tilde{M}_{Q,a} < \infty$ 
and 
\begin{align*}
	|f(\theta ) - a |
	\leq 
	\tilde{M}_{Q,a} \|\nabla f(\theta ) \|^{\tilde{\mu}_{Q,a} }
\end{align*}
for all $\theta \in Q$ satisfying 
$|f(\theta ) - a | \leq \tilde{\delta}_{Q,a}$. 
Hence, the claim holds.

\refstepcounter{appendixcounter}\label{appendix2}
\section*{Appendix \arabic{appendixcounter} }

In this section, a global version of the results of Section \ref{section1} is presented. 
It is also demonstrated how the results of Section \ref{section1} can be extended to
the randomly projected stochastic gradient search. 

First, the stability and the global asymptotic behavior of algorithm (\ref{1.1}) 
are considered. 
To analyze these properties, we introduce the following two assumptions. 

\begin{assumptionappendix}\label{aa2.1} 
$\nabla f(\cdot )$ is (globally) Lipschitz continuous. 
Moreover, $\liminf_{\|\theta\|\rightarrow\infty } \|\nabla f(\theta ) \| > 0$
and 
$\inf_{\theta\in\mathbb{R}^{d_{\theta } } } f(\theta )>-\infty $. 
\end{assumptionappendix}

\begin{assumptionappendix}\label{aa2.2} 
There exists a real number $r\in(1,\infty )$ sucht that 
\begin{align}\label{aa2.2.1}
	\limsup_{n\rightarrow\infty } 
	g(\theta_{n} ) 
	\max_{n\leq j\leq a(n,t) } 
	\left\| 
	\sum_{i=n}^{j} \alpha_{i} \gamma_{i}^{r}\xi_{i} 
	\right\|
	<
	\infty
\end{align}
almost surely for any $t\in(0,\infty )$. 
Here, $g: \mathbb{R}^{d_{\theta } }\rightarrow(0,\infty )$ is the (scaling) function defined by 
\begin{align*}
	g(\theta ) 
	=
	(\|\nabla f(\theta ) \| + 1 )^{-1}
\end{align*}
for $\theta\in\mathbb{R}^{d_{\theta } }$. 
\end{assumptionappendix}

Assumption \ref{aa2.1} is a stability condition. 
In this or a similar form, it is an inevitable ingredient of the stability analysis of 
stochastic gradient search and stochastic approximation 
(see e.g., \cite{benveniste}, \cite{borkar}, \cite{chen} and references cited therein). 
Assumption \ref{aa2.1} is rather restrictive, since it requires $\nabla^{2}f(\cdot )$ to be 
uniformly bounded. 
It also requires $\nabla f(\cdot )$ to grow at most linearly as $\theta\rightarrow\infty$. 
Using the random projections (see (\ref{appendix2.1})), 
these restrictive requirements can considerably be relaxed. 

Assumption \ref{aa2.2} is a noise condition. 
Basically, Assumption \ref{aa2.2} requires 
the effect of the gradient estimator error $\{\xi_{n} \}_{n\geq 0}$
to be compensated by the gradient of the objective function $f(\cdot )$
(i.e., by the stability of the ODE $d\theta/dt=-\nabla f(\cdot )$). 
Assumption \ref{aa2.2} holds whenever (\ref{a1.2.1}) is satisfied w.p.1. 
It also holds for stochastic gradient search with Markovian dynamics
(see Theorem \ref{theorema3.1}, Appendix \ref{appendix3}). 
Assumption \ref{aa2.2} and the results based on it (Theorem \ref{theorema2.1}, below) 
are inspired by the scaled ODE approach 
to the stability analysis of stochastic approximation \cite{borkar&meyn}.\footnote
{The main difference between \cite{borkar&meyn} and the results presented here 
is the choice of the scaling functions. 
The scaling adopted in \cite{borkar&meyn} is (asymptotically) proportional to 
$\|\theta\|$. 
In this paper, the scaling is (asymptotically) proportional to $\|\nabla f(\theta )\|$. } 

Our results on the stability of algorithm (\ref{1.1}) are provided in the next theorem. 

\begin{theoremappendix}\label{theorema2.1}
Let Assumptions \ref{a1.1}, \ref{aa2.1} and \ref{aa2.2} hold. 
Then, there exists a real number $\rho\in[1,\infty )$
such that $\limsup_{n\rightarrow\infty } \|\theta_{n} \| < \rho$ almost surely. 
If in addition, Assumption \ref{a1.3} holds, 
then the following is also true: 
\begin{compactenum}[(i)]
\item
$\hat{\theta}=\lim_{n\rightarrow\infty } \theta_{n}$
exists and satisfies $\nabla f(\hat{\theta} ) = 0$ almost surely. 
\item
There exists a random variable $\hat{K}$ 
such that $0<\hat{K}<\infty$ everywhere and such that (\ref{t1.1.1*}) -- (\ref{t1.1.5*}) 
holds almost surely. 
\end{compactenum}
\end{theoremappendix}

\begin{sproof}
Owing to Assumption \ref{aa2.1}, there exist real numbers
$c\in(0,1)$,  
$\rho_{1}, \tilde{C}_{1} \in [1,\infty )$ 
such that the following is true: 
(i) $\|\nabla f(\theta ) \|\geq c$ for each $\theta\in\mathbb{R}^{d_{\theta } }$
satisfying $\|\theta \|\geq\rho_{1}$,  
(ii) $f(\theta )\leq \tilde{C}_{1}$ 
for any $\theta\in\mathbb{R}^{d_{\theta } }$ satisfying $\|\theta \|\leq \rho_{1}+1$, 
and 
(iii) $f(\theta )> -\tilde{C}_{1}$ for all $\theta\in\mathbb{R}^{d_{\theta } }$. 
On the other side, due to Assumption \ref{aa2.2}, 
there also exists an event $N_{0}\in {\cal F}$ with the following properties: 
(i) $P(N_{0} ) = 0$,
and 
(ii) (\ref{aa2.2.1}) holds on $N_{0}^{c}$ for all $t\in(0,\infty )$. 
Then, relying on the same arguments as in the proof of Lemma \ref{lemma1.1}, 
we conclude
\begin{align}\label{ta2.1.1701}
	\lim_{n\rightarrow\infty } 
	g(\theta_{n} ) 
	\max_{n\leq j < a(n,t) } 
	\left\|\sum_{i=n}^{n} \alpha_{i}\xi_{i} \right\|
	=
	0
\end{align}
on $N_{0}^{c}$ for any $t\in(0,\infty )$. 

Let 
$T=16\tilde{C}_{1}c^{-2}$. 
Moreover, let $\phi:[0,\infty )\rightarrow[0,\infty )$ be the function defined by 
\begin{align*}
	\phi(z)
	=
	\sup\{\|\nabla f(\theta ) \|: \theta\in\mathbb{R}^{d_{\theta } }, \|\theta\|\leq z \}
\end{align*}
for $z\in [0,\infty )$. 
As $\nabla f(\cdot )$ is locally Lipschitz continuous, 
$\phi(\cdot )$ is locally Lipschitz continuous, too. 
$\phi(\cdot )$ is also non-negative and satisfies $\|\nabla f(\theta ) \|\leq \phi(\|\theta \| )$
for all $\theta\in\mathbb{R}^{d_{\theta } }$. 

For $z\in [0,\infty )$, let $\lambda(\cdot\:; z )$ be the solution to 
the ODE $dz/dt=\phi(z)$ satisfying $\lambda(0;z)=z$. 
As $\phi(\cdot )$ is non-negative and locally Lipschitz continuous, 
$\lambda(\cdot\; ;\cdot )$ is well-defined and locally Lipschitz continuous 
(in both arguments) on $[0,\infty )\times [0,\infty )$. 
We also have 
\begin{align}\label{ta2.1.701}
	\lambda(t;z)
	=
	z
	+
	\int_{0}^{t} \phi(\lambda(s;z) ) ds
\end{align}
for all $t,z\in[0,\infty )$. 
Then, there exists $\rho_{2}\in [1,\infty )$ such that 
$\rho_{2}\geq\rho_{1}+1$ and such that 
$|\lambda(t;z) |\leq\rho_{2}$ for all $t\in [0,T]$, $z\in [0,\rho_{1}+1]$. 

Let $\rho=\rho_{2}+1$, $Q=\{\theta\in\mathbb{R}^{d_{\theta } }: \|\theta \|\leq \rho \}$, 
while $\Lambda_{1}$, $\Lambda_{2}$ are the events defined by 
\begin{align*}
	\Lambda_{1}
	=
	\limsup_{n\rightarrow\infty } \{\|\theta_{n} \|<\rho_{1} \} 
	=
	\bigcap_{m=0}^{\infty } \bigcup_{n=m}^{\infty } \{\|\theta_{n} \|<\rho_{1} \}, 
	\;\;\;\;\; 
	\Lambda_{2}
	=
	\liminf_{n\rightarrow\infty } \{\|\theta_{n} \|<\rho \} 
	=
	\bigcup_{m=0}^{\infty } \bigcap_{n=m}^{\infty } \{\|\theta_{n} \|<\rho \}. 
\end{align*}
On the other side, 
let $\tilde{C}_{2}\in [1,\infty )$ stand for a (global) Lipschitz constant of $\nabla f(\cdot )$
and for an upper bound of $\|\nabla f(\cdot )\|$ on $Q$. 
Moreover, let  
$\tilde{C}_{3} = \tilde{C}_{2} \exp(\tilde{C}_{2} )$, 
$\tilde{C}_{4} = 15\tilde{C}_{1}\tilde{C}_{2}\tilde{C}_{3}$, 
while $\tau=8^{-1}\tilde{C}_{4}^{-1}c^{2}$.  

In order to prove the theorem's assertion, it is sufficient to show $N_{0}^{c}\subseteq\Lambda_{2}$
(i.e., to demonstrate that on $N_{0}^{c}$, $\|\theta_{n}\|<\rho$ for all, 
but finitely many $n$).\footnote
{Assumption \ref{a1.2} is a consequence of Assumption \ref{aa2.2}, 
and therefore,  
Parts (i), (ii) of the theorem directly follow from Theorems \ref{theorem1.1}, \ref{theorem1.2}. } 
To demonstrate this, we use contradiction. 
We assume that 
$\|\theta_{n} \|\geq\rho$ for infinitely many $n$ and some $\omega\in N_{0}^{c}$. 
Notice that all formulas which follow in the proof correspond to $\omega$. 


Owing to (\ref{ta2.1.1701}), there exists an integer $k_{1}\geq 0$ (depending on $\omega$)
such that 
\begin{align}\label{ta2.1.1}
	g(\theta_{n} ) 
	\max_{n\leq j\leq a(n,T) } 
	\left\|
	\sum_{i=n}^{j} \alpha_{i} \xi_{i} 
	\right\| 
	\leq 
	\tau^{2}
\end{align}
for $n\geq k_{1}$. 
Hence
$
	\lim_{n\rightarrow\infty } g(\theta_{n} ) \|\alpha_{n}\xi_{n} \|
	=
	0
$. 
Since 
\begin{align*}
	&
	g(\theta_{n} ) \|\theta_{n+1} - \theta_{n} \|
	\leq
	\alpha_{n}
	+
	g(\theta_{n} ) \|\alpha_{n}\xi_{n} \|
\end{align*}
for $n\geq 0$, we get 
$\lim_{n\rightarrow\infty } g(\theta_{n} ) \|\theta_{n+1} - \theta_{n} \| = 0$. 
As
$\lim_{n\rightarrow\infty } \sum_{i=n}^{a(n,\tau)-1} \alpha_{i} = \tau$,  
there exists an integer $k_{2}\geq 0$ (depending on $\omega$)
such that 
\begin{align}\label{ta2.1.3001}
	\sum_{i=n}^{a(n,\tau)-1} \alpha_{i}
	\geq
	\tau/2,
	\;\;\;\;\; 
	g(\theta_{n} ) \|\theta_{n+1} - \theta_{n} \|
	\leq
	\tau
\end{align}
for $n\geq k_{2}$. 

Let $k_{0}=\max\{k_{1}, k_{2} \}$. 
Moreover, let $l_{0}, m_{0}, n_{0}$ be the integers defined as follows. 
If $\omega\in\Lambda_{1}$
(i.e., if $\|\theta_{n} \|<\rho_{1}$ for infinitely many $n$), let 
\begin{align}\label{ta2.1.3005}
	l_{0}
	=
	\min\{n>k_{0}: \|\theta_{n-1} \|<\rho_{1} \}, 
	\;\;\;\: 
	m_{0}
	=
	\min\{n>l_{0}: 
	\|\theta_{n} \|\geq\rho \}, 
	\;\;\;\;\: 
	n_{0}
	=
	\max\{n\leq m_{0}: \|\theta_{n-1} \|<\rho_{1} \}. 
\end{align}
Otherwise, if $\omega\in\Lambda_{1}^{c}$
(i.e., if $\|\theta_{n} \|<\rho_{1}$ for finitely many $n$), let 
\begin{align*}
	l_{0}
	=
	\max\{n>0: \|\theta_{n-1} \|<\rho_{1} \}, 
	\;\;\;\;\; 
	m_{0}
	=
	\infty,
	\;\;\;\;\; 
	n_{0}
	=
	\max\{k_{0}, l_{0} \}. 
\end{align*}
Then, we have $k_{0}<n_{0}\leq m_{0}$
and $\|\theta_{n} \|\geq\rho_{1}$ for $n_{0}\leq n<m_{0}$. 

Let $\zeta_{n,k}$, $\zeta'_{n,k}$, $\zeta''_{n,k}$, 
$\phi_{n,k}$, $\phi'_{n,k}$, $\phi''_{n,k}$
have the same meaning as in Section \ref{section1*}, 
while $\phi_{n}=\phi_{n,a(n,\tau) }$, $\phi'_{n}=\phi'_{n,a(n,\tau) }$, $\phi''_{n}=\phi''_{n,a(n,\tau) }$. 
Now, the asymptotic properties of $\phi_{n}$ are analyzed. 
This analysis follows the same arguments as the proof of Lemma \ref{lemma1.2}. 
Due to (\ref{1.1*}), we have 
\begin{align}\label{ta2.1.1703}
	f(\theta_{a(n,\tau) } ) - f(\theta_{n} )
	=
	-
	(\gamma_{a(n,\tau) } - \gamma_{n} ) \|\nabla f(\theta_{n} ) \|^{2} 
	-
	\phi_{n}
\end{align}
for $n\geq 0$. 
On the other side, (\ref{ta2.1.1}) implies 
\begin{align}\label{ta2.1.51}
	\left\|
	\zeta'_{n,j}
	\right\|
	=
	\left\|
	\sum_{i=n}^{j-1} \alpha_{i}\xi_{i} 
	\right\|
	\leq 
	\tau^{2} g^{-1}(\theta_{n} ) 
\end{align}
for $n_{0}\leq n\leq j\leq \min\{m_{0}, a(n,T)+1 \}$. 
Therefore, 
\begin{align*}
	\|\nabla f(\theta_{j} ) \|
	\leq &
	\|\nabla f(\theta_{n} ) \|
	+
	\|\nabla f(\theta_{j} ) - \nabla f(\theta_{n} ) \|
	\nonumber\\
	\leq &
	\|\nabla f(\theta_{n} ) \| 
	+
	\tilde{C}_{2} \|\theta_{j} - \theta_{n} \| 
	\nonumber\\
	\leq &
	\|\nabla f(\theta_{n} ) \|
	+
	\tilde{C}_{2} \sum_{i=n}^{j-1} \alpha_{i} \|\nabla f(\theta_{i} ) \| 
	+
	\tilde{C}_{2} 
	\left\|
	\sum_{i=n}^{j-1} \alpha_{i} \xi_{i} 
	\right\| 
	\nonumber\\
	\leq &
	\|\nabla f(\theta_{n} ) \| 
	+ 
	\tilde{C}_{2} \tau^{2} g^{-1}(\theta_{n} ) 
	+
	\tilde{C}_{2} \sum_{i=n}^{j-1} \alpha_{i} \|\nabla f(\theta_{i} ) \| 
\end{align*}
for $n_{0}\leq n< j \leq \min\{m_{0}, a(n,\tau)+1 \}$
(notice that $\tau<T$). 
Combining this with Bellman-Gronwall inequality
(see e.g. \cite[Appendix B]{borkar}), we deduce
\begin{align*}
	\|\nabla f(\theta_{j} ) \|
	\leq &
	\left(\|\nabla f(\theta_{n} ) \| + \tilde{C}_{2}\tau^{2} g^{-1}(\theta_{n} ) \right) 
	\exp\left(
	\tilde{C}_{2} \sum_{i=n}^{j-1} \alpha_{i}  
	\right) 
	\nonumber\\
	\leq &
	\left(\|\nabla f(\theta_{n} ) \| + \tilde{C}_{2}\tau^{2} g^{-1}(\theta_{n} ) \right) 
	(1 + \tilde{C}_{3}\tau )
	\nonumber \\
	\leq &
	\|\nabla f(\theta_{n} ) \| 
	+
	(\tilde{C}_{3}\tau + \tilde{C}_{2}\tau^{2} + \tilde{C}_{2}\tilde{C}_{3}\tau^{3} )
	g^{-1}(\theta_{n} ) 
	\nonumber\\
	\leq &
	\|\nabla f(\theta_{n} ) \| 
	+
	\tilde{C}_{4}\tau g^{-1}(\theta_{n} ) 
\end{align*}
for $n_{0} \leq n\leq j \leq \min\{m_{0}, a(n,\tau)+1 \}$.\footnote
{Notice that $\sum_{i=n}^{j-1} \alpha_{i} \leq \tau<1$
when $n \leq j \leq a(n,\tau)$. 
Notice also that 
$\exp(\tilde{C}_{2}\tau )\leq \tilde{C}_{2}\tau\exp(\tilde{C}_{2}\tau )\leq \tilde{C}_{3}\tau$.} 
Consequently, 
\begin{align}\label{ta2.1.5}
	\|\theta_{j} - \theta_{n} \|
	\leq &
	\sum_{i=n}^{j-1} \alpha_{i} \|\nabla f(\theta_{i} ) \| 
	+ 
	\left\|
	\sum_{i=n}^{j-1} \alpha_{i} \xi_{i} 
	\right\| 
	\leq 
	\left(
	\|\nabla f(\theta_{n} ) \|
	+
	\tilde{C}_{4}\tau g^{-1}(\theta_{n} ) 
	\right)
	\left(
	\sum_{i=n}^{j} \alpha_{i} 
	+
	\tau
	\right)
	\leq 
	3 \tau g^{-1}(\theta_{n} )
\end{align}
for $n_{0} \leq n\leq j < \min\{m_{0}, a(n,\tau)+1 \}$
(notice that $\tilde{C}_{4}\tau\leq 1/2$). 
Hence, 
\begin{align*}
	\left\|\zeta''_{n,j} \right\|
	=
	\left\|
	\sum_{i=n}^{j-1} \alpha_{i} (\nabla f(\theta_{i} ) - \nabla f(\theta_{n} ) )
	\right\|
	\leq 
	\tilde{C}_{2} 
	\sum_{i=n}^{j-1} \alpha_{i} \|\theta_{i} - \theta_{n} \|
	\leq
	3\tilde{C}_{2}\tau g^{-1}(\theta_{n} ) \sum_{i=n}^{j-1} \alpha_{i}
	\leq 
	3\tilde{C}_{2}\tau^{2} g^{-1}(\theta_{n} ) 
\end{align*}
for $n_{0} \leq n\leq j \leq \min\{m_{0}, a(n,\tau)+1 \}$. 
Combining this with (\ref{ta2.1.51}), we get
\begin{align*}
	\|\zeta_{n,j} \|
	\leq 
	\|\zeta'_{n,j} \| + \|\zeta''_{n,j} \|
	\leq 
	4\tilde{C}_{2}\tau^{2} g^{-1}(\theta_{n} )
\end{align*}
for $n_{0} \leq n\leq j \leq \min\{m_{0}, a(n,\tau)+1 \}$. 
Therefore, 
\begin{align*}
	|\phi'_{n} |
	\leq &
	\|\nabla f(\theta_{n} ) \|\: \|\zeta_{n,a(n,\tau ) } \|
	\leq 
	4\tilde{C}_{2} \tau^{2} 
	g^{-2}(\theta_{n} )
\end{align*}
for $n\geq n_{0}$ satisfying $a(n,\tau )<m_{0}$. 
We also have 
\begin{align*}
	|\phi''_{n} |
	\leq &
	\tilde{C}_{2} \|\theta_{a(n,\tau)} - \theta_{n} \|^{2}
	\leq 
	9\tilde{C}_{2} \tau^{2} 
	g^{-2}(\theta_{n} )
\end{align*}
for $n\geq n_{0}$ satisfying $a(n,\tau )<m_{0}$. 
Thus, 
\begin{align}\label{ta2.1.307}
	|\phi_{n} |
	\leq &
	\tilde{C}_{4} \tau^{2} 
	g^{-2}(\theta_{n} )
\end{align}
when $n\geq n_{0}$, $a(n,\tau )\leq m_{0}$. 
Then, (\ref{ta2.1.3001}), (\ref{ta2.1.1703}), (\ref{ta2.1.307}) imply 
\begin{align}\label{ta2.1.7}
	f(\theta_{a(n,\tau) } ) - f(\theta_{n} ) 
	\leq&
	-
	\tau\|\nabla f(\theta_{n} ) \|^{2}/2
	+
	\tilde{C}_{4}\tau^{2} g^{-2}(\theta_{n} ) 
	\leq 
	-\tau c^{2}/8
\end{align}
for $n\geq n_{0}$ satisfying $a(n,\tau )<m_{0}$.\footnote{
Notice that 
$\|\nabla f(\theta_{n} ) \|^{2}/2 \geq
\|\nabla f(\theta_{n} ) \|^{2}/4 + c^{2}/4 \geq 
\|\nabla f(\theta_{n} ) \|^{2}/8 + \tilde{C}_{4}\tau g^{-2}(\theta_{n} )$ for $n_{0}\leq n<m_{0}$. }

Let $\{n_{k} \}_{k\geq 0}$ be the sequence recursively defined by 
$n_{k+1}=a(n_{k},\tau )$ for $k\geq 0$. 
Now, we show by contradiction $\omega\in\Lambda_{1}$ 
(i.e., $\|\theta_{n} \|<\rho_{1}$ for infinitely many $n$). 
We assume the opposite. 
Then, $m_{0}=\infty$ and $\|\theta_{n} \|\geq\rho_{1}$ for $n\geq n_{0}$, 
while (\ref{ta2.1.7}) implies 
$f(\theta_{n_{k+1} } ) - f(\theta_{n_{k} } ) \leq - \tau c^{2}/8$ 
for $k\geq 0$. 
Hence, $\lim_{k\rightarrow\infty } f(\theta_{n_{k} } ) = -\infty$.  
However, this is impossible due to Assumption \ref{aa2.1}. 
Thus, $\omega\in\Lambda_{1}$ (i.e., $\|\theta_{n} \|<\rho_{1}$ for infinitely many $n$). 
Therefore, $m_{0}, n_{0}$ are defined through (\ref{ta2.1.3005}), 
while $\|\theta_{n_{0}-1} \|<\rho_{1}$, $\|\theta_{m_{0} } \|\geq\rho$. 
Combining this with (\ref{ta2.1.3001}), we conclude
\begin{align*}
	&
	\|\theta_{n_{0} } - \theta_{n_{0}-1} \| 
	\leq 
	\tau g^{-1}(\theta_{n_{0}-1} ) 
	\leq 
	\tau(\tilde{C}_{2} + 1 ) 
	\leq 
	1/2
\end{align*}
(notice that $\|\nabla f(\theta_{n_{0}-1} ) \|\leq \tilde{C}_{2}$, 
$\tilde{C}_{2}\tau\leq 1/4$). 
Consequently, 
\begin{align}\label{ta2.1.3009}
	&
	\|\theta_{n_{0} } \|
	\leq 
	\|\theta_{n_{0}-1} \|
	+
	\|\theta_{n_{0} } - \theta_{n_{0}-1} \|
	\leq \rho_{1}+1/2
	<\rho. 
\end{align}
Hence, $n_{0}<m_{0}$, 
$f(\theta_{n_{0} } ) \leq \tilde{C}_{1}$. 

Let $i_{0}, j_{0}$ be the integers defined by 
$j_{0} = \max\{j\geq 0: n_{j}< m_{0} \}$, 
$i_{0} = n_{j_{0} }$. 
Then, we have $n_{0}\leq i_{0}=n_{j_{0} }< m_{0}\leq a(i_{0},\tau ) = n_{j_{0}+1}$. 
As a result of this and (\ref{ta2.1.5}), we get
\begin{align*}
	&
	\|\theta_{m_{0} } - \theta_{i_{0} } \| 
	\leq 
	3 \tau g^{-1}(\theta_{i_{0} } ) 
	\leq 
	3 \tau (\tilde{C}_{2} + 1 )
	\leq 
	1/2  
\end{align*}
(notice that $\|\nabla f(\theta_{i_{0} } ) \|\leq \tilde{C}_{2}$, 
$\tilde{C}_{2}\tau\leq 1/12$). 
Consequently, 
\begin{align}
	&\label{ta2.1.703}
	\|\theta_{i_{0} } \|
	\geq 
	\|\theta_{m_{0} } \| 
	-
	\|\theta_{m_{0} } - \theta_{i_{0} } \|
	\geq 
	\rho - 1/2
	>
	\rho_{2}. 
\end{align}

Let $\theta_{0}(\cdot )$ be the stochastic processes defined by 
$\theta_{0}(t)=\theta_{n}$ for $t\in[\gamma_{n},\gamma_{n+1} )$, $n\geq 0$
($\{\gamma_{n} \}_{n\geq 0}$ is defined in Section \ref{section1}). 
Now, we show by contradiction that $\gamma_{i_{0} } - \gamma_{n_{0} }\geq T$. 
We assume the opposite. 
Then, (\ref{ta2.1.51}), (\ref{ta2.1.3009}) imply 
\begin{align}\label{ta2.1.705}
	\|\theta_{0}(t) \|
	=
	\|\theta_{j} \|
	\leq&
	\|\theta_{n_{0} } \| 
	+
	\sum_{i=n_{0} }^{j-1} \alpha_{i} \|\nabla f(\theta_{i} ) \| 
	+
	\left\|
	\sum_{i=n_{0} }^{j-1} \alpha_{i} \xi_{i} 
	\right\|
	\nonumber\\
	\leq&
	\|\theta_{n_{0} } \| 
	+
	\tau^{2} g^{-1}(\theta_{n_{0} } )
	+
	\sum_{i=n_{0} }^{j-1} \alpha_{i} \|\nabla f(\theta_{i} ) \|   
	\nonumber\\
	\leq&
	\rho_{1}+1
	+
	\sum_{i=n_{0} }^{j-1} \alpha_{i} \phi(\|\theta_{i} \| )  
	\nonumber\\
	\leq&
	\rho_{1}+1
	+
	\int_{\gamma_{n_{0} } }^{t} \phi(\|\theta_{0}(s) \| ) ds
\end{align}
for $t\in [\gamma_{j}, \gamma_{j+1} )$, $n_{0}\leq j\leq i_{0}$.\footnote
{As $j\leq i_{0} < m_{0}$, 
we have  
$\gamma_{j}-\gamma_{n_{0} }\leq \gamma_{i_{0} }-\gamma_{n_{0} } \leq T$
and $j<\min\{m_{0},a(n_{0},T)+1\}$. 
We also have $\tau^{2} g^{-1}(\theta_{n_{0} } ) \leq \tau^{2} (\tilde{C}_{2} + 1 ) \leq 1/2$.}
Applying the comparison principle (see \cite[Section 3.4]{khalil}) to
(\ref{ta2.1.701}), (\ref{ta2.1.705}), we conclude 
$\|\theta_{0}(t) \| \leq 
\lambda(t-\gamma_{n_{0} };\rho_{1}+1) \leq \rho_{2}$
for all $t\in[\gamma_{{n}_{0} }, \gamma_{i_{0} } ]$. 
Thus, $\|\theta_{i_{0} } \|=\|\theta_{0}(\gamma_{i_{0} } ) \|\leq \rho_{2}$. 
However, this is impossible, due to (\ref{ta2.1.703}). 
Hence, $\gamma_{i_{0} }-\gamma_{n_{0} }\geq T$. 
Consequently, 
\begin{align}\label{ta2.1.707}
	T
	\leq 
	\gamma_{i_{0} } - \gamma_{n_{0} }
	=
	\sum_{j=0}^{j_{0}-1} (\gamma_{n_{j+1} } - \gamma_{n_{j} } )
	\leq 
	j_{0}\tau 
\end{align}
(notice that $n_{j_{0} } = i_{0}$, 
$\gamma_{n_{j+1} }-\gamma_{n_{j} } = \sum_{i=n_{j} }^{n_{j+1}-1} \alpha_{i} \leq \tau$). 

Owing to (\ref{ta2.1.7}), we have 
$f(\theta_{n_{j+1} } ) - f(\theta_{n_{j} } ) \leq -\tau c^{2}/8$
for $0\leq j \leq j_{0}$. 
Combining this with (\ref{ta2.1.707}), we get 
\begin{align*}
	f(\theta_{i_{0} } )
	=
	f(\theta_{n_{j_{0} } } ) 
	\leq 
	f(\theta_{n_{0} } ) - j_{0}\tau c^{2}/8
	\leq 
	\tilde{C}_{1} - c^{2} T/8
	\leq
	-\tilde{C}_{1}. 
\end{align*}
However, this is impossible, since $f(\theta ) > -\tilde{C}_{1}$ for all 
$\theta\in\mathbb{R}^{d_{\theta } }$. 
Hence, $\|\theta_{n} \|\geq\rho$ for finitely many $n$. 
\end{sproof}

In the rest of section, the results of Section \ref{section1} are extended to 
randomly projected stochastic gradient algorithms. 
These algorithms are defined by the following difference equations: 
\begin{align}\label{appendix2.1}
	&
	\vartheta_{n}
	=
	\theta_{n}
	-
	\alpha_{n} (\nabla f(\theta_{n} ) + \xi_{n} ), 
	\nonumber\\
	&
	\theta_{n+1}
	=
	\vartheta_{n} 
	I_{ \{\|\vartheta_{n} \| \leq \beta_{\sigma_{n} } \} }
	+
	\theta_{0} 
	I_{ \{\|\vartheta_{n} \| > \beta_{\sigma_{n} } \} }, 
	\nonumber\\
	&
	\sigma_{n+1}
	=
	\sigma_{n}
	+
	I_{ \{\|\vartheta_{n} \| > \beta_{\sigma_{n} } \} }, 
	\;\;\;\;\; 
	n\geq 0. 
\end{align}
Here, $\nabla f(\cdot )$, $\{\alpha_{n} \}_{n\geq 0}$, $\{\xi_{n} \}_{n\geq 0}$ 
have the same meaning as in Section \ref{section1}, 
while $\{\beta_{n} \}_{n\geq 0}$ is an increasing sequence of positive real numbers satisfying
$\lim_{n\rightarrow\infty } \beta_{n} = \infty$. 
$\theta_{0}\in\mathbb{R}^{d}$ is a (deterministic) vector satisfying $\|\theta_{0} \|\leq \beta_{0}$, 
while $\sigma_{0}=0$. 
For further details on randomly projected stochastic gradient search and 
stochastic approximation, see \cite{chen}, \cite{tadic11} and references cited therein. 

To study the asymptotic behavior of (\ref{appendix2.1}), we introduce the following two assumption. 

\begin{assumptionappendix}\label{aa2.3} 
$\liminf_{\|\theta\|\rightarrow\infty } \|\nabla f(\theta ) \| > 0$
and 
$\inf_{\theta\in\mathbb{R}^{d_{\theta } } } f(\theta )>-\infty $. 
\end{assumptionappendix}

\begin{assumptionappendix}\label{aa2.4} 
There exists a real number $r\in(1,\infty )$ such that 
\begin{align}\label{aa2.4.1}
	\limsup_{n\rightarrow\infty } 
	\max_{n\leq j\leq a(n,t) } 
	\left\| 
	\sum_{i=n}^{j} \alpha_{i} \gamma_{i}^{r} \xi_{i} 
	\right\|
	I_{ \{\tau_{Q,n}>j \} }
	<\infty 
\end{align}
almost surely for all $t\in(0,\infty )$ and any compact set $Q\subset\mathbb{R}^{d_{\theta } }$. 
Here, $\tau_{Q,n}$ is the stopping time defined by 
\begin{align*}
	\tau_{Q,n}
	=
	\inf\left(\{j\geq n: \theta_{j}\neq\vartheta_{j-1} \text{ or } \theta_{j}\not\in Q \} \cup \{\infty \} \right) 
\end{align*}
for $n\geq 0$. 
\end{assumptionappendix}

Assumption \ref{aa2.3} is a stability condition. 
It is one of the weakest conditions under which the stability of the ODE
$d\theta/dt=-\nabla f(\theta )$ can be demonstrated.  
On the other side, Assumption \ref{aa2.4} is a noise condition. 
It can be considered as a version of the noise conditions adopted in \cite{tadic11}. 

Our results on the asymptotic behavior of algorithm (\ref{appendix2.1}) are provided in the next theorem. 

\begin{theoremappendix}\label{theorema2.2}
Let $\{\theta_{n} \}_{n\geq 0}$ be generated by (\ref{appendix2.1}). 
Moreover, let Assumptions \ref{a1.1}, \ref{aa2.3} and \ref{aa2.4} hold. 
Then, there exists a real number $\rho\in[1,\infty )$
such that $\limsup_{n\rightarrow\infty } \|\theta_{n} \| < \rho$ almost surely. 
If in addition, Assumption \ref{a1.3} holds, 
then the following is also true: 
\begin{compactenum}[(i)]
\item
$\hat{\theta}=\lim_{n\rightarrow\infty } \theta_{n}$
exists and satisfies $\nabla f(\hat{\theta} ) = 0$ almost surely. 
\item
There exists a random variable $\hat{K}$ 
such that $0<\hat{K}<\infty$ everywhere and such that (\ref{t1.1.1*}) -- (\ref{t1.1.5*}) 
holds almost surely. 
\end{compactenum}

\end{theoremappendix}

\begin{sproof}
Due to Assumption \ref{aa2.3}, 
there exist real numbers 
$c\in(0,1)$, $\rho_{1},\tilde{C}_{1} \in [1,\infty )$ 
such that the following is true: 
(i) $\|\nabla f(\theta ) \|\geq c$ for each $\theta\in\mathbb{R}^{d_{\theta } }$
satisfying $\|\theta \|\geq\rho_{1}$, 
(ii) $f(\theta )\leq \tilde{C}_{1}$
for any $\theta\in\mathbb{R}^{d_{\theta } }$ satisfying $\|\theta \|\leq \rho_{1}+1$, 
and 
(iii) $f(\theta )> -\tilde{C}_{1}$ for all $\theta\in\mathbb{R}^{d_{\theta } }$. 
Without loss of generality, it can also be assumed $\|\theta_{0} \|<\rho_{1}$. 
On the other side, owing to Assumption \ref{aa2.4}, 
there exists an event $N_{0}\in {\cal F}$ with the following properties: 
(i) $P(N_{0} ) = 0$,
and 
(ii) (\ref{aa2.4.1}) holds on $N_{0}^{c}$ for all $t\in(0,\infty )$
and any compact set $Q\subset\mathbb{R}^{d_{\theta } }$. 
Then, relying on the same arguments as in the proof of Lemma \ref{lemma1.1}, we conclude 
\begin{align}\label{ta2.2.1701}
	\lim_{n\rightarrow\infty } 
	\max_{n\leq j< a(n,t) } 
	\left\|
	\sum_{i=n}^{j} \alpha_{i}\xi_{i} 
	\right\|
	I_{ \{\tau_{Q,n}>j \} }
	=
	0
\end{align}
on $N_{0}^{c}$ for all $t\in(0,\infty )$ and any compact set 
$Q\subset\mathbb{R}^{d_{\theta } }$. 

Let $T=16\tilde{C}_{1}c^{-2}$, 
while $\phi(\cdot )$, $\lambda(\cdot\:;\cdot )$ 
have the same meaning as in the proof of Theorem \ref{theorema2.1}. 
Then, there exists $\rho_{2}\in [1,\infty )$ such that 
$\rho_{2}\geq\rho_{1}+1$ and such that 
$|\lambda(t;z) |\leq\rho_{2}$ for all $t\in [0,T]$, $z\in [0,\rho_{1}+1]$. 
Moreover, (\ref{ta2.1.701}) holds for all $t,z\in [0,\infty )$. 

Let $\rho=\rho_{2}+1$, $Q=\{\theta\in\mathbb{R}^{d_{\theta } }: \|\theta \|\leq \rho \}$. 
Moreover, let $\sigma = \lim_{n\rightarrow\infty } \sigma_{n}$ and 
\begin{align*}
	\Lambda_{1}
	=
	\limsup_{n\rightarrow\infty } \{\|\theta_{n} \|<\rho_{1} \} 
	=
	\bigcap_{m=0}^{\infty } \bigcup_{n=m}^{\infty } \{\|\theta_{n} \|<\rho_{1} \},  
	\;\;\;\;\; 
	\Lambda_{2}
	=
	\liminf_{n\rightarrow\infty } \{\|\theta_{n} \|<\rho \} 
	=
	\bigcup_{m=0}^{\infty } \bigcap_{n=m}^{\infty } \{\|\theta_{n} \|<\rho \},  
\end{align*}
while 
$\tilde{\rho}=\rho I_{\Lambda_{1} } + \beta(\sigma) I_{\Lambda_{1}^{c} }$, 
$\tilde{Q}=\{\theta\in\mathbb{R}^{d_{\theta } }: \|\theta \|\leq \tilde{\rho} \}$. 
As $\sigma<\infty$ on $\Lambda_{1}^{c}$, 
we have $\theta_{n}, \vartheta_{n}\in \tilde{Q}$ for $n\geq 0$ on the same event. 
We also have $\tilde{\rho}<\infty$ everywhere. 
Consequently, 
\begin{align}\label{ta2.2.5001}
	\lim_{n\rightarrow\infty } 
	\max_{n\leq j\leq a(n,t) }
	\left\|
	\sum_{i=n}^{j} \alpha_{i}\xi_{i}
	\right\|
	I_{ \{\tau_{\tilde{Q},n}>j \} }
	=0
\end{align}
on $N_{0}^{c}$ for all $t\in(0,\infty )$. 

Let $\tilde{C}_{2}\in [1,\infty )$ stand for a local Lipschitz constant of $\nabla f(\cdot )$
on $\tilde{Q}$ and for an upper bound of $\|\nabla f(\cdot )\|$ on the same set. 
In addition to this, let  
$\tilde{C}_{3} = \tilde{C}_{2} \exp(\tilde{C}_{2} )$, 
$\tilde{C}_{4} = 21\tilde{C}_{1}\tilde{C}_{3}^{3}$, 
while $\tau=8^{-1}\tilde{C}_{4}^{-1}c^{-2}$.  

In order to prove the theorem's assertion, it is sufficient to show $N_{0}^{c}\subseteq\Lambda_{2}$
(i.e., to demonstrate that on $N_{0}^{c}$, $\|\vartheta_{n} \|<\rho$ for all, 
but finitely many $n$).\footnote
{On $\Lambda_{2}$, the following holds: 
$\sigma<\infty$ and $\theta_{n}=\vartheta_{n}$, $\tau_{Q,n}=\infty$
for $n>\sigma$. 
Hence, algorithm (\ref{appendix2.1}) asymptotically reduces to (\ref{1.1}) on $\Lambda_{2}$, 
while (\ref{a1.2.1}) holds almost surely on the same event. 
Therefore, Parts (i), (ii) of the theorem directly follow 
from Theorems \ref{theorem1.1}, \ref{theorem1.2}. } 
To demonstrate this, we use contradiction. 
We assume that 
$\|\vartheta_{n} \|\geq\rho$ for infinitely many $n$ and some $\omega\in N_{0}^{c}$. 
Notice that all formulas which follow in the proof correspond to $\omega$. 

As $\{\beta(\sigma_{n} ) \}_{n\geq 0}$ is non-decreasing, 
we have $\beta(\sigma_{n} ) > \rho$ for all, but finitely many $n$.\footnote
{If $\sigma<\infty$, 
then $\rho<\|\theta_{n} \| = \|\vartheta_{n-1} \| \leq \beta(\sigma_{n-1} )$ 
for all, but finitely many $n$. 
On the other side, if $\sigma=\infty$, 
then $\lim_{n\rightarrow\infty } \beta(\sigma_{n} ) = \infty$. } 
Hence, there exists an integer $k_{1}$ (depending on $\omega$) such that 
$\beta(\sigma_{n} )> \rho$ for $n\geq k_{1}$. 
On the other side, 
due to (\ref{ta2.2.5001}), there exists an integer $k_{2}\geq 0$ (depending on $\omega$)
such that 
\begin{align}\label{ta2.2.1}
	\max_{n\leq j\leq a(n,T) } 
	\left\|
	\sum_{i=n}^{j} \alpha_{i} \xi_{i} 
	\right\| 
	I_{ \{\tau_{\tilde{Q},n}>j\} }
	\leq 
	\tau^{2} 
\end{align}
for $n\geq k_{2}$. 
Hence, 
$
	\lim_{n\rightarrow\infty } \|\alpha_{n}\xi_{n} \| I_{ \{\theta_{n}\in\tilde{Q} \} } 
	=
	0
$. 
Since 
\begin{align*}
	&
	\|\vartheta_{n} - \theta_{n} \| I_{ \{\theta_{n}\in\tilde{Q} \} } 
	\leq
	\left(
	\tilde{C}_{2} \alpha_{n}
	+
	\|\alpha_{n}\xi_{n} \|
	\right)
	I_{ \{\theta_{n}\in\tilde{Q} \} } 
\end{align*}
for $n\geq 0$, we get 
$\lim_{n\rightarrow\infty } \|\vartheta_{n} - \theta_{n} \| I_{ \{\theta_{n}\in\tilde{Q} \} } = 0$. 
Then, there exists an integer $k_{3}\geq 0$ (depending on $\omega$)
such that 
\begin{align}\label{ta2.2.3001}
	\sum_{i=n}^{a(n,\tau)-1} \alpha_{i}
	\geq
	\tau/2,
	\;\;\;\;\; 
	\|\vartheta_{n} - \theta_{n} \| I_{ \{\theta_{n}\in\tilde{Q} \} } 
	\leq
	\tau
\end{align}
for $n\geq k_{3}$
(notice that $\lim_{n\rightarrow\infty } \sum_{i=n}^{a(n,\tau)-1} \alpha_{i} = \tau$). 

Let $k_{0}=\max\{k_{1}, k_{2}, k_{3} \}$. 
Moreover, let $l_{0}, m_{0}, n_{0}$ be the integers defined as follows. 
If $\omega\in\Lambda_{1}$
(i.e., if $\|\theta_{n} \|<\rho_{1}$ for infinitely many $n$), let 
\begin{align}\label{ta2.2.3005}
	l_{0}
	=
	\min\{n>k_{0}: \|\theta_{n-1} \|<\rho_{1} \}, 
	\;\;\; 
	m_{0}
	=
	\min\{n>l_{0}: 
	\|\vartheta_{n-1} \|\geq\rho \}, 
	\;\;\;\; 
	n_{0}
	=
	\max\{n\leq m_{0}: \|\theta_{n-1} \|<\rho_{1} \}. 
\end{align}
Otherwise, if $\omega\in\Lambda_{1}^{c}$
(i.e., if $\|\theta_{n} \|<\rho_{1}$ for finitely many $n$), let 
\begin{align*}
	l_{0}
	=
	\max\{n>0: \|\theta_{n-1} \|<\rho_{1} \}, 
	\;\;\;\;\; 
	m_{0}
	=
	\infty,
	\;\;\;\;\; 
	n_{0}
	=
	\max\{k_{0}, l_{0} \}. 
\end{align*}
Then, we have $k_{0}<n_{0}\leq m_{0}$
and $\theta_{n}=\vartheta_{n-1}$, $\rho_{1}\leq\|\theta_{n} \|\leq\tilde{\rho}$ 
for $n_{0}\leq n<m_{0}$.\footnote
{If $\theta_{n}\neq\vartheta_{n-1}$, we have $\|\theta_{n} \|=\|\theta_{0} \|<\rho_{1}$. 
On the other side, if $\omega\in\Lambda_{1}$, 
then $\|\theta_{n} \|=\|\vartheta_{n-1} \|\leq\rho=\tilde{\rho}$ 
for $n_{0}\leq n<m_{0}$. 
Moreover, if $\omega\in\Lambda_{1}^{c}$, 
then $\|\theta_{n} \|\leq\beta(\sigma_{n-1} ) \leq \beta(\sigma )=\tilde{\rho}$
for $n>0$. }
Therefore, 
\begin{align}\label{ta2.2.3091}
	c\leq \|\nabla f(\theta_{n} ) \|\leq \tilde{C}_{2}, 
	\;\;\;\;\; 
	\theta_{n}\in\tilde{Q}, 
	\;\;\;\;\; 
	\tau_{n,\tilde{Q} }\geq m_{0}
\end{align}
for $n_{0}\leq n<m_{0}$, while 
\begin{align}\label{ta2.2.3071}
	\theta_{j}
	=
	\theta_{n}
	-
	\sum_{i=n}^{j-1} \alpha_{i} \nabla f(\theta_{i} ) 
	-
	\sum_{i=n}^{j-1} \alpha_{i} \xi_{i}
\end{align}
for $n_{0}\leq n<j<m_{0}$. 

Let $\zeta_{n,k}$, $\zeta'_{n,k}$, $\zeta''_{n,k}$, 
$\phi_{n,k}$, $\phi'_{n,k}$, $\phi''_{n,k}$
have the same meaning as in Section \ref{section1*}, 
while $\phi_{n}=\phi_{n,a(n,\tau) }$, $\phi'_{n}=\phi'_{n,a(n,\tau) }$, $\phi''_{n}=\phi''_{n,a(n,\tau) }$. 
Now, the asymptotic properties of $\phi_{n}$ are analyzed
using similar arguments as in the proof of Theorem \ref{theorema2.1}. 
Owing to (\ref{1.1*}), we have 
\begin{align}\label{ta2.2.1703}
	f(\theta_{a(n,\tau) } ) - f(\theta_{n} )
	=
	-
	(\gamma_{a(n,\tau) } - \gamma_{n} ) \|\nabla f(\theta_{n} ) \|^{2} 
	-
	\phi_{n}
\end{align}
for $n\geq n_{0}$ satisfying $a(n,\tau)<m_{0}$. 
On the other side, due to (\ref{ta2.2.1}), (\ref{ta2.2.3091}), we have
\begin{align}\label{ta2.2.51}
	\left\|\zeta'_{n,j} \right\|
	\leq 
	\left\|
	\sum_{i=n}^{j-1} \alpha_{i} \xi_{i}
	\right\|
	\leq 
	\tau^{2} 
\end{align}
for $n_{0}\leq n\leq j\leq \min\{m_{0}, a(n,T)+1 \}$.
Using (\ref{ta2.2.3071}), (\ref{ta2.2.51}), we deduce 
\begin{align*}
	\|\nabla f(\theta_{j} ) \|
	\leq &
	\|\nabla f(\theta_{n} ) \|
	+
	\|\nabla f(\theta_{j} ) - \nabla f(\theta_{n} ) \|
	\nonumber\\
	\leq &
	\|\nabla f(\theta_{n} ) \| 
	+
	\tilde{C}_{2} \|\theta_{j} - \theta_{n} \| 
	\nonumber\\
	\leq &
	\|\nabla f(\theta_{n} ) \|
	+
	\tilde{C}_{2} \sum_{i=n}^{j-1} \alpha_{i} \|\nabla f(\theta_{i} ) \| 
	+
	\tilde{C}_{2} 
	\left\|
	\sum_{i=n}^{j-1} \alpha_{i} \xi_{i} 
	\right\| 
	\nonumber\\
	\leq &
	\|\nabla f(\theta_{n} ) \| 
	+ 
	\tilde{C}_{2} \tau^{2} 
	+
	\tilde{C}_{2} \sum_{i=n}^{j-1} \alpha_{i} \|\nabla f(\theta_{i} ) \| 
\end{align*}
for $n_{0}\leq n< j < \min\{m_{0}, a(n,\tau)+1 \}$
(notice that $\tau<T$ and $\theta_{n},\theta_{j}\in\tilde{Q}$ for $n_{0}\leq n< j <m_{0}$). 
Then, Bellman-Gronwall inequality (see e.g., \cite[Appendix B]{borkar}) and (\ref{ta2.2.3091}) imply
\begin{align*}
	\|\nabla f(\theta_{j} ) \|
	\leq &
	\left(\|\nabla f(\theta_{n} ) \| + \tilde{C}_{2}\tau^{2} \right) 
	\exp\left(
	\tilde{C}_{2} \sum_{i=n}^{j-1} \alpha_{i}  
	\right) 
	\nonumber\\
	\leq &
	\left(\|\nabla f(\theta_{n} ) \| + \tilde{C}_{2} \tau^{2} \right) 
	(1 + \tilde{C}_{3}\tau )
	\nonumber \\
	\leq &
	\|\nabla f(\theta_{n} ) \| 
	+
	\tilde{C}_{2}\tilde{C}_{3}\tau + \tilde{C}_{2}\tau^{2} + \tilde{C}_{2}\tilde{C}_{3}\tau^{3} 
	\nonumber\\
	\leq &
	\|\nabla f(\theta_{n} ) \| 
	+
	\tilde{C}_{4}\tau 
\end{align*}
for $n_{0} \leq n\leq j < \min\{m_{0}, a(n,\tau)+1 \}$.\footnote
{Notice that $\sum_{i=n}^{j-1} \alpha_{i} \leq \tau<1$ when $n \leq j \leq a(n,\tau)$. 
Notice also that 
$\exp(\tilde{C}_{2}\tau )\leq \tilde{C}_{2}\tau\exp(\tilde{C}_{2}\tau )\leq \tilde{C}_{3}\tau$. } 
Owing to  (\ref{ta2.2.3091}), (\ref{ta2.2.3071}), (\ref{ta2.2.51}), we have 
\begin{align}\label{ta2.2.5}
	\|\theta_{j} - \theta_{n} \|
	\leq &
	\sum_{i=n}^{j-1} \alpha_{i} \|\nabla f(\theta_{i} ) \| 
	+ 
	\left\|
	\sum_{i=n}^{j-1} \alpha_{i} \xi_{i} 
	\right\| 
	\leq 
	\left(
	\|\nabla f(\theta_{n} ) \|
	+
	\tilde{C}_{4}\tau  
	\right)
	\left(
	\sum_{i=n}^{j} \alpha_{i} 
	+
	\tau
	\right)
	\leq 
	4\tilde{C}_{2} \tau 
\end{align}
for $n_{0} \leq n\leq j < \min\{m_{0}, a(n,\tau)+1 \}$
(notice that $\tilde{C}_{4}\tau\leq 1$). 
Hence, 
\begin{align*}
	\left\|\zeta''_{n,j} \right\|
	=
	\left\|
	\sum_{i=n}^{j-1} \alpha_{i} (\nabla f(\theta_{i} ) - \nabla f(\theta_{n} ) )
	\right\|
	\leq 
	\tilde{C}_{2} 
	\sum_{i=n}^{j-1} \alpha_{i} \|\theta_{i} - \theta_{n} \|
	\leq
	4\tilde{C}_{2}\tau \sum_{i=n}^{j-1} \alpha_{i}
	\leq 
	4\tilde{C}_{2}\tau^{2}  
\end{align*}
for $n_{0} \leq n\leq j < \min\{m_{0}, a(n,\tau)+1 \}$. 
Combining this with (\ref{ta2.1.51}), we get
\begin{align*}
	\|\zeta_{n,j} \|
	\leq 
	\|\zeta'_{n,j} \| + \|\zeta''_{n,j} \|
	\leq 
	5\tilde{C}_{2}\tau^{2} 
\end{align*}
for $n_{0} \leq n\leq j < \min\{m_{0}, a(n,\tau)+1 \}$. 
Consequently, 
\begin{align*}
	|\phi'_{n} |
	\leq &
	\|\nabla f(\theta_{n} ) \|\: \|\zeta_{n,a(n,\tau) } \|
	\leq
	5\tilde{C}_{2}^{3} \tau^{2} 
\end{align*}
for $n\geq n_{0}$ satisfying $a(n,\tau )<m_{0}$
(notice that $\|\nabla f(\theta_{n} ) \|\leq\tilde{C}_{2}$ 
for $n_{0}\leq n< m_{0}$). 
We also have 
\begin{align*}
	|\phi''_{n} |
	\leq &
	\tilde{C}_{2} \|\theta_{a(n,\tau)} - \theta_{n} \|^{2}
	\leq 
	16\tilde{C}_{2}^{3} \tau^{2} 
\end{align*}
for $n\geq n_{0}$ satisfying $a(n,\tau )<m_{0}$
(notice that 
and $\theta_{n}, \theta_{a(n,\tau ) } \in \tilde{Q}$
when $n\geq n_{0}$, $a(n,\tau)<m_{0}$). 
Hence, 
$|\phi_{n} |\leq \tilde{C}_{4} \tau^{2}$
when $n\geq n_{0}$, $a(n,\tau )<m_{0}$. 
Then, (\ref{ta2.2.3001}), (\ref{ta2.2.3091}), (\ref{ta2.2.1703}) imply 
\begin{align}\label{ta2.2.7}
	f(\theta_{a(n,\tau) } ) - f(\theta_{n} ) 
	\leq&
	-
	\tau\|\nabla f(\theta_{n} ) \|^{2}/2
	+
	\tilde{C}_{4}\tau^{2}  
	\leq 
	-
	\tau c^{2}/8
\end{align}
for $n\geq n_{0}$ satisfying $a(n,\tau )<m_{0}$.\footnote
{Notice that 
$\tau\|\nabla f(\theta_{n} ) \|^{2}/2\geq 
\tau\|\nabla f(\theta_{n} ) \|^{2}/4 + \tau c^{2}/4 \geq 
\tau\|\nabla f(\theta_{n} ) \|^{2}/8 + \tilde{C}_{4}\tau$ when $n_{0}\leq n<m_{0}$. } 

Let $\{n_{k} \}_{k\geq 0}$ be the sequence recursively defined by 
$n_{k+1}=a(n_{k},\tau )$ for $k\geq 0$. 
As in the proof of Theorem \ref{theorema2.1}, 
we now show by contradiction $\omega\in\Lambda_{1}$ 
(i.e., $\|\theta_{n} \|<\rho_{1}$ for infinitely many $n$). 
We assume the opposite. 
Then, $m_{0}=\infty$ and $\theta_{n}=\vartheta_{n}$, while (\ref{ta2.2.7}) yields  
$
	f(\theta_{n_{k+1} } ) - f(\theta_{n_{k} } )
	\leq 
	-\tau c^{2}/8
$
for $k\geq 0$. 
Hence, $\lim_{k\rightarrow\infty } f(\theta_{n_{k} } ) = -\infty$.  
However, this is impossible due to Assumption \ref{aa2.3}. 
Thus, $\omega\in\Lambda_{1}$ (i.e., $\|\theta_{n} \|<\rho_{1}$ for infinitely many $n$). 
Therefore, $m_{0}, n_{0}$ are defined through (\ref{ta2.2.3005}), 
while $\|\theta_{n_{0}-1} \|<\rho_{1}$, $\|\vartheta_{m_{0}-1} \|\geq\rho$. 
Combining this with (\ref{ta2.2.3001}), we conclude
$	
	\|\vartheta_{n_{0}-1} - \theta_{n_{0}-1} \|
	\leq 
	\tau 
	\leq 
	1/2
$.
Consequently, 
\begin{align}\label{ta2.2.3009}
	&
	\|\vartheta_{n_{0}-1} \|
	\leq 
	\|\theta_{n_{0}-1} \|
	+
	\|\vartheta_{n_{0}-1} - \theta_{n_{0}-1} \|
	\leq \rho_{1}+1/2
	<\rho. 
\end{align}
Hence, $n_{0}<m_{0}$, 
$f(\theta_{n_{0} } ) \leq \tilde{C}_{1}$
(notice that $\|\theta_{n_{0} } \| = \|\vartheta_{n_{0}-1} \|\leq\rho_{1}+1$). 

Let $i_{0}, j_{0}$ be the integers defined by 
$
	j_{0}
	=
	\max\{j\geq 0: n_{j}< m_{0} \}
$, 
$ 
	i_{0} 
	=
	n_{j_{0} }
$. 
Then, we have $n_{0}\leq i_{0}=n_{j_{0} }< m_{0}\leq a(i_{0},\tau ) = n_{j_{0}+1}$. 
Combining this with (\ref{ta2.2.3001}), (\ref{ta2.2.5}), we get
\begin{align*}
	&
	\|\vartheta_{m_{0}-1} - \theta_{m_{0}-1} \|
	\leq
	\tau
	\leq
	1/2, 
	\;\;\;\;\; 
	\|\theta_{i_{0} }-  \theta_{m_{0}-1} \| 
	\leq 
	4\tilde{C}_{2} \tau 
	\leq 
	1/2.  
\end{align*}
Therefore, 
\begin{align}
	&\label{ta2.2.703}
	\|\theta_{i_{0} } \|
	\geq 
	\|\vartheta_{m_{0}-1} \| 
	-
	\|\vartheta_{m_{0} } - \theta_{m_{0}-1} \|
	-
	\|\theta_{i_{0} }-  \theta_{m_{0}-1} \| 
	>
	\rho - 1
	=
	\rho_{2}. 
\end{align}

Let $\theta_{0}(\cdot )$ be the stochastic processes defined by 
$\theta_{0}(t)=\theta_{n}$ for $t\in[\gamma_{n},\gamma_{n+1} )$, $n\geq 0$
($\{\gamma_{n} \}_{n\geq 0}$ is defined in Section \ref{section1}). 
As in the proof of Theorem \ref{theorema2.1}, 
we now show by contradiction that $\gamma_{i_{0} } - \gamma_{n_{0} }\geq T$. 
We assume the opposite. 
Then, (\ref{ta2.2.3071}), (\ref{ta2.2.51}), (\ref{ta2.2.3009}) yield  
\begin{align}\label{ta2.2.705}
	\|\theta_{0}(t) \|
	=
	\|\theta_{j} \|
	\leq&
	\|\theta_{n_{0} } \| 
	+
	\sum_{i=n_{0} }^{j-1} \alpha_{i} \|\nabla f(\theta_{i} ) \| 
	+
	\left\|
	\sum_{i=n_{0} }^{j-1} \alpha_{i} \xi_{i} 
	\right\|
	\nonumber\\
	\leq&
	\|\theta_{n_{0} } \| 
	+
	\tau^{2} 
	+
	\sum_{i=n_{0} }^{j-1} \alpha_{i} \|\nabla f(\theta_{i} ) \|   
	\nonumber\\
	\leq&
	\rho_{1}+1
	+
	\sum_{i=n_{0} }^{j-1} \alpha_{i} \phi(\|\theta_{i} \| )  
	\nonumber\\
	\leq&
	\rho_{1}+1
	+
	\int_{\gamma_{n_{0} } }^{t} \phi(\|\theta_{0}(s) \| ) ds
\end{align}
for $t\in [\gamma_{j}, \gamma_{j+1} )$, $n_{0}\leq j\leq i_{0}$.\footnote
{Since $j\leq i_{0}<m_{0}$, we have 
$\gamma_{j}-\gamma_{n_{0} }\leq \gamma_{i_{0} }-\gamma_{n_{0} } \leq T$
and $j<\min\{m_{0},a(n_{0},T)+1\}$. 
We also have $\tau^{2}\leq 1/2$. }
Owing to the comparison principle (see \cite[Section 3.4]{khalil}) and 
(\ref{ta2.1.701}), (\ref{ta2.2.705}), we have  
$\|\theta_{0}(t) \| \leq 
\lambda(t-\gamma_{n_{0} };\rho_{1}+1) \leq \rho_{2}$
for all $t\in[\gamma_{{n}_{0} }, \gamma_{i_{0} } ]$. 
Thus, $\|\theta_{i_{0} } \|=\|\theta_{0}(\gamma_{i_{0} } ) \|\leq \rho_{2}$. 
However, this is impossible, due to (\ref{ta2.2.703}). 
Hence, $\gamma_{i_{0} }-\gamma_{n_{0} }\geq T$. 
Consequently, 
\begin{align}\label{ta2.2.707}
	T
	\leq 
	\gamma_{i_{0} } - \gamma_{n_{0} }
	=
	\sum_{j=0}^{j_{0}-1} (\gamma_{n_{j+1} } - \gamma_{n_{j} } )
	\leq 
	j_{0}\tau 
\end{align}
(notice that $n_{j_{0} } = i_{0}$, 
$\gamma_{n_{j+1} }-\gamma_{n_{j} } = \sum_{i=n_{j} }^{n_{j+1}-1} \alpha_{i} \leq \tau$). 

Due to (\ref{ta2.2.7}), 
we have 
$
	f(\theta_{n_{j+1} } ) - f(\theta_{n_{j} } ) 
	\leq 
	-\tau c^{2}/8
$
for $0\leq j \leq j_{0}$. 
Then, (\ref{ta2.2.707}) implies 
\begin{align*}
	f(\theta_{i_{0} } )
	=
	f(\theta_{n_{j_{0} } } ) 
	\leq 
	f(\theta_{n_{0} } ) - j_{0}\tau c^{2}/8
	\leq 
	\tilde{C}_{1} - c^{2} T/8
	\leq 
	-\tilde{C}_{1}. 
\end{align*}
However, this is impossible, since $f(\theta ) > -\tilde{C}_{1}$ for all 
$\theta\in\mathbb{R}^{d_{\theta } }$. 
Hence, $\|\theta_{n} \|\geq\rho$ for finitely many $n$. 
\end{sproof}

\refstepcounter{appendixcounter}\label{appendix3}
\section*{Appendix \arabic{appendixcounter} }

In this section, a global version of the results of Section \ref{section2} is presented.
It is also shown how the results of Section \ref{section2} can be extended 
to the randomly projected stochastic gradient search with Markovian dynamics. 
The results provided in this section can be considered as 
a combination of Theorems \ref{theorema2.1}, \ref{theorema2.2}
(Appendix \ref{appendix2}) with Theorem \ref{theorem2.1} (Section \ref{section2}). 

First, the stability and the global asymptotic behavior of algorithm (\ref{2.1}) are studied.
To analyze these properties, we use the following assumption. 

\begin{assumptionappendix}\label{aa3.1} 
There exists a Borel-measurable function $\varphi:\mathbb{R}^{d_{z} }\rightarrow[1,\infty )$
such that 
\begin{align*}
	&
	\max\{\|F(\theta,z) \|, \|\tilde{F}(\theta,z) \|, \|(\Pi\tilde{F} )(\theta,z) \| \}
	\leq
	\varphi(z) (\|\nabla f(\theta ) \| + 1 ), 
	\\
	&
	\|(\Pi\tilde{F} )(\theta',z) - (\Pi\tilde{F} )(\theta'',z) \|
	\leq 
	\varphi(z) \|\theta' - \theta'' \|
\end{align*}
for all $\theta,\theta',\theta''\in\mathbb{R}^{d_{\theta } }$, 
$z\in\mathbb{R}^{d_{z} }$. 
In addition to this, 
\begin{align*}
	\sup_{n\geq 0} 
	E(\varphi^{2}(Z_{n} )|\theta_{0}=\theta,Z_{0}=z)
	<\infty
\end{align*}
for all $\theta\in\mathbb{R}^{d_{\theta } }$, 
$z\in\mathbb{R}^{d_{z} }$. 
\end{assumptionappendix}

Assumption \ref{aa3.1} is a global version of Assumption \ref{a2.3}.
In a similar form, it is involved in the stability analysis of 
stochastic approximation carried out in \cite[Section II.1.9]{benveniste}. 

Our results on the stability of algorithm (\ref{2.1}) are provided in the next theorem. 

\begin{theoremappendix}\label{theorema3.1}
Let Assumptions \ref{a2.1}, \ref{a2.2}, \ref{aa2.1} and \ref{aa3.1} hold. 
Then, there exists a real number $\rho\in[1,\infty )$
such that $\limsup_{n\rightarrow\infty } \|\theta_{n} \| < \rho$ almost surely. 
If in addition, Assumption \ref{a1.3} holds, 
then the following is also true: 
\begin{compactenum}[(i)]
\item
$\hat{\theta}=\lim_{n\rightarrow\infty } \theta_{n}$
exists and satisfies $\nabla f(\hat{\theta} ) = 0$ almost surely. 
\item
$\|\nabla f(\theta_{n} ) \|^{2} = 
o\big(\gamma_{n}^{-\hat{p} } \big)$,
$|f(\theta_{n} ) - f(\hat{\theta} ) |	=
o\big(\gamma_{n}^{-\hat{p} } \big)$
and
$\|\theta_{n} - \hat{\theta} \| =
o\big(\gamma_{n}^{-\hat{q} } \big)$ 
almost surely on 
$\{\hat{r} > r\}$. 
\item 
$\|\nabla f(\theta_{n} ) \|^{2}	=
O\big(\gamma_{n}^{-\hat{p} } \big)$,
$|f(\theta_{n} ) - f(\hat{\theta} ) | =
O\big(\gamma_{n}^{-\hat{p} } \big)$ 
and 
$\|\theta_{n} - \hat{\theta} \| =
O\big(\gamma_{n}^{-\hat{q} } \big)$
almost surely on 
$\{\hat{r} \leq r\}$. 
\item
$\|\nabla f(\theta_{n} ) \|^{2} = o(\gamma_{n}^{-p} )$ 
and 
$|f(\theta_{n} ) - f(\hat{\theta} ) | = o(\gamma_{n}^{-p} )$ 
almost surely. 
\end{compactenum}
\end{theoremappendix}

\begin{vremark}
$p$, $\hat{p}$, $\hat{q}$, $\hat{r}$ are defined in 
Theorem \ref{theorem2.1} and Corollary \ref{corollary1.1}. 
\end{vremark}

\begin{proof}
Let $g(\cdot )$ be the function defined in Assumption \ref{aa2.2}, 
while $C\in[1,\infty )$ stands for a (global) Lipschitz constant of $\nabla f(\cdot )$. 
Moreover, let $\tau=1/(9C)$. 
On the other side, let $\{\xi_{n} \}_{n\geq 0}$, 
$\{\xi_{1,n} \}_{n\geq 0}$, $\{\xi_{2,n} \}_{n\geq 0}$, $\{\xi_{3,n} \}_{n\geq 0}$
have the same meaning as in the proof of Theorem \ref{theorem2.1}, 
while $\tau_{n}$ is the stopping time defined by 
\begin{align*}
	\tau_{n}
	=
	\min\left(
	\left\{
	j\geq n: g(\theta_{n} ) g^{-1}(\theta_{j} ) > 3 
	\right\}
	\cup
	\{\infty \}
	\right)
\end{align*}
for $n\geq 0$. 
In addition to this, for $\theta\in\mathbb{R}^{d_{\theta } }$, $z\in\mathbb{R}^{d_{z} }$, 
let $E_{\theta,z}(\cdot )$ denote the conditional expectation given $\theta_{0}=\theta$, $Z_{0}=z$. 

As a direct consequence of Assumptions \ref{a2.1}, \ref{aa3.1}, we get
\begin{align*}
	E_{\theta,z}\left(\sum_{n=0}^{\infty } \alpha_{n}^{2} \gamma_{n}^{2r} \varphi^{2}(Z_{n+1} ) \right)
	<\infty
\end{align*}
for all $\theta\in\mathbb{R}^{d_{\theta } }$, $z\in\mathbb{R}^{d_{z} }$. 
We also have 
\begin{align*}
	&
	g(\theta_{n} ) \|\xi_{n} \| 
	\leq 
	g(\theta_{n} ) 
	(\|F(\theta_{n}, Z_{n+1} ) \| + \|\nabla f(\theta_{n} ) \| )
	\leq
	\varphi(Z_{n+1} ) + 1 
	\leq 
	2\varphi(Z_{n+1} ) 
\end{align*}
for $n\geq 0$. 
Consequently, 
\begin{align}\label{ta3.1.1}
	\lim_{n\rightarrow\infty } 
	\alpha_{n}\gamma_{n}^{r}\varphi(Z_{n+1} )
	=
	\lim_{n\rightarrow\infty } 
	\alpha_{n} \gamma_{n}^{r} g(\theta_{n} ) \|\xi_{n} \|
	=0
\end{align}
w.p.1. 

Let $\{m_{k} \}_{k\geq 0}$ be the sequence recursively defined by 
$m_{0}=0$ and $m_{k+1}=a(m_{k},\tau)$ for $k\geq 0$.  
Moreover, let ${\cal F}_{n}=\sigma\{\theta_{0},Z_{0},\dots,\theta_{n},Z_{n} \}$ for $n\geq 0$. 
Due to Assumption \ref{a2.2}, we have
\begin{align*}
	E_{\theta,z}
	\left(
	g(\theta_{n} ) \xi_{1,j} 
	I_{ \{\tau_{n} > j \} }
	|
	{\cal F}_{j} 
	\right)
	=
	g(\theta_{n} ) 
	\left(
	E_{\theta,z}
	(
	\tilde{F}(\theta_{j}, Z_{j+1} )
	|
	{\cal F}_{j} 
	)
	-
	(\Pi\tilde{F} )(\theta_{j},Z_{j} )
	\right) 
	I_{ \{\tau_{n} > j \} }
	=
	0
\end{align*}
w.p.1 for each $\theta\in\mathbb{R}^{d_{\theta } }$, $z\in\mathbb{R}^{d_{z} }$, $0\leq n\leq j$
(notice that $\{\tau_{n}>j\}$ is measurable with respect to ${\cal F}_{j}$). 
On the other side, Assumption \ref{aa3.1} implies 
\begin{align*}
	g(\theta_{n} ) \|\xi_{1,j} \| I_{ \{\tau_{n}>j \} }
	\leq 
	g(\theta_{n} ) g^{-1}(\theta_{j} ) (\varphi(Z_{j} ) + \varphi(Z_{j+1} ) ) I_{ \{\tau_{n}>j \} }
	\leq 
	3 (\varphi(Z_{j} ) + \varphi(Z_{j+1} ) ) 
\end{align*}
for $0\leq n\leq j$. 
Then, as a result of Doob inequality, we get 
\begin{align*}
	E_{\theta,z}\left(
	\max_{n<j<a(n,\tau ) }
	\left\|\sum_{i=n+1}^{j}\alpha_{i}\gamma_{i}^{r}g(\theta_{n} ) \xi_{1,i} \right\|^{2} I_{ \{\tau_{n}>j\} } 
	\right)
	\leq&
	E_{\theta,z}\left(
	\max_{n<j<a(n,\tau ) }
	\left\|\sum_{i=n+1}^{j}\alpha_{i}\gamma_{i}^{r}g(\theta_{n} ) \xi_{1,i} I_{ \{\tau_{n}>i\} } \right\|^{2}
	\right)
	\nonumber\\
	\leq&
	4E_{\theta,z}\left(
	\sum_{i=n+1}^{a(n,\tau )-1} \alpha_{i}^{2}\gamma_{i}^{2r}  
	g^{2}(\theta_{n} ) \|\xi_{1,i} \|^{2} I_{ \{\tau_{n}>i\} } 
	\right)
	\nonumber\\
	\leq& 
	72E_{\theta,z}\left(\sum_{i=n+1}^{a(n,\tau) } \alpha_{i}^{2}\gamma_{i}^{2r}  
	\left(\varphi^{2}(Z_{i} ) + \varphi^{2}(Z_{i+1} ) \right)
	\right)
\end{align*}
for all $\theta\in\mathbb{R}^{d_{\theta } }$, $z\in\mathbb{R}^{d_{z} }$, $n\geq 0$. 
Combining this with Assumptions \ref{a2.1}, \ref{aa3.1}, we deduce 
\begin{align*}
	E_{\theta,z}\!\left(
	\sum_{k=0}^{\infty } 
	g^{2}(\theta_{m_{k} } ) 
	\max_{m_{k}<j<m_{k+1} }
	\left\|\sum_{i=m_{k} }^{j}\alpha_{i}\gamma_{i}^{r}\xi_{1,i} \right\|^{2} \!\! I_{ \{\tau_{m_{k} }>j\} } 
	\right)
	&
	\leq 
	72 E_{\theta,z}\!\left(\sum_{n=1}^{\infty } 
	(\alpha_{i-1}^{2}\gamma_{i-1}^{2r} + \alpha_{i}^{2}\gamma_{i}^{2r} ) \varphi^{2}(Z_{i} ) \right)
	<
	\infty
\end{align*}
for each $\theta\in\mathbb{R}^{d_{\theta } }$, $z\in\mathbb{R}^{d_{z} }$, $n\geq 0$. 
Therefore, 
\begin{align}\label{ta3.1.5}
	\lim_{k\rightarrow\infty } 
	g(\theta_{m_{k} } ) 
	\max_{m_{k}<j<m_{k+1} }
	\left\|\sum_{i=m_{k} }^{j}\alpha_{i}\gamma_{i}^{r}\xi_{1,i} \right\| I_{ \{\tau_{m_{k} }>j\} } 
	=
	0
\end{align}
w.p.1. 

Since $\alpha_{n-1}\alpha_{n}\gamma_{n}^{r} = O(\alpha_{n}^{2}\gamma_{n}^{r} )$, 
$\alpha_{n}\gamma_{n}^{r}-\alpha_{n+1}\gamma_{n+1}^{r} = O(\alpha_{n}^{2}\gamma_{n}^{r} )$ 
for $n\rightarrow\infty$
(see the proof of Theorem \ref{theorem2.1}), 
Assumptions \ref{a2.1}, \ref{aa3.1} yield
\begin{align*}
	E_{\theta,z}
	\left(
	\sum_{n=0}^{\infty } 
	\alpha_{n}\alpha_{n+1}\gamma_{n+1}^{r} \varphi^{2}(Z_{n+1} )
	\right)
	<\infty, 
	\;\;\;\;\;
	E_{\theta,z}
	\left(
	\sum_{n=0}^{\infty } |\alpha_{n}\gamma_{n}^{r} - \alpha_{n+1}\gamma_{n+1}^{r} | \varphi^{2}(Z_{n+1} ) 
	\right)
	<\infty
\end{align*}
for all $\theta\in\mathbb{R}^{d_{\theta } }$, $z\in\mathbb{R}^{d_{z} }$. 
On the other side, 
due to Assumption \ref{aa3.1}, we have 
\begin{align*}
	g(\theta_{n} ) \|\xi_{2,j} \| I_{ \{\tau_{n}>j \} }
	\leq &
	g(\theta_{n} ) \varphi(Z_{j} ) \|\theta_{j} - \theta_{j-1} \| I_{ \{\tau_{n}>j-1 \} } 
	\\
	\leq &
	\alpha_{j-1} g(\theta_{n} ) \varphi(Z_{j} ) 
	\|F(\theta_{j-1},Z_{j} ) \| I_{ \{\tau_{n}>j \} } 
	\\
	\leq &
	\alpha_{j-1} g(\theta_{n} ) g^{-1}(\theta_{j-1} ) \varphi^{2}(Z_{j} ) I_{ \{\tau_{n}>j \} } 
	\\
	\leq &
	3C\alpha_{j-1} \varphi^{2}(Z_{j} ) 
\end{align*}
for $0\leq n<j$.
Consequently, 
\begin{align*}
	&
	g(\theta_{n} ) 
	\left\|
	\sum_{i=n+1}^{j} \alpha_{i}\gamma_{i}^{r} \xi_{2,i} 
	\right\|
	I_{ \{\tau_{n}>j \} }
	\leq 
	\sum_{i=n+1}^{j} \alpha_{i}\gamma_{i}^{r} g(\theta_{n} ) \|\xi_{2,i} \| I_{ \{\tau_{n}>i \} }
	\leq 
	3\sum_{i=n}^{j} \alpha_{i} \alpha_{i+1} \gamma_{i+1}^{r} \varphi^{2}(Z_{i+1} ), 
\end{align*}
for $0\leq n<j$. 
We also have 
\begin{align*}
	&
	g(\theta_{n} ) \|\xi_{3,j} \| I_{ \{\tau_{n}>j \} } 
	\leq 
	g(\theta_{n} ) g^{-1}(\theta_{j} ) \varphi(Z_{j+1} ) I_{ \{\tau_{n}>j \} } 
	\leq 
	3 \varphi(Z_{j+1} ) 
	\leq 
	3 \varphi^{2}(Z_{j+1} )
\end{align*}
for $0\leq n\leq j$. 
Therefore, 
\begin{align*}
	g(\theta_{n} ) 
	\left\|
	\sum_{i=n+1}^{j} (\alpha_{i}\gamma_{i}^{r} - \alpha_{i+1}\gamma_{i+1}^{r} ) \xi_{3,i} 
	\right\|
	I_{ \{\tau_{n}>j \} }
	\leq &
	\sum_{i=n+1}^{j} |\alpha_{i}\gamma_{i}^{r} - \alpha_{i+1}\gamma_{i+1}^{r} | g(\theta_{n} ) \|\xi_{3,i} \| 
	I_{ \{\tau_{n}>i \} }
	\\
	\leq &
	3 \sum_{i=n+1}^{j} |\alpha_{i}\gamma_{i}^{r} - \alpha_{i+1}\gamma_{i+1}^{r} | \varphi^{2}(Z_{i+1} )
\end{align*}
for $0\leq n<j$. 
Hence, 
\begin{align}\label{ta3.1.21}
	&
	\lim_{n\rightarrow\infty }
	g(\theta_{n} ) 
	\max_{j>n} 
	\left\|
	\sum_{i=n+1}^{j} \alpha_{i}\gamma_{i}^{r} \xi_{2,i} 
	\right\|
	I_{ \{\tau_{n}>j \} }
	=
	\lim_{n\rightarrow\infty } 
	g(\theta_{n} ) 
	\max_{j>n} 
	\left\|
	\sum_{i=n+1}^{j} (\alpha_{i}\gamma_{i}^{r} - \alpha_{i+1}\gamma_{i+1}^{r} ) \xi_{3,i} 
	\right\|
	I_{ \{\tau_{n}>j \} }
	=
	0
\end{align}
w.p.1. 
On the other side, (\ref{ta3.1.1}) yields 
\begin{align}\label{ta3.1.23}
	\lim_{n\rightarrow\infty } 
	g(\theta_{n} ) 
	\max_{j>n} \alpha_{j}\gamma_{j}^{r} \|\xi_{3,j-1} \| I_{ \{ \tau_{n}>j \} }
	=
	0
\end{align}
w.p.1. 
Combining (\ref{ta3.1.1}) -- (\ref{ta3.1.23}) with (\ref{2.3*}), 
we deduce 
\begin{align*}
	\lim_{k\rightarrow\infty } 
	g(\theta_{n_{k} } ) 
	\max_{m_{k}\leq j < m_{k+1} }
	\left\|
	\sum_{i=m_{k} }^{j} \alpha_{i}\gamma_{i}^{r} \xi_{i}
	\right\|
	I_{ \{\tau_{m_{k} }>j \} }
	=
	0
\end{align*}
w.p.1. 
Then, relying on the same arguments as in the proof of Lemma \ref{lemma1.1}, 
we conclude
\begin{align}\label{ta3.1.25}
	\lim_{k\rightarrow\infty } 
	g(\theta_{n_{k} } ) 
	\max_{m_{k}\leq j < m_{k+1} }
	\left\|
	\sum_{i=m_{k} }^{j} \alpha_{i} \xi_{i}
	\right\|
	I_{ \{\tau_{m_{k} }>j \} }
	=
	0
\end{align}
w.p.1. 

Owing to Assumption \ref{aa2.1}, we have 
\begin{align*}
	g^{-1}(\theta_{j+1} ) I_{ \{\tau_{n}>j \} }
	\leq &
	g^{-1}(\theta_{n} ) 
	+
	\|\nabla f(\theta_{j+1} ) - \nabla f(\theta_{n} ) \| I_{ \{\tau_{n}>j \} }
	\nonumber\\
	\leq &
	g^{-1}(\theta_{n} ) 
	+
	C\sum_{i=n}^{j} \alpha_{i} \|\nabla f(\theta_{i} ) \| I_{ \{\tau_{n}>j \} }
	+
	C\left\|\sum_{i=n}^{j}\alpha_{i}\xi_{i} \right\| I_{ \{\tau_{n}>j \} }
	\\
	\leq &
	g^{-1}(\theta_{n} ) 
	+
	C\left\|\sum_{i=n}^{j}\alpha_{i}\xi_{i} \right\| I_{ \{\tau_{n}>j \} }
	+
	C \sum_{i=n}^{j-1} \alpha_{i} g^{-1}(\theta_{i} ) I_{ \{\tau_{n}>j \} }
\end{align*}
for $0\leq n\leq j$. 
Combining this with Bellman-Gronwall inequality (see e.g., \cite[Appendix B]{borkar}), we conclude 
\begin{align*}
	g^{-1}(\theta_{j+1} ) I_{ \{\tau_{n}>j \} }
	\leq &
	\left(
	g^{-1}(\theta_{n} ) 
	+
	C \max_{n\leq j<a(n,\tau ) } 
	\left\|\sum_{i=n}^{j} \alpha_{i}\xi_{i} \right\| I_{ \{\tau_{n}>j \} }
	\right)
	\exp\left(C \sum_{i=n}^{j-1} \alpha_{i} \right)
	\\
	\leq &
	2 g^{-1}(\theta_{n} ) 
	\left(
	1
	+
	C g(\theta_{n} ) 
	\max_{n\leq j<a(n,\tau ) } 
	\left\|\sum_{i=n}^{j} \alpha_{i}\xi_{i} \right\| I_{ \{\tau_{n}>j \} }
	\right)
\end{align*}
for $0\leq n\leq j\leq a(n,\tau )$.\footnote
{Notice that $\sum_{i=n}^{j-1}\alpha_{i}\leq\tau$ for $n\leq j\leq a(n,\tau)$. 
Notice also that $\exp(C\tau )\leq\exp(1/2)\leq 2$. }
Then, (\ref{ta3.1.25}) yields
\begin{align}\label{ta3.1.27}
	\limsup_{k\rightarrow\infty } 
	g(\theta_{m_{k} } )
	\max_{m_{k}\leq j<m_{k+1} } g^{-1}(\theta_{j+1} ) I_{ \{\tau_{m_{k} }>j \} }
	\leq 
	2
\end{align}
w.p.1. 

Let $N_{0}$ be the event where (\ref{ta3.1.25}) or (\ref{ta3.1.27}) does not hold. 
Then, in order to prove the theorem's assertion, it is sufficient to show that 
(\ref{aa2.2.1}) is satisfied on $N_{0}^{c}$ for any $t\in(0,\infty )$. 
Let $\omega$ be any sample in $N_{0}$, while $t\in(0,\infty )$ is any real number.  
Notice that all formula which follow in the proof correspond to $\omega$. 

Due to Assumption \ref{a2.1} and (\ref{ta3.1.27}), there exists 
an integer $k_{0}\geq 0$ (depending on $\omega$)
such that 
\begin{align}\label{ta3.1.71}
	\sum_{i=m_{k} }^{m_{k+1}-1} \alpha_{i}
	\geq
	\tau/2, 
	\;\;\;\;\;
	g(\theta_{m_{k} } ) 
	\left\|\sum_{i=m_{k} }^{j} \alpha_{i} \xi_{i} \right\| I_{ \{\tau_{m_{k} }>j \} }
	\leq
	\tau, 
	\;\;\;\;\; 
	g(\theta_{m_{k} } ) g^{-1}(\theta_{j+1} ) I_{ \{\tau_{m_{k} }>j \} }
	\leq 3
\end{align}
for $k\geq k_{0}$, $m_{k}\leq j<m_{k+1}$
(notice that $\lim_{n\rightarrow\infty } \sum_{i=n}^{a(n,\tau)-1} \alpha_{i} =\tau$). 
As $\tau_{n}>n$ for $n\geq 0$, we conclude 
$\tau_{m_{k} }>m_{k+1}$ for $k\geq k_{0}$.\footnote
{If $\tau_{m_{k} }\leq m_{k+1}$, 
then $\tau_{m_{k} }=j$ and 
$g(\theta_{m_{k} } ) g^{-1}(\theta_{j} ) I_{ \{\tau_{m_{k} }>j-1 \} } =
g(\theta_{m_{k} } ) g^{-1}(\theta_{j} ) > 3$ for some $j$
satisfying $m_{k}<j\leq m_{k+1}$. } 
Consequently, 
$I_{ \{\tau_{m_{k} }>j \} } = 1$ for $k\geq k_{0}$, $m_{k}\leq j\leq m_{k+1}$. 
Combining this with (\ref{ta3.1.71}), we get
\begin{align}\label{ta3.1.43}
	g^{-1}(\theta_{j+1} ) 
	\geq &
	g^{-1}(\theta_{m_{k} } ) 
	-
	\|\nabla f(\theta_{j+1} ) - \nabla f(\theta_{n} ) \|
	\nonumber\\
	\geq &
	g^{-1}(\theta_{m_{k} } ) 
	-
	C\sum_{i=m_{k} }^{j} \alpha_{i} \|\nabla f(\theta_{i} ) \|
	-
	C\left\|\sum_{i=m_{k} }^{j} \alpha_{i}\xi_{i} \right\|
	\nonumber\\
	\geq &
	g^{-1}(\theta_{m_{k} } ) 
	-
	C \sum_{i=m_{k} }^{j} \alpha_{i} g^{-1}(\theta_{i} ) 
	-
	C\left\|\sum_{i=m_{k} }^{j} \alpha_{i}\xi_{i} \right\|
	\nonumber\\
	\geq &
	g^{-1}(\theta_{m_{k} } ) 
	(1-3C\tau-C\tau )
	\nonumber\\
	\geq &
	3^{-1} g^{-1}(\theta_{m_{k} } ) 
\end{align}
for $k\geq k_{0}$, $m_{k}\leq j<m_{k+1}$.\footnote
{Notice that $g^{-1}(\theta_{i} )\leq 3g^{-1}(\theta_{m_{k} } )$, 
$\sum_{m_{k} }^{m_{k+1}\!-\!1}\!\! \alpha_{i}\leq\tau$ when $k\geq k_{0}$, $m_{k}\leq i<m_{k+1}$.
Notice also that $C\tau=1/9$. } 

Let $n_{0}=m_{k_{0} }$, while $k(n)=\max\{k\geq 0: m_{k}\leq n \}$, 
$m(n)=m_{k(n) }$ for $n\geq 0$. 
Then, (\ref{ta3.1.43}) implies 
$g(\theta_{n} )\leq 3g(\theta_{m(n) } )$, 
$g(\theta_{m_{k+1} } )\leq 3g(\theta_{m_{k} } )$
for $n\geq n_{0}$, $k\geq k_{0}$
(notice that $k(n)\geq k_{0}$, $m_{k(n) }\leq n<m_{k(n)+1}$ when $n\geq n_{0}$). 
Hence, $g(\theta_{n} )\leq C_{n,k} \:g(\theta_{m_{k} } )$ for $n\geq n_{0}$, $k\geq m(n)$, 
where $C_{n,k}=3^{k-k(n)+1}$. 
Since 
\begin{align*}
	2^{-1} (k(j)-k(n) ) \tau
	\leq 
	\sum_{k=k(n)+1}^{k(j) } \sum_{i=m_{k} }^{m_{k+1}-1} \alpha_{i}
	\leq 
	\sum_{i=n}^{j} \alpha_{i}
	\leq 
	t 
\end{align*}
for $n_{0}\leq n\leq j\leq a(n,\tau )$, 
we conclude $k(j)-k(n)\leq 2t/\tau$ for the same $n,j$. 
Consequently, 
\begin{align*}
	g(\theta_{n} ) 
	\left\|
	\sum_{i=n}^{j} \alpha_{i} \xi_{i}
	\right\|
	=&
	g(\theta_{n} ) 
	\left\|
	\sum_{k=k(n)}^{k(j) }
	\sum_{i=m_{k} }^{m_{k+1}-1} \alpha_{i}\xi_{i}
	-
	\sum_{i=m(n)}^{n-1} \alpha_{i}\xi_{i} 
	+
	\sum_{i=m(j)}^{j} \alpha_{i}\xi_{i} 
	\right\|
	\\
	\leq&
	\sum_{k=k(n)}^{k(j)-1}
	C_{n,k} \: g(\theta_{m_{k} } ) 
	\left\|
	\sum_{i=m_{k} }^{m_{k+1}-1} \alpha_{i}\xi_{i}
	\right\|
	+
	C_{n,k(n) } \: g(\theta_{m(n) } ) 
	\left\|
	\sum_{i=m(n)}^{n-1} \alpha_{i}\xi_{i} 
	\right\|
	\\
	&+
	C_{n,k(j) } \: g(\theta_{m(j) } ) 
	\left\|
	\sum_{i=m(j)}^{j} \alpha_{i}\xi_{i} 
	\right\|
	\\
	\leq &
	C(t) 
	\max_{\stackrel{\scriptstyle m_{k}\leq l<m_{k+1} }{\scriptstyle k(n)\leq k } } 
	g(\theta_{m_{k} } ) 
	\left\|\sum_{i=m_{k} }^{l} \alpha_{i}\xi_{i} \right\|
\end{align*}
for $n_{0}\leq n\leq j\leq a(n,t)$,\footnote
{Here, the following convention is used: If the lower limit of a sum is (strictly) greater than 
the upper limit, then the sum is zero. } 
where $C(t)=(2t/\tau+3) 3^{2t/\tau+3}$.
Since $\tau_{m_{k} }>m_{k+1}$ for $k\geq k_{0}$
(i.e., $I_{ \{\tau_{m_{k} }>j\} } = 1$ for $k\geq k_{0}$, $m_{k}\leq j\leq m_{k+1}$), 
(\ref{ta3.1.25}) implies 
\begin{align*}
	\lim_{n\rightarrow\infty } 
	g(\theta_{n} ) 
	\max_{n\leq j\leq a(n,t) } 
	\left\|\sum_{i=n}^{j} \alpha_{i}\xi_{i} \right\|
	=0
\end{align*}
(notice that $\lim_{n\rightarrow\infty } k(n) =\infty$). 
Hence, (\ref{aa2.2.1}) holds. 
\end{proof}

In the rest of the section, the results of Section \ref{section2} 
are extended to randomly projected stochastic gradient algorithms with Markovian dynamics.
These algorithms are defined by the following difference equations: 
\begin{align}\label{appendix3.1}
	&
	\vartheta_{n}
	=
	\theta_{n}
	-
	\alpha_{n} F(\theta_{n}, Z_{n+1} ), 
	\nonumber\\
	&
	\theta_{n+1}
	=
	\vartheta_{n} 
	I_{ \{\|\vartheta_{n} \| \leq \beta_{\sigma_{n} } \} }
	+
	\theta_{0} 
	I_{ \{\|\vartheta_{n} \| > \beta_{\sigma_{n} } \} }, 
	\nonumber\\
	&
	\sigma_{n+1}
	=
	\sigma_{n}
	+
	I_{ \{\|\vartheta_{n} \| > \beta_{\sigma_{n} } \} }, 
	\;\;\;\;\; 
	n\geq 0. 
\end{align}
Here, $H(\cdot,\cdot)$, $\{\alpha_{n} \}_{n\geq 0}$, $\{Z_{n} \}_{n\geq 0}$
have the same meaning as in Section \ref{section2}, 
while $\theta_{0}$, $\{\beta_{n} \}_{n\geq 0}$ have the same meaning 
as in the case of recursion (\ref{appendix2.1}). 

To analyze the asymptotic behavior of (\ref{appendix3.1}), we use the following assumption. 

\begin{assumptionappendix}\label{aa3.3} 
For any compact set $Q\subset \mathbb{R}^{d_{\theta } }$, 
there exists a Borel-measurable function 
$\varphi_{Q}: \mathbb{R}^{d_{z} } \rightarrow [1,\infty )$ such that 
\begin{align*}
	&
	\max\{
	\|F(\theta,z ) \|, \|\tilde{F}(\theta,z ) \|, \|(\Pi \tilde{F} )(\theta,z ) \|
	\}
	\leq 
	\varphi_{Q}(z ), 
	\\
	&
	\|(\Pi \tilde{F} )(\theta',z ) - (\Pi \tilde{F} )(\theta'',z ) \|
	\leq 
	\varphi_{Q}(z) \|\theta' - \theta'' \| 
\end{align*}
for all 
$\theta, \theta', \theta'' \in Q$, $z \in \mathbb{R}^{d_{z} }$.  
In addition to this, 
\begin{align*} 
	\sup_{n\geq 0}
	E\left(
	\varphi_{Q}^{2}(Z_{n} ) 
	|\theta_{0}=\theta, Z_{0}=z 
	\right)
	< 
	\infty
\end{align*}
for all $\theta \in \mathbb{R}^{d_{\theta } }$, $z \in \mathbb{R}^{d_{z} }$. 
\end{assumptionappendix}

In a similar form, Assumptions \ref{aa3.3} is involved in the analysis
of randomly projected stochastic approximation carried out in \cite{tadic11}. 

Our result on the asymptotic behavior of algorithm (\ref{appendix3.1}) are provided in the next theorem. 

\begin{theoremappendix}\label{theorema3.2}
Let $\{\theta_{n} \}_{n\geq 0}$ be generated by (\ref{appendix3.1}). 
Moreover, let Assumptions \ref{a2.1}, \ref{a2.2}, \ref{aa2.3} and \ref{aa3.3} hold. 
Then, there exists a real number $\rho\in[1,\infty )$
such that $\limsup_{n\rightarrow\infty } \|\theta_{n} \| < \rho$ almost surely. 
If in addition, Assumption \ref{a1.3} holds, 
then the following is also true: 
\begin{compactenum}[(i)]
\item
$\hat{\theta}=\lim_{n\rightarrow\infty } \theta_{n}$
exists and satisfies $\nabla f(\hat{\theta} ) = 0$ almost surely. 
\item
$\|\nabla f(\theta_{n} ) \|^{2} = 
o\big(\gamma_{n}^{-\hat{p} } \big)$,
$|f(\theta_{n} ) - f(\hat{\theta} ) |	=
o\big(\gamma_{n}^{-\hat{p} } \big)$
and
$\|\theta_{n} - \hat{\theta} \| =
o\big(\gamma_{n}^{-\hat{q} } \big)$ 
almost surely on 
$\{\hat{r} > r\}$. 
\item 
$\|\nabla f(\theta_{n} ) \|^{2}	=
O\big(\gamma_{n}^{-\hat{p} } \big)$,
$|f(\theta_{n} ) - f(\hat{\theta} ) | =
O\big(\gamma_{n}^{-\hat{p} } \big)$ 
and 
$\|\theta_{n} - \hat{\theta} \| =
O\big(\gamma_{n}^{-\hat{q} } \big)$
almost surely on 
$\{\hat{r} \leq r\}$. 
\item
$\|\nabla f(\theta_{n} ) \|^{2} = o(\gamma_{n}^{-p} )$ 
and 
$|f(\theta_{n} ) - f(\hat{\theta} ) | = o(\gamma_{n}^{-p} )$ 
almost surely. 
\end{compactenum}
\end{theoremappendix}

\begin{vremark}
$p$, $\hat{p}$, $\hat{q}$, $\hat{r}$ are defined in 
Theorem \ref{theorem2.1} and Corollary \ref{corollary1.1}. 
\end{vremark}

\begin{proof}
Let $Q\subset\mathbb{R}^{d_{\theta } }$ be any compact set, 
while $t\in(0,\infty )$ is any real number. 
Moreover, let $C_{Q}\in[1,\infty )$ be an upper bound of 
$\|\nabla f(\cdot )\|$ on $Q$. 
In order to prove the theorem's assertion, it is sufficient to show that 
(\ref{aa2.4.1}), (\ref{aa2.4.1}) hold w.p.1. 

Due to Assumptions \ref{a2.1}, \ref{aa3.3}, we have 
\begin{align*}
	E_{\theta,z}\left(\sum_{n=0}^{\infty } \alpha_{n}^{2}\gamma_{n}^{2r} \varphi_{Q}^{2}(Z_{n+1} ) \right)
	<\infty
\end{align*}
for all $\theta\in\mathbb{R}^{d_{\theta } }$, $z\in\mathbb{R}^{d_{z} }$. 
Assumption \ref{aa3.3} also yields 
\begin{align*}
	&
	\|\xi_{n} \| I_{ \{\theta_{n}\in Q \} }
	\leq 
	(\|F(\theta_{n}, Z_{n+1} ) \| + \|\nabla f(\theta_{n} ) \| ) I_{ \{\theta_{n}\in Q \} }
	\leq 
	\varphi_{Q}(Z_{n+1} ) + C_{Q} 
	\leq 
	2C_{Q} \varphi_{Q}(Z_{n+1} ) 
\end{align*}
for $n\geq 0$. 
Consequently, 
\begin{align}\label{ta3.2.1}
	\lim_{n\rightarrow\infty } 
	\alpha_{n}\gamma_{n}^{r}\varphi_{Q}(Z_{n+1} )
	=
	\lim_{n\rightarrow\infty } 
	\alpha_{n}\gamma_{n}^{r} \|\xi_{n} \| I_{ \{\theta_{n}\in Q \} }
	=0
\end{align}
w.p.1. 

Let ${\cal F}_{n}=\sigma\{\theta_{0},Z_{0},\dots,\theta_{n},Z_{n} \}$ for $n\geq 0$. 
Owing to Assumption \ref{a2.2}, we have
\begin{align*}
	E_{\theta,z}
	\left(
	\xi_{1,n} 
	I_{ \{\theta_{n}\in Q \} }
	|
	{\cal F}_{n} 
	\right)
	=
	\left(
	E_{\theta,z}
	(
	\tilde{F}(\theta_{n}, Z_{n+1} )
	|
	{\cal F}_{n} 
	)
	-
	(\Pi\tilde{F} )(\theta_{n},Z_{n} )
	\right) 
	I_{ \{\theta_{n}\in Q \} }
	=
	0
\end{align*}
w.p.1 for each $\theta\in\mathbb{R}^{d_{\theta } }$, $z\in\mathbb{R}^{d_{z} }$, $n\geq 0$. 
On the other side, Assumption \ref{aa3.3} implies 
\begin{align*}
	\|\xi_{1,n} \| I_{ \{\theta_{n} \in Q \} }
	\leq 
	\varphi_{Q}(Z_{n} ) + \varphi_{Q}(Z_{n+1} )   
\end{align*}
for $n\geq 0$. 
Combining this with Assumptions \ref{a2.1}, \ref{aa3.3}, we get 
\begin{align*}
	E_{\theta,z}\left(
	\sum_{n=0}^{\infty } \alpha_{n}^{2}\gamma_{n}^{2r} \|\xi_{1,n} \|^{2} I_{ \{\theta_{n}\in Q \} }
	\right)
	\leq
	2E_{\theta,z}\left(
	\sum_{n=0}^{\infty } (\alpha_{n}^{2}\gamma_{n}^{2r} + \alpha_{n+1}^{2}\gamma_{n+1}^{2r} ) 
	\varphi_{Q}^{2}(Z_{n+1} )
	\right)
	<\infty
\end{align*}
for all $\theta\in\mathbb{R}^{d_{\theta } }$, $z\in\mathbb{R}^{d_{z} }$. 
Then, using Doob theorem, we conclude that 
$\sum_{n=0}^{\infty } \alpha_{n}\gamma_{n}^{r} \xi_{1,n} I_{ \{\theta_{n} \in Q \} }$
converges w.p.1. 
Since
\begin{align*}
	\left\|\sum_{i=n}^{j}\alpha_{i}\gamma_{i}^{r}\xi_{1,i} I_{ \{\theta_{i}\in Q \} } \right\| 
	\leq
	\left\|\sum_{i=n}^{j}\alpha_{i}\gamma_{i}^{r}\xi_{1,i} \right\| I_{ \{\tau_{Q,n}>j \} } 
\end{align*}
for $0\leq n\leq j$ (notice that $\theta_{i}\in Q$ for $n\leq i<\tau_{Q,n}$), 
we deduce  
\begin{align}\label{ta3.2.5}
	\lim_{n\rightarrow\infty }  
	\max_{j\geq n}
	\left\|\sum_{i=n}^{j}\alpha_{i}\gamma_{i}^{r}\xi_{1,i} \right\| I_{ \{\tau_{Q,n}>j \} } 
	=
	0
\end{align}
w.p.1. 

As $\alpha_{n-1}\alpha_{n}\gamma_{n}^{r} = O(\alpha_{n}^{2}\gamma_{n}^{r} )$, 
$\alpha_{n}\gamma_{n}^{r}-\alpha_{n+1}\gamma_{n+1}^{r} = O(\alpha_{n}^{2}\gamma_{n}^{r} )$ 
for $n\rightarrow\infty$
(see the proof of Theorem \ref{theorem2.1}), 
Assumptions \ref{a2.1}, \ref{aa3.3} yield
\begin{align*}
	E_{\theta,z}
	\left(
	\sum_{n=0}^{\infty } 
	\alpha_{n}\alpha_{n+1}\gamma_{n+1}^{r} \varphi_{Q}^{2}(Z_{n+1} )
	\right)
	<\infty, 
	\;\;\;\;\;
	E_{\theta,z}
	\left(
	\sum_{n=0}^{\infty } |\alpha_{n}\gamma_{n}^{r} - \alpha_{n+1}\gamma_{n+1}^{r} | 
	\varphi_{Q}^{2}(Z_{n+1} ) 
	\right)
	<\infty
\end{align*}
for all $\theta\in\mathbb{R}^{d_{\theta } }$, $z\in\mathbb{R}^{d_{z} }$. 
On the other side, 
owing to Assumptions \ref{aa3.3}, we have 
\begin{align*}
	\|\xi_{2,j} \| I_{ \{\tau_{Q,n}>j \} }
	\leq &
	\varphi_{Q}(Z_{j} ) \|\theta_{j} - \theta_{j-1} \| I_{ \{\tau_{Q,n}>j\} } 
	\\
	\leq &
	\alpha_{j-1} \varphi(Z_{j} ) 
	\|F(\theta_{j-1},Z_{j} ) \| I_{ \{\theta_{j-1}\in Q \} } 
	\\
	\leq &
	\alpha_{j-1} \varphi_{Q}^{2}(Z_{j} ) 
\end{align*}
for $0\leq n<j$. 
Consequently, 
\begin{align*}
	&
	\left\|
	\sum_{i=n+1}^{j} \alpha_{i}\gamma_{i}^{r} \xi_{2,i} 
	\right\|
	I_{ \{\tau_{Q,n}>j \} }
	\leq 
	\sum_{i=n+1}^{j} \alpha_{i}\gamma_{i}^{r} \|\xi_{2,i} \| I_{ \{\tau_{Q,n}>i \} }
	\leq 
	\sum_{i=n}^{j} \alpha_{i} \alpha_{i+1} \gamma_{i+1}^{r} \varphi_{Q}^{2}(Z_{i+1} ), 
\end{align*}
for $0\leq n<j$. 
We also have 
\begin{align*}
	&
	\|\xi_{3,n} \| I_{ \{\theta_{n}\in Q \} } 
	\leq 
	\varphi_{Q}(Z_{n+1} )  
	\leq 
	\varphi_{Q}^{2}(Z_{n+1} )
\end{align*}
for $n\geq 0$. 
Therefore, 
\begin{align*}
	\left\|
	\sum_{i=n+1}^{j} (\alpha_{i}\gamma_{i}^{r} - \alpha_{i+1}\gamma_{i+1}^{r} ) \xi_{3,i} 
	\right\|
	I_{ \{\tau_{Q,n}>j \} }
	\leq &
	\sum_{i=n+1}^{j} |\alpha_{i}\gamma_{i}^{r} - \alpha_{i+1}\gamma_{i+1}^{r} | \: \|\xi_{3,i} \| 
	I_{ \{\theta_{i}\in Q \} }
	\\
	\leq &
	\sum_{i=n+1}^{j} |\alpha_{i}\gamma_{i}^{r} - \alpha_{i+1}\gamma_{i+1}^{r} | \varphi_{Q}^{2}(Z_{i+1} )
\end{align*}
for $0\leq n<j$. 
Hence, 
\begin{align}\label{ta3.2.21}
	&
	\lim_{n\rightarrow\infty }
	\max_{j>n} 
	\left\|
	\sum_{i=n+1}^{j} \alpha_{i}\gamma_{i}^{r} \xi_{2,i} 
	\right\|
	I_{ \{\tau_{Q,n}>j \} }
	=
	\lim_{n\rightarrow\infty } 
	\max_{j>n} 
	\left\|
	\sum_{i=n+1}^{j} (\alpha_{i}\gamma_{i}^{r} - \alpha_{i+1}\gamma_{i+1}^{r} ) \xi_{3,i} 
	\right\|
	I_{ \{\tau_{Q,n}>j \} }
	=
	0
\end{align}
w.p.1. 
On the other side, (\ref{ta3.2.1}) yields 
\begin{align}\label{ta3.2.23}
	\lim_{n\rightarrow\infty } 
	\alpha_{n+1}\gamma_{n+1}^{r} \|\xi_{3,n} \| I_{ \{ \theta_{n}\in Q \} }
	=
	0
\end{align}
w.p.1. 

Since $\theta_{i}=\vartheta_{i-1}$ for $n\leq i<\tau_{Q,n}$, Assumption \ref{a2.2} and (\ref{2.3*}) yield 
\begin{align*}
	&
	\left\|
	\sum_{i=n+1}^{j}\! \alpha_{i}\gamma_{i}^{r}\xi_{i}
	\right\|
	I_{ \{\tau_{Q,n}>j \} }
	\\
	&=
	\left\|
	\sum_{i=n+1}^{j}\! 
	\alpha_{i}\gamma_{i}^{r} \xi_{1,i} 
	+\!\!
	\sum_{i=n+1}^{j}\! 
	\alpha_{i}\gamma_{i}^{r} \xi_{2,i} 
	-\!\!
	\sum_{i=n+1}^{j}\!
	(\alpha_{i}\gamma_{i}^{r} - \alpha_{i+1}\gamma_{i+1}^{r} ) \xi_{3,i} 
	-
	\alpha_{j+1}\gamma_{j+1}^{r} \xi_{3,j} 
	+
	\alpha_{n+1}\gamma_{n+1}^{r} \xi_{3,n} 
	\right\|
	I_{ \{\tau_{Q,n}>j \} }
	\\
	&\leq 
	\begin{aligned}[t]
	&
	\left\|
	\sum_{i=n+1}^{j}\! 
	\alpha_{i}\gamma_{i}^{r} \xi_{1,i} 
	\right\|
	I_{ \{\tau_{Q,n}>j \} }
	+
	\left\|
	\sum_{i=n+1}^{j}\! 
	\alpha_{i}\gamma_{i}^{r} \xi_{2,i} 
	\right\|
	I_{ \{\tau_{Q,n}>j \} }
	+
	\left\|
	\sum_{i=n+1}^{j}\! 
	(\alpha_{i}\gamma_{i}^{r} - \alpha_{i+1}\gamma_{i+1}^{r} ) \xi_{3,i} 
	\right\|
	I_{ \{\tau_{Q,n}>j \} }
	\\
	&
	+
	\alpha_{j+1}\gamma_{j+1}^{r} \|\xi_{3,j} \| I_{ \{\theta_{j}\in Q \} }
	+
	\alpha_{n+1}\gamma_{n+1}^{r} \|\xi_{3,n} \| I_{ \{\theta_{n}\in Q \} }
	\end{aligned}
\end{align*}
for $0\leq n<j$. 
Combining this with (\ref{ta3.2.1}) -- (\ref{ta3.2.23}), 
we deduce 
\begin{align*}
	\lim_{n\rightarrow\infty } 
	\max_{n\leq j \leq a(n,t) }
	\left\|
	\sum_{i=n}^{j} \alpha_{i}\gamma_{i}^{r}
	\xi_{i}
	\right\|
	I_{ \{\tau_{Q,}>j \} }
	=
	0
\end{align*}
w.p.1. 
Thus, (\ref{aa2.4.1}) holds w.p.1. 
\end{proof}

\refstepcounter{appendixcounter}\label{appendix4} 
\section*{Appendix \arabic{appendixcounter} }

In this section, we rely on the following notation. 
$d, d_{\theta }, d_{v}, d_{w} \geq 1$ are integers. 
$\Theta\subseteq\mathbb{R}^{d_{\theta } }$ is an open set, 
while ${\cal W}\subset\mathbb{R}^{d_{w} }$ is a compact set. 
$A_{\theta }$, $B_{\theta }(w)$, $F(\theta,z)$ are measurable functions 
mapping $\theta\in\Theta$, $w\in{\cal W}$, $z\in\mathbb{R}^{d_{v} }\times{\cal W}$ to 
$\mathbb{R}^{d_{v}\times d_{v} }$, $\mathbb{R}^{d_{v} }$, $\mathbb{R}^{d}$ 
(respectively). 
$\{W_{n} \}_{n\geq 0}$ is a ${\cal W}$-valued Markov chain defined on a probability space 
$(\Omega, {\cal F}, P )$, 
while $P(\cdot, \cdot )$ is its transition kernel. 
$\{V_{n}^{\theta } \}_{n\geq 0}$ is a stochastic processes defined by 
\begin{align*}
	V_{n+1}^{\theta } 
	=
	A_{\theta } V_{n}^{\theta } 
	+ 
	B_{\theta }(W_{n+1} )
\end{align*}
for $\theta\in\Theta$, $n\geq 0$, 
where $V_{0}^{\theta }\in{\cal W}$ is an arbitrary vector. 
$\{Z_{n}^{\theta } \}_{n\geq 0}$ is a Markov chain defined by 
$Z_{n}^{\theta } = [(V_{n}^{\theta } )^{T} \; W_{n}^{T} ]^{T}$
for $\theta\in\Theta$, $n\geq 0$, 
while $\Pi_{\theta }(\cdot,\cdot )$ is its transition kernel. 

\begin{theoremappendix}\label{theorema4.1}
Suppose that the following holds. 
\begin{compactenum}[(i)]
\item
$\{W_{n} \}_{n\geq 0}$ has a unique invariant probability measure
$\pi(\cdot )$. 
\item
There exist real numbers $\rho\in (0,1)$, $C\in [1,\infty )$ such that 
\begin{align*}
	|P^{n}(w,B) - \pi(B) |
	\leq 
	C\rho^{n}
\end{align*}
for all $w\in {\cal W}$, $n\geq 0$ and any measurable set $B\subseteq{\cal W}$. 
\item
For any compact set $Q\subset\Theta$, 
there exist real numbers $\varepsilon_{Q} \in (0,1)$, $K_{1,Q}\in [1,\infty )$ 
such that
$\|A_{\theta }^{n} \| \leq K_{1,Q} \varepsilon_{Q}^{n}$, 
$\|B_{\theta }(w) \| \leq K_{1,Q}$ and 
\begin{align*}
	&
	\max\{\|A_{\theta'} - A_{\theta''} \|, \|B_{\theta'}(w) - B_{\theta''}(w) \| \}
	\leq 
	K_{1,Q} \|\theta' - \theta'' \|
\end{align*}
for all $\theta,\theta',\theta''\in Q$, $w\in {\cal W}$.
\item
There exists a real number $p\in [1,\infty )$ and for any compact set $Q\subset\Theta$, 
there exists another real number $K_{2,Q}\in [1,\infty )$ such that 
\begin{align}
	&\label{ta4.1.1*}
	\|F(\theta,z) \|
	\leq 
	K_{2,Q} (1 + \|z\|^{p+1} ), 
	\\
	&\label{ta4.1.3*}
	\|F(\theta',z) - F(\theta'',z) \|
	\leq 
	K_{2,Q} \|\theta'-\theta''\| (1 + \|z\|^{p+1} ), 
	\\
	&\label{ta4.1.5*}
	\|F(\theta,z') - F(\theta,z'') \|
	\leq 
	K_{2,Q} \|z'-z''\| (1 + \|z'\|^{p} + \|z''\|^{p} )
\end{align}
for all $\theta,\theta',\theta''\in Q$, $z,z',z''\in \mathbb{R}^{d_{v} }\times{\cal W}$.
\end{compactenum}
Then, there exist measurable functions 
$g(\theta )$, $\tilde{F}(\theta,z)$ 
which map $\theta\in\Theta$, $z\in \mathbb{R}^{d_{v} } \times {\cal W}$ to $\mathbb{R}^{d}$
and
which have the following two properties: 
\begin{compactenum}[(i)]
\item
$g(\theta ) = \lim_{n\rightarrow\infty } (\Pi F)(\theta,z )$ and 
\begin{align*}
	F(\theta,z) - g(\theta )
	=
	\tilde{F}(\theta,z) - (\Pi\tilde{F} )(\theta,z)
\end{align*}
for all $\theta\in\Theta$, $z\in \mathbb{R}^{d_{v} } \times {\cal W}$, 
where 
$(\Pi\tilde{F} )(\theta,z) = \int \tilde{F}(\theta,z') \Pi_{\theta}(z,dz')$. 
\item
For any compact set $Q\subset\Theta$ and any real number $s\in (0,1)$ , 
there exists a real number $L_{Q,s}\in [1,\infty )$ 
such that 
\begin{align*}
	&
	\max\{\|\tilde{F}(\theta,z) \|, \|(\Pi\tilde{F} )(\theta,z) \| \} 
	\leq 
	L_{Q,s} (1 + \|z\|^{p+1} ), 
	\\
	&
	\|(\Pi\tilde{F} )(\theta',z) - (\Pi\tilde{F} )(\theta'',z) \|
	\leq 
	L_{Q,s} \|\theta'-\theta''\|^{s} (1 + \|z\|^{p+1} )
\end{align*}
for all $\theta,\theta',\theta''\in Q$, $z\in \mathbb{R}^{d_{v} } \times {\cal W}$. 
\end{compactenum}
\end{theoremappendix}

\begin{sproof}
Let $Q\subset\Theta$ be an arbitrary compact set. 
Moreover, let $G:\mathbb{R}^{d_{v} }\times{\cal W}\rightarrow\mathbb{R}$ 
be any function satisfying 
\begin{align}
	&\label{ta4.1.1}
	|G(z) |
	\leq 
	K(1 + \|z\|^{p+1} ), 
	\\
	&\label{ta4.1.3}
	|G(z') - G(z'') |
	\leq 
	K\|z'-z''\| (1 + \|z'\|^{p} + \|z''\|^{p} )
\end{align}
for all $z,z',z''\in\mathbb{R}^{d_{v} } \times {\cal W}$ 
and some constant $K\in [1,\infty )$. 
On the other side, for $\theta\in\Theta$, $w\in\mathbb{R}^{d_{w} }$, 
let $\tilde{B}_{\theta }(w) = [B_{\theta }^{T}(w) \; w^{T}]^{T}$. 
For the same $\theta$, 
let $\tilde{A}_{\theta }$ be the block-diagonal matrix defined as 
$\tilde{A}_{\theta } = \text{diag}\{A_{\theta},{\boldsymbol 0} \}$, 
where ${\boldsymbol 0}$ denotes $d_{w}\times d_{w}$ zero matrix. 
Then, it is straightforward to verify 
\begin{align}\label{ta4.1.5}
	(\Pi^{n} G)(\theta,z ) 
	=& 
	\int\cdots\int 
	G\left(
	\tilde{A}_{\theta }^{n} z 
	+ 
	\sum_{i=1}^{n} \tilde{A}_{\theta }^{n-i} \tilde{B}_{\theta }(w_{i} )
	\right)
	P(w_{n-1}, dw_{n} ) \cdots P(w_{0},dw_{1} )
	\nonumber\\
	=&
	\begin{aligned}[t]
	\int\cdots\int 
	&
	\left(
	G\left(
	\tilde{A}_{\theta }^{n} z 
	+ 
	\sum_{i=1}^{n} \tilde{A}_{\theta }^{n-i} \tilde{B}_{\theta }(w_{i} )
	\right)
	-
	G\left(
	\sum_{i=k}^{n} \tilde{A}_{\theta }^{n-i} \tilde{B}_{\theta }(w_{i} )
	\right)
	\right)
	\\
	&
	\cdot 
	P(w_{n-1}, dw_{n} ) \cdots P(w_{0},dw_{1} )
	\end{aligned}
	\nonumber\\
	&
	\begin{aligned}[b]
	+
	\int\cdots\int 
	&
	G\left(
	\sum_{i=k}^{n} \tilde{A}_{\theta }^{n-i} \tilde{B}_{\theta }(w_{i} )
	\right)
	P(w_{n-1}, dw_{n} ) \cdots P(w_{k},dw_{k+1} ) 
	\\
	&
	\cdot 
	(P^{k}-\pi )(w_{0},dw_{k} )
	\end{aligned} 
	\nonumber\\
	&
	\begin{aligned}[b]
	+
	\int\cdots\int 
	&
	G\left(
	\sum_{i=k}^{n} \tilde{A}_{\theta }^{n-i} \tilde{B}_{\theta }(w_{i} )
	\right)
	\!P(w_{n-1}, dw_{n} ) \cdots P(w_{k},dw_{k+1} ) 
	\pi(dw_{k} )
	\end{aligned} 
\end{align}
for all $\theta\in\Theta$, $v\in\mathbb{R}^{d_{v} }$, $w_{0}\in{\cal W}$, 
$z=[v^{T} \; w_{0}^{T}]^{T}$, $n\geq k\geq 1$. 
Using condition (iii), it is also easy to show 
\begin{align*}
	\|\tilde{A}_{\theta'}^{n+1} - \tilde{A}_{\theta''}^{n+1} \|
	=&
	\left\|
	\sum_{k=0}^{n} \tilde{A}_{\theta'}^{k} 
	(\tilde{A}_{\theta'} - \tilde{A}_{\theta''} ) \tilde{A}_{\theta''}^{n-k}
	\right\|
	\\
	\leq & 
	\sum_{k=0}^{n} \|\tilde{A}_{\theta'}^{k} \| 
	\|\tilde{A}_{\theta'} - \tilde{A}_{\theta''} \| \|\tilde{A}_{\theta''}^{n-k} \| 
	\\
	\leq &
	K_{1,Q}^{3} n \varepsilon_{Q}^{n} \|\theta'-\theta''\|
\end{align*}
for each $\theta',\theta''\in Q$, $n\geq 0$. 
Thus, there exist real numbers 
$\delta_{1,Q}\in (0,1)$, $\tilde{K}_{1,Q}\in [1,\infty )$ such that 
$\|\tilde{A}_{\theta'}^{n} - \tilde{A}_{\theta''}^{n} \| 
\leq \tilde{K}_{1,Q} \delta_{1,Q}^{n} \|\theta'-\theta''\|$ 
for any $\theta',\theta''\in Q$, $n\geq 1$. 
Consequently, condition (iii) implies that 
there exists another real number $\tilde{K}_{2,Q}\in [1,\infty )$ such that 
\begin{align}\label{ta4.1.7}
	&
	\left\|
	\left(
	\tilde{A}_{\theta'}^{n} z 
	+ 
	\sum_{i=1}^{n} \tilde{A}_{\theta'}^{n-i} \tilde{B}_{\theta'}(w_{i} )
	\right)
	-
	\left(
	\tilde{A}_{\theta''}^{n} z 
	+ 
	\sum_{i=1}^{n} \tilde{A}_{\theta''}^{n-i} \tilde{B}_{\theta''}(w_{i} )
	\right)
	\right\|
	\nonumber\\
	&
	\begin{aligned}[t]
	\leq &
	\|\tilde{A}_{\theta'}^{n} - \tilde{A}_{\theta''}^{n} \| \|z\| 
	+ 
	\sum_{i=1}^{n} 
	\|\tilde{A}_{\theta'}^{n-i} - \tilde{A}_{\theta''}^{n-i} \| 
	\|\tilde{B}_{\theta'}(w_{i} ) \|
	+
	\sum_{i=1}^{n} 
	\|\tilde{A}_{\theta''}^{n-i} \| \|\tilde{B}_{\theta'}(w_{i} ) - \tilde{B}_{\theta''}(w_{i} ) \|
	\end{aligned}
	\nonumber\\
	&
	\leq 
	\tilde{K}_{2,Q} \|\theta'-\theta'' \| (1 + \|z\| )
\end{align}
for all $\theta',\theta''\in Q$, 
$z\in\mathbb{R}^{d_{v} }\times {\cal W}$, 
$w_{1},\dots,w_{n} \in {\cal W}$, 
$n\geq 1$. 
Due to the same reasons, there also exists a real number 
$\tilde{K}_{3,Q} \in [1,\infty )$ such that 
\begin{align}\label{ta4.1.9}
	\left\|
	\tilde{A}_{\theta}^{n} z 
	+ 
	\sum_{i=k}^{l} \tilde{A}_{\theta}^{n-i} \tilde{B}_{\theta}(w_{i} )
	\right\|
	\leq 
	\|\tilde{A}_{\theta}^{n} \| \|z\| 
	+ 
	\sum_{i=k}^{l} 
	\|\tilde{A}_{\theta}^{n-i} \| 
	\|\tilde{B}_{\theta}(w_{i} ) \|
	\leq 
	\tilde{K}_{3,Q} \varepsilon_{Q}^{n-l} (1 + \|z\| )
\end{align}
for each $\theta\in Q$, 
$z\in\mathbb{R}^{d_{v} }\times{\cal W}$, 
$w_{1},\dots,w_{n} \in {\cal W}$, 
$n\geq l \geq k\geq 1$. 
Then, owing to (\ref{ta4.1.1}), (\ref{ta4.1.5}), we have 
\begin{align*}
	|(\Pi^{n} G)(\theta,z ) |
	\leq 
	2^{p+1} K \tilde{K}_{3,Q}^{p+1} 
	(1 + \|z\|^{p} )
\end{align*}
for any $\theta\in Q$, $z\in\mathbb{R}^{d_{v} }\times{\cal W}$, 
$n\geq 1$. 
On the other side, combining (\ref{ta4.1.3}) -- (\ref{ta4.1.9}), we get 
\begin{align*}
	&
	|(\Pi^{n} G)(\theta',z) - (\Pi^{n} G)(\theta'',z) |
	\nonumber\\
	&
	\leq 
	\begin{aligned}[t]
	\int\cdots\int 
	&
	\left|
	G\left(
	\tilde{A}_{\theta'}^{n} z 
	+ 
	\sum_{i=1}^{n} \tilde{A}_{\theta'}^{n-i} \tilde{B}_{\theta'}(w_{i} )
	\right)
	-
	G\left(
	\tilde{A}_{\theta''}^{n} z 
	+ 
	\sum_{i=1}^{n} \tilde{A}_{\theta''}^{n-i} \tilde{B}_{\theta''}(w_{i} )
	\right)
	\right|
	\\
	&
	\cdot 
	P(w_{n-1}, dw_{n} ) \cdots P(w_{0},dw_{1} )
	\end{aligned}
	\nonumber\\
	&
	\leq 
	3^{p+1}K \tilde{K}_{2,Q} \tilde{K}_{3,Q}^{p} 
	\|\theta' - \theta'' \| 
	(1 + \|z\|^{p+1} )
\end{align*}
for all $\theta',\theta''\in Q$, $v\in\mathbb{R}^{d_{v} }$, 
$w_{0}\in {\cal W}$, $z=[v^{T} \; w_{0}^{T} ]^{T}$, 
$n\geq 1$. 
Similarly, using 
(\ref{ta4.1.1}) -- (\ref{ta4.1.5}), (\ref{ta4.1.9}), 
we obtain
\begin{align}\label{ta4.1.25}
	&
	|(\Pi^{n} G)(\theta,z' ) - (\Pi^{n} G)(\theta,z'' ) |
	\nonumber\\
	&
	\leq 
	\begin{aligned}[t]
	&
	\begin{aligned}[t]
	\int\cdots\int 
	&
	\left|
	G\left(
	\tilde{A}_{\theta }^{n} z' 
	+ 
	\sum_{i=1}^{n} \tilde{A}_{\theta }^{n-i} \tilde{B}_{\theta }(w_{i} )
	\right)
	-
	G\left(
	\sum_{i=k}^{n} \tilde{A}_{\theta }^{n-i} \tilde{B}_{\theta }(w_{i} )
	\right)
	\right| 
	\\
	&
	\cdot 
	P(w_{n-1}, dw_{n} ) \cdots P(w_{1},dw_{2} )
	P(w'_{0},dw_{1} ) 
	\end{aligned}
	\\
	&
	+
	\begin{aligned}[t]
	\int\cdots\int 
	&
	\left|
	G\left(
	\tilde{A}_{\theta }^{n} z'' 
	+ 
	\sum_{i=1}^{n} \tilde{A}_{\theta }^{n-i} \tilde{B}_{\theta }(w_{i} )
	\right)
	-
	G\left(
	\sum_{i=k}^{n} \tilde{A}_{\theta }^{n-i} \tilde{B}_{\theta }(w_{i} )
	\right)
	\right| 
	\\
	&
	\cdot 
	P(w_{n-1}, dw_{n} ) \cdots P(w_{1},dw_{2} )
	P(w''_{0},dw_{1} ) 
	\end{aligned}
	\\
	&
	+
	\begin{aligned}[t]
	\int\cdots\int 
	&
	\left|
	G\left(
	\sum_{i=k}^{n} \tilde{A}_{\theta }^{n-i} \tilde{B}_{\theta }(w_{i} )
	\right)
	\right|
	P(w_{n-1}, dw_{n} ) \cdots P(w_{k},dw_{k+1} ) 
	\\
	&
	\cdot 
	\left(|P^{k}-\pi |(w'_{0},dw_{k} ) + |P^{k}-\pi |(w''_{0},dw_{k} ) \right)
	\end{aligned} 
	\end{aligned}
	\nonumber\\
	&
	\leq 
	3^{p+2} K \tilde{K}_{3,Q}^{p+1} \varepsilon_{Q}^{n-k} 
	(1 + \|z'\|^{p+1} + \|z''\|^{p+1} )
	+
	4CK \tilde{K}_{3,Q}^{p+1}\rho^{k} 
\end{align}
for each $\theta\in Q$, $v',v''\in\mathbb{R}^{d_{v} }$, 
$w'_{0},w''_{0}\in {\cal W}$, 
$z'=[(v' )^{T} \; (w'_{0} )^{T} ]^{T}$, $z''=[(v'' )^{T} \; (w''_{0} )^{T} ]^{T}$, 
$n\geq k\geq 1$. 
Then, setting $k=\lfloor{n/2 \rfloor}$ in (\ref{ta4.1.25}), 
we conclude that there exist real numbers 
$\delta_{2,Q} \in (0,1)$, $\tilde{K}_{4,Q} \in [1,\infty )$ such that 
\begin{align*}
	&
	|(\Pi^{n} G)(\theta, z) |
	\leq 
	\tilde{K}_{4,Q} (1 + \|z\|^{p+1} ), 
	\\
	&
	|(\Pi^{n} G)(\theta', z) - (\Pi^{n} G)(\theta'', z) |
	\leq 
	\tilde{K}_{4,Q} \|\theta' - \theta'' \| (1 + \|z\|^{p+1} ), 
	\\
	&
	|(\Pi^{n} G)(\theta, z') - (\Pi^{n} G)(\theta, z'') |
	\leq 
	\tilde{K}_{4,Q} \delta_{2,Q}^{n} (1 + \|z'\|^{p+1} + \|z''\|^{p+1} ) 
\end{align*}
for all $\theta,\theta',\theta'' \in Q$, 
$z,z',z'' \in \mathbb{R}^{d_{v} } \times {\cal W}$, $n\geq 1$. 
Combining this with the results of \cite[Section II.2.2]{benveniste}, 
we deduce that there exist functions $g(\cdot )$, $F(\cdot,\cdot )$
with the properties specified in the statement of the theorem. 
\end{sproof}

\end{document}